\documentclass[11pt,a4paper]{article}
\usepackage{comment}
\usepackage[T1]{fontenc}
\usepackage[utf8]{inputenc}
\usepackage[margin=2.5cm]{geometry}
\usepackage{amsmath,amssymb,amsthm}
\usepackage{microtype}
\usepackage[dvipsnames]{xcolor}
\usepackage{tikz}
\usetikzlibrary{positioning,calc,arrows.meta,fit,backgrounds}

\usepackage[normalem]{ulem}

\usepackage[colorlinks=true,linkcolor=NavyBlue,citecolor=NavyBlue,urlcolor=NavyBlue]{hyperref}
\usepackage{enumitem}
\usepackage{booktabs,longtable,array}

\theoremstyle{plain}
\newtheorem{theorem}{Theorem}[section]
\newtheorem{proposition}[theorem]{Proposition}
\newtheorem{lemma}[theorem]{Lemma}
\newtheorem{corollary}[theorem]{Corollary}
\theoremstyle{definition}
\newtheorem{definition}[theorem]{Definition}
\newtheorem{assumption}{Assumption}

\theoremstyle{remark}

\newcommand{\E}{\mathbb{E}}
\newcommand{\N}{\mathbb{N}}
\newcommand{\He}{\mathrm{He}}
\newcommand{\Var}{\operatorname{Var}}
\newcommand{\Cov}{\operatorname{Cov}}
\newcommand{\Aut}{\mathrm{Aut}}
\newcommand{\CC}{\mathbf{CC}}
\newcommand{\Pcal}{\mathcal{P}}
\newcommand{\Gcup}{\Gamma_{\cup}}
\newcommand{\Gomi}{\Gamma_{\ominus}}
\newcommand{\Gcap}{\Gamma_{\cap}}
\newcommand{\WS}{W_{\mathrm S}^{(r)}}
\newcommand{\WP}{W_{\mathrm P}^{(r)}}
\newcommand{\WR}{W_{\mathrm R}^{(r)}}
\newcommand{\WD}{W_{\mathrm D}^{(r)}}
\newcommand{\Vomi}[1]{V_{\ominus}^{(#1)}}
\newcommand{\PcalSL}{\Pcal^{(\mathrm{SL})}}
\newcommand{\Ppin}{\Pcal^{(\mathrm{pin})}}
\newcommand{\Ppins}{\Pcal^{(\mathrm{pin})*}}
\newcommand{\comp}{\mathrm{comp}_{\mathbf s,H}(G)}
\newcommand{\miCC}{\mathrm{miCC}_{\mathbf s,H}(G)}
\newcommand{\sdrop}{\big(|E_G|+|V_G|\big)-\big(|E_H|+|V_H|\big)}
\newcommand{\sdropb}{\left[\big(|E_G|+|V_G|\big)-\big(|E_H|+|V_H|\big)\right]}
\newcommand{\Pai}{\mathbf{P}}
\newcommand{\one}{\mathbf{1}}
\newcommand{\Gs}{G\setminus V'}
\newcommand{\Psh}{P^{\sharp}}
\newcommand{\Corr}{\mathrm{Corr}}
\newcommand{\ciso}{c_{\mathrm{iso}}}
\newcommand{\upp}{\upsilon_{p}}
\newcommand{\upd}{\upsilon_{\delta}}
\newcommand{\Dcomp}{\Delta_{\mathrm{comp}}}
\DeclareMathOperator{\Span}{span}

\title{\bf The low-degree sample complexity of clustering Gaussian mixtures with more groups than dimensions
}
\author{Alexandra Carpentier\thanks{Institut f\"ur Mathematik, Universit\"at Potsdam, Germany.}
\and Nicolas Verzelen\thanks{INRAE, MISTEA, Univ.\ Montpellier, Montpellier, France.}}
\date{\today}

\begin{document}
\maketitle

\begin{abstract}
We study the computational limits of clustering $n$ observations drawn from a mixture of $K$
isotropic Gaussian distributions in $\mathbb R^d$, when the number $K$ of groups is much larger than
the dimension $d$. In this regime, distance-based and spectral procedures require a squared
separation $\Delta^2$ of order $\sqrt d$. In contrast,  tensor-based procedures cluster at a polylogarithmic
separation but need a sample size that is a large polynomial of $dK$, and the low-degree lower bounds
known for random centers only show that $n$ must exceed $dK^2$. 

Focusing on normally distributed centers, we characterize the minimal separation for clustering with low-degree polynomials in the regime $K>d$, through a low-degree lower bound and a matching upper bound, up to subpolynomial factors. In particular, for $K=d^{p/2}$ with $p> 2$ a real number, we show that a sample size of the order of $K^3 d^{(\lceil p\rceil-p)/2}$ is required to recover the clustering with subpolynomial separation, and that it suffices when $p$ is not an integer.

\end{abstract}

\section{Introduction}\label{sec:intro}

Clustering $n$ observations into $K$ groups is one of the most basic tasks of unsupervised
learning. Its most emblematic statistical formulation is the isotropic Gaussian mixture model,
which has attracted a lot of attention both in statistics and in theoretical computer science
\cite{Dasgupta99,VEMPALA2004,LuZhou2016,Regev2017,fei2018hidden,giraud2019partial,Kwon20,loffler2021optimality,LiuLi2022,Even24}. In this paper, we are interested in
the computational limits of clustering when the number $K$ of groups is much larger than the
dimension $d$. 

\subsection{Clustering Gaussian mixtures}\label{sec:intro-mixtures}

 We observe $n$ independent vectors
$Y_i\sim\mathcal N(\theta_{k^*(i)},I_d)$, $i=1,\dots,n$, where both the centers
$\theta_1,\dots,\theta_K\in\mathbb R^d$ and the labels $k^*:[n]\to[K]$ are unknown, and where each
group $(k^*)^{-1}(\{k\})$ contains of the order of $n/K$ observations. Our aim is to estimate the
partition of $[n]$ given by $k^*$ more accurately than a random guess would. The key parameter is the
separation $\Delta_*:=\min_{k\ne k'}\|\theta_k-\theta_{k'}\|$. Without any constraint on the running
time, clustering is feasible once $\Delta_*^2$ is large compared to
\begin{equation}\label{eq:Dstat}
  \Delta_{\mathrm{stat}}^2\ :=\ \max\Big\{\sqrt{\frac{dK\log K}{n}}\,,\ \log K\Big\} ,
\end{equation}
see \cite{Regev2017,Kwon20,Even24}. The exact $K$-means estimator attains this threshold
\cite{Even24}, but computing it is NP-hard in the worst case \cite{aloise2009np}, and so is
approximating it \cite{awasthi2015hardness}. This raises the question of whether some
polynomial-time method clusters at a separation close to $\Delta_{\mathrm{stat}}$, or whether
there is a statistical-computational gap.

\subsection{Polynomial-time algorithms and low-degree lower bounds}\label{sec:intro-ld}

The worst-case hardness of computing the $K$-means estimator does not settle this question:
clustering is an average-case problem, so
worst-case hardness results say little about it, and computational barriers are instead studied
for restricted families of algorithms. The low-degree polynomial framework
\cite{Hopkins17,KuniskyWeinBandeira,SurveyWein2025} considers estimators that are multivariate
polynomials of degree at most $D$ in the data. For many problems, polynomials of degree of order
$\log n$ perform as well as the best polynomial-time algorithms known, and their failure is
therefore regarded as strong evidence that the problem is computationally hard
\cite{KuniskyWeinBandeira,SurveyWein2025}. Hopkins \cite{hopkins2018statistical} turned this
heuristic into the low-degree conjecture: for testing problems whose null distribution is a product
measure and whose planted distribution is invariant under permutations, if polynomials of degree
slightly above $\log n$ fail, then so does every polynomial-time algorithm, provided a little
independent noise is added. Some variants of this conjecture were refuted
\cite{holmgren2021counterexamples,buhai2025quasi}, and its standard polynomial-time version was
recently disproved as well \cite{mao2026polynomial}. Still, since low-degree polynomials match the
best known polynomial-time algorithms on a broad range of problems, their failure remains one of
the most convincing indications of a computational barrier.

Following \cite{SchrammWein22}, low-degree lower bounds for estimation problems are formulated through the
correlation between a functional of the hidden partition, here $x:=\one\{k^*(1)=k^*(2)\}$, and the
polynomials of degree at most $D$ in the observations $Y=(Y_{ij})_{i\in[n],j\in[d]}$. We use it in
centered form,
\begin{equation}\label{eq:corrintro}
  \Corr_{\le D}\ :=\ \sup_{\deg P\le D,\ P\ne0}\ \frac{\E\big[(x-\E x)\,P(Y)\big]}{\sqrt{\E[P(Y)^2]}} .
\end{equation}
As $\E x=1/K$, the variance of $x$ is $K^{-1}(1-K^{-1})$, and $\Corr_{\le D}^2/\Var(x)$ is the
largest fraction of the variance of $x$ that a polynomial of degree at most $D$ explains. The
uncentered correlation $\sup\{\E[xP]:\deg P\le D,\ \E[P^2]\le1\}$ of \cite{SchrammWein22,carpentier2026small} equals
$(K^{-2}+\Corr_{\le D}^2)^{1/2}$, the constant polynomial alone achieving $1/K$. So a bound
$\Corr_{\le D}\le c/K$, with $c$ a constant, means that low-degree polynomials do not improve over a
random guess by more than a constant factor.

\paragraph{High-dimensional regime ($n\le dK$).} When the dimension is large, the difficulty of
clustering is driven by dimension reduction. Lesieur et al.\ \cite{lesieur2016phase} conjectured that, for a
fixed number of groups, clustering is hard below the BBP threshold \cite{BBP05}, where $\Delta_*^2$
is of order $\sqrt{dK^2/n}$. Even, Giraud and Verzelen \cite{Even24,Even25} established that
low-degree polynomials fail as long as $\Delta_*^2$ is smaller, by a polylogarithmic factor, than
\begin{equation}\label{eq:DHD}
  \Delta^2_{\mathrm{comp,HD}}\ :=\ \max\Big\{1,\ \min\Big\{\sqrt d\,,\ \sqrt{\frac{dK^2}{n}}\Big\}\Big\} .
\end{equation}
This bound is matched, up to polylogarithmic factors and when $n\le dK$, by a combination of two
procedures \cite{VEMPALA2004,Even25}. Hierarchical clustering succeeds as soon as $\Delta_*^2$ exceeds
$\sqrt{d\log n}$ by a large constant factor, as observations of the same group are then closer than
observations of distinct groups. Above $\sqrt{dK^2/n}$, a spectral projection onto a
$\min\{K,d\}$-dimensional subspace, followed by hierarchical clustering, succeeds. Overall, among
distance-based and spectral procedures, including SDP relaxations and Lloyd's algorithm
\cite{LuZhou2016,fei2018hidden,giraud2019partial,loffler2021optimality,ndaoud2022sharp}, the best available guarantees require $\Delta_*^2$ larger than
$\min\{\max\{\sqrt{dK^2/n},\sqrt K\},\sqrt d\}$, up to polylogarithmic factors \cite{Even25}.

\paragraph{Large sample size ($n\ge dK$).} In this
regime, much less is known. The lower bound \eqref{eq:DHD} still holds, but spectral methods no
longer match it: they only succeed when $\Delta_*^2$ exceeds $\min\{\sqrt K,\sqrt d\}$, up to
polylogarithmic factors. Another line of work uses moment tensors of high order, either through
the method of moments and tensor decompositions \cite{hsu2013learning} or through the
sum-of-squares hierarchy
\cite{DBLP:journals/corr/abs-1711-07454,kothari2018robust,diakonikolas2018list}. The latter
procedures succeed at separation $\Delta_*\ge K^\gamma$ for any fixed $\gamma>0$,
with sample size and running time polynomial in $d$ and $K$, but with an exponent
that increases with $1/\gamma$. Liu and Li \cite{LiuLi2022} designed a polynomial-time method, based on
iterated projections of tensors, that recovers the partition once $\Delta_*^2$ is a constant factor
above $(\log K)^{1+\epsilon}$, for some $\epsilon>0$, which nearly matches
$\Delta^2_{\mathrm{stat}}$. It requires, however, $n\ge(dK)^{C/\epsilon}$, where $C$ is a large
numerical constant that is not made explicit.
One reason why spectral methods fail and higher-order tensors are needed is that second moments
may carry no information on the partition. Given the centers and the labels, the second-moment
matrix of the observations is
\[
\E\Big[\frac1n\sum_{i=1}^nY_iY_i^{\top}\Big]\ =\ I_d+\frac1n\sum_{k=1}^K\big|(k^*)^{-1}(\{k\})\big|\,\theta_k\theta_k^{\top} ,
\]
that is $I_d+\frac1K\sum_{k\le K}\theta_k\theta_k^{\top}$ for groups of equal sizes, and for some
configurations of the centers this matrix is proportional to $I_d$ as soon as $K>d$, for
instance when the centers are the vectors of an orthonormal basis of $\mathbb R^d$, multiplied by
$\Delta$, and their opposites, so that $K=2d$, as in the construction of \cite{carpentier2026low}.

\paragraph{Low-degree lower bounds for large $K$ ($K$ possibly much larger than $d$).} When $n\ge dK$
and $\Delta_*^2\ge\Delta^2_{\mathrm{stat}}$, which is then of order $\log K$, the data can be
projected on a space of dimension $\min\{d,K\}$ without distorting the centers~\cite{Even25}. So,
for lower bounds, one may assume $K\ge d$.
In~\cite{carpentier2026low} and~\cite{carpentier2026small}, we introduced prior distributions on
the centers that enforce algebraic conditions either on the second-moment matrix above or on
higher-order moment tensors. For $K=2d$, the centers of~\cite{carpentier2026low} are pairwise
orthogonal or opposite, so that the second-moment matrix is proportional to the identity: this
construction makes spectral approaches useless.
Low-degree polynomials fail when $\Delta_*^2$ is below $\min\{K/n^{1/4},\sqrt K\}$, and above this
threshold an efficient non-spectral procedure recovers the functional $x$ under this prior.

As several of the procedures above are based on moment tensors, \cite{carpentier2026small}
extends this cancellation to tensors of higher order: there $K=(d+1)^m$ for an integer $m$, and the observations are drawn from a mixture of $2K$ product
Bernoulli distributions whose centers are built from a dual BCH code
\cite[Chapter~9]{macwilliams1977theory}. This construction makes the moment tensors of the
observations coincide with those of a Rademacher vector up to the order $2m+1$: these tensors carry
no information on the partition, and a procedure has to use tensors of higher order, which makes the
problem harder. Low-degree
polynomials fail as long as $\Delta^2$ is small compared to $\min\{(K^2d^{3/2}/n)^{1/(2m+2)},\sqrt
d\}$, so that $n\ge K^2d^{3/2}$ is needed at a subpolynomial separation, and evidence is given that
this bound is tight for these centers, up to a possible factor $d^{1/2}$. A main message of
these two constructions is that the more centers there are, compared to the dimension, the more
moment tensors they can cancel.
Although these configurations are appealing in terms of tensor matching, other configurations of
centers could lead to harder instances. Moreover, the configurations of~\cite{carpentier2026small}
are rigid, the construction being algebraic, and the threshold could be drastically different for
random centers.

In summary, despite some recent advances, these two general questions remain largely open: 
\begin{itemize}[itemsep=2pt,topsep=3pt]
\item When the number $K$ of groups is large compared to the dimension $d$, what is the minimal
separation $\Delta_*$, as a function of $n$, $K$ and $d$, for which clustering is possible in
polynomial time?
\item More specifically, in this regime, how large must the sample size $n$ be compared to $K$ and
$d$ for clustering in polynomial time with only a polylogarithmic separation?
\end{itemize}

\subsection{Our contribution}\label{sec:contribution}

As in \cite{lesieur2016phase,banks2018information,Even24,Even25}, we focus on the arguably emblematic setting where we draw the centers at random according to a Gaussian distribution,
$\theta_k=\frac{\Delta}{\sqrt d}(X_{k,1},\dots,X_{k,d})$ with i.i.d.\ standard Gaussian entries
$X_{k,j}$, and the labels $k^*(1),\dots,k^*(n)$ i.i.d.\ uniform on $[K]$ (Section~\ref{sec:model}).
As $\|X_k-X_{k'}\|^2$ is twice a chi-square variable with $d$ degrees of freedom, a union bound over
the $K(K-1)/2$ pairs shows that $\Delta_*^2/(2\Delta^2)$ is close to $1$ with high probability when
$(\log K)/d$ is small, which is the case here. So the separation is
$\sqrt2\,\Delta$ up to a factor close to $1$, and our results can be read in terms of $\Delta_*$.

\medskip 

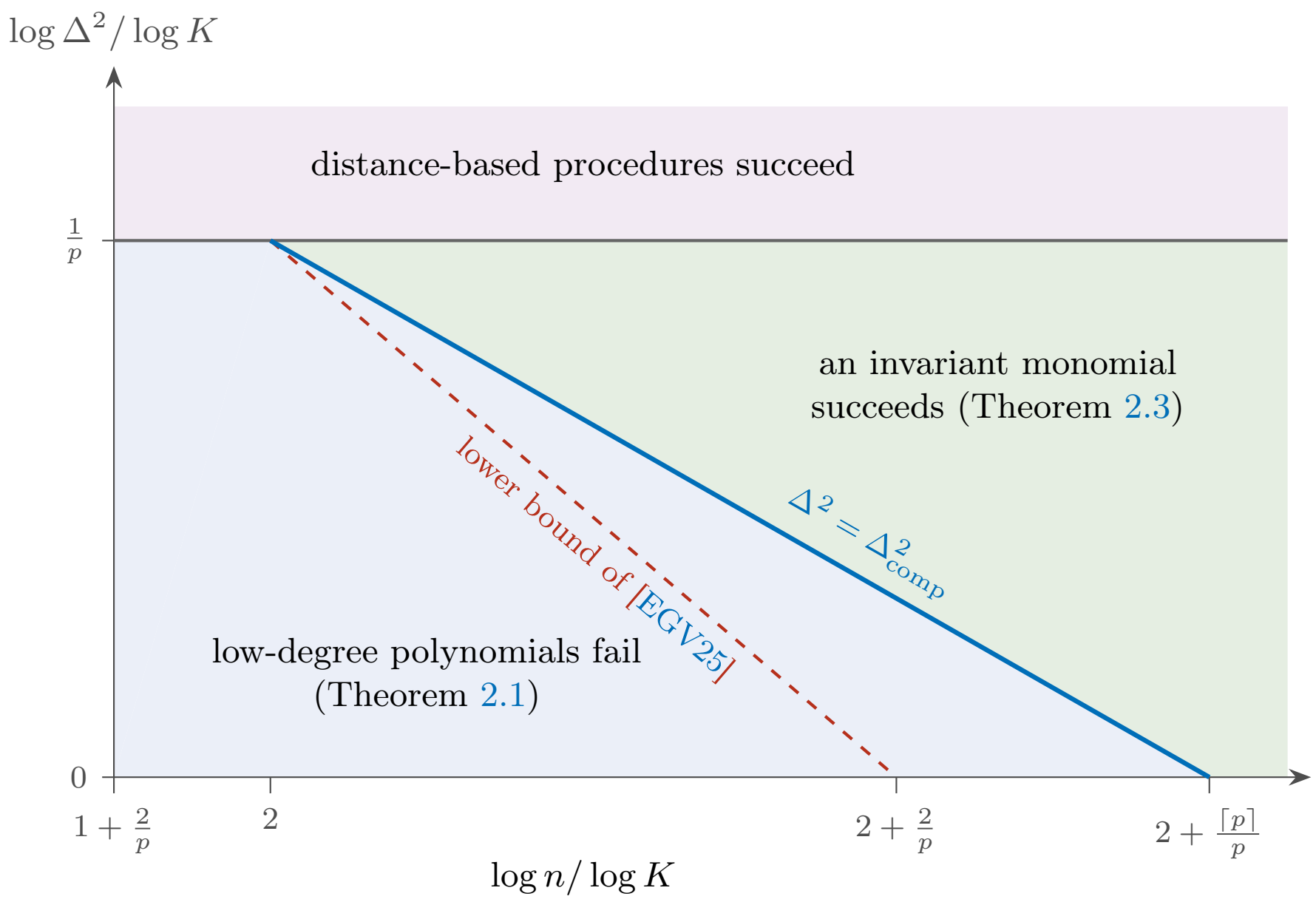
\begin{figure}[ht]
\centering
\begin{tikzpicture}[x=7cm,y=12cm,
  ax/.style={-{Stealth[length=2mm]},black!70},
  th/.style={NavyBlue,very thick},
  hd/.style={BrickRed,thick,dashed}]
  \fill[NavyBlue!7] (1.8,0) -- (2,0.4) -- (3.2,0) -- cycle;
  \fill[NavyBlue!7] (1.8,0) -- (1.8,0.4) -- (2,0.4) -- cycle;
  \fill[OliveGreen!10] (2,0.4) -- (3.3,0.4) -- (3.3,0) -- (3.2,0) -- cycle;
  \fill[Plum!8] (1.8,0.4) -- (3.3,0.4) -- (3.3,0.5) -- (1.8,0.5) -- cycle;
  \draw[ax] (1.8,0) -- (3.33,0);
  \node[font=\footnotesize] at (2.4,-0.075) {$\log n/\log K$};
  \draw[ax] (1.8,0) -- (1.8,0.53) node[above,font=\footnotesize] {$\log\Delta^2/\log K$};
  \foreach \v/\l in {1.8/{$1+\frac{2}{p}$},2/{$2$},2.8/{$2+\frac{2}{p}$},3.2/{$2+\frac{\lceil p\rceil}{p}$}}{
    \draw[black!70] (\v,0) -- (\v,-0.012) node[below,font=\scriptsize] {\l};}
  \draw[black!70] (1.8,0.4) -- (1.785,0.4) node[left,font=\scriptsize] {$\frac1p$};
  \draw[black!70] (1.8,0) -- (1.785,0) node[left,font=\scriptsize] {$0$};
  \draw[black!60,thick] (1.8,0.4) -- (3.3,0.4);
  \draw[hd] (2,0.4) -- node[sloped,below,font=\scriptsize,pos=0.55] {lower bound of \cite{Even25}} (2.8,0);
  \draw[th] (2,0.4) -- node[sloped,above,font=\scriptsize,pos=0.62] {$\Delta^2=\Dcomp^2$} (3.2,0);
  \node[font=\footnotesize,align=center] at (2.4,0.455) {distance-based procedures succeed};
  \node[font=\footnotesize,align=center] at (2.2,0.075) {low-degree polynomials fail\\(Theorem~\ref{thm:LB})};
  \node[font=\footnotesize,align=center] at (2.93,0.29) {an invariant monomial\\succeeds (Theorem~\ref{thm:UB})};
\end{tikzpicture}
\caption{Phase diagram for low-degree polynomials at $d=K^{2/p}$, drawn for $p=5/2$, so that
$\lceil p\rceil=3$ and $\sqrt d=K^{1/p}$, in the coordinates $\log n/\log K$ and $\log\Delta^2/\log
K$ and up to subpolynomial factors. Below the solid line, polynomials of
degree polylogarithmic in $K$ fail (Theorem~\ref{thm:LB}). Between this line and $\Delta^2=\sqrt d$,
the invariant monomial of a graph of degree $\lceil p\rceil$ decides $x$ (Theorem~\ref{thm:UB}), on a subsample
when $n>K^2d^{\lceil p\rceil/2}$. Above
$\sqrt d$, distance-based procedures succeed. The dashed line is the lower bound
$\sqrt{dK^2/n}$ of \cite{Even25}, the value of $\Dcomp^2$ at $\lceil p\rceil=2$. The figure starts at $n=dK$, where
the statistical threshold $\Delta^2_{\mathrm{stat}}$ becomes polylogarithmic, at the bottom of the
figure.}
\label{fig:phase}
\end{figure}

For random (Gaussian) centers, we answer both questions for low-degree polynomials, up to subpolynomial
factors, whenever $2\log K/\log d$ is not an integer\footnote{When $p=2\log K/\log d$ is an
integer, Theorem~\ref{thm:LB} still gives the lower bound, with $\lceil p\rceil=p$. But the invariant monomial of a
template of degree $p$ then fails, as $d^{p/2}=K$ exactly, and the upper bound must use a template of
degree $p+1$, which requires a polynomially larger separation. The threshold is therefore only
known to lie between two values in a subpolynomial window of values of $d$ around $K^{2/p}$
(Section~\ref{sec:integer}).}. In this introduction, we write $K=d^{p/2}$, with a real $p>2$, so
that $p=2\log K/\log d$, and we call $\lceil p\rceil$ the \textbf{critical order}: it is the
smallest integer $q$ such that $d^{q/2}\ge K$, and, in effect, the order of the first informative
moment tensor of the centers (Section~\ref{sec:bch}). The separation at which low-degree polynomials start
to cluster is then, up to subpolynomial factors,
\begin{equation}\label{eq:Dcomp}
  \Dcomp^2\ :=\ \sqrt d\,\min\Big\{1,\ \Big(\frac{K^2}{n}\Big)^{1/\lceil p\rceil}\Big\} .
\end{equation}
In Theorem~\ref{thm:LB}, we establish that polynomials of degree of order $\log K$ do not estimate
the functional $x=\one\{k^*(1)=k^*(2)\}$ significantly better than a random guess when $\Delta^2$ is
smaller than $\Dcomp^2$ by a polylogarithmic factor (Corollary~\ref{cor:LBpoly}). Conversely, for
$n\le K^2$, hierarchical clustering~\cite{VEMPALA2004,Even24} clusters once $\Delta^2$ is large
compared to $\sqrt d=\Dcomp^2$. For $n>K^2$ and $p$ not an integer, we deduce from
Theorem~\ref{thm:UB} that a polynomial of degree of order $\log K$ decides whether $k^*(1)=k^*(2)$
with high probability as soon as $\Delta^2$ exceeds $\max\{\Dcomp^2,1\}$ by a factor
$e^{c(\log K)^{1/2}}$ (Corollary~\ref{cor:UBpoly} and Section~\ref{sec:resUB}). This gives the
complete picture for low-degree polynomials when $p$ is not an integer (Figure~\ref{fig:phase}).

\paragraph{Comparison with information bounds and previous results for Gaussian centers.}
For $K\le d$, that is $p\le2$, the formula \eqref{eq:Dcomp} with $\lceil p\rceil$ replaced by $2$
is the bound $\min\{\sqrt d,\sqrt{dK^2/n}\}$ of \eqref{eq:DHD}, so the new regime is $K>d$.
At a separation $\Delta_*$ subpolynomial in $K$, our results imply that low-degree polynomials need
$n=K^2d^{\lceil p\rceil/2}=K^3 d^{(\lceil p\rceil-p)/2}$ observations, and that this
is enough when $p$ is not an integer. For $p>2$, this exceeds the sample size $dK^2$ below which
\cite{Even25} shows that low-degree polynomials fail, by the factor $d^{\lceil p\rceil/2-1}\ge\sqrt d$.
It also exceeds by the factor $Kd^{\lceil p\rceil/2-1}\log K$ the sample size
$dK/\log K$ at which the optimal separation $\Delta_{\mathrm{stat}}^2$, without any computational constraint,  becomes $\log K$.  Figure~\ref{fig:phase} draws the phase
diagram of recovery in our model. In it, the region where low-degree polynomials fail lies above the
statistical threshold, which is polylogarithmic once $n\ge dK$: it is a statistical-computational gap. 

\paragraph{Comparison with algebraic constructions of the centers.} Gaussian centers match no
moment exactly, unlike the algebraic constructions of \cite{carpentier2026low,carpentier2026small}, and
they are not always easier to cluster. For $K=2d$, the orthogonal construction of
\cite{carpentier2026low} requires $n\ge K^4$ at a polylogarithmic separation, up to polylogarithmic
factors, while Gaussian centers need $K^2d^{3/2}$, of order $K^{7/2}$. But for $K$ of order $d^m$, with
an integer $m\ge2$, $\Dcomp^2$ is polynomially larger than the threshold at which
\cite{carpentier2026small} conjectures that its BCH instances become easy, so that, if this conjecture
holds, exact moment matching is not what makes clustering hard (Section~\ref{sec:bch}). The algebraic
construction makes the first informative tensor of the BCH centers extremely spiked: it differs from its
average on a few large entries only, which are easy to find, while that of Gaussian centers differs
slightly from its average on all its entries, in a way that is much harder to detect. Symmetrizing
the Gaussian centers makes them harder still, and leads us to conjecture that $n\ge K^3d$ suffices at
a subpolynomial separation for every configuration of the centers (Section~\ref{sec:symcenters}).

\subsection{Techniques and Intuition}\label{sec:intuition}

 Let $G=(V,W,E)$ be a bipartite graph, with $E\subseteq V\times W$,
whose nodes in $V$, the \emph{individual nodes}, stand for observations and whose nodes in $W$, the
\emph{variable nodes}, stand for coordinates, and let $v_1\ne v_2$ be two distinguished nodes of $V$. A labeling of $G$ is a pair $\pi=(\pi_1,\pi_2)$ of
injective maps $\pi_1:V\to[n]$ and $\pi_2:W\to[d]$ with $\pi_1(v_1)=1$ and $\pi_1(v_2)=2$: it places
each $v\in V$ on the observation $\pi_1(v)$ and each $w\in W$ on the coordinate $\pi_2(w)$. Suppose
in this section that every $w\in W$ has degree $2$, so that $w$ may be seen as an edge between its two
neighbors, and $G$ as a graph on $V$ whose edges carry distinct coordinates. Define the invariant monomial of $G$ by
\begin{equation}\label{eq:TGintro}
  T_G(Y)\ :=\ \sum_{\pi}\ T_{G,\pi}(Y),\qquad T_{G,\pi}(Y)\ :=\ \prod_{(v,w)\in E}Y_{\pi_1(v),\pi_2(w)} ,
\end{equation}
where the sum runs over the labelings $\pi$ of $G$. It is invariant under the permutations of the
observations other than $1,2$ and of the coordinates, and it counts the copies of $G$ in the data, each
weighted by the product of its entries. Up to normalization, $T_G$ is the polynomial $P_G$ of
Definition~\ref{def:PG}. Conditionally on $k^*$, the columns of $Y$ are independent and
$\E[Y_{ij}Y_{i'j}\mid k^*]=\frac{\Delta^2}{d}\one\{k^*(i)=k^*(i')\}$ for $i\ne i'$. So, when $G$ is
connected, the term of a labeling $\pi$ has conditional mean $(\Delta^2/d)^{|W|}$ if all the
observations $\pi_1(V)$ lie in the same group, and $0$ otherwise. Hence $\E[T_G\mid x=0]=0$, and $\E[T_G\mid x=1]$ is $(\Delta^2/d)^{|W|}K^{-(|V|-2)}$ times the number of
labelings, which is close to $n^{|V|-2}d^{|W|}$. As in \cite{Hopkins17,carpentier2025phaseII}, our upper bound amounts to choosing
a suitable graph $G$ and to showing, by a second moment method, that the fluctuations of $T_G$ are
small compared to this gap $\E[T_G\mid x=1]-\E[T_G\mid x=0]$, so that comparing $T_G$ with its two conditional means decides $x$.

Conversely, the law of $(x,Y)$ is invariant under permutations of the observations $3,\dots,n$,
under permutations of the coordinates and under sign changes of the coordinates. The polynomial
achieving the supremum in \eqref{eq:corrintro} inherits these invariances, and is therefore a
linear combination of $1$ and of invariant monomials $T_G$, where $G$ is now a bipartite multigraph with at
most $D$ edges whose variable nodes have even degree and where, in \eqref{eq:TGintro}, the product of the entries
$Y_{\pi_1(v),\pi_2(w)}$ over $m$ parallel edges between $v$ and $w$ is replaced by a Hermite
polynomial of degree $m$ of this entry (Proposition~\ref{prop:sym}). Our lower bound therefore
amounts to showing that every such linear combination is poorly correlated with $x$.

\begin{figure}[ht]
\centering
\begin{tikzpicture}[x=1cm,y=1cm,
  di/.style={circle,draw=black,fill=white,inner sep=1pt,font=\scriptsize},
  ob/.style={circle,fill=black,inner sep=1.4pt},
  lab/.style={font=\scriptsize,inner sep=1.5pt},
  co/.style={font=\tiny,inner sep=1pt},
  cp/.style={NavyBlue,thick},
  cq/.style={Plum,thick,densely dashed}]
\begin{scope}
  \node[di] (a1) at (0,0) {$1$};
  \node[ob,label={[lab]above:$3$}] (aa) at (1.2,0) {};
  \node[ob,label={[lab]above:$4$}] (ab) at (2.4,0) {};
  \node[di] (a2) at (3.6,0) {$2$};
  \draw[cp] (a1) to[bend left=28] node[co,above]{$j_1$} (aa);
  \draw[cp] (aa) to[bend left=28] node[co,above]{$j_2$} (ab);
  \draw[cp] (ab) to[bend left=28] node[co,above]{$j_3$} (a2);
  \draw[cq] (a1) to[bend right=28] node[co,below]{$j_1$} (aa);
  \draw[cq] (aa) to[bend right=28] node[co,below]{$j_2$} (ab);
  \draw[cq] (ab) to[bend right=28] node[co,below]{$j_3$} (a2);
  \node[font=\footnotesize] at (1.8,-1.4) {$\mathcal V_{\mathrm{noise}}$: $\pi'=\pi$};
\end{scope}
\begin{scope}[xshift=5.4cm]
  \begin{scope}[on background layer]
    \fill[OliveGreen!15,rounded corners=6pt] (0.9,-0.95) rectangle (1.5,0.95);
    \fill[Orange!20,rounded corners=6pt] (2.1,-0.95) rectangle (2.7,0.95);
  \end{scope}
  \node[di] (b1) at (0,0) {$1$};
  \node[ob,label={[lab]above:$3$}] (ba) at (1.2,0.5) {};
  \node[ob,label={[lab]above:$4$}] (bb) at (2.4,0.5) {};
  \node[ob,label={[lab]below:$5$}] (bc) at (1.2,-0.5) {};
  \node[ob,label={[lab]below:$6$}] (bd) at (2.4,-0.5) {};
  \node[di] (b2) at (3.6,0) {$2$};
  \draw[cp] (b1) -- node[co,above]{$j_1$} (ba) -- node[co,above]{$j_2$} (bb) -- node[co,above]{$j_3$} (b2);
  \draw[cq] (b1) -- node[co,below]{$j_1$} (bc) -- node[co,below]{$j_2$} (bd) -- node[co,below]{$j_3$} (b2);
  \node[font=\footnotesize] at (1.8,-1.4) {$\mathcal V_{\mathrm{centers}}$: $\pi'_2=\pi_2$, same groups};
\end{scope}
\end{tikzpicture}
\caption{The pairs of labelings behind $\mathcal V_{\mathrm{noise}}$ and $\mathcal V_{\mathrm{centers}}$, for the path
$v_1,v_3,v_4,v_2$. The copy placed by $\pi$ is solid, the copy placed by $\pi'$ dashed. Left:
$\mathcal V_{\mathrm{noise}}$, the part of the variance due to the noise, where the two labelings coincide: the two
copies sit on the same observations and on the same coordinates. Right: $\mathcal V_{\mathrm{centers}}$, the part due to the random centers, where the two labelings place the links on the same coordinates $j_1,j_2,j_3$, and the nodes $v_3$ and $v_4$ on the distinct observations $3$ and
$5$, and $4$ and $6$, forced into the same groups (shaded), so that
each coordinate carries the product of the variances of the coordinates of two centers.}
\label{fig:twoterms}
\end{figure}
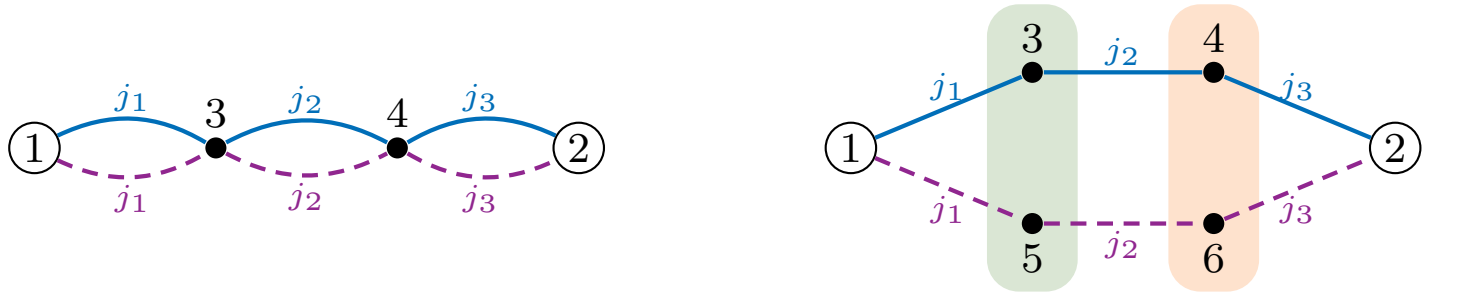

\paragraph{Where the critical order $\lceil p\rceil$ and the separation $\Dcomp$ come from.} Let $G$
be connected, with every $w\in W$ of degree $2$ and every $v\in V$ other than $v_1,v_2$ of degree
$q$, and neglect below the nodes of $W$ adjacent to $v_1$ or $v_2$, so that $|W|=q(|V|-2)/2$. As each
entry $Y_{ij}=\frac{\Delta}{\sqrt d}X_{k^*(i),j}+\epsilon_{ij}$ is the sum of a signal and a noise, the
variance of $T_G$ is made of several terms, one for each way of combining them: in a term, each
individual node other than $v_1,v_2$ contributes either through the noise of its observation or
through the center of its group. Conditionally on $x$, $\E[T_G^2\mid x]$ is the sum, over the pairs
$(\pi,\pi')$ of labelings, of the nonnegative numbers $\E[T_{G,\pi}T_{G,\pi'}\mid x]$
(Section~\ref{sec:moments}), and the pairs that share no observation other than $1,2$ and no
coordinate make up $\E[T_G\mid x]^2$, up to a factor close to $1$. So, up to a term small compared
with $\E[T_G\mid x=1]^2$, the variance of $T_G$ is the sum of the other terms, and two families of
pairs show where the threshold comes from (Figure~\ref{fig:twoterms}). They are the two extreme cases:
in the first, every individual node contributes through the noise of its observation, and, in the
second, through the center of its group, which gives the largest cross terms between the two
labelings. 

\noindent\uline{$\mathcal V_{\mathrm{noise}}$: the two labelings coincide.} The pairs with $\pi'=\pi$, in which
every node contributes through the noise of its observation, give the part of the variance due to the
noise. As each entry has variance close to $1$, it is close to the number of labelings, where,
here and below, $a\asymp b$ means that $a/b$ lies between two positive constants that do not depend on
$n$, $d$, $K$ and $\Delta$:
\[
  \mathcal V_{\mathrm{noise}}\ :=\ \sum_\pi\E\big[T_{G,\pi}^2\big]\ \asymp\ n^{|V|-2}d^{|W|}\ =\ \big(n\,d^{q/2}\big)^{|V|-2} .
\]

\noindent\uline{$\mathcal V_{\mathrm{centers}}$: same coordinates, distinct observations in the same groups.} The pairs that
place the variable nodes on the same coordinates, $\pi'_2=\pi_2$, and every individual node other than
$v_1,v_2$ on two distinct observations give the part of the variance due to the random centers: each
such node then contributes through the center of its group, and only $v_1$ and $v_2$ tie the term to
the observations $1$ and $2$. These pairs contribute because they force the two observations of each
node into the same group, which happens with
probability of order $K^{-(|V|-2)}$, for about $n^{2(|V|-2)}d^{|W|}$ pairs. Each coordinate
$\pi_2(w)$, $w=\{v,v'\}$, then contributes the product $(\Delta^2/d)^2$ of the variances of the
coordinates of the two centers it involves, which gives
\[
  \mathcal V_{\mathrm{centers}}\ :=\ n^{2(|V|-2)}d^{|W|}\,K^{-(|V|-2)}\Big(\frac{\Delta^2}d\Big)^{2|W|}\ =\ \Big(\frac{n^2\Delta^{2q}}{K\,d^{q/2}}\Big)^{|V|-2} .
\]
Comparing the squared gap $\E[T_G\mid x=1]^2\asymp(n^2\Delta^{2q}/K^2)^{|V|-2}$ with $\mathcal V_{\mathrm{noise}}$ and with $\mathcal V_{\mathrm{centers}}$, we see that
$T_G$ decides $x$ with high probability only if
\begin{align}\label{eq:condition:Deltacomp}
  n\Delta^{2q}\ \ge\ K^2d^{q/2}\qquad\text{and}\qquad d^{q/2}\ \ge\ K .
\end{align}
Both comparisons are made node by node: in the ratio of the squared gap to $\mathcal
V_{\mathrm{noise}}$, each individual node contributes the factor $n\Delta^{2q}/(K^2d^{q/2})$, and, in its ratio
to $\mathcal V_{\mathrm{centers}}$, the factor $d^{q/2}/K$. For nodes of high degree, $d^{q/2}\ge K$, and only the
comparison with $\mathcal V_{\mathrm{noise}}$ matters: the invariant monomial succeeds once $n$ is large enough. For
nodes of small degree, $d^{q/2}<K$, and $\mathcal V_{\mathrm{centers}}$ exceeds the squared gap whatever $n$: the
fluctuations due to the random centers of their groups exceed the signal, and polynomials containing
such nodes are useless. This is the
difficulty of the regime $K>d$, and the reason why spectral procedures, which only use second
moments, that is nodes of degree $2$, need $\Delta^2$ of order $\sqrt d$ when $K>d$, whatever $n$. 

As we are interested in the regime where $\Delta^2<\sqrt d$, the first condition in~\eqref{eq:condition:Deltacomp} is the weakest at the
smallest $q$, and the second one asks that the integer $q$ be at least $p$. The best choice is
$q=\lceil p\rceil$, which gives the threshold $\Dcomp$.

\paragraph{Techniques.} Both proofs follow this heuristic. For the upper bound, we use the second
moment method on the invariant monomial of a template graph, as for stochastic block models
\cite{Hopkins17,carpentier2025phaseII}. The template is made of a $\lceil p\rceil$-regular graph
with expansion, called the body, and of $2L$ nodes of the body far apart from one another, joined to
$v_1$ or $v_2$, with $L$ of order $(\log K)^{1/2}$. Its isoperimetric properties ensure that no pair
of labelings other than those of $\mathcal V_{\mathrm{noise}}$ and $\mathcal V_{\mathrm{centers}}$ makes the variance large. The body
of \cite{carpentier2025phaseII}, a blown-up cycle, does not have the required expansion, and we use
a cycle together with random perfect matchings. For the lower bound, we follow the approach of
\cite{CGGV25,carpentier2025phase,carpentier2026low,carpentier2026small}, which compares the
polynomials with a family whose Gram matrix is controlled, with two differences. First, as
anticipated in \cite[Section~4]{carpentier2026small} for Gaussian mixtures, the polynomials $T_G$ are
built from Hermite polynomials of the entries $Y_{ij}$. Second, instead of showing that these
polynomials are almost orthonormal, we project them onto an exactly orthonormal family of
polynomials of the latent variables, in the spirit of~\cite{SohnWein25}. In this family, an
individual node of degree smaller than $p$ is replaced by a center, and the others
keep their noise: this is the comparison with $\mathcal V_{\mathrm{centers}}$ and $\mathcal V_{\mathrm{noise}}$ above, made node by
node. After a recursive centering, which plays the role of the centering by connected components of
\cite{montanari2024equivalence,CGGV25,carpentier2025phase,carpentier2026low,carpentier2026small},
the matrix of these projections is diagonally dominant once rescaled, and this gives the bound.
Section~\ref{sec:ideas} explains both proofs in more detail.

\subsection{Other related works}\label{sec:related}

\paragraph{Proof techniques for low-degree lower bounds.} In detection problems,
proving a low-degree lower bound reduces to controlling the norm of the projection of the
likelihood ratio onto low-degree polynomials, and this ratio has an explicit form in many models
\cite{Hopkins17,KuniskyWeinBandeira}. The low-degree correlation \eqref{eq:corrintro} admits no
such closed form. Schramm and Wein \cite{SchrammWein22} control it through the joint cumulants of
$x$ and of the entries of $Y$. Their method is flexible and has been used for submatrix estimation
\cite{SchrammWein22}, graphon estimation \cite{luo2023computational} and the clustering of
Gaussian mixtures \cite{Even24}, among other problems. However, the computations can become heavy,
and the resulting thresholds are sometimes tight only up to polylogarithmic factors. In latent
variable models, conditioning on the latent variables makes the method simpler and sharper
\cite{Even25}, and Sohn and Wein \cite{SohnWein25} extended it so as to recover sharp thresholds,
for instance the BBP transition for the spiked Wigner model and the Kesten--Stigum threshold for
the stochastic block model. Another approach, introduced by Montanari and Wein
\cite{montanari2024equivalence} for rank-one matrix estimation and developed in general in
\cite{CGGV25}, constructs polynomials that are nearly orthonormal under the law of the data. Then
$\Corr_{\le D}^2$ is at most, up to a factor close to one, the sum over this family of the squared
correlations with $x$. This approach has been used for stochastic block models
\cite{carpentier2025phase}, clustering \cite{carpentier2026low,carpentier2026small} and ranking
\cite{carpentier2025rank}. Our proof is a variant of it, in which the orthonormal family is made of
polynomials of the latent variables. Last, Li \cite{li2026algorithmic} has recently suggested an
argument that depends little on the model: a cross-validation step converts any efficient recovery
procedure into a test, and low-degree bounds for that test then conclude. The lower bounds obtained
this way, however, rely on the low-degree conjecture for detection, while the approaches above,
ours included, bound $\Corr_{\le D}$ without any conjecture.

\paragraph{Tensor decompositions.} Given the centers, the moment tensor of order $q$ of an
observation involves $\frac1K\sum_{k\le K}X_k^{\otimes q}$, a symmetric tensor of order $q$ in
dimension $d$ with $K$ random components. For random tensors of order $3$ with $r$ components in
dimension $d$, polynomial-time algorithms recover the components for $r$ up to nearly $d^{3/2}$
\cite{ma2016polynomial}, and Wein \cite{wein2022average} showed that low-degree polynomials succeed when $r$ is much
smaller than $d^{3/2}$ and fail when it is much larger, and indicated the threshold $d^{q/2}$ for the number of components for a  tensor of order $q$. The condition $d^{q/2}\ge K$ of Section~\ref{sec:contribution} is the analogue,
for the tensor of order $q$ of the centers, of the condition that $r$ be smaller than $d^{q/2}$, and
the critical order is the smallest order at which it holds. 

Decompositions of moment tensors of a fixed order could be directly applied to
learn mixtures with many more components than dimensions. For generic or smoothed components, a
setting that includes Gaussian centers, tensors of order $2\ell$ can be decomposed in polynomial
time up to a rank $d^{\ell}$ \cite{bhaskara2019smoothed}, see also
\cite{hopkins2019robust} for $\ell=2$, in line with the condition $d^{q/2}\ge K$. To our
knowledge, however, the available guarantees for mixtures remain far from the threshold of
Section~\ref{sec:contribution}: they require either an inverse-polynomial error on the moment
tensors, so that the sample size is a large polynomial in $d$ and $K$
\cite{bhaskara2014smoothed,bhaskara2019smoothed}, or, when $K$ is much larger than $d$, an
initialization close to the centers \cite{anandkumar2015learning}.

\paragraph{Anisotropic Gaussian mixtures.}
Statistical query lower bounds \cite{SQclustering,pmlr-v195-diakonikolas23b}
also point to statistical-computational gaps in clustering, but only for non-isotropic mixtures,
where the barrier is far stronger: when the covariance matrix is unknown and only assumed to be
bounded by $I_d$, every SQ procedure has complexity of order $d^{c/\gamma}$, for some constant
$c>0$, at separation $\Delta_*=K^\gamma$ \cite{pmlr-v195-diakonikolas23b}. However, non-isotropic mixtures are of a very different nature.

\subsection{Organization of the manuscript}

Section~\ref{sec:setting} states the model and the two theorems, with their corollaries at degree of
order $\log K$. Section~\ref{sec:discussion} compares the results with those of
\cite{carpentier2026small} for the BCH instances, and discusses integer values of $p$, the running
time, the recovery of the whole partition and symmetrized Gaussian centers.
Section~\ref{sec:ideas} explains the ideas of the two proofs: the first two moments of the
invariant monomials, then the lower bound and the upper bound. Section~\ref{sec:graphs} introduces the
graph notation of the proofs. Section~\ref{part:lb} proves Theorem~\ref{thm:LB},
Section~\ref{part:ub} proves Theorem~\ref{thm:UB}, and Section~\ref{ss:tech1} proves the
corollaries. The appendix builds the body of the template, collects technical lemmas, and proves the
recovery of the partition and the lower bound for symmetrized centers.

\paragraph{Notation.} For a positive integer $b$, $[b]:=\{1,\dots,b\}$, and $|A|$ is the
cardinality of a finite set $A$. For a real $p$, $\lceil p\rceil$ is the smallest integer at least
$p$ and $\lfloor p\rfloor$ the largest integer at most $p$. The logarithm $\log$ is natural.
$\|\cdot\|$ is the Euclidean norm of a vector of $\mathbb R^d$ and, for a square-integrable random
variable $P$, $\|P\|:=\E[P^2]^{1/2}$. In the heuristic passages only, for positive $a$ and $b$,
$a\lesssim b$ (or $b\gtrsim a$) means that $a\le Cb$ for a constant $C$
that does not depend on $n$, $d$, $K$ and $\Delta$, and $a\asymp b$
means $a\lesssim b\lesssim a$.

\section{Model and main results}\label{sec:setting}

\subsection{The model}\label{sec:model}

Let $K,d,n\in\N^*$ and $\Delta>0$. We consider a clustering problem of $n$ observations in dimension
$d$. In the proofs, the observations are called \textbf{individuals}, indexed by $i\in[n]$, and the
coordinates are called \textbf{variables}, indexed by $j\in[d]$. Draw the
group labels, or \textbf{colors} (the name used in Section~\ref{part:lb}), $k^*(1),\dots,k^*(n)$ of the individuals i.i.d.\ uniform on
$[K]$, the center coordinates $X_{k,j}$ ($k\in[K]$, $j\in[d]$) i.i.d.\ $\mathcal N(0,1)$, and the
noise $\epsilon_{i,j}$ i.i.d.\ $\mathcal N(0,1)$, independent of everything, and set
\begin{equation}\label{eq:model}
Y_{i,j}\;=\;\frac{\Delta}{\sqrt d}\,X_{k^*(i),j}+\epsilon_{i,j},\qquad
\Var(Y_{ij})=1+\frac{\Delta^2}d .
\end{equation}
The center of group $k$ is $\tfrac{\Delta}{\sqrt d}(X_{k,1},\dots,X_{k,d})$, of mean squared norm
$\Delta^2$. Given the centers, the rows $Y_{i,\cdot}$ are i.i.d.\ from the mixture of the $K$
distributions $\mathcal N(\tfrac{\Delta}{\sqrt d}X_{k,\cdot},I_d)$ with equal weights, whose
separation is $\sqrt2\,\Delta$ up to a factor close to $1$ (Section~\ref{sec:contribution}).
Throughout, $p:=2\log K/\log d$, so that $K=d^{p/2}$, and $\lceil p\rceil$ is the critical order of
Section~\ref{sec:contribution}.

Write $i\sim i'$ for $k^*(i)=k^*(i')$, and
\begin{equation}\label{eq:x}
x\;:=\;\one\{k^*(1)=k^*(2)\}
\end{equation}
for the indicator that the first two individuals lie in the same group. Conditionally on $k^*$, the
entries in distinct coordinates are independent, and for each coordinate $j$,
\begin{equation}\label{eq:cov}
\E[Y_{i,j}Y_{i',j}\mid k^*]=\frac{\Delta^2}{d}\one\{i\sim i'\}\quad(i\ne i'),\qquad
\E[Y_{i,j}^2\mid k^*]=1+\frac{\Delta^2}{d} .
\end{equation}
We bound, for the polynomials $P$ of degree at most $D$ in the $nd$ variables $(Y_{ij})$, the
correlation \eqref{eq:corrintro} with $x$, restated here:
\begin{equation}\label{eq:corr}
\Corr_{\le D}\;=\;\sup_{\deg P\le D,\ P\ne0}
\frac{\E\big[(x-\E x)\,P(Y)\big]}{\sqrt{\E[P^2]}}.
\end{equation}
All expectations are under the joint law of $(k^*,X,\epsilon)$, and throughout
$\langle P,P'\rangle:=\E[PP']$ and $\|P\|:=\sqrt{\E[P^2]}$. As $\E x=1/K$ and
$\|x-\E x\|^2=K^{-1}(1-K^{-1})$, the correlation is always at most $K^{-1/2}$.
We recall the classical identity\footnote{In \cite{SchrammWein22}, the right-hand side is written as $\E[x^2]$ minus the square of the
uncentered correlation, the latter being equal to $\big((\E x)^2+\Corr_{\le D}^2\big)^{1/2}$.}~\cite{SchrammWein22}, which relates this
correlation to the mean squared error of low-degree estimators of $x$: since the constants are
polynomials of degree at most $D$,
\begin{equation}\label{eq:MMSE}
\inf_{\deg P\le D}\E\big[(P(Y)-x)^2\big]\;=\;\Var(x)-\Corr_{\le D}^2= \frac{1}{K}\left(1-\frac{1}{K}\right) -  \Corr_{\le D}^2.
\end{equation}

\subsection{The lower bound}\label{sec:resLB}

The lower bound holds under the following assumption. The
integer $D$ is the degree of the polynomials.

\begin{assumption}\label{ass:A}
$K=d^{p/2}$ with $p\ge2$, $D\ge2$ an integer, $n\ge16D^2$, and
\[
\text{(i) } \frac{\Delta}{d^{1/4}}\le D^{-200},\qquad
\text{(ii) } \frac{\Delta^{2\lceil p\rceil}\,n}{K^{2}d^{\,\lceil p\rceil/2}}\le D^{-400\lceil p\rceil},\qquad
\text{(iii) } d^{\,p-\lceil p\rceil+1}\ \ge\ D^{1600\lceil p\rceil} .
\]
\end{assumption}

\begin{theorem}\label{thm:LB}
Under Assumption~\ref{ass:A}, $\ \Corr_{\le D}\le4/K$.
\end{theorem}

The proof is given in Section~\ref{sec:solve}.
By \eqref{eq:MMSE}, Theorem~\ref{thm:LB} states that every polynomial $P$ of degree at most $D$
has mean squared error
\[
\E\big[(P(Y)-x)^2\big]\ \ge\ \frac1K\Big(1-\frac{17}{K}\Big),
\]
whereas the constant estimator $\E x=1/K$, which ignores the data, has mean squared error
$\Var(x)=\frac1K\big(1-\frac1K\big)$. Hence low-degree polynomials reduce the error of this trivial
estimator by a fraction at most $16/(K-1)$.
In Assumption~\ref{ass:A}, condition \textup{(ii)} is the main one, and it
reads $\Delta^2\le D^{-400}\Dcomp^2$ when $n\ge K^2$. Condition \textup{(i)} keeps $\Delta^2$ below
the separation $\sqrt d$ at which distance-based procedures succeed, and when $n\le K^2$ it implies
\textup{(ii)}, as the left-hand side of \textup{(ii)} is then at most $(\Delta^2/\sqrt d)^{\lceil
p\rceil}$. Condition \textup{(iii)} reads $K\ge d^{(\lceil p\rceil-1)/2}D^{800\lceil p\rceil}$: the
order $q=\lceil p\rceil-1$ fails the condition $d^{q/2}\ge K$ by the factor $D^{800\lceil p\rceil}$ at
least. It excludes the values of $p$ slightly above an integer, and it bounds $D$ by a small power of
$d$. At degree of order $\log K$, which is of order $\log n$ in the relevant regimes, the assumption
takes the following form, in which the two displayed conditions say together that $\Delta^2\le(C_{\mathrm{poly}}\log K)^{-400}\Dcomp^2$.

\begin{corollary}\label{cor:LBpoly}
Fix a real $p\ge2$ and a constant $C_{\mathrm{poly}}>0$, and put
$D:=\lfloor C_{\mathrm{poly}}\log K\rfloor$. For all $K$ large enough in terms of $p$ and
$C_{\mathrm{poly}}$: whenever $d=K^{2/p}$, $n\ge(\log K)^{3}$,
\[
  \frac{\Delta^{2\lceil p\rceil}\,n}{K^{2}d^{\,\lceil p\rceil/2}}
  \ \le\ \big(C_{\mathrm{poly}}\log K\big)^{-400\lceil p\rceil}
  \qquad\text{and}\qquad
  \Delta\ \le\ d^{1/4}\big(C_{\mathrm{poly}}\log K\big)^{-200} ,
\]
one has $\Corr_{\le D}\le4/K$.
\end{corollary}
The proof is given in Section~\ref{sss:LBpoly}.

\subsection{The upper bound}\label{sec:resUB}

The upper bound counts the copies, in the data, of a template graph built in
Section~\ref{sec:construction}. It holds under the following assumption, in which $L$ sets the size
of the template and the integer $q$ its degree. Condition \textup{(i)} is the condition $d^{q/2}\ge
K$ of Section~\ref{sec:intuition}, up to the factor $K^{3201q/L}$.

\begin{assumption}\label{ass:B}
$q\ge3$ an integer, $K\ge d$, $K^2<n\le K^2d^{q/2}$, $L\ge(400q)^{600q^{2}}$ an integer, and
\[
\begin{gathered}
\text{(i) } d^{\,q/2}\ \ge\ K^{1+3201q/L},\qquad
\text{(ii) } \frac{\Delta^{2q}\,n}{K^{2}d^{\,q/2}}\ \ge\ K^{5q/L},
\qquad
\text{(iii) } \Delta^{2}<\sqrt d .
\end{gathered}
\]
\end{assumption}

\begin{theorem}\label{thm:UB}
Under Assumption~\ref{ass:B} there is one polynomial $P$ in the $nd$ variables $(Y_{ij})$, of degree at most
$2(q+2)L^{2}$, such that the test
\[
  \hat x\ :=\ \one\Big\{\big|P(Y)-\E[P\mid x=1]\big|<\big|P(Y)-\E[P\mid x=0]\big|\Big\} ,
\]
which declares $x=1$ when $P(Y)$ is closer to $\E[P\mid x=1]$ than to $\E[P\mid x=0]$, satisfies,
for $\xi\in\{0,1\}$,
\[
  \Pr\big[\hat x\ne x\ \big|\ x=\xi\big]\ \le\ 4608q^{2}L^{4}K^{-1/q}\ +\ 7776q^{4}L^{8}K^{-2/L} .
\]
\end{theorem}

We are interested in the regime where $q$ is a constant and $L$ is large, but at most of order
$(\log K)^{1/2}$, so that the degree of $P$ is at most of order $\log K$. Then the error probability
of Theorem~\ref{thm:UB} is small. The constants of Assumption~\ref{ass:B} and
Theorem~\ref{thm:UB} are explicit but not optimized. The condition $K^2<n$ of Assumption~\ref{ass:B} is not a
restriction of its own: it follows from \textup{(ii)} and \textup{(iii)}, as $\Delta^{2q}<d^{q/2}$.
We obtain the following corollary by taking $L$ of order $(\log K)^{1/2}$.
\begin{corollary}\label{cor:UBpoly}
Fix an integer $q\ge3$ and a constant
$C_{\mathrm{poly}}>0$, and put
$D:=\lfloor C_{\mathrm{poly}}\log K\rfloor$. There is a $c\ge1$, depending on $q$ and
$C_{\mathrm{poly}}$ alone, such that for all $K$ large enough in terms of $q$ and $C_{\mathrm{poly}}$: whenever
$K\ge d$, $K^2<n\le K^2d^{q/2}$,
\[
  d\ \ge\ K^{2/q}e^{\,c\,(\log K)^{1/2}}
\]
and
\[
  \frac{\Delta^{2q}\,n}{K^{2}d^{\,q/2}}\ \ge\ e^{\,c\,(\log K)^{1/2}}
  \qquad\text{and}\qquad \Delta^{2}<\sqrt d ,
\]
the test $\hat x$ of Theorem~\ref{thm:UB}, built on a polynomial of degree at most $D$ in the $nd$
data, satisfies $\Pr[\hat x\ne x\mid x=\xi]\le e^{-(\log K)^{1/2}/c}$ for $\xi\in\{0,1\}$.
\end{corollary}
The proof is given in Section~\ref{sss:UBpoly}.
To read this corollary, let $K=d^{p/2}$ with $p$ fixed and not an integer, and take
$q=\lceil p\rceil$. Then the condition $d\ge K^{2/q}e^{c(\log K)^{1/2}}$ holds for $K$ large enough,
and, when $n>K^2$, the signal condition reads $\Delta^2\ge\Dcomp^2e^{c(\log K)^{1/2}/q}$: it matches
$\Dcomp^2$ up to a subpolynomial factor. The condition $\Delta^2<\sqrt d$ is no restriction: when
$\Delta^2$ is large compared to $\sqrt d$, the squared distance $\|Y_1-Y_2\|^2$ decides $x$ with a
small error.

With $q=\lceil p\rceil$, the range of $n$ in Assumption~\ref{ass:B} is $K^2<n\le
K^2d^{\lceil p\rceil/2}$, a mild restriction. For $n\le K^2$, $\Dcomp^2=\sqrt d$, and the squared
distance $\|Y_1-Y_2\|^2$ decides $x$ once $\Delta^2$ is large compared to $\sqrt d$. For
$n>K^2d^{\lceil p\rceil/2}$, the test can be computed on a subsample of size
$K^2d^{\lceil p\rceil/2}$, and Theorem~\ref{thm:UB} then only requires $\Delta^{2\lceil p\rceil}\ge
K^{5\lceil p\rceil/L}$, which, for $L$ of order $(\log K)^{1/2}$, is the statistical threshold
$\Delta^2_{\mathrm{stat}}$ up to a subpolynomial factor. In summary, low-degree polynomials decide
$x$ once $\Delta^2$ exceeds $\max\{\Dcomp^2,\Delta^2_{\mathrm{stat}}\}$ by a subpolynomial factor,
for every $n$, when $p$ is not an integer. Integer values of $p$ are discussed in Section~\ref{sec:integer}.

\section{Discussion}\label{sec:discussion}

\subsection{Gaussian centers versus BCH centers}\label{sec:bch}

In \cite{carpentier2026small}, $K=(d+1)^m$ with $m$ a positive integer, the observations are
Bernoulli, and the $2K$ centers $\pm X_k$ are, up to a random permutation and random signs of the
coordinates, the codewords of a dual BCH code rescaled by $\Delta/\sqrt d$. This construction
ensures that the moment tensors of the observations coincide with those of a Rademacher vector up
to the order $2m+1$. Low-degree polynomials fail as long as $\Delta^2$ is small compared to
$\min\{(K^2d^{3/2}/n)^{1/(2m+2)},\sqrt d\}$ \cite[Theorem~1]{carpentier2026small}. Conversely, the
entries of the moment tensor of order $2m+2$ can be estimated one by one once $n\gtrsim
K^2d^2/\Delta^{4m+4}$, and \cite[Section~4]{carpentier2026small} conjectures that the clusters can
then be recovered in polynomial time if $\Delta^2\gtrsim m\log d$.

\paragraph{Comparison at $K\asymp d^m$.} Let $m\ge2$, and let $K$ be larger than $d^m$ by a large
enough polylogarithmic factor, so that the critical order is $\lceil p\rceil=2m+1$ and
Assumption~\ref{ass:A}\textup{(iii)} holds (for $K=d^m$ exactly, see Section~\ref{sec:integer}). Up to subpolynomial
factors, the thresholds are as follows.
\begin{center}\small
\begin{tabular}{@{}lcc@{}}
\toprule
& threshold on $\Delta^2$ ($n\ge K^2$) & threshold on $n$ (polylog.\ $\Delta^2$)\\
\midrule
Gaussian centers (Thms.~\ref{thm:LB}, \ref{thm:UB}) & $(K^2d^{m+1/2}/n)^{1/(2m+1)}$ & $K^2d^{m+1/2}$\\
BCH, lower bound & $(K^2d^{3/2}/n)^{1/(2m+2)}$ & $K^2d^{3/2}$\\
BCH, conjectured upper bound & $\max\{(K^2d^{2}/n)^{1/(2m+2)},\,m\log d\}$ & $K^2d^{2}$\\
\bottomrule
\end{tabular}
\end{center}
The first row is $\Dcomp^2$. On the whole range $K^2\le n<K^2d^{m+1/2}$, it is larger than the
conjectured threshold of the BCH centers (by the factor $d^{(m-1)/(2m+2)}$ at $n=K^2$), and, at a
polylogarithmic separation, Gaussian centers require $d^{m-3/2}$ times more observations than the
conjectured sample size $K^2d^2$ of the BCH centers. So, if the
conjecture of \cite{carpentier2026small} holds, Gaussian centers are harder to cluster than BCH
centers, although they match no moment exactly.

\paragraph{The role of the moment tensors.} Moment matching alone would suggest the opposite. The
BCH centers match exactly the tensors up to the order $2m+1$, so that the first informative tensor
has order $2m+2$. Gaussian centers match no moment exactly but, for $K\asymp d^m$, graphs whose
individual nodes have degree at most $2m$ are irrelevant: a template of degree $q$ is informative only
if $d^{q/2}\ge K$ (Section~\ref{sec:intuition}). In effect,
the first informative tensor of Gaussian centers has order $2m+1$, one less than for the BCH centers.

The difference lies in the structure of this first informative tensor. Given the centers, the
moment tensor of order $q$ of an observation has, at distinct indices $j_1,\dots,j_q$, the entry
$(\Delta/\sqrt d)^q\frac1K\sum_{k\le K}X_{k,j_1}\cdots X_{k,j_q}$. For the BCH centers and
$q=2m+2$, this entry is $\pm(\Delta/\sqrt d)^{2m+2}$ on the tuples of coordinates that
support a parity check of the code, a fraction of order $1/K$ of all tuples, and $0$ elsewhere. This sparse
tensor with large entries is easy to detect: an empirical average estimates each entry once
$n(\Delta^2/d)^{2m+2}$ exceeds a logarithmic factor, that is, $n\gtrsim K^2d^2/\Delta^{4m+4}$. Its
support and its signs reflect the hidden permutation and signs of the code, that is, the structure
of the whole configuration of the centers. In both cases, the squared entries of the tensor of
order $q$ sum to about $\Delta^{2q}/K$, but, for Gaussian centers and $q=2m+1$, this total is spread
evenly over all the tuples: every entry is of order $(\Delta/\sqrt d)^qK^{-1/2}$. No entry can be estimated on its own before $n\gtrsim
K(d/\Delta^2)^{2m+1}\asymp K^2d^{m+1}/\Delta^{4m+2}$, and aggregating the entries along a template
only gains a factor $\sqrt d$, which leads to $\Dcomp$. In short, the BCH construction gains one
order of exact moment matching, but at the price of a tensor of order $2m+2$ with a rigid structure
that is easy to detect and reveals the centers. Gaussian centers have no such weakness, which
explains why $\Dcomp^2$ is larger than the conjectured threshold of the BCH centers. A
configuration of centers that would match exactly the moment tensors of $\mathcal N(0,I_d)$ up to the
order $2m+1$, as the Gaussian designs discussed in \cite[Section~4]{carpentier2026small}, while
having a tensor of order $2m+2$ without such a rigid structure, could be harder to cluster than
Gaussian centers. Whether such configurations exist is open.

\subsection{Integer values of $p$}\label{sec:integer}

The threshold $\Dcomp^2$ is discontinuous at the integer values of $p$. Fix $K$ and $n$ with
$K^2<n\le K^3$, and an integer $q\ge2$ (at $d=K^{2/q}$, $K^3=K^2d^{q/2}$, the largest sample size
allowed by Assumption~\ref{ass:B}). As $d$ decreases through $K^{2/q}$, that is, as $p$ increases
through $q$, the factor $\sqrt d$ in \eqref{eq:Dcomp} varies continuously, but the critical order
jumps from $q$ to $q+1$. Hence $\Dcomp^2$ jumps from $\sqrt d\,(K^2/n)^{1/q}$ to
$\sqrt d\,(K^2/n)^{1/(q+1)}$, that is, it is multiplied by $(n/K^2)^{1/(q(q+1))}$.

This jump comes from the condition $d^{q/2}\ge K$ of Section~\ref{sec:intuition}: a template of
degree $q$ is informative only if $d^{q/2}\ge K$, and below $d=K^{2/q}$ one has to rely on a
template of degree $q+1$, which requires a larger separation. At $d=K^{2/q}$ itself, the variance
term $\mathcal V_{\mathrm{centers}}$ is no longer small compared with the squared gap, so that the template of degree
$q$ is not expected to succeed.

At $d=K^{2/q}$ exactly, that is at $p=q$, the critical order is still $\lceil p\rceil=q$. Hence
Theorem~\ref{thm:LB} applies:
the polynomials of degree at most $D$ fail as long as $\Delta^2\le D^{-400}\Dcomp^2$, where
$\Dcomp^2=\sqrt d\,(K^2/n)^{1/q}$ for $n\geq K^2$. In contrast,
Corollary~\ref{cor:UBpoly} requires the choice $q+1$, as its condition $d\ge K^{2/q}e^{c(\log
K)^{1/2}}$ fails for $q$. With $q+1$ in place of $q$, its conditions on $d$ and $n$ hold for $K$
large enough, and its signal condition reads
\[
  \Delta^2\ \ge\ \sqrt d\,\Big(\frac{K^2}{n}\Big)^{1/(q+1)}e^{\,c(\log K)^{1/2}/(q+1)} .
\]
At $d=K^{2/q}$, our lower and
upper bounds therefore differ by the factor $(n/K^2)^{1/(q(q+1))}$ of the jump.
Our results thus leave open a subpolynomial window of values
of $d$ around $K^{2/q}$, in which the threshold is only known to lie between the two one-sided
limits.

\subsection{Running time}\label{sec:algo}

The polynomial of Theorem~\ref{thm:UB} is a sum over the labelings of a template with at most
$3q L^2$ individual nodes and at most $3q L^2$ variable nodes, so evaluating
it term by term takes a number of operations of order $(nd)^{3q L^2}$, quasi-polynomial
when $L^2$ is of order $\log K$. This leaves open whether polynomial-time algorithms cluster
at the separation $\Dcomp$, and whether tensor-based procedures
such as \cite{LiuLi2022} cluster at a polylogarithmic separation with a sample size of order
$K^2d^{\lceil p\rceil/2}$. Simple local algorithms are other candidates: for tensor PCA, the recovery of one unknown vector
from a noisy observation of one of its tensor powers of order at least $3$, Das and Schramm recently
showed that local Markov chains, Glauber and Langevin dynamics with annealing and averaging, reach in
polynomial time the algorithmic threshold of spectral and sum-of-squares methods
\cite{das2026averagedI,das2026averagedII}. Whether algorithms of this kind reach the threshold
$\Dcomp$ for mixtures is open.

\subsection{Reconstructing the full partition with high probability}\label{sec:reconstruction}

Corollary~\ref{cor:UBpoly} decides whether the individuals $1$ and $2$ belong to the same group
with an error probability that tends to $0$, but more slowly than $n^{-2}$, so that a union bound
over the pairs of individuals does not recover the partition. Following \cite{carpentier2025phase},
we amplify the probability of success by a median-of-means step. Given a pair of individuals, we
split the other individuals, and the variables, into $\beta$ batches, with $\beta$ of order $\log n$. The
data of the pair and of the individuals of a batch, restricted to the variables of the same batch,
follow the model~\eqref{eq:model} with smaller values of $n$, $d$ and $\Delta$, and these $\beta$ data
sets are independent conditionally on $x=\one\{k^*(1)=k^*(2)\}$. The majority vote of the $\beta$
resulting tests is then wrong with probability at most $n^{-3}$, and a union bound over the pairs
recovers the partition.

Let us describe the procedure for the pair formed by individuals $1$ and $2$. By symmetry of the
model, an identical construction, with $1,2$ replaced by $i,i'$, yields an estimator $\hat
x_{ii'}$ of $\one\{i\sim i'\}$ for every other pair $i<i'$. Put $\beta:=\lceil24\log n\rceil$ and, when
$n\ge \beta+2$ and $d\ge \beta$,
\[
  n_0:=\Big\lfloor\frac{n-2}{\beta}\Big\rfloor+2,\qquad
  d_0:=\Big\lfloor\frac d\beta\Big\rfloor,\qquad
  \Delta_0:=\Delta\Big(\frac{d_0}{d}\Big)^{1/2} .
\]
Listing the elements of $[n]\setminus\{1,2\}$ in increasing order, let $I_b$ be the $b$-th run of
$n_0-2$ consecutive elements of this list, and let $J_b:=\{(b-1)d_0+1,\dots,bd_0\}$, for
$b\in[\beta]$. The remaining individuals and variables are not used. Let $Y^{(b)}$ be the $n_0\times d_0$ matrix
whose rows are the restrictions to the variables $J_b$ of $Y_{1,\cdot}$, of $Y_{2,\cdot}$, and then
of the $Y_{u,\cdot}$, $u\in I_b$, in increasing order of $u$. Given a test $\psi:\mathbb
R^{n_0\times d_0}\to\{0,1\}$, put
\[
  \hat x\ :=\ \one\Big\{\sum_{b=1}^{\beta}\psi\big(Y^{(b)}\big)>\frac \beta2\Big\} ,
\]
the majority vote of the $\beta$ values $\psi(Y^{(b)})$. Combined with Corollary~\ref{cor:UBpoly}, this
majority vote recovers the partition under
the conditions of Corollary~\ref{cor:UBpoly}, up to logarithmic factors.

\begin{corollary}\label{cor:partition}
Fix an integer $q\ge3$ and a constant
$C_{\mathrm{poly}}>0$, and let $c$ be given by Corollary~\ref{cor:UBpoly}. For all $K$ large enough
in terms of $q$ and $C_{\mathrm{poly}}$: whenever $K\ge d$,
\[
  K^2\ <\ n\ \le\ K^2\Big(\frac d{50\log n}\Big)^{q/2}\qquad\text{and}\qquad d\ \ge\ K^{2/q}e^{\,2c\,(\log K)^{1/2}}
\]
and
\[
  \frac{\Delta^{2q}\,n}{K^{2}d^{\,q/2}}\ \ge\ e^{\,2c\,(\log K)^{1/2}}
  \qquad\text{and}\qquad \Delta^{2}<\sqrt d ,
\]
the estimator $\hat x$, built with $\psi$ the test of Corollary~\ref{cor:UBpoly} in the model with
parameters $(K,n_0,d_0,\Delta_0)$, a test based on a polynomial of degree at most
$C_{\mathrm{poly}}\log K$, satisfies $\Pr[\hat x=x]\ge1-n^{-3}$.
\end{corollary}
The proof is given in Appendix~\ref{app:MOM}.
Up to the factor $50\log n$ in the upper bound on $n$ and to the replacement of $c$ by $2c$,
the conditions of Corollary~\ref{cor:partition} are those of Corollary~\ref{cor:UBpoly}: the
median-of-means step only changes the conditions by logarithmic factors.

By the same argument, for every pair of individuals $i<i'$ the analogous estimator $\hat x_{ii'}$ of
$x_{ii'}:=\one\{i\sim i'\}$ satisfies $\Pr[\hat x_{ii'}=x_{ii'}]\ge1-n^{-3}$, so that a union bound over
the $\binom n2$ pairs gives $\hat x_{ii'}=x_{ii'}$ for all $i<i'$ with probability at least
$1-\frac1{2n}$. As the partition is determined by the values of the $x_{ii'}$, it is then recovered
with the same probability.

\subsection{Symmetrized Gaussian centers}\label{sec:symcenters}

Section~\ref{sec:bch} contrasts two
features. The set of the BCH centers of \cite{carpentier2026small} is symmetric, so that all its
moment tensors of odd order vanish exactly, but its first informative tensor has a rigid structure
that is easy to detect. Gaussian centers have no such structure, but none of their moment tensors
vanishes. The following model combines the two features: we sample $K$ Gaussian centers and add
their opposites.

Draw $k^*$, $X$ and $\epsilon$ as in Section~\ref{sec:model} and, independently of them, signs
$z_1,\dots,z_n$ i.i.d.\ uniform on $\{-1,1\}$, and set
\begin{equation}\label{eq:modelsym}
  Y_{i,j}\;=\;z_i\,\frac{\Delta}{\sqrt d}\,X_{k^*(i),j}+\epsilon_{i,j} .
\end{equation}
Given the centers, the rows $Y_{i,\cdot}$ are i.i.d.\ from the mixture, with equal weights, of the
$2K$ distributions $\mathcal N(\pm\tfrac{\Delta}{\sqrt d}X_{k,\cdot},I_d)$, $k\in[K]$. Two opposite
centers are at distance close to $2\Delta$ and two other distinct centers at distance close to
$\sqrt2\,\Delta$, so that the separation is still $\sqrt2\,\Delta$ up to a factor close to $1$. As
in \cite{carpentier2026small}, the functional of interest remains $x=\one\{k^*(1)=k^*(2)\}$
\eqref{eq:x}, the indicator that the centers of the first two individuals are equal or opposite, so
that $\E x=1/K$. We keep $p=2\log K/\log d$, $K$ being now the number of pairs of opposite centers.

In this model, the law of $(x,Y)$ is also invariant under the sign change of any row $Y_{i,\cdot}$,
$i\in[n]$, the rows $1$ and $2$ included: $-Y_{i,\cdot}$ is obtained by replacing
$(z_i,\epsilon_{i,\cdot})$ by $(-z_i,-\epsilon_{i,\cdot})$, which has the same law and leaves $x$
unchanged. Hence, in the analogue of Proposition~\ref{prop:sym}, the polynomials $P_G$ may be
restricted to the graphs all of whose nodes, \emph{individual nodes included}, have even degree.
Accordingly, write $\lceil p\rceil_{\mathrm{even}}$ for the smallest even integer larger than or
equal to $p$, that is $\lceil p\rceil$ if $\lceil p\rceil$ is even and $\lceil p\rceil+1$ otherwise.
It will play the role of the critical order $\lceil p\rceil$.

\begin{theorem}\label{thm:LBsym}
Let $K=d^{p/2}$ with $p\ge2$, $D\ge2$ an integer, $n\ge16D^2$, and
\[
\text{(i) } \frac{\Delta}{d^{1/4}}\le D^{-200},\qquad
\text{(ii) } \frac{\Delta^{2\lceil p\rceil_{\mathrm{even}}}\,n}{K^{2}d^{\,\lceil p\rceil_{\mathrm{even}}/2}}\le
D^{-400\lceil p\rceil_{\mathrm{even}}},\qquad
\text{(iii) } d^{\,p-\lceil p\rceil_{\mathrm{even}}+2}\ \ge\ D^{1600\lceil p\rceil_{\mathrm{even}}} .
\]
Then the correlation \eqref{eq:corr}, computed for the observations \eqref{eq:modelsym}, is at
most $4/K$.
\end{theorem}
The proof is given in Appendix~\ref{app:LBsym}.

Conditions \textup{(i)} and \textup{(ii)} are those of Assumption~\ref{ass:A} with
$\lceil p\rceil_{\mathrm{even}}$ in place of $\lceil p\rceil$, and condition \textup{(iii)} only
excludes the values of $p$ slightly above an even integer. When $\lceil p\rceil$ is even, the
conditions are those of Assumption~\ref{ass:A}, with a weaker \textup{(iii)}. When $\lceil p\rceil$
is odd, $\lceil p\rceil_{\mathrm{even}}=\lceil p\rceil+1$: up to polylogarithmic factors,
low-degree polynomials fail as long as $\Delta^2$ is smaller than $\sqrt
d\min\{1,(K^2/n)^{1/(\lceil p\rceil+1)}\}$, which is larger than $\Dcomp^2$ when $n>K^2$, and, at a
polylogarithmic separation, they need $n\ge K^2d^{(\lceil p\rceil+1)/2}$ instead of
$K^2d^{\lceil p\rceil/2}$. For instance, for $K=2d$, Gaussian centers need $K^2d^{3/2}$, of order
$K^{7/2}$, while the orthogonal construction of \cite{carpentier2026low}, whose centers come with their
opposites, requires $K^4$. Symmetrized Gaussian centers recover this bound: for $K$ larger than $d$ by a
large enough polylogarithmic factor, the critical order becomes $4$, and the sample size needed at a
polylogarithmic separation becomes $K^2d^2$, which is $K^4$ up to polylogarithmic factors. In particular, as $\lceil p\rceil_{\mathrm{even}}<p+2$ and $d^{p/2}=K$,
the sample size $K^2d^{\lceil p\rceil_{\mathrm{even}}/2}$ required at a polylogarithmic separation is
always smaller than $K^3d$, and it is of this order, up to subpolynomial factors, when $p$ is
slightly above an even integer. For Gaussian centers, $K^2d^{\lceil p\rceil/2}$ is always smaller
than $K^3\sqrt d$. This could lead us to conjecture that $n\ge K^3d$ suffices to cluster at a
subpolynomial separation, for arbitrary configurations of the centers. The orthogonal construction of
\cite{carpentier2026low}, at $K=2d$, attains this order: at a polylogarithmic separation it requires
$n$ of order $K^4=2dK^3$.

\section{Proof ideas}\label{sec:ideas}

Section~\ref{sec:intuition} gave some intuition on the computational threshold $\Dcomp$ by comparing the squared
gap $\E[T_G\mid x=1]^2$ of an invariant monomial with two parts of its variance, $\mathcal V_{\mathrm{noise}}$, due to
the noise, and $\mathcal V_{\mathrm{centers}}$, due to the random centers. The two proofs show that this comparison
indeed gives the threshold. Below $\Dcomp$, the lower bound shows that $\Corr_{\le D}$ is at most of order $1/K$: low-degree
polynomials do not improve over a constant (Section~\ref{sec:intro-ld}). Above $\Dcomp$, the upper bound exhibits one
invariant monomial that decides $x$.

For the lower bound, the difficulty is that a polynomial combines many invariant monomials, and that these turn out to be
correlated. For the upper bound, consider a
graph whose individual nodes other than $v_1,v_2$ have degree $\lceil p\rceil$. Above the threshold $\Dcomp$, the variance terms $\mathcal
V_{\mathrm{noise}}$ and $\mathcal V_{\mathrm{centers}}$ of Section~\ref{sec:intuition} are small enough, but the variance contains many other terms, and the difficulty is to find a template that makes all of them small enough.

\subsection{The first two moments of an invariant monomial}\label{sec:moments}

Both proofs rest on the first two moments of the invariant monomials. We compute them for a connected graph
$G=(V,W,E)$ whose variable nodes are \emph{links}: each is joined to two distinct individual nodes $v,v'$, and
written $w=\{v,v'\}$. Such a graph has no parallel edges. The template of the upper bound is such a graph, and these
graphs are the most informative ones for the lower bound (Section~\ref{sec:othersketch}). The path
$v_1,v_3,v_4,v_2$ of Figure~\ref{fig:twoterms} serves as an example. Recall that the invariant monomial $T_G$ of such a graph is defined in \eqref{eq:TGintro}.

\paragraph{The first moment.} By Section~\ref{sec:intuition}, $\E[T_G\mid x=0]=0$, and $\E[T_G\mid x=1]$ is of
order $n^{|V|-2}\Delta^{2|W|}K^{-(|V|-2)}$. As
$\E x=1/K$,
\begin{equation}\label{eq:corrTG}
  \E\big[(x-\E x)\,T_G\big]\ =\ \frac1K\Big(1-\frac1K\Big)\E[T_G\mid x=1] .
\end{equation}

\paragraph{The second moment.} The second moment is a sum over the pairs $(\pi,\pi')$ of labelings, $\E[T_G^2\mid
x]=\sum_{\pi,\pi'}\E[T_{G,\pi}T_{G,\pi'}\mid x]$. Call a pair \emph{disjoint} when $\pi_1(V)\cap\pi'_1(V)=\{1,2\}$
and $\pi_2(W)\cap\pi'_2(W)=\varnothing$. Given the groups, the two terms of a disjoint pair carry distinct
coordinates, hence are independent, and their observations other than $1,2$ are distinct. Hence, for such a disjoint pair, we have 
\begin{align}\label{eq:first_moment_intuition}
\E[T_{G,\pi}T_{G,\pi'}\mid x]= \E[T_{G,\pi}\mid x]\E[T_{G,\pi'}\mid x]\enspace .
\end{align}
Write $\rho$ for the proportion of disjoint pairs among all pairs $(\pi,\pi')$. We have $1-\rho\lesssim|V|^2/n+|W|^2/d$ and 
\begin{equation}\label{eq:varsplit}
  \Var(T_G\mid x)\ =\ \sum_{(\pi,\pi')\ \text{not disjoint}}\E\big[T_{G,\pi}T_{G,\pi'}\mid x\big]\ -\
  (1-\rho)\,\E[T_G\mid x]^2 ,
\end{equation}
and the variance comes from the pairs that share an observation other than $1,2$, or a coordinate.

To sort these pairs, write each entry as the sum of a signal and a noise,
$Y_{ij}=\frac{\Delta}{\sqrt d}X_{k^*(i),j}+\epsilon_{ij}$, and expand $T_{G,\pi}$: each edge $(v,w)$ carries the
{\bf signal} part $\frac{\Delta}{\sqrt d}X_{k^*(\pi_1(v)),\pi_2(w)}$ of its entry, or its {\bf noise} part
$\epsilon_{\pi_1(v),\pi_2(w)}$. For a set $S$ of individual nodes other than $v_1,v_2$, keep the terms in which the
edges at the nodes of $S$ carry signal parts, the other edges noise parts, and the nodes of $S$ lie in distinct
groups:
\begin{equation}\label{eq:VS}
\begin{aligned}
  T^S_G\ &:=\ \sum_{\pi}\one\big\{k^*\circ\pi_1\ \text{is one-to-one on}\ S\big\}
  \prod_{v\in S}\ \prod_{w:\,(v,w)\in E}\frac{\Delta}{\sqrt d}\,X_{k^*(\pi_1(v)),\pi_2(w)}\\
  &\qquad\qquad\times\prod_{v\in V\setminus S}\ \prod_{w:\,(v,w)\in E}\epsilon_{\pi_1(v),\pi_2(w)} ,\\
  \mathcal V(S)\ &:=\ \E\big[(T^S_G)^2\big] .
\end{aligned}
\end{equation}
In $T^S_G$, each node of $S$ is replaced by the center of its group, its \emph{color}, and the other nodes keep their
noise. It turns out that  $\Var(T_G)\geq \sum_{S\subset V\setminus\{v_1,v_2\}} \mathcal V(S)$.  Within the decomposition of $\Var(T_G)$ as a sum over the pairs of labelings $\pi$ and $\pi'$, we interpret
$\mathcal V(S)$ as the sum over  the pairs $(\pi,\pi')$ of labelings such that  $\pi_2(w)=\pi'_2(w)$ for $w\in W$ and $\pi_1(v)=\pi'_1(v)$ for $v\notin S$ (as the noises of distinct observations are independent), and such that each node of $S$
has the same color under $\pi$ and $\pi'$ ($k^*(\pi_1(v))=k^*(\pi'_1(v))$ for $v\in S$). The two parts of Section~\ref{sec:intuition} are the extreme cases,
$\mathcal V_{\mathrm{noise}}\asymp\mathcal V(\varnothing)$ and $\mathcal V_{\mathrm{centers}}\asymp\mathcal
V(V\setminus\{v_1,v_2\})$.

The variance $\Var(T_G)$ contains other terms, and their analysis will be important for the upper bound  (Section~\ref{sec:construction}). For the sketch of the lower bound, however, the inequality $\Var(T_G)\geq \sum_{S\subset V\setminus\{v_1,v_2\}} \mathcal V(S)$ will be sufficient. 

Let us evaluate $\mathcal V(S)$. In light of the previous discussion, there are of the order $d^{|W|}n^{|V|-2+|S|}$ pairs $(\pi,\pi')$ that contribute to $\mathcal{V}(S)$. The probability that each node of $S$ has the same color under $\pi$ and $\pi'$ is $K^{-|S|}$. Taking into account the signal term  
$(\Delta^2/d)^{\sum_{v\in S} \deg v}$,
and together with \eqref{eq:first_moment_intuition}, we obtain
\begin{equation}\label{eq:nodefactors}
  \frac{\E[T_G\mid x=1]^2}{\mathcal V(S)}\ \asymp\ \prod_{v\in S}\frac{d^{\deg v/2}}{K}\
  \prod_{v\notin S\cup\{v_1,v_2\}}\frac{n\,\Delta^{2\deg v}}{K^2d^{\deg v/2}} .
\end{equation}
Each node brings its factor of Section~\ref{sec:intuition}, and $T_G$ can only decide $x$ if \eqref{eq:nodefactors}
is at least of order $1$ for every $S$.

\subsection{The lower bound: projection on an orthonormal family, and centering of the polynomials}\label{sec:overview}

For the lower bound, we must show that $\Corr_{\le D}$ is at most of order $1/K$ (Theorem~\ref{thm:LB}).

\subsubsection{Reduction to invariant polynomials}\label{sec:redsketch}

The supremum \eqref{eq:corrintro} runs over all the polynomials of degree at most $D$. As the law of $(x,Y)$ is
invariant under the permutations of the observations $3,\dots,n$, under the permutations of the coordinates and
under their sign changes (Section~\ref{sec:intuition}), the supremum is unchanged when it runs over the invariant polynomials
\begin{equation}\label{eq:span}
  P\ =\ u_0+\sum_Gu_G\,P_G ,
\end{equation}
with real coefficients $u_0$ and $u_G$ (Proposition~\ref{prop:sym}), where $G$ runs over the bipartite multigraphs
with at most $D$ edges whose variable nodes have even degree. To simplify the discussion, we first restrict ourselves to graphs $G$ that are connected and whose variable nodes are links, so that $P_G$ is a multiple of $T_G$ (see Section~\ref{sec:othersketch} for the extension to all $G$).

\subsubsection{The blocks}\label{sec:blocktypes}

\begin{figure}[ht]
\centering
\begin{tikzpicture}[x=1cm,y=1cm,
  di/.style={circle,draw=black,fill=white,inner sep=1pt,font=\scriptsize},
  nz/.style={circle,fill=black,inner sep=1.6pt},
  sg/.style={rectangle,fill=OliveGreen!80!black,inner sep=2.2pt},
  lab/.style={font=\scriptsize,inner sep=1.5pt}]
  \node[di] (v1) at (0,0) {$v_1$};
  \node[nz,label={[lab]above:$v_3$}] (a) at (1.3,0) {};
  \node[sg,label={[lab]above:$v_4$}] (b) at (2.5,0.65) {};
  \node[sg,label={[lab]below:$v_5$}] (c) at (2.5,-0.65) {};
  \node[nz,label={[lab]above:$v_6$}] (e) at (3.7,0) {};
  \node[di] (v2) at (5,0) {$v_2$};
  \draw[thick] (v1)--(a) (a)--(b) (a)--(c) (b)--(e) (c)--(e) (e)--(v2);
  \node[nz] at (6,0.45) {};
  \node[lab,anchor=west] at (6.2,0.45) {keeps its noise: $\epsilon_{\pi_1(v),\pi_2(w)}$ at each edge $(v,w)$};
  \node[sg] at (6,-0.45) {};
  \node[lab,anchor=west] at (6.2,-0.45) {replaced by a center: $X_{\sigma(v),\pi_2(w)}$ at each edge $(v,w)$};
\end{tikzpicture}
\caption{A block of type $S=\{v_4,v_5\}$ (Definition~\ref{def:block}). Each line is a link. The
nodes $v_4$ and $v_5$ are replaced by the centers of their colors $\sigma(v_4)\ne\sigma(v_5)$, and the other
nodes keep the noise of their observations, $v_1$ and $v_2$ being placed on the observations $1$ and $2$. At
$\lceil p\rceil=3$, this is the type kept by the lower bound (Section~\ref{sec:blocksketch}): $v_4$ and $v_5$
have degree $2<p$, and $v_3$ and $v_6$ degree $3\ge p$. }
\label{fig:blocks}
\end{figure}
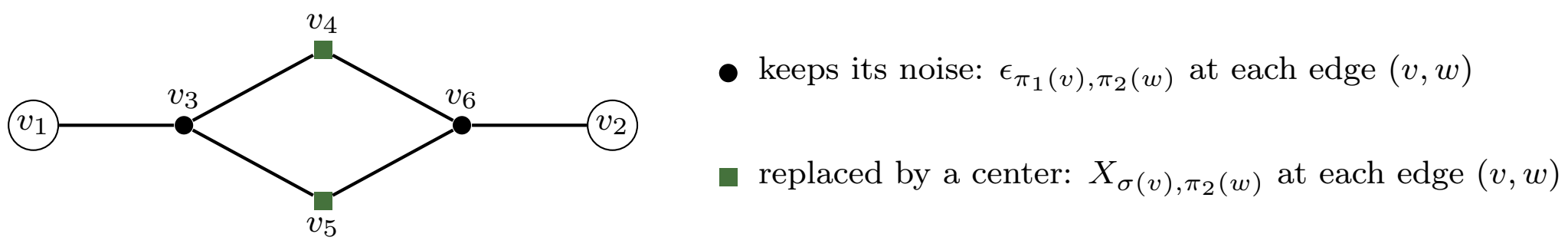

Grouping the terms of $T^S_G$ in \eqref{eq:VS} by the colors $\sigma(v)\in[K]$ of the nodes of $S$ gives the
products
\begin{equation}\label{eq:blockS}
  Q^S_{G,\pi_1,\pi_2,\sigma}\ :=\ \prod_{v\in S}\ \prod_{w:\,(v,w)\in E}X_{\sigma(v),\pi_2(w)}\
  \prod_{v\in V\setminus S}\ \prod_{w:\,(v,w)\in E}\epsilon_{\pi_1(v),\pi_2(w)} ,
\end{equation}
indexed by the distinct colors $\sigma(v)$ of the nodes of $S$, the observations $\pi_1(v)$ of the other individual
nodes, and the coordinates $\pi_2(w)$ of the variable nodes (Figure~\ref{fig:blocks}). These \emph{blocks of $G$ of type
$S$} (Definition~\ref{def:block}) are products of distinct independent standard Gaussian variables, hence
orthonormal, and the blocks of distinct types are orthogonal. Define the block of type $S$ as $Q^S_G= \sum_{\sigma, \pi_1,\pi_2}Q^S_{G,\pi_1,\pi_2,\sigma}$. It turns out that $\mathcal V(S)$ is of the order of the squared norm of the projection of $T_G$ on the blocks $Q^S_G$.

\subsubsection{One type of blocks for each graph}\label{sec:blocksketch}

Among the blocks of $G$, a single type is of interest. As announced in Section~\ref{sec:intuition}, it is the
type $S_G$ in which the low-degree nodes
are replaced by centers and the high-degree nodes keep their noise. More precisely, for a graph $G$, we define $S_G:=\{v\in V\setminus\{v_1,v_2\}:\deg_G(v)<p\}$. 

This choice is the crucial one: at the threshold $\Delta=\Dcomp$, it corresponds to the largest term $\mathcal V(S)$, as it takes, node by node, the smaller of the two factors of
\eqref{eq:nodefactors}. Below the threshold ($\Delta<\Dcomp$),  $\mathcal V(S_G)$ is not necessarily the largest term, but it still plays a special role.

The \emph{invariant block} $\bar Q_G$ of $G$ is the normalized version of the block $Q_G^{S_G}$. As explained above, the invariant blocks of distinct graphs are orthonormal (Lemma~\ref{lem:ortho}). The
projection of $T_G$ on the blocks of type $S_G$ is $\langle T_G,\bar Q_G\rangle\bar Q_G$, with $\langle T_G,\bar
Q_G\rangle^2\asymp\mathcal V(S_G)$, and \eqref{eq:nodefactors} reads
\begin{equation}\label{eq:nodefactorsG}
  \frac{\E[T_G\mid x=1]^2}{\langle T_G,\bar Q_G\rangle^2}\ \asymp\
  \prod_{\substack{v\ne v_1,v_2\\ \deg_G(v)<p}}\frac{d^{\deg_G(v)/2}}{K}\
  \prod_{\substack{v\ne v_1,v_2\\ \deg_G(v)\ge p}}\frac{n\,\Delta^{2\deg_G(v)}}{K^2d^{\deg_G(v)/2}} .
\end{equation}
Every factor is small below the threshold: a signal node has degree at most $\lceil p\rceil-1<p$, so that
$d^{\deg_G(v)/2}<K$, and a noise node has degree at least $\lceil p\rceil$, so that its factor is at most the one at
degree $\lceil p\rceil$, small below the threshold (Assumption~\ref{ass:A}). By \eqref{eq:corrTG}, the correlation of
$T_G$ with $x$, divided by $|\langle T_G,\bar Q_G\rangle|$, is therefore at most $1/K$ times a small factor for each
individual node (Lemma~\ref{lem:master}). For the path of Figure~\ref{fig:twoterms} and $p>2$, both $v_3$ and $v_4$
are signal nodes, and the invariant block of the path carries $\mathcal V_{\mathrm{centers}}$.

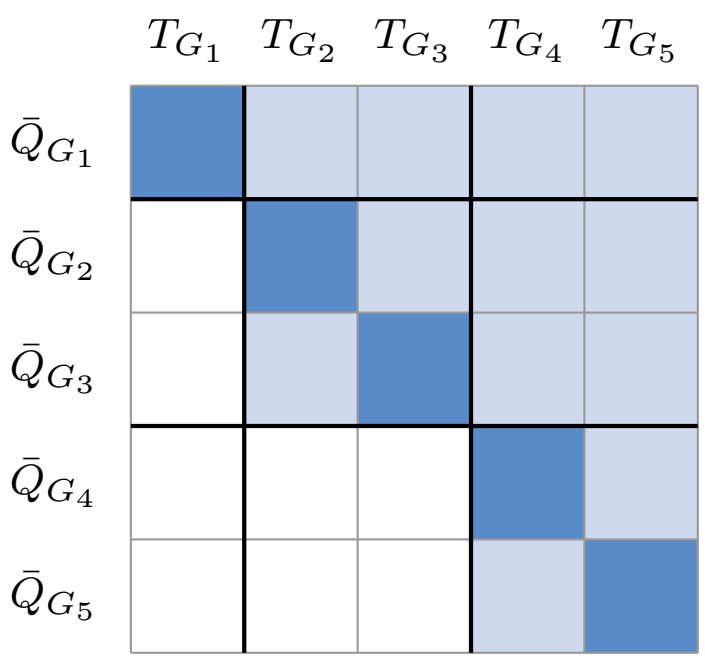
\begin{figure}[ht]
\centering
\begin{tikzpicture}[x=0.75cm,y=0.75cm]
  \foreach \r in {1,...,5}{
    \foreach \c in {1,...,5}{
      \ifnum\r=\c\relax\fill[NavyBlue!70] (\c-1,-\r) rectangle (\c,1-\r);\fi
      \ifnum\r<\c\relax\fill[NavyBlue!18] (\c-1,-\r) rectangle (\c,1-\r);\fi
    }
  }
  \fill[NavyBlue!18] (1,-3) rectangle (2,-2);
  \fill[NavyBlue!18] (3,-5) rectangle (4,-4);
  \foreach \i in {0,...,5} \draw[black!40] (\i,0)--(\i,-5) (0,-\i)--(5,-\i);
  \draw[black,thick] (1,0)--(1,-5) (3,0)--(3,-5) (0,-1)--(5,-1) (0,-3)--(5,-3);
  \foreach \c in {1,...,5} \node[font=\scriptsize] at (\c-0.5,0.4) {$T_{G_\c}$};
  \foreach \r in {1,...,5} \node[font=\scriptsize,anchor=east] at (-0.1,0.5-\r) {$\bar Q_{G_\r}$};
\end{tikzpicture}
\caption{The projections $\langle T_G,\bar Q_H\rangle$ of the invariant monomials (columns) on the invariant
blocks (rows), the graphs being ordered by their numbers of edges ($G_2$ and $G_3$ have the same number of edges, as
$G_4$ and $G_5$). Each invariant monomial projects on its own invariant block (dark), possibly on the invariant
blocks of graphs with at most as many edges (light), and on no other (white). Weighted by the correlations with $x$,
the matrix is diagonally dominant (Section~\ref{sec:combsketch}).}
\label{fig:matrix}
\end{figure}

\subsubsection{From invariant monomials to invariant polynomials}\label{sec:combsketch}

Here, our approach follows the ideas of Sohn and Wein~\cite{SohnWein25}: we bound the variance of a polynomial
$P=\sum_Gu_GT_G$ from below by its projections onto the orthonormal family formed by the invariant blocks $\bar Q_H$, that
is, we use $\Var P\geq\sum_H\langle P,\bar Q_H\rangle^2$. However, we depart from it in the way we control the ratio $\E\left[(x-\E x)P\right]^2/\sum_H\langle P,\bar Q_H\rangle^2$: we do not construct a certificate as in~\cite{SohnWein25} but we rather rely on a more direct approach.

If each invariant monomial $T_G$ were orthogonal to $\bar Q_H$ for all $H$ non-isomorphic to $G$, then we would easily bound the low-degree correlation. Indeed, consider a polynomial of the form $P= \sum_{G}u_G T_G$. Then, by Bessel's inequality and by the orthogonality of the $\bar{Q}_H$'s, 
\begin{align*}
\Var{P}\geq \sum_{H}\langle P,  \bar{Q}_H \rangle^2\asymp \sum_{H} u_H^2 \mathcal V(S_{H})\ ,
\end{align*}
so that, by \eqref{eq:corrTG} and the Cauchy--Schwarz inequality, we would get
\begin{align*}
\frac{\E\left[(x-\E x)P\right]^2}{\Var{P}}\lesssim \frac{1}{K^2} \frac{\left[\sum_{G} u_G\E[T_G\mid x=1]\right]^2 }{\sum_{G} u_G^2 \mathcal V(S_G)}\leq \frac{1}{K^2}\sum_{G} \frac{\E[T_G\mid x=1]^2}{\mathcal V(S_G)}\ .  
\end{align*}
In light of~\eqref{eq:nodefactorsG}, and since the number of non-isomorphic graphs with a given number of edges is controlled (Lemma~\ref{lem:count}), the right-hand side would be $O(1/K^2)$ for $\Delta$ small compared to $\Dcomp$.

Unfortunately, the monomial $T_G$ is not orthogonal to $\bar{Q}_H$ for $G\neq H$ (see Figure~\ref{fig:matrix}). In particular, we may have $\langle T_G, \bar{Q}_H \rangle \neq 0$ when $H$ has at most as many edges as $G$. Nor does the matrix $(\langle T_G, \bar{Q}_H \rangle)_{G,H}$ satisfy the nice property of being diagonally dominant. However, it turns out that a slight variant of the weighted matrix $B$ defined by  
\[
B_{H,G}= \left(\frac{\mathbb{E}[T_H\mid x=1]}{\mathbb{E}[T_G\mid x=1] [\mathcal V(S_{H})]^{1/2}}\langle T_G, \bar{Q}_H \rangle\right)
\]
is diagonally dominant, with a diagonal close to one (see Corollary~\ref{cor:diag} and Lemma~\ref{lem:colsum}). In other words, the projection of $T_G$ onto $\bar{Q}_H$, rescaled by the correlation with $x$, is extremely small for $H\neq G$.
This diagonal dominance is instrumental in the proof. Indeed, for any polynomial $P= \sum_{G}u_G T_G$, it implies (see Proposition~\ref{prop:reduce}) that 
\[
\sum_{G}  \frac{\mathbb{E}[T_G\mid x=1]}{[\mathcal V(S_G)]^{1/2}} \big|\langle P,\bar Q_G\rangle\big| \gtrsim  \sum_G \mathbb{E}[T_G\mid x=1]|u_G|\enspace . 
\]
Together with \eqref{eq:corrTG}, the orthonormality of the invariant blocks and the Cauchy--Schwarz inequality, this gives
\begin{align*}
  \E\big[(x-\E x)P\big]\ & \le\ \frac{1}{K}\sum_G \mathbb{E}[T_G\mid x=1]|u_G|\\ &  \lesssim \frac{1}{K}\sum_G \frac{\mathbb{E}[T_G\mid x=1]}{[\mathcal V(S_G)]^{1/2}} \big|\langle P,\bar Q_G\rangle\big|
  \ \le  \frac{1}{K}\Big(\sum_G \frac{\mathbb{E}[T_G\mid x=1]^2}{\mathcal V(S_G)}\Big)^{1/2}\|P\| .
\end{align*}
We conclude since $\sum_G \frac{\mathbb{E}[T_G\mid x=1]^2}{\mathcal V(S_G)}$ is bounded for $\Delta$ small compared to $\Dcomp$.

\subsubsection{The other graphs, and the centering}\label{sec:othersketch}

So far, the graphs were connected and their variable nodes links. The invariant polynomials \eqref{eq:span} also
involve graphs with parallel edges, with several connected components, and with variable nodes of degree larger than $2$. 

Their
invariant monomials are centered (Definition~\ref{def:center}): multiples of the invariant monomials of smaller
graphs are subtracted, as in the centering by connected components of \cite{montanari2024equivalence,CGGV25}. The
centered polynomials carry no correlation with $x$, or a correlation that is small relative to their projection on
their invariant block, and they project on the invariant blocks as the invariant monomials of
Section~\ref{sec:combsketch} do (Theorem~\ref{thm:cancel}), so that the bound of Section~\ref{sec:combsketch}
applies to all the invariant polynomials. Three mechanisms are at work.

\noindent\uline{(a) Parallel edges.} The $m$ entries carried by $m$ parallel edges enter through a Hermite polynomial
of the entry, which has mean zero given the groups.

\noindent\uline{(b) Several connected components.} Given the groups, distinct components are independent. If $v_1$
and $v_2$ lie in distinct components, the invariant monomial is uncorrelated with $x$, by symmetry.
If $G$ is the disjoint union of two graphs $C$ and $G_1$, where $C$ contains neither $v_1$ nor $v_2$, then, by conditional independence and recentering, we have, up to the labelings in which $G_1$ and $C$ overlap,
\[
  \E\big[(x-\E x)\,T_{G_1}\big(T_C-\E[T_C]\big)\big]\ =\ 0 .
\]

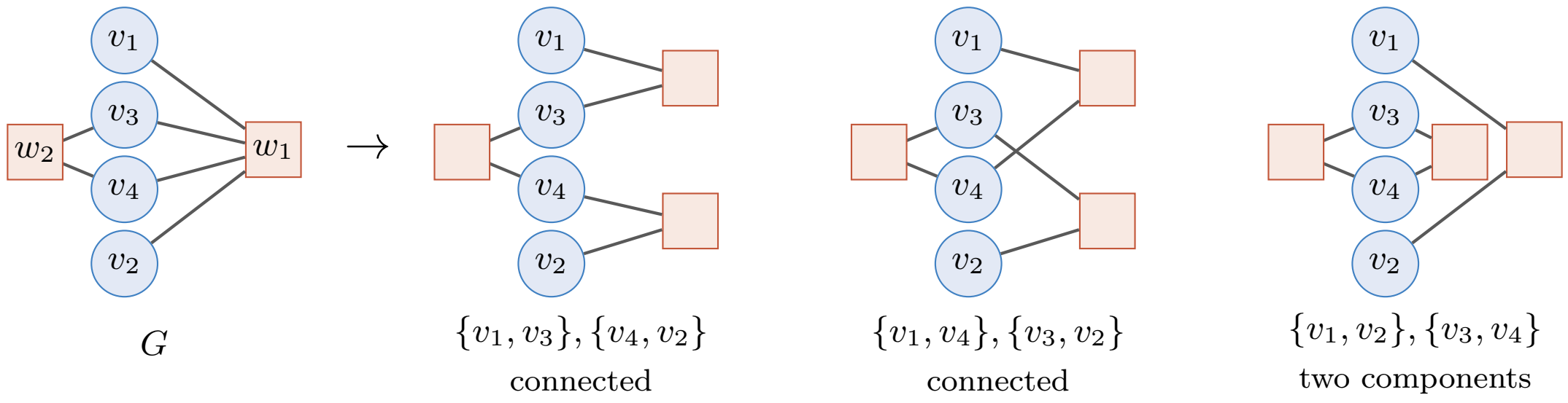
\begin{figure}[ht]
\centering
\begin{tikzpicture}[
  ind/.style={circle,draw=NavyBlue!80,fill=NavyBlue!10,inner sep=1pt,minimum size=6mm,font=\footnotesize},
  var/.style={rectangle,draw=BrickRed!80,fill=BrickRed!8,inner sep=1.5pt,minimum size=5mm,font=\footnotesize},
  e/.style={black!65,thick},scale=0.9]
\newcommand{\fournodes}{%
  \node[ind] (n1) at (0,1.1) {$v_1$}; \node[ind] (n2) at (0,0.35) {$v_3$};
  \node[ind] (n3) at (0,-0.4) {$v_4$}; \node[ind] (n4) at (0,-1.15) {$v_2$};
  \draw[e] (n2)--(-0.9,-0.025)--(n3);}
\begin{scope}
  \fournodes
  \node[var] at (-0.9,-0.025) {$w_2$};
  \node[var] (w) at (1.5,0) {$w_1$};
  \foreach \n in {n1,n2,n3,n4}{\draw[e] (\n)--(w);}
  \node at (0.3,-1.95) {\footnotesize $G$};
\end{scope}
\node at (2.45,0) {\large $\to$};
\begin{scope}[xshift=4.3cm]
  \fournodes
  \node[var] at (-0.9,-0.025) {};
  \node[var] (p) at (1.4,0.72) {}; \node[var] (q) at (1.4,-0.72) {};
  \draw[e] (n1)--(p); \draw[e] (n2)--(p); \draw[e] (n3)--(q); \draw[e] (n4)--(q);
  \node[align=center] at (0.3,-2.05) {\footnotesize $\{v_1,v_3\},\{v_4,v_2\}$\\[-1pt]\scriptsize connected};
\end{scope}
\begin{scope}[xshift=8.5cm]
  \fournodes
  \node[var] at (-0.9,-0.025) {};
  \node[var] (p) at (1.4,0.72) {}; \node[var] (q) at (1.4,-0.72) {};
  \draw[e] (n1)--(p); \draw[e] (n3)--(p); \draw[e] (n2)--(q); \draw[e] (n4)--(q);
  \node[align=center] at (0.3,-2.05) {\footnotesize $\{v_1,v_4\},\{v_3,v_2\}$\\[-1pt]\scriptsize connected};
\end{scope}
\begin{scope}[xshift=12.7cm]
  \fournodes
  \node[var] at (-0.9,-0.025) {};
  \node[var] (p) at (1.5,0) {}; \node[var] (q) at (0.75,-0.025) {};
  \draw[e] (n1)--(p); \draw[e] (n4)--(p); \draw[e] (n2)--(q); \draw[e] (n3)--(q);
  \node[align=center] at (0.3,-2.05) {\footnotesize $\{v_1,v_2\},\{v_3,v_4\}$\\[-1pt]\scriptsize two components};
\end{scope}
\end{tikzpicture}
\caption{The three pairings of a variable node $w_1$ joined to $v_1$, $v_3$, $v_4$ and $v_2$, in a graph
with a link $w_2=\{v_3,v_4\}$ (Proposition~\ref{prop:pairid}). Each pairing replaces $w_1$ by two links on the
coordinate of $w_1$. The first two give the paths $v_1,v_3,v_4,v_2$ and $v_1,v_4,v_3,v_2$, connected graphs of
links whose first and last links share a coordinate, hence a signal smaller by a factor of order $\sqrt d$.
The third leaves $w_2$ away from $v_1$ and $v_2$, and the centering removes it.}
\label{fig:centering}
\end{figure}

\noindent\uline{(c) Variable nodes of degree larger than $2$.} 
Here, we rely on Wick's formula \cite[Chapter~3]{janson1997gaussian} (Lemma~\ref{lem:wick}) to control both the first and the second moments.\footnote{Actually, Isserlis's formula (a special case of Wick's formula) was already used to control the second moments in the previous part of the sketch.}
Given the groups, a cross-moment of the entries is a sum over pairings. For instance, we have 
\begin{equation}\label{eq:wickfour}
\begin{aligned}
  \E\big[Y_{i_1j}Y_{i_2j}Y_{i_3j}Y_{i_4j}\mid k^*\big]
&  =\E[Y_{i_1j}Y_{i_2j}\mid k^*]\,\E[Y_{i_3j}Y_{i_4j}\mid k^*]
  +\E[Y_{i_1j}Y_{i_3j}\mid k^*]\,\E[Y_{i_2j}Y_{i_4j}\mid k^*]\\ &\quad
  +\E[Y_{i_1j}Y_{i_4j}\mid k^*]\,\E[Y_{i_2j}Y_{i_3j}\mid k^*] ,
\end{aligned}
\end{equation}
so that we can express the first moment (and the second moment) as a sum, over the pairings of the edges, of the moments of new graphs in which each variable node has degree $2$. For most pairings, the new graph is disconnected, so that, by (b),
this pairing carries no correlation with $x$ after the centering. But some pairings may leave the graph
connected, as in Figure~\ref{fig:centering}, where the first two pairings give the paths $v_1,v_3,v_4,v_2$ and
$v_1,v_4,v_3,v_2$, whose first and last links share a coordinate.
Nevertheless, the original graph has fewer variable nodes than the
connected graphs obtained from the pairings, and its correlation with $x$, relative to its projection on its invariant block, is smaller than theirs by a factor $\sqrt d$ for each coordinate saved (Lemma~\ref{lem:master}).

\subsection{The upper bound: construction of the template}\label{sec:construction}

The upper bound needs one invariant monomial $T_G$ whose standard deviation, given each value of $x$, is small
compared with its gap $\E[T_G\mid x=1]$: comparing $T_G$ with its two conditional means then decides $x$
(Lemma~\ref{lem:criterion}). A connected graph of links whose individual nodes other than $v_1,v_2$ have degree
$q=\lceil p\rceil$ has small parts $\mathcal V_{\mathrm{noise}}$ and $\mathcal V_{\mathrm{centers}}$ above the
threshold (Section~\ref{sec:intuition}), but the other parts of its variance may be large, and the template must be
chosen so that all of them are small. As two copies of a graph can overlap in very many ways, this is the delicate part of the upper bound.

The purpose of this subsection is to work out the main terms of $\Var(T_G)$. In turn, this will allow us to identify the key properties that our template $G$ needs to satisfy.

\subsubsection{The second moment of the template}\label{sec:secondsketch}

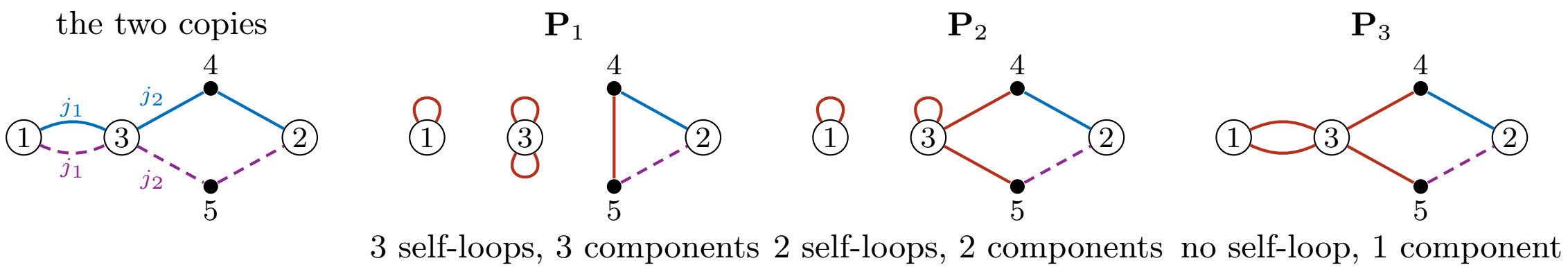
\begin{figure}[ht]
\centering
\begin{tikzpicture}[x=0.85cm,y=0.85cm,
  ob/.style={circle,fill=black,inner sep=1.4pt},
  di/.style={circle,draw=black,fill=white,inner sep=1pt,font=\scriptsize},
  lab/.style={font=\scriptsize,inner sep=1.5pt},
  cp/.style={NavyBlue,thick},
  cq/.style={Plum,thick,densely dashed},
  pr/.style={BrickRed,thick},
  sl/.style={BrickRed,thick,min distance=3.5mm},
  co/.style={font=\tiny,inner sep=1pt}]
  \newcommand{\obs}[1]{%
    \node[di] (o1#1) at (0,0.5) {$1$};
    \node[di] (o3#1) at (1.1,0.5) {$3$};
    \node[ob,label={[lab]above:$4$}] (o4#1) at (2.1,1.05) {};
    \node[ob,label={[lab]below:$5$}] (o5#1) at (2.1,-0.05) {};
    \node[di] (o2#1) at (3.1,0.5) {$2$};}
  \begin{scope}
    \obs{a}
    \draw[cp] (o1a) to[bend left=25] node[co,above] {$j_1$} (o3a);
    \draw[cq] (o1a) to[bend right=25] node[co,below] {$j_1$} (o3a);
    \draw[cp] (o3a) -- node[co,above left] {$j_2$} (o4a);
    \draw[cq] (o3a) -- node[co,below left] {$j_2$} (o5a);
    \draw[cp] (o4a) -- (o2a);
    \draw[cq] (o5a) -- (o2a);
    \node[font=\footnotesize] at (1.55,1.75) {the two copies};
  \end{scope}
  \begin{scope}[xshift=3.85cm]
    \obs{b}
    \draw[sl] (o1b) edge[loop,out=120,in=60] (o1b);
    \draw[sl] (o3b) edge[loop,out=120,in=60] (o3b);
    \draw[sl] (o3b) edge[loop,out=-60,in=-120] (o3b);
    \draw[pr] (o4b) -- (o5b);
    \draw[cp] (o4b) -- (o2b);
    \draw[cq] (o5b) -- (o2b);
    \node[font=\footnotesize] at (1.55,1.75) {$\Pai_1$};
    \node[font=\footnotesize] at (1.55,-0.75) {3 self-loops, 3 components};
  \end{scope}
  \begin{scope}[xshift=7.7cm]
    \obs{c}
    \draw[sl] (o1c) edge[loop,out=120,in=60] (o1c);
    \draw[sl] (o3c) edge[loop,out=120,in=60] (o3c);
    \draw[pr] (o3c) -- (o4c);
    \draw[pr] (o3c) -- (o5c);
    \draw[cp] (o4c) -- (o2c);
    \draw[cq] (o5c) -- (o2c);
    \node[font=\footnotesize] at (1.55,1.75) {$\Pai_2$};
    \node[font=\footnotesize] at (1.55,-0.75) {2 self-loops, 2 components};
  \end{scope}
  \begin{scope}[xshift=11.55cm]
    \obs{d}
    \draw[pr] (o1d) to[bend left=25] (o3d);
    \draw[pr] (o1d) to[bend right=25] (o3d);
    \draw[pr] (o3d) -- (o4d);
    \draw[pr] (o3d) -- (o5d);
    \draw[cp] (o4d) -- (o2d);
    \draw[cq] (o5d) -- (o2d);
    \node[font=\footnotesize] at (1.55,1.75) {$\Pai_3$};
    \node[font=\footnotesize] at (1.55,-0.75) {no self-loop, 1 component};
  \end{scope}
\end{tikzpicture}
\caption{Pairings, for two terms of the path $v_1,v_3,v_4,v_2$ of Figure~\ref{fig:twoterms}. Left: the two
copies, on the observations $1,3,4,2$ (solid) and $1,3,5,2$ (dashed), share the observation $3$ and the
coordinates $j_1$ and $j_2$ of their first two links. At $j_1$, the four entries sit at $1,3$ and $1,3$, and at
$j_2$, at $3,4$ and $3,5$. Right: $\Gomi$ and the self-loops for three pairings, the pairs at $j_1,j_2$ in red.
$\Pai_1$ pairs each entry of one copy with the corresponding entry of the other copy, at $j_1$ and at $j_2$. The
pairing $\Pai_2$ does so at $j_1$ only and pairs the entries within each copy at $j_2$, and $\Pai_3$ pairs them
within each copy at both. The third pairing at $j_1$ or $j_2$ gives the same graph as the pairing within
each copy. Given $x=0$, the term of $\Pai_3$ vanishes, as $1$ and $2$ lie in one component of $\Gomi$.}
\label{fig:pairings}
\end{figure}

Let $G=(V,W,E)$ be connected, with every variable node a link, read as a multigraph on $V$ as in
Section~\ref{sec:moments}. Two labelings $\pi,\pi'$ place two copies of $G$ on the data, and their overlap is
recorded by the \emph{matching} $\mathbf M=(\mathbf M_V,\mathbf M_W)$ of the two copies (Section~\ref{sec:graphs}),
\[
  \mathbf M_V\ =\ \big\{(v,v'):\ \pi_1(v)=\pi'_1(v')\big\},\qquad
  \mathbf M_W\ =\ \big\{(w,w'):\ \pi_2(w)=\pi'_2(w')\big\} ,
\]
which always contains $(v_1,v_1)$ and $(v_2,v_2)$. The observations and coordinates of the matched nodes are
\emph{shared}. About $n^{2(|V|-2)-(|\mathbf M_V|-2)}d^{\,2|W|-|\mathbf M_W|}$ pairs of labelings follow a given
matching, against about $n^{2(|V|-2)}d^{\,2|W|}$ pairs in all. Given the groups, a coordinate that carries one link
of one copy (that is, not shared), at two observations $i\ne i'$, contributes $\E[Y_{ij}Y_{i'j}\mid k^*]$, and a shared coordinate carries
one link of each copy and contributes the sum \eqref{eq:wickfour} over the three ways of pairing its four entries.
Expanding the product over the coordinates, $\E[T_{G,\pi}T_{G,\pi'}\mid k^*]$ is a sum over the \emph{pairings}
$\Pai$, which choose one of these ways at each shared coordinate (Figure~\ref{fig:pairings}). In the term of $\Pai$,
a pair of entries $Y_{ij},Y_{i'j}$ contributes $\frac{\Delta^2}d\one\{i\sim i'\}$ if $i\ne i'$, an \emph{edge}
between $i$ and $i'$, and $1+\frac{\Delta^2}d$ if $i=i'$, a \emph{self-loop}. The self-loops form a set $E_\cap$, and
the edges a graph $\Gomi$ on the observations $\pi_1(V)\cup\pi'_1(V)$, whose components, an isolated observation
counting as one, must each lie in one group (Lemma~\ref{lem:mono}).

Write $\CC_{\Gomi}(\{1,2\})$ for the set of the
components of $\Gomi$ that hold neither $1$ nor $2$. Up to factors close to $1$, a matching and a pairing contribute
to $\E[T_G^2\mid x=1]$ the fraction
\begin{equation}\label{eq:Rratio}
  R\ :=\ \Big(\frac Kn\Big)^{|\mathbf M_V|-2}\,\frac1{d^{\,|\mathbf M_W|}}\,\Big(\frac d{\Delta^2}\Big)^{|E_\cap|}\,
  K^{\,|\CC_{\Gomi}(\{1,2\})|}
\end{equation}
of $\E[T_G\mid x=1]^2$, and to $\E[T_G^2\mid x=0]$ the same fraction if $1$ and $2$ lie in distinct components of
$\Gomi$, and $0$ otherwise. Each shared observation multiplies $R$ by $K/n$ and each shared coordinate by $1/d$. Each
self-loop multiplies it by $d/\Delta^2$, as it replaces a covariance $\Delta^2/d$ by a variance close to $1$, and each
component of $\Gomi$ holding neither $1$ nor $2$ by $K$, the free choice of its group. When only
$v_1$ and $v_2$ are matched, the pairs are the disjoint pairs of Section~\ref{sec:moments}: there is no self-loop,
$\Gomi$ is the union of the two copies, glued at $1$ and $2$, and $R=1$. Hence, by \eqref{eq:varsplit}, both
conditional variances are at most $\E[T_G\mid x=1]^2$ times the sum of $R$ over the other matchings and pairings, up
to factors close to $1$, and $T_G$ decides $x$ when this sum of $R$'s is small. In Section~\ref{part:ub}, $R$ is the factor
$K^{\mathcal E}$ of Proposition~\ref{prop:exponent}.

\subsubsection{The shape of the template}\label{sec:shapesketch}

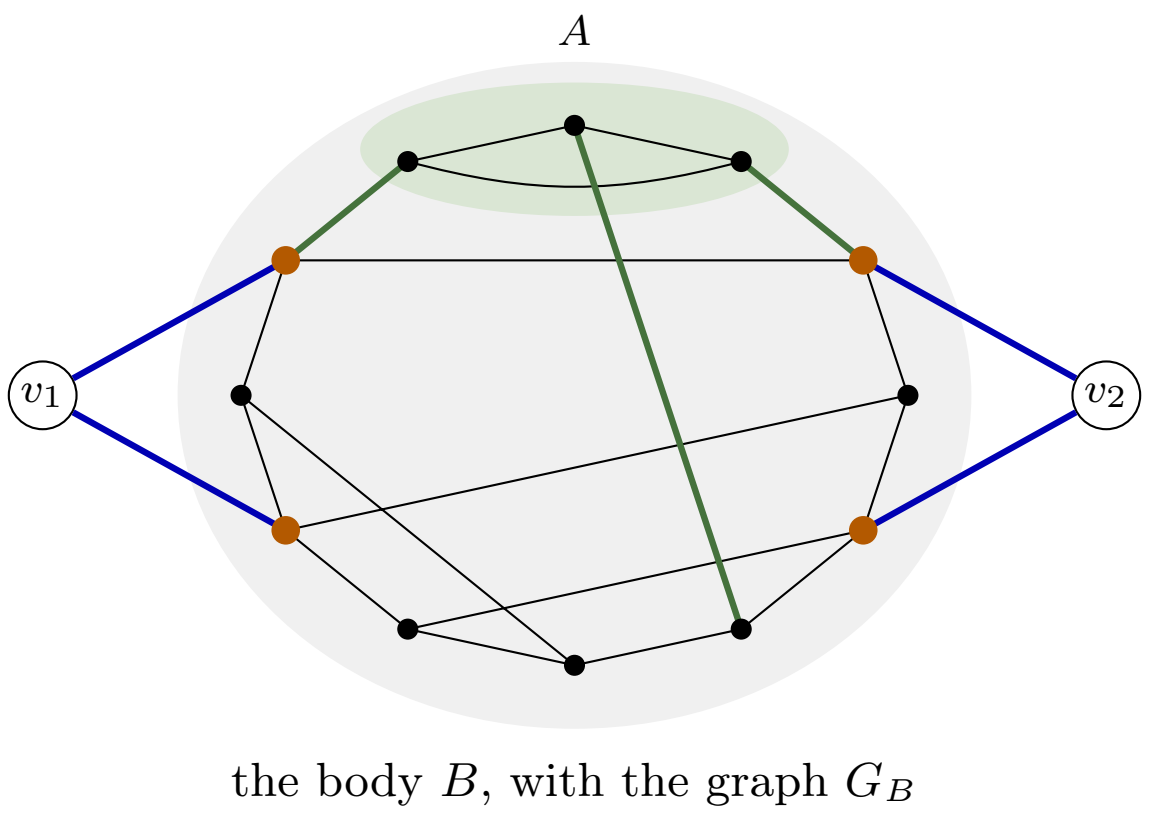
\begin{figure}[ht]
\centering
\begin{tikzpicture}[scale=1.05,
  ob/.style={circle,fill=black,inner sep=1.4pt},
  fs/.style={circle,fill=orange!70!black,inner sep=1.9pt},
  di/.style={circle,draw=black,fill=white,inner sep=1.2pt,font=\scriptsize},
  lab/.style={font=\scriptsize,inner sep=1.5pt},
  fl/.style={blue!70!black,very thick},
  bd/.style={OliveGreen!80!black,very thick}]
  \foreach \k in {0,...,11} \coordinate (n\k) at ({2.1*cos(30*\k)},{1.7*sin(30*\k)});
  \begin{scope}[on background layer]
    \fill[black!6] (0,0) ellipse (2.5 and 2.1);
    \fill[OliveGreen!15] (0,1.55) ellipse (1.35 and 0.42);
  \end{scope}
  \foreach \k in {0,...,11}{\pgfmathtruncatemacro{\l}{mod(\k+1,12)}\draw[black] (n\k)--(n\l);}
  \draw[black] (n0)--(n7) (n1)--(n5) (n6)--(n9) (n8)--(n11);
  \draw[black] (n2) to[bend left=15] (n4);
  \draw[bd] (n1)--(n2) (n4)--(n5) (n3)--(n10);
  \node[di] (v1) at (-3.35,0) {$v_1$};
  \node[di] (v2) at (3.35,0) {$v_2$};
  \draw[fl] (v1)--(n5) (v1)--(n7) (v2)--(n1) (v2)--(n11);
  \foreach \k in {0,2,3,4,6,8,9,10} \node[ob] at (n\k) {};
  \foreach \k in {1,5,7,11} \node[fs] at (n\k) {};
  \node[lab] at (0,2.3) {$A$};
  \node[font=\footnotesize] at (0,-2.45) {the body $B$, with the graph $G_B$};
\end{tikzpicture}
\caption{The shape of the template, at $q=3$. Each line is a link, that is, a variable node
joined to its two endpoints. The body $B$ (shaded) carries the graph $G_B$ (black), in which every
node has degree $3$. Here $G_B$ is a cycle together with a perfect matching, as in
Appendix~\ref{app:body}. The four fasteners (orange) are each joined by a fastener link (blue) to
$v_1$ or to $v_2$, so that $|W|=3|B|/2+|\Phi|=22$. The set $A$ (green) has
$\partial_B(A)=3$: three edges of $G_B$ (thick green) have exactly one endpoint
in $A$.}
\label{fig:shape}
\end{figure}

We look for $G$ of the following shape, for an integer $q\ge3$. The individual nodes other than $v_1,v_2$ form the
\emph{body} $B$, and the links between them form the subgraph $G_B$ of $G$ induced on $B$, a multigraph in which
every node has degree $q$. In addition, each node of a set $\Phi\subseteq B$ of \emph{fasteners} is joined by one link
to $v_1$ or to $v_2$, so that $|W|=q|B|/2+|\Phi|$ (Figure~\ref{fig:shape}). For $A\subseteq B$, let $\partial_B(A)$ be
the number of edges of $G_B$ with exactly one endpoint in $A$. In Figure~\ref{fig:shape}, $\partial_B(A)=3=q$, as
(T3) below requires.

\subsubsection{The matchings that matter}\label{sec:overlapsketch}

For each matching discussed below, we compute $R$ for the pairing that pairs together the two occurrences of each entry $Y_{ij}$ that appears in both copies. This pairing has the most self-loops and, in these matchings, the most
components of $\Gomi$, so that its $R$ is the largest among the pairings of the matching. We recall that we look for a template $G$ such that $R$ is small for all matchings and pairings.
 Each such condition
translates into a property of the template.

\noindent\uline{(a) The two extreme matchings.} They are those of $\mathcal V_{\mathrm{noise}}$ and $\mathcal
V_{\mathrm{centers}}$, the two extreme cases of Section~\ref{sec:intuition} (Figure~\ref{fig:twoterms}). When
$\pi'=\pi$, every pair is a self-loop and $\Gomi$ has no edge. When $\pi'_2=\pi_2$ and $\pi'_1(V)\cap\pi_1(V)=\{1,2\}$,
the pairing joins $\pi_1(v)$ to $\pi'_1(v)$ for each $v\in B$ and puts the self-loops at $1$ and $2$. By
\eqref{eq:Rratio}, these two matchings have
\[
  R\ \asymp\ \Big(\frac{K^2d^{\,q/2}}{n\Delta^{2q}}\Big)^{|B|}\Big(\frac d{\Delta^4}\Big)^{|\Phi|}
  \qquad\text{and}\qquad
  R\ \asymp\ \Big(\frac K{d^{\,q/2}}\Big)^{|B|}\,\Delta^{-2|\Phi|} ,
\]
the inverses of \eqref{eq:nodefactors} with $S=\varnothing$ and $S=B$, with the factors of the fastener links. Both are
small when $n\Delta^{2q}>K^2d^{q/2}$ and $d^{q/2}>K$, that is, $q>p$, which leads to $q=\lceil p\rceil$ when $p$ is not
an integer, and when the body is large compared with the number of fasteners.

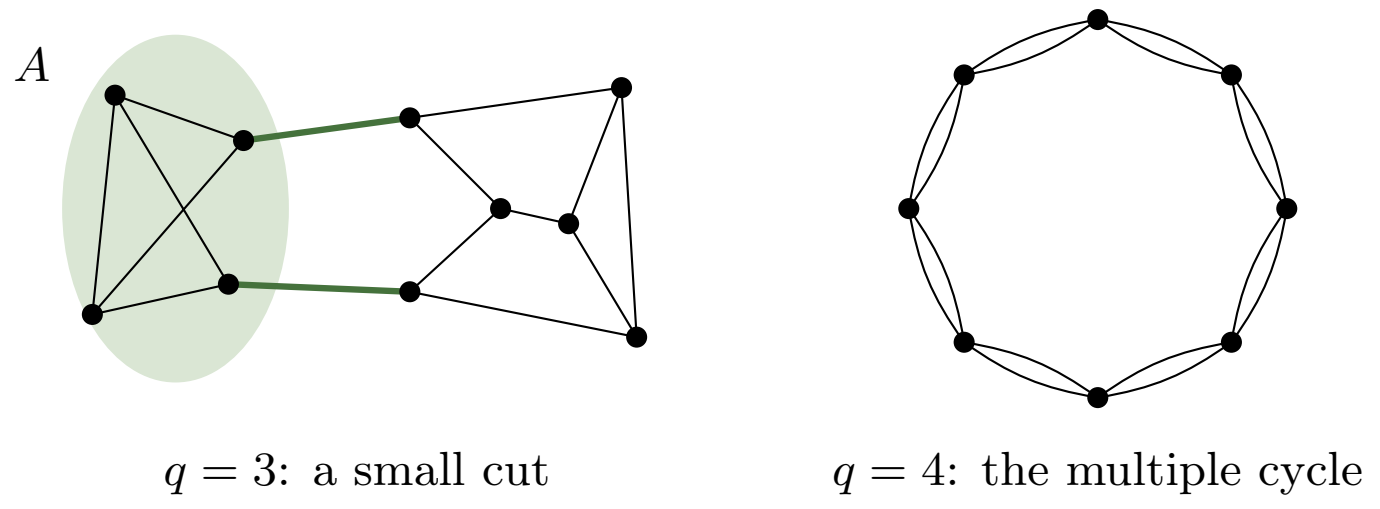
\begin{figure}[ht]
\centering
\begin{tikzpicture}[
  ob/.style={circle,fill=black,inner sep=1.4pt},
  lk/.style={black},
  bd/.style={OliveGreen!80!black,very thick}]
\begin{scope}[xshift=-1.6cm]
  \coordinate (a) at (0,0.75); \coordinate (b) at (-0.15,-0.7);
  \coordinate (c) at (0.85,0.45); \coordinate (e) at (0.75,-0.5);
  \coordinate (r4) at (1.95,0.6); \coordinate (r6) at (1.95,-0.55); \coordinate (r5) at (2.55,0);
  \coordinate (r1) at (3.35,0.8); \coordinate (r2) at (3.0,-0.1); \coordinate (r3) at (3.45,-0.85);
  \begin{scope}[on background layer]
    \fill[OliveGreen!15] (0.4,0) ellipse (0.75 and 1.15);
  \end{scope}
  \draw[lk] (a)--(b) (a)--(c) (a)--(e) (b)--(c) (b)--(e);
  \draw[lk] (r1)--(r2) (r1)--(r3) (r1)--(r4) (r2)--(r3) (r2)--(r5) (r3)--(r6) (r4)--(r5) (r5)--(r6);
  \draw[bd] (c)--(r4) (e)--(r6);
  \foreach \n in {a,b,c,e,r1,r2,r3,r4,r5,r6} \node[ob] at (\n) {};
  \node[font=\footnotesize] at (-0.55,0.95) {$A$};
  \node[font=\footnotesize] at (1.6,-1.75) {$q=3$: a small cut};
\end{scope}
\begin{scope}[xshift=4.9cm]
  \foreach \k in {0,...,7} \coordinate (b\k) at ({1.25*cos(45*\k)},{1.25*sin(45*\k)});
  \foreach \k in {0,...,7}{\pgfmathtruncatemacro{\l}{mod(\k+1,8)}
    \draw[lk] (b\k) to[bend left=13] (b\l); \draw[lk] (b\k) to[bend right=13] (b\l);}
  \foreach \k in {0,...,7} \node[ob] at (b\k) {};
  \node[font=\footnotesize] at (0,-1.75) {$q=4$: the multiple cycle};
\end{scope}
\end{tikzpicture}
\caption{Small cuts. Left: a graph in which every node has degree $3$, but only $2$ edges (green) leave the
set $A$ (shaded), a small cut which makes a part of the variance large (Section~\ref{sec:overlapsketch}). Right: at
$q=4$, the multiple cycle, every link doubled: every arc of the cycle has $4$ edges leaving it, so that this graph
has no small cut.}
\label{fig:smallcuts}
\end{figure}

\begin{figure}[ht]
\centering
\begin{tikzpicture}[scale=1.25,
  ob/.style={circle,fill=black,inner sep=1.4pt},
  di/.style={circle,draw=black,fill=white,inner sep=1.2pt,font=\scriptsize},
  lab/.style={font=\scriptsize,inner sep=1.5pt},
  own/.style={black},
  sh/.style={BrickRed,thick,densely dashed},
  pr/.style={BrickRed,thick},
  sl/.style={BrickRed,thick,min distance=4mm},
  fs/.style={circle,fill=orange!70!black,inner sep=1.9pt}]
  \begin{scope}
  \coordinate (La) at (0,1); \coordinate (Lb) at (0.75,0.5); \coordinate (Lc) at (0,0);
  \coordinate (Ld) at (2.25,1); \coordinate (Le) at (1.5,0.5); \coordinate (Lf) at (2.25,0);
  \begin{scope}[on background layer]
    \fill[OliveGreen!12,rounded corners=8pt] (-0.25,-0.25) rectangle (1,1.25);
    \fill[OliveGreen!12,rounded corners=8pt] (1.25,-0.25) rectangle (2.5,1.25);
  \end{scope}
  \node[lab] at (0.5,-0.45) {$A_1$};
  \node[lab] at (1.75,-0.45) {$A_2$};
  \node[di] (Lv1) at (-0.9,1) {$v_1$};
  \node[di] (Lv2) at (3.15,0) {$v_2$};
  \draw[own] (La)--(Lb)--(Lc)--(La) (Ld)--(Le)--(Lf)--(Ld);
  \draw[sh] (La)--(Ld) (Lb)--(Le) (Lc)--(Lf) (Lv1)--(La) (Lf)--(Lv2);
  \foreach \n/\m/\pos in {a/3/above left,b/4/left,c/5/below left,d/6/above right,e/7/right,f/8/below right}
    \node[ob,label={[lab]\pos:$v_{\m}$}] at (L\n) {};
  \foreach \n in {a,f} \node[fs] at (L\n) {};
  \node[font=\footnotesize] at (1.125,1.6) {the template $G$};
  \end{scope}
  \begin{scope}[xshift=5cm]
  \foreach \n/\x in {a/0,b/0.6,c/1.2,d/2.2,e/2.8,f/3.4}{
    \coordinate (T\n) at (\x,1); \coordinate (B\n) at (\x,0);}
  \begin{scope}[on background layer]
    \fill[OliveGreen!12,rounded corners=8pt] (-0.3,-0.5) rectangle (1.5,1.5);
    \fill[OliveGreen!12,rounded corners=8pt] (1.9,-0.5) rectangle (3.7,1.5);
  \end{scope}
  \draw[own] (Ta)--(Tb)--(Tc) (Ta) to[out=50,in=130] (Tc);
  \draw[own] (Ba)--(Bb)--(Bc) (Ba) to[out=-50,in=-130] (Bc);
  \draw[own] (Td)--(Te)--(Tf) (Td) to[out=50,in=130] (Tf);
  \draw[own] (Bd)--(Be)--(Bf) (Bd) to[out=-50,in=-130] (Bf);
  \foreach \n in {a,b,c,d,e,f} \draw[pr] (T\n)--(B\n);
  \draw[pr] (Ta) to[bend right=35] (Ba);
  \draw[pr] (Tf) to[bend left=35] (Bf);
  \foreach \n/\t/\b in {a/3/9,b/4/10,c/5/11,d/6/12,e/7/13}{
    \node[ob,label={[lab]below right:$\t$}] at (T\n) {};
    \node[ob,label={[lab]above right:$\b$}] at (B\n) {};}
  \node[ob,label={[lab]above right:$8$}] at (Tf) {};
  \node[ob,label={[lab]below right:$14$}] at (Bf) {};
  \node[di] (one) at (-0.95,0.5) {$1$};
  \node[di] (two) at (4.35,0.5) {$2$};
  \draw[sl] (one) edge[loop,out=120,in=60] (one);
  \draw[sl] (two) edge[loop,out=120,in=60] (two);
  \node[font=\footnotesize] at (1.7,1.85) {the graph $\Gomi$ and the self-loops};
  \end{scope}
\end{tikzpicture}

\vspace{4mm}
\begin{tikzpicture}[scale=1.25,
  ob/.style={circle,fill=black,inner sep=1.4pt},
  di/.style={circle,draw=black,fill=white,inner sep=1.2pt,font=\scriptsize},
  lab/.style={font=\scriptsize,inner sep=1.5pt},
  own/.style={black},
  sh/.style={BrickRed,thick,densely dashed},
  pr/.style={BrickRed,thick},
  sl/.style={BrickRed,thick,min distance=4mm},
  fs/.style={circle,fill=orange!70!black,inner sep=1.9pt}]
  \begin{scope}
  \coordinate (La) at (0,1); \coordinate (Lb) at (0.75,0.5); \coordinate (Lc) at (0,0);
  \coordinate (Ld) at (2.25,1); \coordinate (Le) at (1.5,0.5); \coordinate (Lf) at (2.25,0);
  \begin{scope}[on background layer]
    \fill[Plum!15] (La) circle (0.22);
    \fill[OliveGreen!12] plot[smooth cycle,tension=0.55] coordinates
      {(-0.28,-0.15) (0.75,0.78) (2.3,1.25) (2.5,0.5) (2.4,-0.25) (1,-0.27)};
  \end{scope}
  \node[lab] at (-0.25,1.3) {$U$};
  \node[lab] at (1.4,-0.5) {$A$};
  \node[di] (Lv1) at (-0.9,1) {$v_1$};
  \node[di] (Lv2) at (3.15,0) {$v_2$};
  \draw[own] (Lb)--(Lc) (Lb)--(Le) (Lc)--(Lf) (Ld)--(Le)--(Lf)--(Ld);
  \draw[sh] (La)--(Lb) (La)--(Lc) (La)--(Ld) (Lv1)--(La) (Lf)--(Lv2);
  \foreach \n/\m/\pos in {a/3/above,b/4/below,c/5/below left,d/6/above right,e/7/right,f/8/below right}
    \node[ob,label={[lab]\pos:$v_{\m}$}] at (L\n) {};
  \foreach \n in {a,f} \node[fs] at (L\n) {};
  \node[font=\footnotesize] at (1.125,1.6) {the template $G$};
  \end{scope}
  \begin{scope}[xshift=5.5cm]
  \foreach \n/\x in {c/0,b/0.6,e/1.2,d/1.8,f/2.4}{
    \coordinate (T\n) at (\x,1); \coordinate (B\n) at (\x,0);}
  \begin{scope}[on background layer]
    \fill[OliveGreen!12,rounded corners=8pt] (-0.3,-0.55) rectangle (2.7,1.55);
  \end{scope}
  \draw[own] (Tc)--(Tb)--(Te)--(Td)--(Tf) (Tc) to[out=40,in=140] (Tf) (Te) to[out=45,in=135] (Tf);
  \draw[own] (Bc)--(Bb)--(Be)--(Bd)--(Bf) (Bc) to[out=-40,in=-140] (Bf) (Be) to[out=-45,in=-135] (Bf);
  \foreach \n in {b,c,d,f} \draw[pr] (T\n)--(B\n);
  \foreach \n/\t/\b in {c/5/11,b/4/10,e/7/13,d/6/12,f/8/14}{
    \node[ob,label={[lab]below right:$\t$}] at (T\n) {};
    \node[ob,label={[lab]above right:$\b$}] at (B\n) {};}
  \node[di] (sa) at (-0.95,0.5) {$3$};
  \foreach \t in {45,135,225,315} \draw[sl] (sa) edge[loop,out=\t+22,in=\t-22] (sa);
  \node[di] (one) at (-2.05,0.5) {$1$};
  \node[di] (two) at (3.3,0.5) {$2$};
  \draw[sl] (one) edge[loop,out=120,in=60] (one);
  \draw[sl] (two) edge[loop,out=120,in=60] (two);
  \node[font=\footnotesize] at (0.6,1.85) {the graph $\Gomi$ and the self-loops};
  \end{scope}
\end{tikzpicture}
\caption{Free components, at $q=3$. Top row: the matching of (b), with the partition $B=A_1\cup A_2$ and the
fasteners $v_3,v_8$ (orange). The red dashed links are placed on shared coordinates, the black ones on distinct
coordinates in the two copies. On the right, $\Gomi$ and the self-loops, when $\pi$ places each node $v_i$ of the
body on the observation $i$ and $\pi'$ places it on $i+6$: the two shaded components are free, and
$R=\Delta^{-4}\big(K/d^{\,3/2}\big)^2$, as $\partial_B(A_1)=\partial_B(A_2)=3$. Bottom row: the matching of (c), with
$U=\{v_3\}$ and $A=B\setminus U$, $\pi'$ placing $v_3$ on $3$: the shared observations $1$, $2$ and $3$ are components
on their own, with four self-loops at $3$, the shaded component is free, and
$R=K\Delta^{-4}\big(K^2d^{\,3/2}/(n\Delta^6)\big)(d/\Delta^2)\,d^{-3/2}$.}
\label{fig:free}
\end{figure}
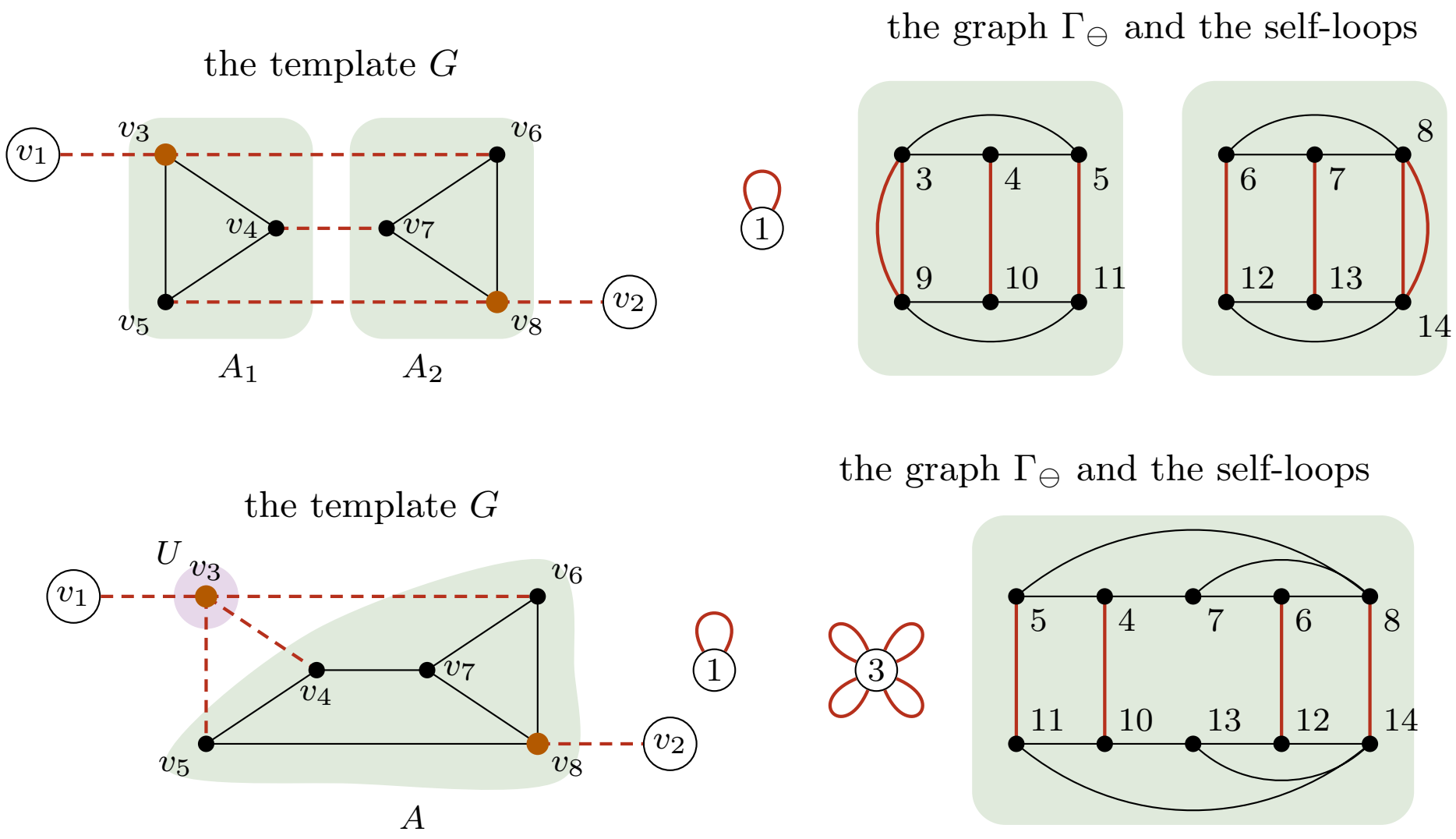

\noindent\uline{(b) Free components and small cuts.} A component of $\Gomi$ is \emph{free} when it holds no shared
observation. In the matching of $\mathcal V_{\mathrm{centers}}$, each of the $|B|$ components
$\{\pi_1(v),\pi'_1(v)\}$ is free, and its factor $K$ is compensated by the $q$ shared coordinates of the links at $v$.
Free components may be larger. Partition $B$ into sets $A_i$, each inducing a connected subgraph of $G_B$, let
$\pi'_1(V)\cap\pi_1(V)=\{1,2\}$ and $\pi'_2(w)=\pi_2(w)$ on the fastener links and on the links between distinct sets $A_i$, and let
the two copies place the links inside the sets on distinct coordinates, that is, the $\pi_2(w)$ and $\pi'_2(w)$ are all distinct for $w$ corresponding to a link within an $A_i$.  Then each $\pi_1(A_i)\cup\pi'_1(A_i)$ is a
free component (Figure~\ref{fig:free}, top), and
\[
  R\ =\ \Delta^{-2|\Phi|}\prod_i\frac K{d^{\,\partial_B(A_i)/2}} .
\]
When $q=\lceil p\rceil$, $d^{(q-1)/2}<K$, so that a
set $A_i$ that fewer than $q$ edges leave, a \emph{small cut}, makes $R$ larger by a power of $K$. In
Figure~\ref{fig:smallcuts} (left), every node has degree $3$, but only $2$ edges leave the set $A$. This leads to the first isoperimetric property of the body,
$\partial_B(A)\ge q$ for every $\varnothing\ne A\subsetneq B$.

\noindent\uline{(c) Large free components and the fasteners.} When $\mathbf M_V=\{(v_1,v_1),(v_2,v_2)\}$ (no other individual node is matched) and $|\mathbf M_W|=|\Phi|$, the matched links being the links between the fasteners and $v_1$ or $v_2$, we get
$R=K\Delta^{-2|\Phi|}$: a free component holding the whole body is
compensated by the fastener links alone, each bringing $1/\Delta^2$.
So $|\Phi|$ must be large compared with $\log
K/\log\Delta^2$. 

The fasteners cannot be gathered in a small region of $G_B$ either. Let the two copies also share the observations of a set $U\subsetneq B$ such that
$A:=B\setminus U$ induces a connected subgraph of $G_B$: $\pi'_1=\pi_1$ on $U\cup\{v_1,v_2\}$,
$\pi'_1(A)\cap\pi_1(V)=\varnothing$, and $\pi'_2=\pi_2$ on all links but those inside $A$. Then
$\pi_1(A)\cup\pi'_1(A)$ is a free component, every shared observation is a component on its own
(Figure~\ref{fig:free}, bottom), and
\[
  R\ =\ \frac K{\Delta^{2|\Phi|}}\Big(\frac{K^2d^{\,q/2}}{n\Delta^{2q}}\Big)^{|U|}\Big(\frac
  d{\Delta^2}\Big)^{|U\cap\Phi|}d^{\,-\partial_B(U)/2} .
\]
Each shared observation of the body brings the factor of $\mathcal V_{\mathrm{noise}}$, smaller than $1$ above the
threshold, each shared fastener the large factor $d/\Delta^2$, as its link to $v_1$ or $v_2$ now carries two
self-loops, and each edge of $G_B$ leaving $U$ the factor $1/\sqrt d$. For $R$ to be small, two edges should leave
$U$ for each fastener in $U$: $\partial_B(U)\ge2|U\cap\Phi|$ whenever $|U|<|B|/2$, a larger $U$ leaving a free
component of at most half of the body, which (b) handles. This fails when many fasteners gather in a region with few
edges leaving it, as on an arc of a multiple cycle. By Lemma~\ref{lem:fast}, it holds when the fasteners are far
apart from one another in $G_B$ and $G_B$ satisfies the linear isoperimetric inequality
$\partial_B(A)\ge\ciso\min\{|A|,|B|-|A|\}$, for a constant $\ciso>0$: a connected set holding at least two
fasteners then has more than $2/\ciso$ nodes per fastener, hence more than two edges leaving it per fastener.

\paragraph{The template.} The matchings above ask the following of the template.
\begin{enumerate}[label=\textup{(T\arabic*)},leftmargin=2.8em,itemsep=1pt,topsep=3pt]
\item \emph{Degree.} Every node of the body has degree $q=\lceil p\rceil$ in $G_B$.
\item \emph{Size of the body.} The body is large compared with the number of fasteners, so that the factor
$(K^2d^{q/2}/(n\Delta^{2q}))^{|B|}$ of (a) beats the factor $(d/\Delta^4)^{|\Phi|}$.
\item \emph{No small cut.} $\partial_B(A)\ge q$ for every $\varnothing\ne A\subsetneq B$.
\item \emph{Many fasteners.} $|\Phi|$ is large compared with $\log K/\log\Delta^2$.
\item \emph{Spread fasteners.} $\partial_B(U)\ge2|U\cap\Phi|$ whenever $|U|<|B|/2$, which holds when the fasteners
are far apart in $G_B$ and $G_B$ satisfies the linear isoperimetric inequality (Lemma~\ref{lem:fast}).
\end{enumerate}
The template of Definition~\ref{def:template} has this shape, with $2L$ fasteners, $L$ of them joined to $v_1$ and $L$
to $v_2$, and a body of $2L^2$ nodes: a cycle through $B$ together with $q-2$ random perfect matchings of $B$, which
satisfies (T3) and the linear isoperimetric inequality (Appendix~\ref{app:body}). Section~\ref{part:ub} shows that
these properties suffice: the factor $K$ of each free component is compensated by the links leaving the set of its
nodes in one copy, through (T3) and the linear isoperimetric inequality when this set holds at most half of the body,
and through the fasteners otherwise (Lemmas~\ref{lem:small} and~\ref{lem:large}).

\section{Individuals, variables, and graphs}\label{sec:graphs}

This section fixes the graph notation shared by Sections~\ref{part:lb} and~\ref{part:ub}.

$G=(V,W,E)$ is a bipartite \textbf{multi}graph, its three sets written $V_G,W_G,E_G$ when the
graph has to be named:
\begin{itemize}
\item \textbf{Individual and variable nodes:} $V$ is the set of \textbf{individual nodes}, $W$ the set of
\textbf{variable nodes}, \emph{individuals} and \emph{variables} for short: a
\textbf{labeling} \eqref{eq:embed} sends them to individuals and to variables, and sends two
\textbf{distinguished} nodes $v_1\ne v_2$ to the individuals $1$ and $2$.
 No node is isolated, except possibly $v_1$ or $v_2$.
\item \textbf{Edges, multiplicities and degrees:} $m_{v,w}$ is the multiplicity of the pair
$(v,w)$, and the \textbf{edges} of $G$ --- the elements of $E$ --- are the triples $(v,\iota,w)$
with $1\le\iota\le m_{v,w}$, so that $|E|=\sum_{v,w}m_{v,w}$. \emph{The $m_{v,w}$ edges joining
one $v$ to one $w$ are distinct}, the index $\iota$ telling them apart. The degrees are
$\deg_G(v)=\sum_wm_{v,w}$ and $\deg_G(w)=\sum_vm_{v,w}$, and $N(w)=\{v:m_{v,w}\ge1\}$.
\end{itemize}
An \textbf{isomorphism} $G^{(1)}\to G^{(2)}$ is a pair of bijections $V^{(1)}\to V^{(2)}$, $W^{(1)}\to W^{(2)}$ preserving the multiplicities and fixing $v_1$ and $v_2$, and
$\Aut(G)$ is the group of automorphisms of $G$. Every isomorphism in this paper is of this kind.

Let $G^{(1)}=(V^{(1)},W^{(1)},E^{(1)})$ and $G^{(2)}=(V^{(2)},W^{(2)},E^{(2)})$ be two graphs. A
\textbf{matching} $\mathbf M$ of $G^{(1)}$ with $G^{(2)}$ is a pair of sets
$\mathbf M_V\subseteq V^{(1)}\times V^{(2)}$ and $\mathbf M_W\subseteq W^{(1)}\times W^{(2)}$ in
which every node occurs at most once. The two nodes of a pair of $\mathbf M$ are matched,
and a node lying in no pair is unmatched. Write $\mathcal M(G^{(1)},G^{(2)})$ for the set
of those matchings that fix $v_1$ and $v_2$, i.e.\ have $(v_1,v_1)$ and $(v_2,v_2)$ in
$\mathbf M_V$.

\paragraph{Labelings.} A \textbf{labeling} of $G$ is a pair $\pi=(\pi_1,\pi_2)$ of injections,
\begin{equation}\label{eq:embed}
\pi_1:V\to[n]\ \text{with}\ \pi_1(v_1)=1,\ \pi_1(v_2)=2,\qquad \pi_2:W\to[d] ,
\end{equation}
 the set of them is $\Pi_G$. As $\pi_1$ sends $v_1,v_2$ to $1,2$ and the other $|V|-2$ individuals
injectively into $[n]\setminus\{1,2\}$, and $\pi_2$ sends $W$ injectively into $[d]$,
\begin{equation}\label{eq:Pisize}
  |\Pi_G|\ =\ \prod_{k=2}^{|V|-1}(n-k)\ \cdot\ \prod_{k=0}^{|W|-1}(d-k) .
\end{equation}
Let $\mathbf M$ be a matching of $G$ with a second graph and let $\pi'=(\pi'_1,\pi'_2)$ place,
injectively, nodes of that graph, among them every node $\mathbf M$ matches. The labelings of
$G$ \textbf{following} $\mathbf M$ are those putting each matched node of $G$ where $\pi'$ puts
its partner and no node of $G$ where $\pi'$ puts a node that is not its partner:
\begin{equation}\label{eq:PiM}
\begin{aligned}
\Pi_G(\mathbf M,\pi')\ :=\ \big\{\pi\in\Pi_G\ :\
  &\mathbf M_V=\big\{(v,v'):\pi_1(v)=\pi'_1(v')\big\}\\
  &\text{and}\ \ \mathbf M_W=\big\{(w,w'):\pi_2(w)=\pi'_2(w')\big\}\big\} ,
\end{aligned}
\end{equation}
$v'$ and $w'$ running over the nodes $\pi'$ places.

\paragraph{Links, self-loops.} A variable node $w$ of degree $2$ is a \textbf{link} when
$N(w)=\{v,v'\}$ with $v\ne v'$, and is then written $w=\{v,v'\}$. It is a \textbf{self-loop} when
$N(w)=\{v\}$.

\paragraph{Pairings.} The edges of $G$ at a variable $w$ are the $\deg_G(w)$
elements $(v,\iota,w)$ of $E$ with $v\in N(w)$. A \textbf{pairing of $w$} is a partition of
those edges into pairs and a \textbf{pairing of $G$} a family $\Pai=(\Pai_w)_{w\in W}$ of pairings $\Pai_w$ of the variables $w$. We identify $\Pai$ with the set of all its pairs $\{e,e'\}$, and $v_e$ denotes the individual at which an edge $e$ sits. The
graph $G[\Pai]$ has the same individual set $V$ and one variable per pair of $\Pai$, joined to
the individuals of the two edges of that pair, so that
\[
|V_{G[\Pai]}|=|V|,\quad |E_{G[\Pai]}|=|E|,\quad |W_{G[\Pai]}|=|E|/2,\quad
\deg_{G[\Pai]}(v)=\deg_G(v)\ \ \forall v .
\]
A pair whose two edges sit at one individual makes a self-loop of $G[\Pai]$ and every other a
link. Write $\PcalSL(G)$ for the set of all pairings of $G$, and $\Pcal(G)\subseteq\PcalSL(G)$ for the set
of those with no self-loop.

\paragraph{Components.} Two nodes, of either kind, lie in the same connected component of a
graph when a path of that graph joins them. For a graph $G$, $\CC_{G}$ is the set of its
connected components, \emph{every component counting, an isolated node included}, and $V_{C}$
is the set of individuals of a component $C$. So $G$ is \textbf{connected} when $|\CC_{G}|=1$.

\paragraph{Merging two graphs.} The \textbf{merged graph}
$\Gcup(G^{(1)},G^{(2)},\mathbf M)$ is the disjoint union
$G^{(1)}\sqcup G^{(2)}$ with the two nodes of each pair of $\mathbf M$ identified. A node of
$\Gcup$ so obtained, carrying one node of $G^{(1)}$ and one of $G^{(2)}$, is again called
\textbf{matched}. The same construction is applied to subgraphs of $G^{(1)}$ and of $G^{(2)}$,
a pair with a node outside them identifying nothing. For a graph $G$  and a
set $\bar V$ of its individuals, write
\[
  \CC_{G}(\bar V)\ :=\ \big\{C\in\CC_{G}\ :\ V_{C}\cap\bar V=\varnothing\big\}
\]
for the components of $G$ that miss $\bar V$. These are the \textbf{free} components and the
others the \textbf{pinned} ones.  $\Gcup[\Pai]$  stands for $\Gcup(G^{(1)},G^{(2)},\mathbf M)[\Pai]$.

\section{Proof of the lower bound}\label{part:lb}

The proof of Theorem~\ref{thm:LB} compares the polynomials of degree at most $D$ with an
orthonormal family of polynomials of the latent variables, and shows that the matrix of their inner
products, once rescaled, is diagonally dominant. It has five steps.

\noindent\uline{Step 1: conditional moments.} Given $k^*$, the columns of $X$ and $\epsilon$ are
independent Gaussian vectors, so $\E[P_{G,\pi}\mid k^*]$ follows from Wick's formula: it is a sum
over the pairings of the edges at each variable node in which no two edges at one individual are paired (Propositions~\ref{prop:master}
and~\ref{prop:pairid}). A pairing turns $G$ into a graph whose variable nodes are links, and
contributes the indicator that $k^*$ is constant along its components
(Corollary~\ref{cor:selfloop}). Sections~\ref{sec:moments} and~\ref{sec:othersketch} sketch this step.

\noindent\uline{Step 2: centering.} The polynomials $\Psh_G$ of Definition~\ref{def:center} subtract
from $P_G$, recursively, the polynomials of the graphs obtained by deleting the deletable sets of
individual nodes, the sets whose edges can be paired among themselves (Definition~\ref{def:delete}),
such as the individual nodes that a pairing separates from $v_1,v_2$. The coefficients are such
that the covariance with $x$ only counts the pairings $\Pai$ in $\Pcal(G)$ for which $G[\Pai]$ is connected:
$\E[(x-\E x)\Psh_G]=\hat s_G(1-\frac1K)$ times their number (Proposition~\ref{prop:num}), $\hat s_G$
being the explicit number \eqref{eq:shat}. Section~\ref{sec:othersketch} sketches this step.

\noindent\uline{Step 3: signal and noise individuals, and the blocks.} Each Hermite factor of an
entry $Y_{ij}=\frac{\Delta}{\sqrt d}X_{k^*(i),j}+\epsilon_{ij}$ splits into products of Hermite
polynomials of the center coordinate and of the noise (Lemma~\ref{lem:split}). An individual node of
$G$ other than $v_1,v_2$ is called signal when its degree is smaller than $p$ and noise otherwise
\eqref{eq:blocksplit}. A block of $G$ (Definition~\ref{def:block}) is the product of the Hermite
polynomials of $\epsilon$ at the noise individual nodes and of $X$ at the signal ones, the latter
being sent injectively to centers $\sigma(v)\in[K]$. Every block has norm $1$, and the normalized
sums $\bar Q_H$ of the blocks of the graphs $H$ form an orthonormal family (Lemma~\ref{lem:ortho}). Sections~\ref{sec:blocktypes} and~\ref{sec:blocksketch} sketch this step.

\noindent\uline{Step 4: the projections on the blocks.} Expanding each Hermite factor of $\Psh_G$ by
Lemma~\ref{lem:split}, the inner products $\langle\Psh_G,\bar Q_H\rangle$ are explicit sums over the
matchings of the noise part of $G$ with that of $H$ and over the pairings of the merged graph
(Lemma~\ref{lem:leak}). The centering of Step~2 cancels exactly the terms in which a component of
the merged graph consists of individual nodes of $G$ all of whose factors are expanded on the center
coordinates (Theorem~\ref{thm:cancel}). With the weights $\eta_G:=\hat s_GD^{50(|E_G|+|V_G|)}$, the matrix
$A_{H,G}:=\frac{\eta_H}{\eta_G}\langle\Psh_G,\bar Q_H\rangle$ has non-negative entries, a diagonal at
least $\frac12$ and off-diagonal column sums at most $\frac15$ (Corollary~\ref{cor:diag} and
Lemma~\ref{lem:colsum}). Hence $\Corr_{\le D}\le\frac{10}3(\sum_G\eta_G^2)^{1/2}$
(Proposition~\ref{prop:reduce}). Sections~\ref{sec:combsketch} and~\ref{sec:othersketch} sketch this step.

\noindent\uline{Step 5: the critical order.} By \eqref{eq:shat}, as $|W|\le|E|/2$ and
$\Delta\le d^{1/4}$,
\[
  \hat s_G\ \le\ \frac1K\prod_{\substack{v\ \text{signal}}}\Big(\frac{d^{\deg_G(v)/2}}{K}\Big)^{1/2}
  \prod_{\substack{v\ \text{noise}\\ v\ne v_1,v_2}}\Big(\frac{\Delta^{2\deg_G(v)}\,n}{K^{2}d^{\deg_G(v)/2}}\Big)^{1/2} .
\]
A signal node has degree at most $\lceil p\rceil-1$, so its factor is at most
$(d^{(\lceil p\rceil-1)/2}/K)^{1/2}$, small by Assumption~\ref{ass:A}\textup{(iii)}. A noise node has
degree at least $\lceil p\rceil$, so its factor is the factor of degree $\lceil p\rceil$, small by
Assumption~\ref{ass:A}\textup{(ii)}, times $(\Delta/d^{1/4})^{\deg_G(v)-\lceil p\rceil}$, at most $1$
by Assumption~\ref{ass:A}\textup{(i)}. This gives $\eta_G\le K^{-1}D^{-2|E_G|}$
(Lemma~\ref{lem:master}), and summing over the graphs gives Theorem~\ref{thm:LB}. Section~\ref{sec:blocksketch} sketches this step.

Assumption~\ref{ass:A} is in force throughout Section~\ref{part:lb}. Here $x$ is the same-cluster indicator
\eqref{eq:x}, so $\E x=\tfrac1K$, and $\Corr_{\le D}$ is \eqref{eq:corr}.

\subsection{Graphs and the Hermite basis}\label{sec:hermite}

\subsubsection{The normalized Hermite basis}\label{sec:sym}

\begin{definition}[Hermite basis]\label{def:hermite}
$\He_m$ is the probabilists' Hermite polynomial of degree $m$, see~\cite[\S18.3, Table~18.3.1]{DLMF}.
Hence the products
\[
  \prod_{i,j}\frac{1}{\sqrt{m_{i,j}!}}\,
  \He_{m_{i,j}}\!\Big(\frac{Y_{i,j}}{\sqrt{1+\Delta^2/d}}\Big),
  \qquad \sum_{i,j}m_{i,j}\ \le\ D ,
\]
form a basis of the polynomials of degree at most $D$ in the $nd$ variables $Y$. Each of
them is indexed by one graph --- the one with multiplicities $m_{i,j}$ --- together with one
labeling of it, and is, up to the factor $(1+\Delta^2/d)^{|E_G|/2}$, the product of the Hermite factors of \eqref{eq:PG}.
\end{definition}

\begin{lemma}[{Splitting property of Hermite polynomials, see~\cite[(18.18.11)]{DLMF}}]\label{lem:split}
For all $i\in[n]$, $j\in[d]$ and integers $m\ge0$, the sum running over the integers $\mathbf s,\mathbf n\ge0$ with $\mathbf s+\mathbf n=m$,
\[
\frac{(1+\Delta^2/d)^{m/2}}{\sqrt{m!}}\He_m\!\Big(\frac{Y_{ij}}{\sqrt{1+\Delta^2/d}}\Big)
\;=\sum_{\mathbf s+\mathbf n=m}\sqrt{\binom m{\mathbf s}}\,\Big(\frac{\Delta}{\sqrt d}\Big)^{\!\mathbf s}
\frac{\He_{\mathbf s}(X_{k^*(i),j})}{\sqrt{\mathbf s!}}\,\frac{\He_{\mathbf n}(\epsilon_{ij})}{\sqrt{\mathbf n!}} .
\]
\end{lemma}
Proved in \S\ref{sss:split}.

\paragraph{Signal and noise.} A
\textbf{splitting} of $G$ is an array
$\mathbf s=(\mathbf s_{v,w})$, $0\le\mathbf s_{v,w}\le m_{v,w}$, making $\mathbf s_{v,w}$ of the
$m_{v,w}$ edges joining $v$ to $w$ \textbf{signal} and $\mathbf n_{v,w}:=m_{v,w}-\mathbf s_{v,w}$ \textbf{noise}. It cuts $G$ into
$G_{\mathbf s}$ and $G_{\mathbf n}$, of multiplicities
$\mathbf s$ and $\mathbf n$, so that $|E|=|E_{G\mathbf s}|+|E_{G\mathbf n}|$. The nodes of $G_{\mathbf s}$ are those carrying a signal edge, and those of $G_{\mathbf n}$
are those carrying a noise edge \textbf{together with $v_1,v_2$}. Let $V_{G\mathbf s}$ and $V_{G\mathbf n}$ be their sets of individuals, so that
\[
V\ =\ V_{G\mathbf s}\cup V_{G\mathbf n},\qquad\text{the \textbf{mixed} individuals being those of }
V_{G\mathbf s}\cap V_{G\mathbf n}.
\]
An individual of $V_{G\mathbf s}\setminus V_{G\mathbf n}$ --- all of whose edges are
signal --- is \textbf{purely signal}.

One specific splitting is used throughout:
\begin{equation}\label{eq:blocksplit}
V_{G\mathbf n^\dagger}\;:=\;\{v_1,v_2\}\ \cup\ \{v:\deg_G(v)\ge p\},\qquad
V_{G\mathbf s^\dagger}\;:=\;V\setminus V_{G\mathbf n^\dagger},
\end{equation}
an edge being signal exactly when its individual lies in $V_{G\mathbf s^\dagger}$: \emph{$v_1$
and $v_2$ are always noise}, and any other individual is signal exactly when its degree is at
most $\lceil p\rceil-1$, and noise exactly when it is at least $\lceil p\rceil$. It has no mixed
individual, so 
$|V|=|V_{G\mathbf s^\dagger}|+|V_{G\mathbf n^\dagger}|$, and, degrees
being preserved by isomorphisms, it is preserved by every automorphism of $G$. It is the
\textbf{block splitting}, the one the blocks of \S\ref{sec:blocks} are built from. As a splitting it is written
\begin{equation}\label{eq:sdagger}
  \mathbf s^\dagger(G)_{v,w}\ :=\ m_{v,w}\,\one\{v\in V_{G\mathbf s^\dagger}\},\qquad
  \mathbf n^\dagger(G)\ :=\ m-\mathbf s^\dagger(G) .
\end{equation}
\emph{The dagger always marks the block splitting, and always that of the graph it is attached
to}: $G_{\mathbf s^\dagger},G_{\mathbf n^\dagger}$ abbreviate the two parts of
$\mathbf s^\dagger(G)$ --- with individual sets $V_{G\mathbf s^\dagger},V_{G\mathbf n^\dagger}$
and edge sets $E_{G\mathbf s^\dagger},E_{G\mathbf n^\dagger}$ --- while
$H_{\mathbf s^\dagger},H_{\mathbf n^\dagger}$ are the two parts of $\mathbf s^\dagger(H)$, and so
on. An undaggered $G_{\mathbf s},G_{\mathbf n}$ means the parts of whatever splitting is in force
(\S\ref{sec:expansion}), which may have mixed individuals.

\paragraph{The number $\hat s_G$.} To each graph attach the explicit number
\begin{equation}\label{eq:shat}
\hat s_G\ :=\ \frac{1}{\sqrt{|\Aut(G)|}}\Big(\frac{\Delta}{\sqrt d}\Big)^{\!|E_{G\mathbf n^\dagger}|}
d^{\,|W|/2}\,n^{\frac{|V_{G\mathbf n^\dagger}|-2}2}\,
K^{\,1-|V_{G\mathbf n^\dagger}|-\frac{|V_{G\mathbf s^\dagger}|}2}\Big(\prod_{v,w}m_{v,w}!\Big)^{-1/2} ,
\end{equation}
\emph{the splitting being the block one \eqref{eq:sdagger} and no other}: $\hat s_G$ is a
function of $G$ alone.

\begin{definition}\label{def:PG}
For $\pi\in\Pi_G$ set
\begin{equation}\label{eq:PG}
\begin{aligned}
P_{G,\pi}(Y)\ :=\ &\hat s_G\,\frac{K^{\,|V_G|-1}}{|\Pi_G|}\Big(\frac{\sqrt d}{\Delta}\Big)^{\!|E_G|}
\Big(\prod_{v,w}m_{v,w}!\Big)^{\!1/2}\\[2pt]
&\times\ \prod_{(v,w):\,m_{v,w}\ge1}\frac{(1+\Delta^2/d)^{m_{v,w}/2}}{\sqrt{m_{v,w}!}}\,
\He_{m_{v,w}}\!\Big(\frac{Y_{\pi_1(v),\pi_2(w)}}{\sqrt{1+\Delta^2/d}}\Big),
\end{aligned}
\end{equation}
and $P_G:=\sum_{\pi\in\Pi_G}P_{G,\pi}$.
\end{definition}

\begin{proposition}[Symmetry]\label{prop:sym}
The supremum \eqref{eq:corr} is unchanged if $P$ is restricted to
\[
  \Span\big(1,\ (P_G)_G\big)
  \qquad\textup{(Definition~\ref{def:PG})} ,
\]
$G$ running over the graphs with $2\le|E_G|\le D$ and every $\deg_G(w)$ \textbf{even}.
\end{proposition}
This is the analogue of \cite[Lemmas~1 and~4]{carpentier2026low}. Proved in \S\ref{sss:sym}. 

\emph{$\mathcal G_{\le D}$ denotes, from here to the end of Section~\ref{part:lb}, a set of representatives of
the isomorphism classes of those graphs}, one per class. Every object attached to a graph below is
defined from the graph alone, so an isomorphism carries it to the corresponding object of the image graph, and
a graph isomorphic to a member of $\mathcal G_{\le D}$ is identified with that member. The variable degrees of these graphs being even and positive,
\begin{equation}\label{eq:evencons}
|E|\ \text{is even},\qquad |W|\ \le\ |E|/2.
\end{equation}

\subsection{The conditional moment, and the pairing identity}\label{sec:engine}

\subsubsection{Wick, and the conditional moment}

For a multiset $\{a_1,\dots,a_r\}$ of positive integers put
$\mu(a_1,\dots,a_r):=\E\big[\prod_{i=1}^r\He_{a_i}(Z)\big]$ with $Z\sim\mathcal N(0,1)$, and
$\mu(\varnothing):=1$.

\begin{lemma}[{Wick's formula, see~\cite[Chapter~3]{janson1997gaussian}}]\label{lem:wick}
Give the $i$th factor $a_i$ points, and call a \textbf{Wick pairing} a partition of those
$\sum_ia_i$ points into pairs, no pair lying inside a single factor. Then $\mu(a_1,\dots,a_r)$ is
the number of Wick pairings. Consequently:
\begin{enumerate}[label=\textup{(\roman*)},leftmargin=*,itemsep=1pt,topsep=3pt]
\item $\mu(a)=0\ \text{ for every }a\ge1$ --- with one factor every pair would lie inside it.
\item $\mu(a,a')=a!\,\one\{a=a'\}$ --- with two factors every point must be paired into the other one.
\end{enumerate}
\end{lemma}
Proved in \S\ref{sss:wick}.

\begin{proposition}\label{prop:master}
For a map $\chi:V\to[K]$, put
$f_G(\chi):=\prod_{w\in W}\prod_{k\in[K]}\mu\big(\{m_{v,w}:v\in N(w),\ \chi(v)=k\}\big)$, each
$\{m_{v,w}:v\in N(w),\ \chi(v)=k\}$ being the multiset with one element $m_{v,w}$ for each $v\in N(w)$ with
$\chi(v)=k$. Then
\[
\E\big[P_{G,\pi}\mid k^*\big]\ =\ \hat s_G\,\frac{K^{\,|V_G|-1}}{|\Pi_G|}\ f_G(k^*\circ\pi_1) .
\]
\end{proposition}

\begin{proof}
Conditionally on $k^*$, distinct variables involve independent columns of $X$ and $\epsilon$, so
the expectation factorizes over the variables. In the factor of one variable,
Lemma~\ref{lem:split} splits each Hermite factor into signal and noise terms, only the term with no
noise survives the expectation over $\epsilon$, and the expectation over $X$ factorizes over the
color classes of the individuals at that variable.

\noindent\uline{Step 1: factorization over the variables.} Grouping the Hermite factors of
\eqref{eq:PG} by variable,
\[
\begin{aligned}
  \E\big[P_{G,\pi}\mid k^*\big]\ =\ &\hat s_G\,\frac{K^{\,|V_G|-1}}{|\Pi_G|}
  \Big(\frac{\sqrt d}{\Delta}\Big)^{\!|E_G|}\Big(\prod_{v,w}m_{v,w}!\Big)^{\!1/2}\\
  &\times\prod_{w\in W}\E\Big[\prod_{v\in N(w)}\frac{(1+\Delta^2/d)^{m_{v,w}/2}}{\sqrt{m_{v,w}!}}\,
  \He_{m_{v,w}}\Big(\frac{Y_{\pi_1(v),\pi_2(w)}}{\sqrt{1+\Delta^2/d}}\Big)\,\Big|\,k^*\Big] .
\end{aligned}
\]
Indeed, by \eqref{eq:model}, the factor of $w$ is a function of $k^*$ and of the columns
$(X_{k,\pi_2(w)})_{k\in[K]}$ and $(\epsilon_{i,\pi_2(w)})_{i\in[n]}$ alone. Distinct variables
use distinct columns, $\pi_2$ being injective, and conditionally on $k^*$ the pairs of columns
$\big((X_{k,j})_{k\in[K]},(\epsilon_{i,j})_{i\in[n]}\big)$, $j\in[d]$, are independent, the
entries of $X$ and $\epsilon$ being i.i.d.\ and independent of $k^*$.

\noindent\uline{Step 2: one variable.} Fix $w$, let $j=\pi_2(w)$, and abbreviate $m_v:=m_{v,w}$ for
$v\in N(w)$, so that $\sum_{v\in N(w)}m_v=\deg_G(w)$. Applying Lemma~\ref{lem:split} to each
of the $|N(w)|$ factors and expanding the product of the resulting sums, one splitting
$\mathbf s_v+\mathbf n_v=m_v$ per individual at $w$,
\begin{align*}
\prod_{v\in N(w)}\frac{(1+\Delta^2/d)^{m_v/2}}{\sqrt{m_v!}}&\He_{m_v}\!\Big(\frac{Y_{\pi_1(v),j}}{\sqrt{1+\Delta^2/d}}\Big)\\[-2pt]
&=\sum_{\substack{\mathbf s_v+\mathbf n_v=m_v\\ (v\in N(w))}}\
\prod_{v\in N(w)}\sqrt{\binom{m_v}{\mathbf s_v}}\Big(\frac{\Delta}{\sqrt d}\Big)^{\!\mathbf s_v}
\frac{\He_{\mathbf s_v}\big(X_{k^*(\pi_1(v)),j}\big)}{\sqrt{\mathbf s_v!}}\,
\frac{\He_{\mathbf n_v}\big(\epsilon_{\pi_1(v),j}\big)}{\sqrt{\mathbf n_v!}} .
\end{align*}
The rows $\pi_1(v)$, $v\in N(w)$, are distinct, $\pi_1$ being injective, so conditionally on
$k^*$ and $(X_{k,j})_{k\in[K]}$ the $\epsilon_{\pi_1(v),j}$ are i.i.d.\ $\mathcal N(0,1)$. As
$\E\He_r(\epsilon_{i,j})=0$ for $r\ge1$ and $\He_0=1$,
\[
  \E\Big[\prod_{v\in N(w)}\frac{\He_{\mathbf n_v}(\epsilon_{\pi_1(v),j})}{\sqrt{\mathbf n_v!}}
  \ \Big|\ k^*,(X_{k,j})_{k\in[K]}\Big]\ =\ \prod_{v\in N(w)}\one\{\mathbf n_v=0\} .
\]
So only the term with $\mathbf s_v=m_v$ for every $v\in N(w)$ survives, with
$\binom{m_v}{m_v}=1$, and
\[
\begin{aligned}
  &\E\Big[\prod_{v\in N(w)}\frac{(1+\Delta^2/d)^{m_v/2}}{\sqrt{m_v!}}
  \He_{m_v}\Big(\frac{Y_{\pi_1(v),j}}{\sqrt{1+\Delta^2/d}}\Big)\,\Big|\,k^*\Big]\\
  &\quad=\ \Big(\frac{\Delta}{\sqrt d}\Big)^{\!\deg_G(w)}
  \frac{1}{\big(\prod_{v\in N(w)}m_v!\big)^{1/2}}\;
  \E\Big[\prod_{v\in N(w)}\He_{m_v}\big(X_{k^*(\pi_1(v)),\,j}\big)\,\Big|\,k^*\Big]\\
  &\quad=\ \Big(\frac{\Delta}{\sqrt d}\Big)^{\!\deg_G(w)}
  \frac{1}{\big(\prod_{v\in N(w)}m_v!\big)^{1/2}}\;
  \prod_{k\in[K]}\mu\big(\{m_v:v\in N(w),\ k^*(\pi_1(v))=k\}\big) ,
\end{aligned}
\]
the last equality by grouping the $v\in N(w)$ by their color $k^*(\pi_1(v))=k$: conditionally
on $k^*$ the $X_{k,j}$, $k\in[K]$, are i.i.d.\ $\mathcal N(0,1)$, and the expectation of the
product over one color class is $\mu$ of that class, $\mu(\varnothing)=1$ for an empty one.

\noindent\uline{Step 3: product over the variables.} Multiplying Step~2 over $w\in W$, and using
$\sum_{w}\deg_G(w)=|E_G|$, $\prod_w\prod_{v\in N(w)}m_{v,w}!=\prod_{v,w}m_{v,w}!$ and the
definition of $f_G$ at $\chi=k^*\circ\pi_1$, the product over $w\in W$ in Step~1 equals
\[
  \Big(\frac{\Delta}{\sqrt d}\Big)^{\!|E_G|}\Big(\prod_{v,w}m_{v,w}!\Big)^{\!-1/2}
  f_G(k^*\circ\pi_1) .
\]
Substituting this into Step~1, the factors $(\sqrt d/\Delta)^{|E_G|}(\prod_{v,w}m_{v,w}!)^{1/2}$
cancel, which leaves the statement.
\end{proof}

\begin{corollary}\label{cor:selfloop}
If every variable of $G$ is a link, then for every $\chi:V\to[K]$
\[
\begin{aligned}
  f_G(\chi)\ &=\ \prod_{w\in W}\ \prod_{v,v'\in N(w)}\one\{\chi(v)=\chi(v')\}\\
  &=\ \one\{\chi\ \text{is constant on each connected component of }G\} .
\end{aligned}
\]
\end{corollary}

\begin{proof}
By Proposition~\ref{prop:master}, $f_G(\chi)$ is the product, over $w\in W$ and $k\in[K]$, of
the factors $\mu\big(\{m_{v,w}:v\in N(w),\ \chi(v)=k\}\big)$. For a link $w=\{v,v'\}$,
$m_{v,w}=m_{v',w}=1$ and the non-empty sets $\chi^{-1}(k)\cap N(w)$ are $\{v,v'\}$ if
$\chi(v)=\chi(v')$, and $\{v\}$ and $\{v'\}$ otherwise, every other $k$ contributing
$\mu(\varnothing)=1$. So the product over $k$ of the factors of $w$ is
\[
  \mu(1,1)=1\ \ \text{if }\chi(v)=\chi(v'),\qquad
  \mu(1)\,\mu(1)=0\ \ \text{if }\chi(v)\ne\chi(v') ,
\]
by Lemma~\ref{lem:wick}\textup{(ii)} and \textup{(i)}, that is
$\prod_{v,v'\in N(w)}\one\{\chi(v)=\chi(v')\}$. Multiplying over $w\in W$ gives the first identity. For the second, the product is $1$ exactly when $\chi$ is constant on $V_C$ for
every component $C$ of $G$, as two individuals of one component are joined by a path of $G$, on
which consecutive individuals are the two individuals of a link, and the two individuals of a
link lie in one component.
\end{proof}

\begin{lemma}\label{lem:colourcount}
Let $C\subseteq[n]$ be partitioned into non-empty parts $C_1,\dots,C_N$, let $0\le b\le N$, let
$k_1,\dots,k_{N-b}\in[K]$, and let
\[
  \mathcal A\ :=\ \bigcap_{a\le N-b}\big\{k^*(i)=k_a\ \ \forall i\in C_a\big\}\ \cap
  \bigcap_{a>N-b}\big\{k^*(i)=k^*(i')\ \ \forall i,i'\in C_a\big\} .
\]
Then $\Pr[\mathcal A]=K^{\,b-|C|}$.
\end{lemma}

\begin{proof}
Write $\mathcal A=\bigcap_{a\le N}\mathcal A_a$, $\mathcal A_a$ being the event of index $a$ in
the definition of $\mathcal A$. The $C_a$ being disjoint and the $k^*(i)$, $i\in C$, i.i.d.\
uniform on $[K]$, the events $\mathcal A_a$ are independent, and
\[
\begin{aligned}
  \Pr[\mathcal A_a]\ &=\ \prod_{i\in C_a}\Pr\big[k^*(i)=k_a\big]\ =\ K^{-|C_a|}
  &&(a\le N-b),\\
  \Pr[\mathcal A_a]\ &=\ \sum_{k\in[K]}\ \prod_{i\in C_a}\Pr\big[k^*(i)=k\big]\ =\ K^{1-|C_a|}
  &&(a>N-b) .
\end{aligned}
\]
Hence, as $\sum_{a\le N}|C_a|=|C|$,
\[
  \Pr[\mathcal A]\ =\ \prod_{a\le N}\Pr[\mathcal A_a]\ =\ K^{\,b-\sum_{a\le N}|C_a|}\ =\ K^{\,b-|C|} .
\]
\end{proof}

\subsubsection{The pairing identity}\label{sec:pair}

\begin{proposition}\label{prop:pairid}
For every $G$ and every $\chi:V\to[K]$,
\[
f_G(\chi)=\sum_{\Pai\in\Pcal(G)}f_{G[\Pai]}(\chi) ,
\]
where every variable of $G[\Pai]$ is a link, so that Corollary~\ref{cor:selfloop} gives
$f_{G[\Pai]}(\chi)$ (the three pairings of a variable of degree $4$ are drawn in
Figure~\ref{fig:centering}).
\end{proposition}
\begin{proof}
At each variable, the factor of $f_G(\chi)$ counts, by Lemma~\ref{lem:wick}, the pairings of the
edges at that variable in which no pair lies inside one individual and the two edges of each pair sit
at individuals of the same color. Multiplying over the variables gives a sum over $\Pcal(G)$, whose
summands are the $f_{G[\Pai]}(\chi)$ by Corollary~\ref{cor:selfloop}.

\noindent\uline{Step 1: one variable.} Fix $w$. By Proposition~\ref{prop:master} its factor in
$f_G(\chi)$ is $\prod_{k\in[K]}\mu(\{m_{v,w}:v\in N(w),\ \chi(v)=k\})$. By
Lemma~\ref{lem:wick}, applied for each $k$ with the individuals $v\in N(w)$ of color
$\chi(v)=k$ as the factors and their $m_{v,w}$ edges at $w$ as the points, the factor of $k$ is
the number of partitions into pairs of the edges at $w$ whose individual has color $k$, no pair
lying inside one individual. A family of such partitions, one for each $k$, is the same as a
partition of all the edges at $w$ into pairs such that no pair lies inside one individual
and the two edges of each pair sit at individuals of the same color. The partitions
with the first property are the pairings $\Pai_w$ of $w$ with no pair inside one individual, and such a
$\Pai_w$ has the second exactly when $\prod_{\{e,e'\}\in\Pai_w}\one\{\chi(v_e)=\chi(v_{e'})\}=1$. Hence
\[
  \prod_{k\in[K]}\mu\big(\{m_{v,w}:v\in N(w),\ \chi(v)=k\}\big)
  \ =\ \sum_{\Pai_w}\ \prod_{\{e,e'\}\in\Pai_w}\one\{\chi(v_e)=\chi(v_{e'})\} ,
\]
the sum running over the pairings $\Pai_w$ of $w$ with no pair inside one individual.

\noindent\uline{Step 2: all variables.} By definition, $\Pcal(G)$ is the set of the families $\Pai=(\Pai_w)_{w\in W}$ in
which each $\Pai_w$ is a pairing of $w$ with no pair inside one individual. Multiplying Step~1
over $w\in W$ and expanding the product of sums,
\[
\begin{aligned}
  f_G(\chi)\ &=\ \prod_{w\in W}\ \sum_{\Pai_w}\ \prod_{\{e,e'\}\in\Pai_w}\one\{\chi(v_e)=\chi(v_{e'})\}\\
  &=\ \sum_{\Pai\in\Pcal(G)}\ \prod_{w\in W}\ \prod_{\{e,e'\}\in\Pai_w}\one\{\chi(v_e)=\chi(v_{e'})\} .
\end{aligned}
\]

\noindent\uline{Step 3: the summands.} Let $\Pai\in\Pcal(G)$. Every variable of $G[\Pai]$ has degree $2$
and none is a self-loop, so each is a link, that of the pair $\{e,e'\}$ joining $v_e$ to $v_{e'}$. By
Corollary~\ref{cor:selfloop},
\[
  f_{G[\Pai]}(\chi)\ =\ \prod_{w\in W}\ \prod_{\{e,e'\}\in\Pai_w}\one\{\chi(v_e)=\chi(v_{e'})\} ,
\]
and with Step~2 this gives the statement.
\end{proof}

\subsection{The centering of the Hermite basis}\label{sec:center}

\subsubsection{Deletable sets}

Given a graph $G'$ and a set $\bar V$ of its individuals, write
\begin{equation}\label{eq:Sdef}
\mathcal V^{(\mathrm{free})}_{\bar V}(G')\ :=\!\!\bigcup_{C\in\CC_{G'}(\bar V)}\!\!V_C ,
\end{equation}
the set of the individuals of the free components of $G'$ (\S\ref{sec:graphs}), those missing $\bar V$. Until \S\ref{sec:cancel}, $G'=G[\Pai]$ and $\bar V=\{v_1,v_2\}$.

\begin{definition}[Deletable sets]\label{def:delete}
For $V'\subseteq V\setminus\{v_1,v_2\}$ let $G_{V'}$ be the graph with individual set $V'$, variable
set $\{w:N(w)\cap V'\ne\varnothing\}$ and only the edges incident to $V'$ (so
$\deg_{G_{V'}}(w)=\sum_{v\in V'}m_{v,w}$) and let $\Gs$ be $G$ with the
individuals of $V'$ and their edges removed (and any variable thereby left with no edge removed).
Call $V'$ \textbf{deletable} if
\[
V'\ne\varnothing\qquad\text{and}\qquad f_{G_{V'}}\neq 0,\ \text{that is, by Proposition~\ref{prop:pairid} and Corollary~\ref{cor:selfloop},}\ \Pcal(G_{V'})\ne\varnothing,
\]
and write $\mathcal D(G)$ for the set of deletable sets.
\end{definition}

\begin{definition}[Centered polynomials]\label{def:center}
For $V'\in\mathcal D(G)$ put
\begin{equation}\label{eq:LambdaS}
\Lambda_{V'}\ :=\ \frac{\hat s_G}{\hat s_{\Gs}}\!\!\sum_{\Pai_{V'}\in\Pcal(G_{V'})}\!\!K^{\,|\CC_{G_{V'}[\Pai_{V'}]}|} .
\end{equation}
Define,
\textbf{by recursion on $|E_G|+|V_G|$},
\begin{equation}\label{eq:Psh}
 \Psh_G\;:=\;P_G\;-\sum_{V'\in\mathcal D(G)}\Lambda_{V'}\,\Psh_{\Gs}.
\end{equation}
\end{definition}

The recursion is well founded, as $V'\ne\varnothing$ gives $|V_{\Gs}|<|V_G|$, while $|E_{\Gs}|\le|E_G|$, and $\Psh_G=P_G$ when
$\mathcal D(G)=\varnothing$. This is the case of the \textbf{empty graph} $\varnothing$: individual set $\{v_1,v_2\}$, no
variable and no edge. Its set of labelings $\Pi_\varnothing$ is a single point, $\hat s_\varnothing=K^{-1}$
and $P_\varnothing=\Psh_\varnothing=1$. It is the only graph of the recursion carrying no edge and
it is not a member of $\mathcal G_{\le D}$.

\begin{lemma}\label{lem:deletable}
\begin{enumerate}[label=\textup{(\alph*)},leftmargin=*,itemsep=2pt,topsep=2pt]
\item For $G\in\mathcal G_{\le D}$ and $V'\in\mathcal D(G)$, every
variable $w$ of $\Gs$ has even degree $\deg_{\Gs}(w)\ge2$, no individual of $\Gs$ other than $v_1,v_2$ is
isolated, and $|E_{\Gs}|<|E_G|$, $|V_{\Gs}|<|V_G|$. Hence $\Gs\in\mathcal G_{\le D}$ if
$|E_{\Gs}|\ge2$, and $\Gs=\varnothing$ otherwise.
\item  $\mathcal V^{(\mathrm{free})}_{\{v_1,v_2\}}(G[\Pai])\in\mathcal D(G)$ for every $\Pai\in\Pcal(G)$
with $G[\Pai]$ disconnected and $v_1,v_2$ in one component of it.
\end{enumerate}
\end{lemma}
\begin{proof}
For \textup{(a)}, a pairing of $G_{V'}$ makes every variable degree of $G_{V'}$ even, so the
variable degrees of $\Gs$ stay even, and $\Gs$ has fewer edges and individuals than $G$ as $V'$ is
non-empty. For \textup{(b)}, with $V'$ the set of the individuals of the components of $G[\Pai]$
missing $\{v_1,v_2\}$, the pairs of $\Pai$ at $V'$ form a pairing of $G_{V'}$ with no self-loop.

\textbf{(a)} \uline{Step 1: the degrees in $G_{V'}$ are even.} A pairing in $\Pcal(G_{V'})\ne\varnothing$
partitions the $\deg_{G_{V'}}(w)$ edges of $G_{V'}$ at each variable $w$ into pairs, so $\deg_{G_{V'}}(w)$
is even.

\noindent\uline{Step 2: the variables of $\Gs$.} Every variable degree of $G$ is even (definition of
$\mathcal G_{\le D}$), so, by Step~1, $\deg_{\Gs}(w)\ =\ \deg_{G}(w)-\deg_{G_{V'}}(w)$ is even for every variable $w$ of $G$. The variables with $\deg_{\Gs}(w)=0$ being removed (Definition~\ref{def:delete}), every
variable $w$ of $\Gs$ has even degree $\deg_{\Gs}(w)\ge2$.

\noindent\uline{Step 3: the individuals of $\Gs$.} An individual $v\notin V'$ keeps all its edges, so
$\deg_{\Gs}(v)=\deg_{G}(v)$, and $\deg_G(v)\ge1$ for $v\ne v_1,v_2$, no node of $G$ other than
$v_1,v_2$ being isolated (\S\ref{sec:graphs}).

\noindent\uline{Step 4: the sizes.} Let $v\in V'$. Then $v\ne v_1,v_2$, so $v$ carries an edge of $G$,
which is an edge of $G_{V'}$. So
$|E_{G_{V'}}|\ge1$, and
\[
  |E_{\Gs}|\ =\ |E_G|-|E_{G_{V'}}|\ <\ |E_G|,\qquad |V_{\Gs}|\ =\ |V_G|-|V'|\ <\ |V_G| .
\]

\noindent\uline{Step 5: conclusion.} If $|E_{\Gs}|\ge2$: by Steps~2--4 every variable degree of $\Gs$ is
even, $2\le|E_{\Gs}|<|E_G|\le D$, and no node of $\Gs$ other than $v_1,v_2$ is isolated, so
$\Gs\in\mathcal G_{\le D}$. If $|E_{\Gs}|<2$: by Step~2 $\Gs$ has no variable, so
$|E_{\Gs}|=0$, and by Step~3 no individual other than $v_1,v_2$ (which are not in $V'$). So
$\Gs=\varnothing$.

\textbf{(b)} Put $V':=\mathcal V^{(\mathrm{free})}_{\{v_1,v_2\}}(G[\Pai])$, the union \eqref{eq:Sdef} of the sets $V_C$
over the components $C$ of $G[\Pai]$ missing $\{v_1,v_2\}$.

\noindent\uline{Step 1: $\varnothing\ne V'\subseteq V\setminus\{v_1,v_2\}$.} The component of $G[\Pai]$
holding $v_1$ and $v_2$ meets $\{v_1,v_2\}$, so $v_1,v_2\notin V'$. As $G[\Pai]$ is
disconnected, it has a second component. This one misses $\{v_1,v_2\}$, and it holds an
individual, every variable of $G[\Pai]$ having degree $2$. So $V'\ne\varnothing$.

\noindent\uline{Step 2: $\Pcal(G_{V'})\ne\varnothing$.} $V'$ is a union of components of $G[\Pai]$, so a pair of
$\Pai$ with one edge at an individual of $V'$ has its other edge there too. Hence the pairs of
$\Pai$ formed by edges at $V'$ make a pairing $\Pai_{V'}$ of $G_{V'}$, with no self-loop as
$\Pai\in\Pcal(G)$: $\Pai_{V'}\in\Pcal(G_{V'})$. With Step~1, this gives $V'\in\mathcal D(G)$.
\end{proof}

\subsection{The correlation of $x$ and rescaled Hermite basis elements}

As $x$ is a function of $k^*$, summing Proposition~\ref{prop:master} over the $|\Pi_G|$ labelings gives
\begin{equation}\label{eq:numG}
\E\big[(x-\E x)P_G\big]\;=\;\hat s_G\,K^{\,|V_G|-1}\,\Cov(x,f_G) ,
\end{equation}
where $\Cov(x,f_G)$ abbreviates $\Cov\big(x,f_G(k^*\circ\pi_1)\big)$: as $\pi_1$ is injective
with $\pi_1(v_1)=1$, $\pi_1(v_2)=2$, the joint law of $\big(x,k^*\circ\pi_1\big)$ is the same for
every labeling, so that number does not depend on $\pi$.

\begin{lemma}\label{lem:numP}
For every $G\in\mathcal G_{\le D}\cup\{\varnothing\}$,
\[
  \E\big[(x-\E x)P_G\big]\ =\ \hat s_G\Big(1-\frac1K\Big)
  \!\!\sum_{\substack{\Pai\in\Pcal(G):\ v_1,v_2\ \text{in one}\\ \text{component of}\ G[\Pai]}}\!\!
  K^{\,|\CC_{G[\Pai]}|-1} .
\]
\end{lemma}

\begin{proof}
By Proposition~\ref{prop:pairid}, $\Cov(x,f_G)$ in \eqref{eq:numG} is the sum over
$\Pai\in\Pcal(G)$ of the covariances of $x$ with $f_{G[\Pai]}$, the indicator that the colors are
constant on each component of $G[\Pai]$. Lemma~\ref{lem:colourcount} computes each covariance, which
vanishes unless $v_1,v_2$ lie in one component of $G[\Pai]$.

\noindent\uline{Step 1: the covariance of one pairing.} Fix $\Pai\in\Pcal(G)$ and put $\chi:=k^*\circ\pi_1$:
as $\pi_1$ is injective, the $\chi(v)$, $v\in V$, are i.i.d.\ uniform on $[K]$, and
$x=\one\{\chi(v_1)=\chi(v_2)\}$ by \eqref{eq:x}, $\pi_1(v_1)=1$, $\pi_1(v_2)=2$. Every variable
of $G[\Pai]$ has degree $2$, so every component of $G[\Pai]$ holds an individual and the sets
$V_C$, $C\in\CC_{G[\Pai]}$, partition $V$. By Corollary~\ref{cor:selfloop},
$f_{G[\Pai]}(\chi)$ is the indicator that $\chi$ is constant on each $V_C$, so by
Lemma~\ref{lem:colourcount} (parts $\pi_1(V_C)$, none prescribed)
\[
  \E\big[f_{G[\Pai]}(\chi)\big]\ =\ K^{\,|\CC_{G[\Pai]}|-|V|} .
\]
\textbf{Case 1: $v_1$ and $v_2$ lie in two distinct components $C_1\ne C_2$.} Then
$x\,f_{G[\Pai]}(\chi)$ is the indicator that $\chi$ is constant on each $V_C$,
$C\ne C_1,C_2$, and on $V_{C_1}\cup V_{C_2}$. By Lemma~\ref{lem:colourcount} with these
$|\CC_{G[\Pai]}|-1$ parts,
\[
  \E\big[x\,f_{G[\Pai]}(\chi)\big]\ =\ K^{\,|\CC_{G[\Pai]}|-1-|V|}\ =\ \E[x]\,\E\big[f_{G[\Pai]}(\chi)\big],
  \qquad\text{so}\qquad \Cov\big(x,f_{G[\Pai]}\big)\ =\ 0 .
\]
\textbf{Case 2: $v_1$ and $v_2$ lie in one component.} Then $f_{G[\Pai]}(\chi)\ne0$ forces
$\chi(v_1)=\chi(v_2)$, so $x\,f_{G[\Pai]}(\chi)=f_{G[\Pai]}(\chi)$ and
\[
  \Cov\big(x,f_{G[\Pai]}\big)\ =\ \E\big[f_{G[\Pai]}(\chi)\big]-\E[x]\,\E\big[f_{G[\Pai]}(\chi)\big]
  \ =\ \Big(1-\frac1K\Big)K^{\,|\CC_{G[\Pai]}|-|V|} .
\]

\noindent\uline{Step 2: summing over the pairings.} By Proposition~\ref{prop:pairid},
$\Cov(x,f_G)=\sum_{\Pai\in\Pcal(G)}\Cov(x,f_{G[\Pai]})$. Inserting Step~1 into \eqref{eq:numG} gives the statement, as
$K^{|V|-1}K^{|\CC_{G[\Pai]}|-|V|}=K^{|\CC_{G[\Pai]}|-1}$.
\end{proof}

\begin{proposition}\label{prop:num}
For every $G\in\mathcal G_{\le D}\cup\{\varnothing\}$,
\[
\E\big[(x-\E x)\,\Psh_G\big]\;=\;\hat s_G\Big(1-\frac1K\Big)
\cdot\#\{\Pai\in\Pcal(G):G[\Pai]\text{ connected}\} .
\]
\end{proposition}

\begin{proof}
By Lemma~\ref{lem:numP}, $\E[(x-\E x)P_G]$ is a sum over the pairings $\Pai\in\Pcal(G)$ with
$v_1,v_2$ in one component of $G[\Pai]$. We group those with $G[\Pai]$ disconnected by the set $V'$ of
the individuals of the components missing $\{v_1,v_2\}$, and show, by induction on $|E_G|+|V_G|$,
that the term of $V'$ in \eqref{eq:Psh} removes exactly their contribution.

\noindent\textbf{Induction hypothesis.} For an integer $m\ge2$, let $(\mathrm H_m)$ be the
statement:
\begin{center}
\emph{the identity of Proposition~\ref{prop:num} holds for every
$G\in\mathcal G_{\le D}\cup\{\varnothing\}$ with $|E_G|+|V_G|\le m$.}
\end{center}
Every graph has $|V_G|\ge2$, as $v_1,v_2\in V_G$. So it suffices to prove $(\mathrm H_m)$ for
every $m\ge2$, by induction on $m$.

\noindent\textbf{Initialization: $(\mathrm H_2)$.} A graph of $\mathcal G_{\le D}$ has
$|E_G|\ge2$, hence $|E_G|+|V_G|\ge4$. So the only graph concerned is $\varnothing$. There
$\Psh_\varnothing=1$, so the left-hand side is $\E[x-\E x]=0$. The only pairing of
$\varnothing$ is the empty one, and $\varnothing$ itself, with the two isolated individuals
$v_1,v_2$, has two components, so the right-hand side is $0$ as well.

\noindent\textbf{Induction step: $(\mathrm H_{m-1})\Rightarrow(\mathrm H_m)$, for $m\ge3$.}
Assume $(\mathrm H_{m-1})$. It covers the graphs with $|E_G|+|V_G|\le m-1$, so let
$G\in\mathcal G_{\le D}$ have $|E_G|+|V_G|=m$. By Lemma~\ref{lem:deletable}(a), for every
$V'\in\mathcal D(G)$ the graph $\Gs$ is $\varnothing$ or a member of $\mathcal G_{\le D}$, with
$|E_{\Gs}|+|V_{\Gs}|\le m-1$. So, by $(\mathrm H_{m-1})$, the identity holds for every $\Gs$,
$V'\in\mathcal D(G)$. Steps~1--4 below prove it for $G$, which gives $(\mathrm H_m)$.

\noindent\uline{Step 1: the connected pairings.} A pairing $\Pai$ with $G[\Pai]$ connected has $|\CC_{G[\Pai]}|=1$, so
Lemma~\ref{lem:numP} reads
\begin{equation}\label{eq:numsplit}
\begin{aligned}
  \E\big[(x-\E x)P_G\big]\ &=\ \hat s_G\Big(1-\frac1K\Big)\#\{\Pai\in\Pcal(G):G[\Pai]\text{ connected}\}\\
  &\qquad+\hat s_G\Big(1-\frac1K\Big)
  \!\!\!\sum_{\substack{\Pai\in\Pcal(G):\ v_1,v_2\ \text{in one component}\\ \text{of}\ G[\Pai],\
  G[\Pai]\ \text{disconnected}}}\!\!\!K^{\,|\CC_{G[\Pai]}|-1} .
\end{aligned}
\end{equation}

\noindent\uline{Step 2: what the recursion subtracts.} For $V'\in\mathcal D(G)$, by the identity
for $\Gs$ given by $(\mathrm H_{m-1})$ and by \eqref{eq:LambdaS}, the factors $\hat s_{\Gs}$ canceling,
\begin{equation}\label{eq:numsub}
\begin{aligned}
  \Lambda_{V'}\,\E\big[(x-\E x)\Psh_{\Gs}\big]\ &=\ \hat s_G\Big(1-\frac1K\Big)\sum_{\Pai_{V'}\in\Pcal(G_{V'})}K^{\,|\CC_{G_{V'}[\Pai_{V'}]}|}\\
  &\qquad\times\#\{\Pai'\in\Pcal(\Gs):\Gs[\Pai']\text{ connected}\} .
\end{aligned}
\end{equation}

\noindent\uline{Step 3: the disconnected pairings, grouped by $V'$.} Let $\Pai\in\Pcal(G)$ have $v_1,v_2$ in
one component of $G[\Pai]$. Then $G[\Pai]$ is disconnected exactly when
$\mathcal V^{(\mathrm{free})}_{\{v_1,v_2\}}(G[\Pai])\ne\varnothing$, and in that case
$\mathcal V^{(\mathrm{free})}_{\{v_1,v_2\}}(G[\Pai])\in\mathcal D(G)$ by Lemma~\ref{lem:deletable}(b). We show that, for
every $V'\in\mathcal D(G)$,
\begin{equation}\label{eq:numclass}
\begin{aligned}
  &\sum_{\substack{\Pai\in\Pcal(G):\ v_1,v_2\ \text{in one}\\ \text{component of}\ G[\Pai],\
  \mathcal V^{(\mathrm{free})}_{\{v_1,v_2\}}(G[\Pai])=V'}}\!\!\!K^{\,|\CC_{G[\Pai]}|-1}\\
  &\qquad=\ \sum_{\Pai_{V'}\in\Pcal(G_{V'})}K^{\,|\CC_{G_{V'}[\Pai_{V'}]}|}\cdot\#\{\Pai'\in\Pcal(\Gs):\Gs[\Pai']\text{ connected}\} .
\end{aligned}
\end{equation}
Let $\Pai$ be a pairing summed on the left. $V'$ is a union of components of $G[\Pai]$, so a pair
of $\Pai$ with one edge at an individual of $V'$ has its other edge there too. Hence $\Pai$ is the
union of a pairing $\Pai_{V'}$ of $G_{V'}$ (the pairs of edges at $V'$) and a pairing $\Pai'$ of $\Gs$
(the other pairs), neither with a self-loop, and $G[\Pai]$ is the disjoint union of
$G_{V'}[\Pai_{V'}]$ and $\Gs[\Pai']$. The components of $G[\Pai]$ missing $\{v_1,v_2\}$ are those of
$G_{V'}[\Pai_{V'}]$, their individuals making up $V'$. So every component of $\Gs[\Pai']$ meets
$\{v_1,v_2\}$, and, $v_1,v_2$ lying in one component, $\Gs[\Pai']$ is connected and
\[
  |\CC_{G[\Pai]}|\ =\ |\CC_{G_{V'}[\Pai_{V'}]}|+1 .
\]
Conversely, for $\Pai_{V'}\in\Pcal(G_{V'})$ and $\Pai'\in\Pcal(\Gs)$ with $\Gs[\Pai']$ connected, their
union $\Pai$ lies in $\Pcal(G)$, and $G[\Pai]$ is the disjoint union of $G_{V'}[\Pai_{V'}]$, whose
components miss $\{v_1,v_2\}$ and have individual sets covering $V'$, and of the connected
$\Gs[\Pai']$, which holds $v_1,v_2$. So $v_1,v_2$ lie in one component of $G[\Pai]$ and
$\mathcal V^{(\mathrm{free})}_{\{v_1,v_2\}}(G[\Pai])=V'$. The two maps are inverse to each other, which gives
\eqref{eq:numclass}.

\noindent\uline{Step 4: conclusion.} By Step~3, summing \eqref{eq:numclass} over $V'\in\mathcal D(G)$ gives
the last sum of \eqref{eq:numsplit}. Hence, by \eqref{eq:Psh}, \eqref{eq:numsplit} and
\eqref{eq:numsub},
\[
\begin{aligned}
  \E\big[(x-\E x)\Psh_G\big]\ &=\ \E\big[(x-\E x)P_G\big]
  -\sum_{V'\in\mathcal D(G)}\Lambda_{V'}\,\E\big[(x-\E x)\Psh_{\Gs}\big]\\
  &=\ \hat s_G\Big(1-\frac1K\Big)\#\{\Pai\in\Pcal(G):G[\Pai]\text{ connected}\} .
\end{aligned}
\]
\end{proof}

\subsection{The blocks, and the projection on the blocks}\label{sec:blocks}

\subsubsection{The blocks}\label{sec:blockdef}

\begin{definition}[Blocks]\label{def:block}
For an injection $\pi_2:W\to[d]$, an injection $\pi_1$ of $V_{G\mathbf n^\dagger}$ into $[n]$ with
$\pi_1(v_1)=1$ and $\pi_1(v_2)=2$, and an \textbf{injection} $\sigma$ of $V_{G\mathbf s^\dagger}$ into
$[K]$, put
\[
Q_{G,\pi_1,\pi_2,\sigma}:=
\prod_{v\in V_{G\mathbf n^\dagger}}\prod_w\frac{\He_{m_{v,w}}(\epsilon_{\pi_1(v),\pi_2(w)})}{\sqrt{m_{v,w}!}}
\ \cdot\prod_{v\in V_{G\mathbf s^\dagger}}\prod_w\frac{\He_{m_{v,w}}(X_{\sigma(v),\pi_2(w)})}{\sqrt{m_{v,w}!}} .
\]
These are the \textbf{blocks} of $G$. They are non-constant, $G$ having an edge. The
\textbf{invariant block} $\bar Q_G$ is the sum of all of
them divided by its norm, $\bar Q_G:=\sum_{\pi_1,\pi_2,\sigma}Q_{G,\pi_1,\pi_2,\sigma}\big/\big\|\sum_{\pi_1,\pi_2,\sigma}Q_{G,\pi_1,\pi_2,\sigma}\big\|$, a unit vector (Lemma~\ref{lem:ortho}).
\end{definition}

\begin{lemma}\label{lem:ortho}
Every block has norm $1$, and blocks of distinct graphs of $\mathcal G_{\le D}$ are
orthogonal. For $G\in\mathcal G_{\le D}$,
writing $\#\{(\pi_1,\pi_2,\sigma)\}$ for the number of triples of Definition~\ref{def:block} for $G$,
\begin{equation}\label{eq:Qbarnorm}
\begin{gathered}
\Big\|\sum_{\pi_1,\pi_2,\sigma}Q_{G,\pi_1,\pi_2,\sigma}\Big\|^2\ =\ |\Aut(G)|\cdot\#\{(\pi_1,\pi_2,\sigma)\},\\
\tfrac12\,n^{|V_{G\mathbf n^\dagger}|-2}d^{|W_G|}K^{|V_{G\mathbf s^\dagger}|}\ \le\ \#\{(\pi_1,\pi_2,\sigma)\}\ \le\
n^{|V_{G\mathbf n^\dagger}|-2}d^{|W_G|}K^{|V_{G\mathbf s^\dagger}|} ,
\end{gathered}
\end{equation}
so that \textbf{the invariant blocks $\bar Q_G$, $G\in\mathcal G_{\le D}$, are orthonormal}.
\end{lemma}

\begin{proof}
A block is a product of normalized Hermite polynomials of distinct entries of $\epsilon$ and
$X$, so, by Lemma~\ref{lem:wick}, two blocks are equal or orthogonal according as their
multi-indices agree. The main step shows that the multi-index determines the graph, and the triple up
to an automorphism of the graph, which gives the orthogonality across graphs and the norm of the sum
of the blocks of one graph.

\noindent\uline{Step 1: products of Hermite factors.} For a \emph{multi-index} $\alpha$, a map
from the entries of $\epsilon$ and $X$ to $\{0,1,2,\dots\}$ non-zero at finitely many of them, put
$\mathcal H_\alpha:=\prod_Z\He_{\alpha(Z)}(Z)/\sqrt{\alpha(Z)!}$. The entries being i.i.d.\
$\mathcal N(0,1)$ and $\He_0=1$, Lemma~\ref{lem:wick}\textup{(i),(ii)} gives
\begin{equation}\label{eq:hermprod}
  \E\big[\mathcal H_\alpha\,\mathcal H_\beta\big]\ =\ \one\{\alpha=\beta\} .
\end{equation}
A block $Q_{G,\pi_1,\pi_2,\sigma}$ is $\mathcal H_\alpha$, its multi-index $\alpha$ being
$m_{v,w}$ at $\epsilon_{\pi_1(v),\pi_2(w)}$ ($v\in V_{G\mathbf n^\dagger}$) and at
$X_{\sigma(v),\pi_2(w)}$ ($v\in V_{G\mathbf s^\dagger}$), entries that are pairwise distinct as
$\pi_1,\pi_2,\sigma$ are injective, and $0$ elsewhere. So every block has norm $1$, and two blocks
are equal or orthogonal according as their multi-indices are equal or not.

\noindent\uline{Step 2: the multi-index determines the graph, and the triple up to an automorphism.} Let
$G,G'\in\mathcal G_{\le D}$, with triples $(\pi_1,\pi_2,\sigma)$ of $G$ and
$(\pi'_1,\pi'_2,\sigma')$ for $G'$ giving blocks with the same multi-index. Send each
$v\in V_{G\mathbf n^\dagger}$ to the row $\pi_1(v)$ of $\epsilon$, each $v\in V_{G\mathbf
s^\dagger}$ to the row $\sigma(v)$ of $X$, each $w\in W_G$ to the column $\pi_2(w)$. This map is
injective. Every node of $G$ other than $v_1,v_2$ carries an edge (\S\ref{sec:graphs}), so its
image is the row or column of an entry with non-zero exponent, while $v_1,v_2$ are sent to the rows $1,2$
of $\epsilon$. So the image is the set made of the rows $1,2$ of $\epsilon$ and of the rows and
columns of the entries with non-zero exponent, which depends on the multi-index only. The same
holds for $G'$. Composing the first map with the inverse of the second gives a bijection $\varphi$
from the nodes of $G$ onto those of $G'$ sending $V_{G\mathbf n^\dagger}$ to $V_{G'\mathbf
n^\dagger}$, $V_{G\mathbf s^\dagger}$ to $V_{G'\mathbf s^\dagger}$, $W_G$ to $W_{G'}$, fixing
$v_1,v_2$, and with $m^{G'}_{\varphi(v),\varphi(w)}=m_{v,w}$, both being the exponent at the
entry of $(v,w)$. So $\varphi$ is an isomorphism $G\to G'$. As $\mathcal G_{\le D}$ holds one
graph per isomorphism class, $G=G'$, $\varphi\in\Aut(G)$, and
$(\pi'_1,\pi'_2,\sigma')\circ\varphi=(\pi_1,\pi_2,\sigma)$. Conversely, for $\varphi\in\Aut(G)$,
the triple $(\pi_1,\pi_2,\sigma)\circ\varphi^{-1}$ is again a triple of $G$ (an automorphism
preserves degrees and fixes $v_1,v_2$, hence preserves $V_{G\mathbf n^\dagger}$ by
\eqref{eq:blocksplit}) and gives the same block, by the substitution
$(v,w)\mapsto(\varphi(v),\varphi(w))$ in Definition~\ref{def:block} and
$m_{\varphi(v),\varphi(w)}=m_{v,w}$. Finally, if
$(\pi_1,\pi_2,\sigma)\circ\varphi=(\pi_1,\pi_2,\sigma)$, injectivity of the three maps gives
$\varphi=\mathrm{id}$. 
With Step~1, this shows that blocks of distinct graphs of $\mathcal G_{\le D}$ are orthogonal, and that two triples
give the same block of $G$ exactly when they differ by an automorphism of $G$, no automorphism
other than the identity fixing a triple.

\noindent\uline{Step 3: the number of triples.} Counting the injections $\pi_1,\pi_2,\sigma$,
\[
  \#\{(\pi_1,\pi_2,\sigma)\}\ =\ \prod_{k=2}^{|V_{G\mathbf n^\dagger}|-1}(n-k)\cdot
  \prod_{k=0}^{|W_G|-1}(d-k)\cdot\prod_{k=0}^{|V_{G\mathbf s^\dagger}|-1}(K-k) ,
\]
and \eqref{eq:Picount} applies to each factor, as $|V_G|\le2D$ and $|W_G|,|V_{G\mathbf s^\dagger}|\le D$
(Lemma~\ref{lem:sizes}). As $(\tfrac78)^3\ge\tfrac12$, this gives the second line of \eqref{eq:Qbarnorm}.

\noindent\uline{Step 4: the norm.} Expanding the square and using Step~1,
\[
  \Big\|\sum_{\pi_1,\pi_2,\sigma}Q_{G,\pi_1,\pi_2,\sigma}\Big\|^2
  \ =\sum_{(\pi_1,\pi_2,\sigma)}\ \sum_{(\pi'_1,\pi'_2,\sigma')}
  \one\big\{Q_{G,\pi_1,\pi_2,\sigma}=Q_{G,\pi'_1,\pi'_2,\sigma'}\big\}
  \ =\ |\Aut(G)|\cdot\#\{(\pi_1,\pi_2,\sigma)\} ,
\]
since by Step~2, for a fixed triple $(\pi_1,\pi_2,\sigma)$, the triples giving the same block are
the $(\pi_1,\pi_2,\sigma)\circ\varphi^{-1}$, $\varphi\in\Aut(G)$, and these are pairwise
distinct. This is the first line of \eqref{eq:Qbarnorm}.

\noindent\uline{Step 5: the invariant blocks.} By Step~4 the sum of the blocks of $G$ has positive norm, so
$\bar Q_G$ is well defined and has norm $1$. For $G\ne G'$ in $\mathcal G_{\le D}$, every block of
$G$ is orthogonal to every block of $G'$ (Steps~1 and~2), so
$\langle\bar Q_G,\bar Q_{G'}\rangle=0$.
\end{proof}

\subsubsection{The projection on the blocks}\label{sec:expansion}

Fix $G$, a second graph $H$, and a block
$Q=Q_{H,\pi^0_1,\pi^0_2,\sigma}$ of $H$ (Definition~\ref{def:block}). Write
$\pi^0:=(\pi^0_1,\pi^0_2)$ and $m^H$ for the multiplicities of $H$. Until \S\ref{sec:solve}, the letter $u$ denotes
an individual of $H_{\mathbf s^\dagger}$.

$P_{G,\pi}$ has one Hermite factor for each pair $(v,w)$ joined by at least one edge, and
Lemma~\ref{lem:split} breaks that factor into $m_{v,w}+1$ terms, one per value of
$\mathbf s_{v,w}$. So multiplying the factors out gives one term per splitting $\mathbf s$ of $G$
(\S\ref{sec:sym}).  Explicitly, the term of
\eqref{eq:PG} indexed by $\mathbf s$ is
\begin{equation}\label{eq:PGs}
\begin{aligned}
P_{G,\pi,\mathbf s}\ :=\ &\hat s_G\,\frac{K^{\,|V_G|-1}}{|\Pi_G|}\Big(\frac{\sqrt d}{\Delta}\Big)^{\!|E_G|}
\Big(\prod_{v,w}m_{v,w}!\Big)^{\!1/2}\\[2pt]
&\times\prod_{(v,w):\,m_{v,w}\ge1}\sqrt{\binom{m_{v,w}}{\mathbf s_{v,w}}}
\Big(\frac{\Delta}{\sqrt d}\Big)^{\!\mathbf s_{v,w}}
\frac{\He_{\mathbf s_{v,w}}\big(X_{k^*(\pi_1(v)),\,\pi_2(w)}\big)}{\sqrt{\mathbf s_{v,w}!}}\,
\frac{\He_{\mathbf n_{v,w}}\big(\epsilon_{\pi_1(v),\,\pi_2(w)}\big)}{\sqrt{\mathbf n_{v,w}!}} ,
\end{aligned}
\end{equation}
the first line being the scalar of \eqref{eq:PG} unchanged, and
$P_{G,\mathbf s}:=\sum_{\pi\in\Pi_G}P_{G,\pi,\mathbf s}$ is the corresponding part of $P_G$, so
that $P_G=\sum_{\mathbf s}P_{G,\mathbf s}$.

The same decomposition applies to $\Psh_G$: running \eqref{eq:Psh} one splitting at a time, set
\begin{equation}\label{eq:Pshs}
\Psh_{G,\mathbf s}\ :=\ P_{G,\mathbf s}\ -\!\!\sum_{\substack{V'\in\mathcal D(G)\\ \text{purely
signal in }\mathbf s}}\!\!\Lambda_{V'}\,\Psh_{\Gs,\,\mathbf s|_{\Gs}} .
\end{equation}

The pairs $(v,w)\in V\times W$ with $m_{v,w}>0$ carried by $\Gs$ are exactly those of $G$ with
$v\notin V'$, so restriction is a bijection from the splittings making $V'$ purely signal onto the
splittings of $\Gs$. Hence, if $\Psh_{\Gs}=\sum_{\mathbf s'}\Psh_{\Gs,\mathbf s'}$ for every
$V'\in\mathcal D(G)$, summing \eqref{eq:Pshs} over $\mathbf s$ returns \eqref{eq:Psh}. As $\Psh_{\varnothing,\mathbf s}=1$
for the one splitting $\mathbf s$ of $\varnothing$, induction on $|E_G|+|V_G|$ gives $\Psh_G=\sum_{\mathbf s}\Psh_{G,\mathbf s}$. 

\paragraph{The matchings a block admits, and the scalar $\tau_{\mathbf s,H}$.} The matchings of $G$ with $H$
that concern us are those in
\begin{equation}\label{eq:Mset}
\begin{aligned}
\mathcal M(\mathbf s,H):=\big\{\mathbf M:\
&\mathbf M|_{G_{\mathbf n}}\ \text{is an isomorphism}\ G_{\mathbf n}\to H_{\mathbf n^\dagger},\ \
\mathbf M_V\subseteq V_{G\mathbf n}\times V_{H\mathbf n^\dagger},\\
&\{w':(w,w')\in\mathbf M_W\}=W_H\big\} .
\end{aligned}
\end{equation}

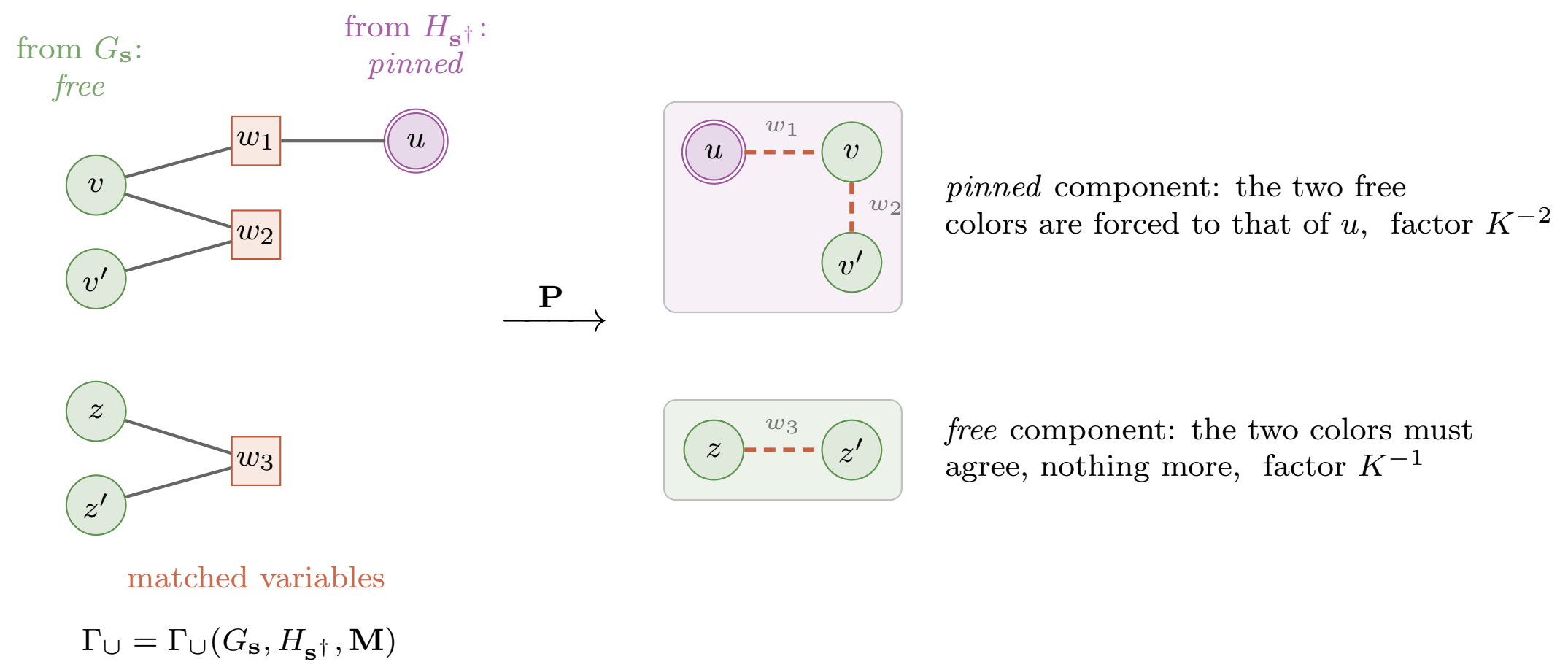
\begin{figure}[htbp]
\centering
\begin{tikzpicture}[
  fr/.style={circle,draw=OliveGreen!80,fill=OliveGreen!12,inner sep=0.4pt,minimum size=5.4mm,font=\scriptsize},
  pn/.style={circle,draw=Plum!80,fill=Plum!15,double,inner sep=0.4pt,minimum size=5.4mm,font=\scriptsize},
  va/.style={rectangle,draw=BrickRed!80,fill=BrickRed!8,inner sep=1pt,minimum size=4.4mm,font=\scriptsize},
  e/.style={black!60,thick},
  pr/.style={BrickRed!75,very thick,densely dashed}]
\begin{scope}
  \node[fr] (v)  at (0,1.30)  {$v$};
  \node[fr] (vp) at (0,0.45)  {$v'$};
  \node[fr] (z)  at (0,-0.75) {$z$};
  \node[fr] (zp) at (0,-1.60) {$z'$};
  \node[va] (w1) at (1.45,1.70)  {$w_1$};
  \node[va] (w2) at (1.45,0.85)  {$w_2$};
  \node[va] (w3) at (1.45,-1.20) {$w_3$};
  \node[pn] (u)  at (2.90,1.70)  {$u$};
  \draw[e] (v)--(w1); \draw[e] (w1)--(u);
  \draw[e] (v)--(w2); \draw[e] (vp)--(w2);
  \draw[e] (z)--(w3); \draw[e] (zp)--(w3);
  \node[OliveGreen!65,font=\scriptsize,align=center] at (-0.15,2.35)
        {from $G_{\mathbf s}$:\\ \emph{free}};
  \node[Plum!75,font=\scriptsize,align=center] at (2.90,2.55)
        {from $H_{\mathbf s^\dagger}$:\\ \emph{pinned}};
  \node[BrickRed!70,font=\scriptsize] at (1.45,-2.25) {matched variables};
  \node[font=\scriptsize] at (1.3,-2.85) {$\Gcup=\Gcup(G_{\mathbf s},H_{\mathbf s^\dagger},\mathbf M)$};
\end{scope}
\node at (4.15,0.2) {$\xrightarrow{\ \ \Pai\ \ }$};
\begin{scope}[xshift=5.6cm]
  \node[pn] (U)  at (0,1.60)  {$u$};
  \node[fr] (V)  at (1.25,1.60) {$v$};
  \node[fr] (VP) at (1.25,0.60) {$v'$};
  \draw[pr] (U)--(V)  node[midway,above,font=\tiny,black!55] {$w_1$};
  \draw[pr] (V)--(VP) node[midway,right,font=\tiny,black!55] {$w_2$};
  \begin{scope}[on background layer]
   \node[draw=black!25,rounded corners=3pt,fill=Plum!5,fit=(U)(V)(VP),inner sep=5pt] (c1) {};
  \end{scope}
  \node[font=\scriptsize,anchor=west,align=left] at ($(c1.east)+(0.25,0)$)
        {\emph{pinned} component: the two free\\ colors are forced to that of $u$,\ \ factor $K^{-2}$};
  \node[fr] (Z)  at (0,-1.10)    {$z$};
  \node[fr] (ZP) at (1.25,-1.10) {$z'$};
  \draw[pr] (Z)--(ZP) node[midway,above,font=\tiny,black!55] {$w_3$};
  \begin{scope}[on background layer]
   \node[draw=black!25,rounded corners=3pt,fill=OliveGreen!7,fit=(Z)(ZP),inner sep=5pt] (c2) {};
  \end{scope}
  \node[font=\scriptsize,anchor=west,align=left] at ($(c2.east)+(0.25,0)$)
        {\emph{free} component: the two colors must\\ agree, nothing more,\ \ factor $K^{-1}$};
\end{scope}
\end{tikzpicture}
\caption{The merged graph $\Gcup$ and one of its pairings (Lemma~\ref{lem:leak}(b)).
\textbf{Left:} the signal part of the splitting and the signal part of the block, merged along the
matched variables. A \emph{free} node carries a random color, a \emph{pinned} node a fixed one.
The pinned colors are pairwise distinct. \textbf{Right:}
the pairing $\Pai$ replaces each variable by links (dashed), and every component it creates must
be monochromatic. A component holding a pinned node forces its $r$ free colors to that fixed
color: factor $K^{-r}$. A \emph{free} component only forces its $r$ colors to agree: factor
$K^{1-r}$, with \emph{no} constraint relative to the pinned colors. Here $|\CC_{\Gcup[\Pai]}(V_{H\mathbf s^\dagger})|=1$ and $|V_{G\mathbf s}|=4$, total
factor $K^{|\CC_{\Gcup[\Pai]}(V_{H\mathbf s^\dagger})|-|V_{G\mathbf s}|}=K^{-3}$.}
\label{fig:gamma}
\end{figure}

For $\mathbf M\in\mathcal M(\mathbf s,H)$ let
$\Gcup:=\Gcup(G_{\mathbf s},H_{\mathbf s^\dagger},\mathbf M)$ be the merged graph of
\S\ref{sec:graphs}, whose matched nodes are variables only, as $\mathbf M_V\subseteq V_{G\mathbf n}\times V_{H\mathbf n^\dagger}$ by \eqref{eq:Mset}. Its individuals of $V_{G\mathbf s}$ are called
\textbf{free} and those of $V_{H\mathbf s^\dagger}$ \textbf{pinned} (Figure~\ref{fig:gamma}).
Let $\Ppin(\mathbf s,H,\mathbf M)$ be the set of the 
$\Pai\in\Pcal(\Gcup)$ with no component of $\Gcup[\Pai]$ holding two pinned nodes. The other pairings of $\Gcup$ contribute $0$ in
Lemma~\ref{lem:leak}(b): a self-loop does not occur in Proposition~\ref{prop:pairid}, and a
component holding two pinned nodes $u\ne u'$ makes the indicator of that proposition vanish, as
$\sigma(u)\ne\sigma(u')$.

Write
\begin{equation}\label{eq:tau}
\tau_{\mathbf s,H}\ :=\ \hat s_G\,\frac{K^{\,|V_G|-1}}{|\Pi_G|}\Big(\frac{\sqrt d}{\Delta}\Big)^{\!|E_{G\mathbf n}|}
\Big(\prod_{v,w}\binom{m_{v,w}}{\mathbf s_{v,w}}\Big)
\left(\frac{\prod_{v,w}\mathbf n_{v,w}!}{\prod_{u,w}m^H_{u,w}!}\right)^{\!1/2} ,
\end{equation}
a number depending on $G,\mathbf s,H$ alone.

\begin{lemma}\label{lem:leak}
Fix a splitting $\mathbf s$ of $G$, a graph $H$, and a block $Q=Q_{H,\pi^0_1,\pi^0_2,\sigma}$ of
$H$.
\begin{enumerate}[label=\textup{(\alph*)},leftmargin=*,itemsep=3pt]
\item  $\big\langle P_{G,\pi,\mathbf s},Q\big\rangle=0$ unless
$\pi\in\Pi_G(\mathbf M,\pi^0)$ for some $\mathbf M\in\mathcal M(\mathbf s,H)$, and if $\mathcal M(\mathbf s,H)\ne\varnothing$ then
\begin{equation}\label{eq:noisematch}
|V_{G\mathbf n}|=|V_{H\mathbf n^\dagger}|,\qquad |E_{G\mathbf n}|=|E_{H\mathbf n^\dagger}|,\qquad
|W_H|\le|W_G| .
\end{equation}

\item  For $\mathbf M\in\mathcal M(\mathbf s,H)$ and $\pi\in\Pi_G(\mathbf M,\pi^0)$,
\begin{equation}\label{eq:leakb}
\big\langle P_{G,\pi,\mathbf s},\,Q\big\rangle\ =\
\tau_{\mathbf s,H}\!\!\sum_{\Pai\in\Ppin(\mathbf s,H,\mathbf M)}\!\!K^{\,|\CC_{\Gcup[\Pai]}(V_{H\mathbf s^\dagger})|-|V_{G\mathbf s}|} ,
\end{equation}
$\Gcup$ being the merged graph of $\mathbf M$. In particular $\langle P_{G,\pi,\mathbf s},Q\rangle$ depends on $\pi$
only through $\mathbf M$.
\end{enumerate}
\end{lemma}

\begin{proof}
As $\epsilon$ is independent of $(X,k^*)$, the inner product is a scalar times an expectation
over $\epsilon$ and an expectation over $(X,k^*)$. The first vanishes unless the matching of $G$ with
$H$ defined by $\pi$ and $\pi^0$ restricts to an isomorphism $G_{\mathbf n}\to H_{\mathbf n^\dagger}$,
and the second vanishes when this matching leaves a variable of $H$ unmatched, which gives
\textup{(a)}. Otherwise the second is the expectation of $f_{\Gcup}$, which
Proposition~\ref{prop:pairid} and Lemma~\ref{lem:colourcount} turn into \eqref{eq:leakb}.

By \eqref{eq:PGs}, $P_{G,\pi,\mathbf s}$ is a scalar times the product of its
\emph{$\epsilon$-part} and its \emph{$X$-part},
\[
  \prod_{v,w}\frac{\He_{\mathbf n_{v,w}}(\epsilon_{\pi_1(v),\pi_2(w)})}{\sqrt{\mathbf n_{v,w}!}}
  \qquad\text{and}\qquad
  \prod_{v,w}\frac{\He_{\mathbf s_{v,w}}(X_{k^*(\pi_1(v)),\pi_2(w)})}{\sqrt{\mathbf s_{v,w}!}} .
\]
By Definition~\ref{def:block}, $Q$ is the product of its $\epsilon$-part and its $X$-part,
\[
  \prod_{v\in V_{H\mathbf n^\dagger}}\prod_w
  \frac{\He_{m^H_{v,w}}(\epsilon_{\pi^0_1(v),\pi^0_2(w)})}{\sqrt{m^H_{v,w}!}}
  \qquad\text{and}\qquad
  \prod_{u,w}\frac{\He_{m^H_{u,w}}(X_{\sigma(u),\pi^0_2(w)})}{\sqrt{m^H_{u,w}!}} .
\]
As $\epsilon$ is independent of $(X,k^*)$,
\begin{equation}\label{eq:leakfact}
  \big\langle P_{G,\pi,\mathbf s},Q\big\rangle\ =\ (\text{scalar})\cdot
  \E\big[(\epsilon\text{-part of }P_{G,\pi,\mathbf s})(\epsilon\text{-part of }Q)\big]\cdot
  \E\big[(X\text{-part of }P_{G,\pi,\mathbf s})(X\text{-part of }Q)\big] .
\end{equation}
Let $\mathbf M$ be the matching of $G$ with $H$ given by
$\mathbf M_V=\{(v,v'):\pi_1(v)=\pi^0_1(v')\}$, $v'\in V_{H\mathbf n^\dagger}$, and
$\mathbf M_W=\{(w,w'):\pi_2(w)=\pi^0_2(w')\}$, $w'\in W_H$. By \eqref{eq:PiM},
$\pi\in\Pi_G(\mathbf M,\pi^0)$, and $\mathbf M$ is the only matching with this property.

\textbf{(a)} \uline{Step 1: the $\epsilon$-parts.} Both $\epsilon$-parts are of the form $\mathcal H_\alpha$ of
\eqref{eq:hermprod}, $\pi_1,\pi_2,\pi^0_1,\pi^0_2$ being injective. So the second factor of
\eqref{eq:leakfact} is $1$ if the two multi-indices are equal, and $0$ otherwise. Assume they
are equal: for every $(v,w)$ with $\mathbf n_{v,w}\ge1$ there are $v'\in V_{H\mathbf n^\dagger}$
and $w'\in W_H$ with $(v,v')\in\mathbf M_V$, $(w,w')\in\mathbf M_W$ and
$m^H_{v',w'}=\mathbf n_{v,w}$, and, conversely, for every $v'\in V_{H\mathbf n^\dagger}$ and
$w'\in W_H$ with $m^H_{v',w'}\ge1$ there are $v,w$ with these properties. The nodes of $G_{\mathbf
n}$ are $v_1,v_2$ and the nodes carrying a noise edge, and those of $H_{\mathbf n^\dagger}$ are
$v_1,v_2$ and the nodes carrying an edge at $V_{H\mathbf n^\dagger}$ (an individual of
$V_{H\mathbf n^\dagger}\setminus\{v_1,v_2\}$ has degree at least $p$). Moreover
$(v_1,v_1),(v_2,v_2)\in\mathbf M_V$, as $\pi_1(v_i)=i=\pi^0_1(v_i)$. So every node of
$G_{\mathbf n}$ is matched with a node of $H_{\mathbf n^\dagger}$ and conversely. Then, $\mathbf M$
being injective, it restricts to a bijection between the two node sets, which preserves the
multiplicities and fixes $v_1,v_2$:
\[
  \mathbf M|_{G_{\mathbf n}}\ \text{is an isomorphism}\ G_{\mathbf n}\to H_{\mathbf n^\dagger},
  \qquad \mathbf M_V\ \subseteq\ V_{G\mathbf n}\times V_{H\mathbf n^\dagger} ,
\]
the inclusion because each $v'\in V_{H\mathbf n^\dagger}$ occurs once in $\mathbf M_V$, with its
partner in $V_{G\mathbf n}$.

\noindent\uline{Step 2: a variable of $H$ left unmatched.} Assume the multi-indices are equal and some
$w'\in W_H$ is unmatched. By Step~1, $w'$ is not a node of $H_{\mathbf n^\dagger}$, so its edges,
at least one, sit at individuals $u\in V_{H\mathbf s^\dagger}$. Its column $\pi^0_2(w')$ is
$\pi_2(w)$ for no $w\in W_G$, so the entries of that column in the two $X$-parts are the
$X_{\sigma(u),\pi^0_2(w')}$, $u\in N(w')$, pairwise distinct ($\sigma$ being injective) and all in
the $X$-part of $Q$. Conditionally on $k^*$ they are independent of the other entries of the two
$X$-parts, so, by Lemma~\ref{lem:wick}\textup{(i)} ($m^H_{u,w'}\ge1$),
\[
  \E\big[(X\text{-part of }P_{G,\pi,\mathbf s})(X\text{-part of }Q)\,\big|\,k^*\big]\ =\
  \E\big[\cdots\,\big|\,k^*\big]\cdot\prod_{u\in N(w')}
  \frac{\E\big[\He_{m^H_{u,w'}}(X_{\sigma(u),\pi^0_2(w')})\big]}{\sqrt{m^H_{u,w'}!}}\ =\ 0 ,
\]
the dots standing for the product of the other factors.

\noindent\uline{Step 3: conclusion of (a).} By \eqref{eq:leakfact} and Steps~1--2,
$\langle P_{G,\pi,\mathbf s},Q\rangle\ne0$ forces the three conditions of \eqref{eq:Mset}, i.e.\
$\mathbf M\in\mathcal M(\mathbf s,H)$, and $\pi\in\Pi_G(\mathbf M,\pi^0)$. For any
$\mathbf M\in\mathcal M(\mathbf s,H)$, the isomorphism $\mathbf M|_{G_{\mathbf n}}$ gives
$|V_{G\mathbf n}|=|V_{H\mathbf n^\dagger}|$ and $|E_{G\mathbf n}|=|E_{H\mathbf n^\dagger}|$, and
$\mathbf M_W$, injective with image $W_H$, gives $|W_H|\le|W_G|$: this is \eqref{eq:noisematch}.

\textbf{(b)} \uline{Step 4: the scalar.} Let $\mathbf M\in\mathcal M(\mathbf s,H)$ and
$\pi\in\Pi_G(\mathbf M,\pi^0)$. By \eqref{eq:PiM} $\pi$ puts each node of $G_{\mathbf n}$ where
$\pi^0$ puts its image under the isomorphism $\mathbf M|_{G_{\mathbf n}}$, so the two
$\epsilon$-parts have the same multi-index and the second factor of \eqref{eq:leakfact} is $1$.
The scalar of \eqref{eq:leakfact} is the first line of \eqref{eq:PGs} times
$\prod_{v,w}\sqrt{\binom{m_{v,w}}{\mathbf s_{v,w}}}(\Delta/\sqrt d)^{\mathbf s_{v,w}}$. With
$\sum_{v,w}\mathbf s_{v,w}=|E_G|-|E_{G\mathbf n}|$ and
$m_{v,w}!=\mathbf s_{v,w}!\,\mathbf n_{v,w}!\binom{m_{v,w}}{\mathbf s_{v,w}}$ it equals
\[
  \hat s_G\,\frac{K^{\,|V_G|-1}}{|\Pi_G|}\Big(\frac{\sqrt d}{\Delta}\Big)^{\!|E_{G\mathbf n}|}
  \Big(\prod_{v,w}\binom{m_{v,w}}{\mathbf s_{v,w}}\Big)
  \Big(\prod_{v,w}\mathbf s_{v,w}!\Big)^{\!1/2}\Big(\prod_{v,w}\mathbf n_{v,w}!\Big)^{\!1/2} .
\]
Moving the normalizations $(\prod_{v,w}\mathbf s_{v,w}!)^{-1/2}$ and
$(\prod_{u,w}m^H_{u,w}!)^{-1/2}$ of the two $X$-parts into it gives $\tau_{\mathbf s,H}$
\eqref{eq:tau}:
\begin{equation}\label{eq:leaktau}
  \big\langle P_{G,\pi,\mathbf s},Q\big\rangle\ =\ \tau_{\mathbf s,H}\;
  \E\Big[\prod_{v,w}\He_{\mathbf s_{v,w}}\big(X_{k^*(\pi_1(v)),\pi_2(w)}\big)
  \prod_{u,w}\He_{m^H_{u,w}}\big(X_{\sigma(u),\pi^0_2(w)}\big)\Big] .
\end{equation}

\noindent\uline{Step 5: the expectation is $\E f_{\Gcup}$.} The individuals of $\Gcup$ are those of
$V_{G\mathbf s}$ and of $V_{H\mathbf s^\dagger}$, none identified as
$\mathbf M_V\subseteq V_{G\mathbf n}\times V_{H\mathbf n^\dagger}$. Color them by
$\chi(v):=k^*(\pi_1(v))$, $v\in V_{G\mathbf s}$, and $\chi(u):=\sigma(u)$,
$u\in V_{H\mathbf s^\dagger}$. By \eqref{eq:PiM}, $\pi_2(w)=\pi^0_2(w')$ exactly when
$(w,w')\in\mathbf M_W$, so the columns occurring in \eqref{eq:leaktau} correspond one to one to
the variables of $\Gcup$, and the factors in the column of a variable are one
$\He$ of the multiplicity in $\Gcup$ for each individual of $\Gcup$ at that variable, in the row
of its color. Grouping the factors of each column by color, as in Step~2 of the proof of
Proposition~\ref{prop:master} (conditionally on $k^*$ the entries of $X$ are i.i.d.\
$\mathcal N(0,1)$),
\[
  \E\Big[\prod_{v,w}\He_{\mathbf s_{v,w}}\big(X_{k^*(\pi_1(v)),\pi_2(w)}\big)
  \prod_{u,w}\He_{m^H_{u,w}}\big(X_{\sigma(u),\pi^0_2(w)}\big)\,\Big|\,k^*\Big]
  \ =\ f_{\Gcup}(\chi) .
\]

\noindent\uline{Step 6: the pairings.} By Proposition~\ref{prop:pairid} and
Corollary~\ref{cor:selfloop}, $f_{\Gcup}(\chi)$ is the sum, over $\Pai\in\Pcal(\Gcup)$, of the indicator that $\chi$ is constant on each component of $\Gcup[\Pai]$. A
component holding two pinned nodes $u\ne u'$ makes that indicator $0$, as
$\sigma(u)\ne\sigma(u')$. So the sum runs over $\Ppin(\mathbf s,H,\mathbf M)$. Fix $\Pai\in\Ppin(\mathbf s,H,\mathbf M)$. Every
individual of $\Gcup$ carries an edge (a free node carries a signal edge, and a pinned node is not $v_1,v_2$, so is not isolated). Hence each component of $\Gcup[\Pai]$ holds an individual of
$V_{G\mathbf s}$ (a pinned node is linked to an individual that is not pinned). The sets
$V_C\cap V_{G\mathbf s}$, $C\in\CC_{\Gcup[\Pai]}$, partition $V_{G\mathbf s}$, and the indicator
says that $k^*\circ\pi_1$ is constant on each of them and equals $\sigma(u)$ on that of a
component holding the pinned node $u$. The $k^*(\pi_1(v))$, $v\in V_{G\mathbf s}$, being i.i.d.\
uniform, Lemma~\ref{lem:colourcount}, whose parts without prescription are those of the
components in $\CC_{\Gcup[\Pai]}(V_{H\mathbf s^\dagger})$, gives
\[
  \E\big[f_{\Gcup[\Pai]}(\chi)\big]\ =\ K^{\,|\CC_{\Gcup[\Pai]}(V_{H\mathbf s^\dagger})|-|V_{G\mathbf s}|} .
\]
Taking expectations in Step~5 and substituting into \eqref{eq:leaktau} gives
\eqref{eq:leakb}. Its right-hand side involves $G$,
$\mathbf s$, $H$ and $\Gcup$ only, and $\Gcup$ is determined by $\mathbf M$.
\end{proof}

\subsection{The centering cancels exactly}\label{sec:cancel}

Here the set $\bar V$ of \eqref{eq:Sdef} consists of the pinned nodes and the mixed individuals:
\[
\bar V\ :=\ V_{H\mathbf s^\dagger}\cup\big(V_{G\mathbf s}\cap V_{G\mathbf n}\big),
\qquad\text{so that}\qquad
\mathcal V^{(\mathrm{free})}_{\bar V}(\Gcup[\Pai])\ =\!\!\bigcup_{\substack{C\in\CC_{\Gcup[\Pai]}(V_{H\mathbf s^\dagger})\\ V_C\cap V_{G\mathbf n}=\varnothing}}
\!\!V_C\ \subseteq\ V_{G\mathbf s} ,
\]
the second form because the individuals of a free component, one missing $V_{H\mathbf s^\dagger}$ (\S\ref{sec:expansion}), lie in $V_{G\mathbf s}$, so that it misses the mixed individuals exactly when it misses $V_{G\mathbf n}$. \textbf{Free} still refers to $V_{H\mathbf s^\dagger}$, not to $\bar V$.
Every individual of $\mathcal V^{(\mathrm{free})}_{\bar V}(\Gcup[\Pai])$ is purely signal: by the second form, the set lies in $V_{G\mathbf s}$ and misses $V_{G\mathbf n}$. In particular $v_1,v_2\notin \mathcal V^{(\mathrm{free})}_{\bar V}(\Gcup[\Pai])$, both lying in $V_{G\mathbf n}$. The set is
a union of components of $\Gcup[\Pai]$, so no pair of $\Pai$ joins it to its complement.  Finally, $\mathcal V^{(\mathrm{free})}_{\bar V}(\Gcup[\Pai])=\varnothing$ says exactly that every free
component of $\Gcup[\Pai]$ holds a mixed individual. Put
\begin{equation}\label{eq:MP}
\Ppins(\mathbf s,H,\mathbf M)\ :=\ \big\{\Pai\in\Ppin(\mathbf s,H,\mathbf M)\ :\ \mathcal V^{(\mathrm{free})}_{\bar V}(\Gcup[\Pai])=\varnothing\big\}
\qquad(\mathbf M\in\mathcal M(\mathbf s,H)) ,
\end{equation}
$\Gcup$ being the merged graph of $\mathbf M$.

\begin{theorem}\label{thm:cancel}
Fix $G$ in $\mathcal G_{\le D}\cup\{\varnothing\}$, a splitting $\mathbf s$ of $G$, a graph $H$ in $\mathcal G_{\le D}$ and a
block $Q=Q_{H,\pi^0_1,\pi^0_2,\sigma}$ of $H$, with $\pi^0=(\pi^0_1,\pi^0_2)$. Then
\begin{equation}\label{eq:cancelQ}
\big\langle\Psh_{G,\mathbf s},\,Q\big\rangle\;=\;\tau_{\mathbf s,H}
\!\!\sum_{\mathbf M\in\mathcal M(\mathbf s,H)}\ \sum_{\Pai\in\Ppins(\mathbf s,H,\mathbf M)}
\!\!|\Pi_G(\mathbf M,\pi^0)|\,K^{\,|\CC_{\Gcup[\Pai]}(V_{H\mathbf s^\dagger})|-|V_{G\mathbf s}|} .
\end{equation}
\end{theorem}
\begin{proof}
The proof is an induction on $|E_G|+|V_G|$, parallel to that of Proposition~\ref{prop:num}. By
Lemma~\ref{lem:leak}, $\langle P_{G,\mathbf s},Q\rangle$ is the right-hand side of \eqref{eq:cancelQ}
with $\Ppins(\mathbf s,H,\mathbf M)$ replaced by $\Ppin(\mathbf s,H,\mathbf M)$. Grouped by the set
$V'=\mathcal V^{(\mathrm{free})}_{\bar V}(\Gcup[\Pai])$, the extra pairs $(\mathbf M,\Pai)$ add up to the
terms $\Lambda_{V'}\langle\Psh_{\Gs,\mathbf s|_{\Gs}},Q\rangle$ that the centering \eqref{eq:Pshs}
subtracts. $H$ and $Q$ stay fixed throughout.

\noindent\textbf{Induction hypothesis.} For an integer $m\ge2$, let $(\mathrm H_m)$ be the
statement:
\begin{center}
\emph{\eqref{eq:cancelQ} holds for every $G\in\mathcal G_{\le D}\cup\{\varnothing\}$ with
$|E_G|+|V_G|\le m$ and every splitting $\mathbf s$ of $G$.}
\end{center}
Every graph has $|V_G|\ge2$, as $v_1,v_2\in V_G$. So it suffices to prove $(\mathrm H_m)$ for
every $m\ge2$, by induction on $m$.

\noindent\textbf{Initialization: $(\mathrm H_2)$.} As in the proof of Proposition~\ref{prop:num},
the only graph concerned is $\varnothing$, with its one splitting. There $\Psh_{\varnothing,\mathbf s}=1$ (\S\ref{sec:expansion}), and
$\E[Q]=0$ by \eqref{eq:hermprod}, $Q$ being non-constant. So the left-hand side of
\eqref{eq:cancelQ} is $0$. $H$ has an edge, hence a variable, and $\varnothing$ has none, so no
matching of $\varnothing$ with $H$ has $W_H$ as image of $\mathbf M_W$:
$\mathcal M(\mathbf s,H)=\varnothing$, and the right-hand side is $0$ as well.

\noindent\textbf{Induction step: $(\mathrm H_{m-1})\Rightarrow(\mathrm H_m)$, for $m\ge3$.}
Assume $(\mathrm H_{m-1})$, and let $G\in\mathcal G_{\le D}$ have $|E_G|+|V_G|=m$ and $\mathbf s$
be a splitting of $G$. By Lemma~\ref{lem:deletable}(a), for every $V'\in\mathcal D(G)$ the graph
$\Gs$ is $\varnothing$ or a member of $\mathcal G_{\le D}$, with $|E_{\Gs}|+|V_{\Gs}|\le m-1$. So,
by $(\mathrm H_{m-1})$, \eqref{eq:cancelQ} holds for $\Gs$ with the splitting $\mathbf s|_{\Gs}$,
for every $V'\in\mathcal D(G)$ purely signal in $\mathbf s$. Steps~1--6 below prove
\eqref{eq:cancelQ} for $G$ and $\mathbf s$, which gives $(\mathrm H_m)$.

\noindent\uline{Step 1: reduction.} By Lemma~\ref{lem:leak}(a), $\langle P_{G,\pi,\mathbf s},Q\rangle=0$
unless $\pi$ lies in one of the sets $\Pi_G(\mathbf M,\pi^0)$, $\mathbf M\in\mathcal M(\mathbf
s,H)$, which are pairwise disjoint, $\mathbf M$ being determined by $\pi$
\eqref{eq:PiM}. By
Lemma~\ref{lem:leak}(b),
\begin{equation}\label{eq:cancelP}
  \big\langle P_{G,\mathbf s},Q\big\rangle\ =\ \sum_{\pi\in\Pi_G}\big\langle P_{G,\pi,\mathbf s},Q\big\rangle
  \ =\ \tau_{\mathbf s,H}\!\!\sum_{\substack{\mathbf M\in\mathcal M(\mathbf s,H)\\ \Pai\in\Ppin(\mathbf s,H,\mathbf M)}}
  \!\!|\Pi_G(\mathbf M,\pi^0)|\,K^{\,|\CC_{\Gcup[\Pai]}(V_{H\mathbf s^\dagger})|-|V_{G\mathbf s}|} .
\end{equation}
Sort the pairs $(\mathbf M,\Pai)$ of \eqref{eq:cancelP} by $\mathcal V^{(\mathrm{free})}_{\bar V}(\Gcup[\Pai])$. Those
with $\mathcal V^{(\mathrm{free})}_{\bar V}(\Gcup[\Pai])=\varnothing$ are those with $\Pai\in\Ppins(\mathbf s,H,\mathbf M)$ \eqref{eq:MP}, and their
part of \eqref{eq:cancelP} is the right-hand side of \eqref{eq:cancelQ}. By \eqref{eq:Pshs} it therefore suffices to prove:
\begin{enumerate}[label=\textup{(\roman*)},leftmargin=*,itemsep=1pt,topsep=3pt]
\item every non-empty $\mathcal V^{(\mathrm{free})}_{\bar V}(\Gcup[\Pai])$, $\mathbf M\in\mathcal M(\mathbf s,H)$,
$\Pai\in\Ppin(\mathbf s,H,\mathbf M)$, lies in $\mathcal D(G)$ and is purely signal in $\mathbf s$,
\item for every $V'\in\mathcal D(G)$ purely signal in $\mathbf s$,
\begin{equation}\label{eq:cancelclass}
  \tau_{\mathbf s,H}\!\!\sum_{\substack{\mathbf M\in\mathcal M(\mathbf s,H),\
  \Pai\in\Ppin(\mathbf s,H,\mathbf M)\\ \mathcal V^{(\mathrm{free})}_{\bar V}(\Gcup[\Pai])=V'}}
  \!\!|\Pi_G(\mathbf M,\pi^0)|\,K^{\,|\CC_{\Gcup[\Pai]}(V_{H\mathbf s^\dagger})|-|V_{G\mathbf s}|}
  \ =\ \Lambda_{V'}\,\big\langle\Psh_{\Gs,\,\mathbf s|_{\Gs}},Q\big\rangle .
\end{equation}
\end{enumerate}
Indeed, by \eqref{eq:Pshs}, \eqref{eq:cancelP}, (i) and (ii), $\langle\Psh_{G,\mathbf s},Q\rangle$
is \eqref{eq:cancelP} minus the sum over $V'$ of the right-hand sides of \eqref{eq:cancelclass},
that is, minus the part of \eqref{eq:cancelP} with $\mathcal V^{(\mathrm{free})}_{\bar V}(\Gcup[\Pai])\ne\varnothing$.

\noindent\uline{Step 2: proof of (i).} Let $V':=\mathcal V^{(\mathrm{free})}_{\bar V}(\Gcup[\Pai])\ne\varnothing$. By the paragraph
before \eqref{eq:MP}, $V'$ is a union of components of $\Gcup[\Pai]$ holding no pinned node, its
individuals are purely signal and $v_1,v_2\notin V'$. So the edges of $G$ at $V'$ are those of
$G_{V'}$, all in $\Gcup$. As in Step~2 of the proof of Lemma~\ref{lem:deletable}(b), the pairs of
$\Pai$ formed by these edges make an element of $\Pcal(G_{V'})$, the two edges of such a pair sitting at one
variable of $G$, as a variable of $\Gcup$ carries at most one variable of $G$. So $V'\in\mathcal D(G)$.

\noindent\uline{Step 3: $G$ and $\Gs$.} Fix $V'\in\mathcal D(G)$ purely signal in $\mathbf s$. As
$V'\cap V_{G\mathbf n}=\varnothing$, every noise edge of $G$ sits at an individual of $\Gs$. So
$G$ with $\mathbf s$ and $\Gs$ with $\mathbf s|_{\Gs}$ have the same noise part $G_{\mathbf n}$,
the same mixed individuals and hence the same $\bar V$, and the signal individuals of $\Gs$
are those of $V_{G\mathbf s}\setminus V'$. The three conditions of \eqref{eq:Mset} involve $G$ only through 
$G_{\mathbf n}$ and $V_{G\mathbf n}$. So a matching of $\Gs$ with $H$ lies in
$\mathcal M(\mathbf s|_{\Gs},H)$ exactly when, as a matching of $G$ with $H$, it lies in
$\mathcal M(\mathbf s,H)$. For such an $\mathbf M$ let $\Gamma'_{\cup}$ be its merged graph for
$\Gs$. The edges of $\Gcup$ are the signal edges of $G$ and the edges of $H_{\mathbf s^\dagger}$,
and those of $\Gamma'_{\cup}$ are the same minus the edges of $G_{V'}$. Two edges not at $V'$ sit at
one variable of $\Gcup$ exactly when they sit at one variable of $\Gamma'_{\cup}$: for two edges of
$G$, or two of $H$, this means one variable of $G$, resp.\ of $H$, and for an edge of $G$ at $w$ and
one of $H$ at $w'$, it means $(w,w')\in\mathbf M_W$ in both graphs, $w$ carrying a signal edge
not at $V'$. Hence, for a pairing $\Pai$ of $\Gcup$ in which no pair joins an edge at $V'$ to an
edge not at $V'$, the pairs formed by edges at $V'$ make a pairing $\Pai_{V'}$ of $G_{V'}$, the
other pairs make a pairing $\Pai'$ of $\Gamma'_{\cup}$, every couple $(\Pai_{V'},\Pai')$ arises
from exactly one such $\Pai$, and
\begin{equation}\label{eq:cancelsplit}
  \Gcup[\Pai]\ =\ G_{V'}[\Pai_{V'}]\ \sqcup\ \Gamma'_{\cup}[\Pai'] .
\end{equation}

\noindent\uline{Step 4: the class of $V'$ is a product.} Let $\mathbf M\in\mathcal M(\mathbf s,H)$ and
$\Pai\in\Ppin(\mathbf s,H,\mathbf M)$ with $\mathcal V^{(\mathrm{free})}_{\bar V}(\Gcup[\Pai])=V'$. First, $\mathbf M$ is a matching of
$\Gs$: its individuals lie in $V_{G\mathbf n}$, which misses $V'$, and every
$(w,w')\in\mathbf M_W$ has $w\in W_{\Gs}$. If $w'$ carries an edge at
$V_{H\mathbf n^\dagger}$, then $w$ carries a noise edge, $\mathbf M|_{G_{\mathbf n}}$ being an
isomorphism onto $H_{\mathbf n^\dagger}$, and that edge is not at $V'$. Otherwise $w'$ carries an
edge $e$ of $H_{\mathbf s^\dagger}$, at a pinned node $u$. Its pair under $\Pai$ is an edge at
the same variable of $\Gcup$, at an individual $v\ne u$ (no self-loop) which is not pinned (no
component with two pinned nodes) and not in $V'$ ($v$ lies in the component of $u$, and $V'$ holds
no pinned node). So that edge is an edge of $G$ at $w$, not at $V'$. In both cases
$w\in W_{\Gs}$. By Step~3, $\mathbf M\in\mathcal M(\mathbf s|_{\Gs},H)$.

Next, as $V'$ is a union
of components of $\Gcup[\Pai]$, no pair of $\Pai$ joins an edge at $V'$ to an edge not at $V'$. So, by Step~3, $\Pai$ is a pairing $\Pai_{V'}\in\Pcal(G_{V'})$ together with a
pairing $\Pai'$ of $\Gamma'_{\cup}$, and \eqref{eq:cancelsplit} holds. $\Pai'$ has no self-loop, as $\Pai$ has none, and the components of $\Gamma'_{\cup}[\Pai']$, being
components of $\Gcup[\Pai]$, hold at most one pinned node, and, as those missing
$\bar V$ make up $V'$, none misses $\bar V$: $\Pai'\in\Ppins(\mathbf s|_{\Gs},H,\mathbf M)$. And
\begin{equation}\label{eq:cancelexp}
  |\CC_{\Gcup[\Pai]}(V_{H\mathbf s^\dagger})|-|V_{G\mathbf s}|
  \ =\ \big(|\CC_{G_{V'}[\Pai_{V'}]}|-|V'|\big)
  +\big(|\CC_{\Gamma'_{\cup}[\Pai']}(V_{H\mathbf s^\dagger})|-|V_{G\mathbf s}\setminus V'|\big) .
\end{equation}
Conversely, by Step~3, $\Pai_{V'}\in\Pcal(G_{V'})$ and
$\Pai'\in\Ppins(\mathbf s|_{\Gs},H,\mathbf M)$, $\mathbf M\in\mathcal M(\mathbf s|_{\Gs},H)$, make a
pairing $\Pai\in\Ppin(\mathbf s,H,\mathbf M)$ with $\mathcal V^{(\mathrm{free})}_{\bar V}(\Gcup[\Pai])=V'$ ($\mathbf M\in\mathcal M(\mathbf s,H)$ by
Step~3, and $G_{V'}[\Pai_{V'}]$ holds no pinned and no mixed individual). The two maps are inverse
to each other.

\noindent\uline{Step 5: the labelings and the scalar.} Let $\mathbf M\in\mathcal M(\mathbf s|_{\Gs},H)$. As
$\mathbf M$ matches all of $V_{H\mathbf n^\dagger}$ and $W_H$ with nodes of $\Gs$ \eqref{eq:Mset}, the elements
of $\Pi_G(\mathbf M,\pi^0)$ are the injective extensions of those of $\Pi_{\Gs}(\mathbf M,\pi^0)$ to $V'$ and
$W_G\setminus W_{\Gs}$, which number $|\Pi_G|/|\Pi_{\Gs}|$ \eqref{eq:Pisize}. So
\begin{equation}\label{eq:cancelcount}
  |\Pi_G(\mathbf M,\pi^0)|\ =\ \frac{|\Pi_G|}{|\Pi_{\Gs}|}\;|\Pi_{\Gs}(\mathbf M,\pi^0)| .
\end{equation}
In \eqref{eq:tau}, the individuals of $V'$, being purely signal, contribute
$\binom{m_{v,w}}{\mathbf s_{v,w}}=\binom{m_{v,w}}{m_{v,w}}=1$ and $\mathbf n_{v,w}!=0!=1$, and
$G$ and $\Gs$ have the same $E_{G\mathbf n}$. With $|V_G|-|V_{\Gs}|=|V'|$,
\begin{equation}\label{eq:cancelratio}
  \tau_{\mathbf s,H}\ =\ \frac{\hat s_G}{\hat s_{\Gs}}\,K^{|V'|}\,\frac{|\Pi_{\Gs}|}{|\Pi_G|}\;
  \tau_{\mathbf s|_{\Gs},H} .
\end{equation}

\noindent\uline{Step 6: proof of (ii).} By Step~4 and \eqref{eq:cancelexp}--\eqref{eq:cancelratio}, the
left-hand side of \eqref{eq:cancelclass} equals
\[
\begin{aligned}
  &\frac{\hat s_G}{\hat s_{\Gs}}K^{|V'|}\!\!\sum_{\Pai_{V'}\in\Pcal(G_{V'})}\!\!K^{\,|\CC_{G_{V'}[\Pai_{V'}]}|-|V'|}\\
  &\qquad\cdot\tau_{\mathbf s|_{\Gs},H}\!\!\!\!\sum_{\mathbf M\in\mathcal M(\mathbf s|_{\Gs},H)}\ \sum_{\Pai'\in\Ppins(\mathbf s|_{\Gs},H,\mathbf M)}
  \!\!\!\!|\Pi_{\Gs}(\mathbf M,\pi^0)|\,
  K^{\,|\CC_{\Gamma'_{\cup}[\Pai']}(V_{H\mathbf s^\dagger})|-|V_{G\mathbf s}\setminus V'|} .
\end{aligned}
\]
The first factor is $\Lambda_{V'}$ by \eqref{eq:LambdaS}. The second is the right-hand side of
\eqref{eq:cancelQ} for $\Gs$ with the splitting $\mathbf s|_{\Gs}$ (whose signal individuals are
$V_{G\mathbf s}\setminus V'$, Step~3), hence equals $\langle\Psh_{\Gs,\mathbf s|_{\Gs}},Q\rangle$
by $(\mathrm H_{m-1})$. This is (ii).

\end{proof}

\subsection{The pairings $\Ppins(\mathbf s,H,\mathbf M)$}\label{sec:comb}

\subsubsection{Constraints on $(\mathbf s,H)$ for pairings in $\Ppins(\mathbf s,H,\mathbf M)$}

\begin{lemma}\label{lem:move}
Let $\mathbf M\in\mathcal M(\mathbf s,H)$ and $\Pai\in\Ppin(\mathbf s,H,\mathbf M)$, not necessarily with
$\mathcal V^{(\mathrm{free})}_{\bar V}(\Gcup[\Pai])=\varnothing$.
\begin{enumerate}[label=\textup{(\alph*)},leftmargin=*,itemsep=2pt]
\item  Every edge of $H_{\mathbf s^\dagger}$ is paired with an
edge of $G_{\mathbf s}$, and 
$|E_G|-|E_H|$ is the number of edges of $G_{\mathbf s}$ paired with an edge
of $G_{\mathbf s}$. In particular $|E_H|\le|E_G|$.
\item  Each of the
$|V_{H\mathbf s^\dagger}|$ pinned nodes lies in a component holding at least one free node and
no other pinned node, each of the
$|\CC_{\Gcup[\Pai]}(V_{H\mathbf s^\dagger})|$ free components holds at least two free nodes, and these components are disjoint. Hence
\[
|V_{H\mathbf s^\dagger}|+2|\CC_{\Gcup[\Pai]}(V_{H\mathbf s^\dagger})|\ \le\ |V_{G\mathbf s}|,\qquad\text{and}\qquad
|V_G|-|V_H|\ =\ |V_{G\mathbf s}|-|V_{G\mathbf s}\cap V_{G\mathbf n}|-|V_{H\mathbf s^\dagger}| ,
\]
so that $|V_G|-|V_H|\ge2|\CC_{\Gcup[\Pai]}(V_{H\mathbf s^\dagger})|-|V_{G\mathbf s}\cap V_{G\mathbf n}|$.
\item Every edge at an individual of a free component of $\Gcup[\Pai]$ is an edge of
$G_{\mathbf s}$ paired with an edge of $G_{\mathbf s}$.
\end{enumerate}
\end{lemma}

\begin{proof}
We first show that an edge of a pinned node is paired with an edge of $G_{\mathbf s}$ at a free
node of the same component. Counting the pairs of $\Pai$ then gives \textup{(a)}, counting the free
nodes of each component of $\Gcup[\Pai]$ gives \textup{(b)}, and \textup{(c)} holds because the two
edges of a pair lie in one component.

\noindent\textbf{Preliminary.} Let $\{e,e'\}$ be a pair of $\Pai$, $e$ being an edge of
$H_{\mathbf s^\dagger}$ at the pinned node $u$. The individual $v$ of $e'$ is not $u$, as $\Pai$
has no self-loop, and is not pinned, as it lies in the component of $u$ and
$\Pai\in\Ppin(\mathbf s,H,\mathbf M)$. So $v\in V_{G\mathbf s}$ and $e'$ is an edge of $G_{\mathbf s}$.

\noindent\textbf{(a)} By the preliminary, a pair of $\Pai$ is made either of an edge of
$H_{\mathbf s^\dagger}$ and an edge of $G_{\mathbf s}$, or of two edges of $G_{\mathbf s}$. Hence,
 as
$|E_{G\mathbf n}|=|E_{H\mathbf n^\dagger}|$ \eqref{eq:noisematch},
\[
  |E_G|-|E_H|\ =\ |E_{G\mathbf s}|-|E_{H\mathbf s^\dagger}|
  \ =\ \#\big\{e\in E_{G\mathbf s}:\ \{e,e'\}\in\Pai\ \text{for some}\ e'\in E_{G\mathbf s}\big\}\ \ge\ 0 .
\]

\noindent\textbf{(b)} Every individual of $V_{G\mathbf s}$ carries a signal edge, and an
individual of $V_{H\mathbf s^\dagger}$ is not $v_1,v_2$, which lie in $V_{H\mathbf n^\dagger}$, so is not isolated. So every individual of $\Gcup$ carries an
edge, and is linked in $\Gcup[\Pai]$ to an individual distinct from it. Hence, by the preliminary
and as $\Pai\in\Ppin(\mathbf s,H,\mathbf M)$, the component of a pinned node holds a free node and
no other pinned node, and a free component holds an individual (every variable of $\Gcup[\Pai]$
having degree $2$), linked to a second one, both free. This gives the inequality. The identity
follows from $|V_H|=|V_{H\mathbf n^\dagger}|+|V_{H\mathbf s^\dagger}|$,
$|V_G|=|V_{G\mathbf n}|+|V_{G\mathbf s}|-|V_{G\mathbf s}\cap V_{G\mathbf n}|$ and
$|V_{G\mathbf n}|=|V_{H\mathbf n^\dagger}|$ \eqref{eq:noisematch}, and the last inequality from
the two.

\noindent\textbf{(c)} An edge at an individual of a free component $C$ is paired with an edge at
an individual of $C$, which is free. So the edges at the individuals of $C$ are edges of
$G_{\mathbf s}$ paired with edges of $G_{\mathbf s}$.
\end{proof}

If $\Ppins(\mathbf s,H,\mathbf M)=\varnothing$ for every $\mathbf M\in\mathcal M(\mathbf s,H)$, \eqref{eq:cancelQ} vanishes and
$(\mathbf s,H)$ may be ignored. Otherwise set
\begin{equation}\label{eq:complexity}
\begin{gathered}
\miCC\ :=\ |V_{G\mathbf s}\cap V_{G\mathbf n}|-\!\!\max_{\mathbf M\in\mathcal M(\mathbf s,H),\ \Pai\in\Ppins(\mathbf s,H,\mathbf M)}\!\!|\CC_{\Gcup[\Pai]}(V_{H\mathbf s^\dagger})| ,\\[2pt]
\comp\ :=\ \big(|E_G|-|E_H|\big)\ +\ \big(|V_G|-|V_H|\big)_+\ +\ \miCC ,
\end{gathered}
\end{equation}
and call $\comp$ the \textbf{complexity} of $(\mathbf s,H)$. All three summands are non-negative
(Lemma~\ref{lem:complexity}(a)), so $|E_G|-|E_H|\le\comp$.

\begin{lemma}\label{lem:complexity}
Let $G,H\in\mathcal G_{\le D}$, and let $\mathbf s$ be a splitting of $G$ with $\Ppins(\mathbf s,H,\mathbf M)\ne\varnothing$ for some $\mathbf M\in\mathcal M(\mathbf s,H)$.
\begin{enumerate}[label=\textup{(\alph*)},leftmargin=*,itemsep=2pt]
\item For every $\mathbf M\in\mathcal M(\mathbf s,H)$ and $\Pai\in\Ppins(\mathbf s,H,\mathbf M)$,
$|\CC_{\Gcup[\Pai]}(V_{H\mathbf s^\dagger})|\le|V_{G\mathbf s}\cap V_{G\mathbf n}|$, so $\miCC\ge0$ and
$\comp\ge0$, and
\[
\begin{gathered}
|V_{G\mathbf s}\cap V_{G\mathbf n}|\ \le\ \big(|E_G|-|E_H|\big)+\miCC\ \le\ \comp,\\[3pt]
\sum_{v\in V_{G\mathbf s}\cap V_{G\mathbf n}}\deg_{G\mathbf s}(v)\ \le\ \comp+\lceil p\rceil \miCC .
\end{gathered}
\]
\item $\big(|E_G|+|V_G|\big)-\big(|E_H|+|V_H|\big)\ \ge\ \comp-2\miCC$, and every
$\mathbf M\in\mathcal M(\mathbf s,H)$ and $\Pai\in\Ppins(\mathbf s,H,\mathbf M)$ satisfy
$|V_{G\mathbf s}|-|V_{H\mathbf s^\dagger}|-|\CC_{\Gcup[\Pai]}(V_{H\mathbf s^\dagger})|\le \comp$.
\item $\comp\ge1$, unless $(\mathbf s,H)$ is the \textbf{diagonal} pair: $H=G$ and
$\mathbf s=\mathbf s^\dagger(G)$.
\end{enumerate}
\end{lemma}
\begin{proof}
The proof rests on Lemma~\ref{lem:move} and on the property that defines
$\Ppins(\mathbf s,H,\mathbf M)$: every free component of $\Gcup[\Pai]$ holds a mixed individual.
Counting the mixed individuals and their signal edges against the pairs of $\Pai$ gives
\textup{(a)}, and the identity of Lemma~\ref{lem:move}\textup{(b)} gives \textup{(b)}. When
$\comp=0$, each component of $\Gcup[\Pai]$ holds one free and one pinned node, which yields an
isomorphism $G\to H$ for \textup{(c)}.

\noindent\textbf{(a)} \uline{Step 1: positivity.} Let $\mathbf M\in\mathcal M(\mathbf s,H)$ and
$\Pai\in\Ppins(\mathbf s,H,\mathbf M)$. Every free component $C$ satisfies $V_C\subseteq V_{G\mathbf s}$ and
$V_C\cap V_{G\mathbf n}\ne\varnothing$ (paragraph before \eqref{eq:MP}), so it holds a mixed individual. The
free components being disjoint,
\[
  |\CC_{\Gcup[\Pai]}(V_{H\mathbf s^\dagger})|\ \le\
  \#\{\text{mixed individuals in free components of }\Gcup[\Pai]\}\ \le\
  |V_{G\mathbf s}\cap V_{G\mathbf n}| .
\]
Hence $\miCC\ge0$, and $\comp\ge0$ by Lemma~\ref{lem:move}(a).

\noindent\uline{Step 2: the number of mixed individuals.} Fix $\mathbf M\in\mathcal M(\mathbf s,H)$ and
$\Pai\in\Ppins(\mathbf s,H,\mathbf M)$ attaining the maximum in \eqref{eq:complexity}, and let $F$ be the set of
mixed individuals in its free components. By Step~1, $|F|\ge|\CC_{\Gcup[\Pai]}(V_{H\mathbf s^\dagger})|$, so
$|(V_{G\mathbf s}\cap V_{G\mathbf n})\setminus F|\le\miCC$. Each $v\in F$ carries a signal edge, which is an
edge of $G_{\mathbf s}$ paired with an edge of $G_{\mathbf s}$ (Lemma~\ref{lem:move}(c)), distinct $v$ giving
distinct edges. So $|F|\le|E_G|-|E_H|$ (Lemma~\ref{lem:move}(a)). Adding,
\[
  |V_{G\mathbf s}\cap V_{G\mathbf n}|\ \le\ \big(|E_G|-|E_H|\big)+\miCC\ \le\ \comp .
\]

\noindent\uline{Step 3: the signal degree of the mixed individuals.} With $(\mathbf M,\Pai)$ and $F$ as in
Step~2, a signal edge at a mixed individual $v$ is paired
\begin{itemize}[itemsep=1pt,topsep=2pt,leftmargin=*]
\item either with an edge of $G_{\mathbf s}$: there are at most $|E_G|-|E_H|$ such edges
(Lemma~\ref{lem:move}(a)),
\item or with an edge of $H_{\mathbf s^\dagger}$ at the pinned node $u$ of the component of $v$: then
$v\notin F$, there are at most $\miCC$ such $v$ (Step~2), and the edges of $v$ of this kind are paired
with distinct edges of $u$, so number at most $\deg_H(u)\le\lceil p\rceil-1$, $u$ being a signal
individual of $H$ \eqref{eq:blocksplit}.
\end{itemize}
Adding,
\[
  \sum_{v\in V_{G\mathbf s}\cap V_{G\mathbf n}}\deg_{G\mathbf s}(v)\ \le\ \big(|E_G|-|E_H|\big)
  +(\lceil p\rceil-1)\,\miCC\ \le\ \comp+\lceil p\rceil\,\miCC .
\]

\noindent\textbf{(b)} Both relations follow from the identity
$|V_G|-|V_H|=|V_{G\mathbf s}|-|V_{G\mathbf s}\cap V_{G\mathbf n}|-|V_{H\mathbf s^\dagger}|$ of
Lemma~\ref{lem:move}(b). By \eqref{eq:complexity},
\[
  \big(|E_G|+|V_G|\big)-\big(|E_H|+|V_H|\big)-\big(\comp-2\miCC\big)\ =\ \miCC-\big(|V_H|-|V_G|\big)_+ ,
\]
so, as $\miCC\ge0$, the first relation amounts to $|V_H|-|V_G|\le\miCC$. At the pair of Step~2 of
\textup{(a)}, Lemma~\ref{lem:move}(b) and the choice of that pair give
\[
\begin{aligned}
  |V_G|-|V_H|\ &\ge\ 2\,|\CC_{\Gcup[\Pai]}(V_{H\mathbf s^\dagger})|-|V_{G\mathbf s}\cap V_{G\mathbf n}|\\
  &=\ |\CC_{\Gcup[\Pai]}(V_{H\mathbf s^\dagger})|-\miCC\ \ge\ -\miCC .
\end{aligned}
\]
For the second, let $\mathbf M\in\mathcal M(\mathbf s,H)$ and $\Pai\in\Ppins(\mathbf s,H,\mathbf M)$. By the
identity, by Step~2 of \textup{(a)}, and as $|\CC_{\Gcup[\Pai]}(V_{H\mathbf s^\dagger})|\ge0$,
\[
\begin{aligned}
  &|V_{G\mathbf s}|-|V_{H\mathbf s^\dagger}|-|\CC_{\Gcup[\Pai]}(V_{H\mathbf s^\dagger})|\\
  &\qquad=\ \big(|V_G|-|V_H|\big)+|V_{G\mathbf s}\cap V_{G\mathbf n}|
  -|\CC_{\Gcup[\Pai]}(V_{H\mathbf s^\dagger})|\\
  &\qquad\le\ \big(|V_G|-|V_H|\big)_++\big(|E_G|-|E_H|\big)+\miCC\ =\ \comp .
\end{aligned}
\]

\noindent\textbf{(c)} Assume $\comp=0$. By \eqref{eq:complexity} and \textup{(a)}, $|E_G|=|E_H|$,
$|V_G|\le|V_H|$ and $\mathbf s$ has no mixed individual. So $G$, like $H$ \eqref{eq:blocksplit}, is
the union of its two parts, which share no individual. Fix $\mathbf M\in\mathcal M(\mathbf s,H)$ and
$\Pai\in\Ppins(\mathbf s,H,\mathbf M)$. By Lemma~\ref{lem:move}(b),
\[
  0\ \ge\ |V_G|-|V_H|\ =\ |V_{G\mathbf s}|-|V_{H\mathbf s^\dagger}|\ \ge\
  2\,|\CC_{\Gcup[\Pai]}(V_{H\mathbf s^\dagger})|\ \ge\ 0 ,
\]
so each component of $\Gcup[\Pai]$ holds exactly one free node and one pinned node. By
Lemma~\ref{lem:move}(a) and $|E_G|=|E_H|$, $\Pai$ pairs the edges of $G_{\mathbf s}$ with those of
$H_{\mathbf s^\dagger}$, hence the edges of a free node at a variable $w$ bijectively with those of
the pinned node of its component at the variable merged with $w$. So $\mathbf M$, extended by sending
each free node to the pinned node of its component, is an isomorphism $G\to H$ (each variable of $G$
is matched: by $\mathbf M|_{G_{\mathbf n}}$ if it carries a noise edge, through $\Pai$ otherwise).
Hence $H=G$, as $\mathcal G_{\le D}$ holds one graph per isomorphism class, and this automorphism of
$G$ maps $V_{G\mathbf n}$ onto $V_{G\mathbf n^\dagger}$ and preserves $V_{G\mathbf n^\dagger}$
\eqref{eq:blocksplit}, so $V_{G\mathbf n}=V_{G\mathbf n^\dagger}$ and $\mathbf s=\mathbf s^\dagger(G)$
\eqref{eq:sdagger}.
\end{proof}

\subsubsection{Bounding the size of some sets and of some factors}

\begin{lemma}\label{lem:count}
For every $e\le D$, the isomorphism classes of graphs with $|E_G|=e$ number at most $(3D^2)^e$.
\end{lemma}
\begin{proof}
As $|V|\le e+2$ and $|W|\le e$, every class contains a graph with individual set in $[D+2]$,
$v_1,v_2$ being $1,2$, and variable set in $[D]$. Such a graph is determined by the multiset of
the pairs $(v,w)$ of its $e$ edges, a multiset of size $e$ in $[D+2]\times[D]$. So the classes
number at most
\[
  \binom{(D+2)D+e-1}{e}\ \le\ \big(D^{2}+3D\big)^{e}\ \le\ \big(3D^{2}\big)^{e} .
\]
\end{proof}

\begin{lemma}\label{lem:counts}
Fix $G\in\mathcal G_{\le D}$, and consider the pairs $(\mathbf s,H)$ of a splitting $\mathbf s$ of $G$ and a graph
$H\in\mathcal G_{\le D}$ such that $\Ppins(\mathbf s,H,\mathbf M)\ne\varnothing$ for some
$\mathbf M\in\mathcal M(\mathbf s,H)$, with $\comp$ and $\miCC$ as in \eqref{eq:complexity}.
\begin{enumerate}[label=\textup{(\alph*)},leftmargin=*,itemsep=2pt]
\item Every such pair satisfies
$\displaystyle\prod_{v,w}\binom{m_{v,w}}{\mathbf s_{v,w}}\ \le\ D^{\,\comp+\lceil p\rceil \miCC}$.
\item Every such pair satisfies
$\displaystyle\frac{\sum_{\mathbf M\in\mathcal M(\mathbf s,H)}\big|\Ppins(\mathbf s,H,\mathbf M)\big|}{|\Aut(H)|\,\prod_{u,w}m^H_{u,w}!}\ \le\ D^{6\comp}$.
\item For all integers $c,m\ge0$, at most $D^{\,14c+2\lceil p\rceil m}$ such pairs have $\comp=c$ and
$\miCC=m$.
\end{enumerate}
\end{lemma}

\begin{proof}
Only mixed individuals contribute to the product in \textup{(a)}, so
Lemma~\ref{lem:complexity}\textup{(a)} bounds it. For \textup{(b)}, each pair $(\mathbf M,\Pai)$ is
encoded by its shadow, defined below: one shadow comes from at most
$|\Aut(H)|\prod_{u,w}m^H_{u,w}!$ pairs, and at most $D^{6\comp}$ shadows occur. For \textup{(c)}, a
shadow determines $H$, and $\mathbf s$ is determined by three data whose numbers of values are
bounded in terms of $\comp$ and $\miCC$. Let $(\mathbf s,H)$ be such a pair.

\noindent\textbf{(a)} If $\mathbf s_{v,w}\in\{0,m_{v,w}\}$ then $\binom{m_{v,w}}{\mathbf s_{v,w}}=1$.
Otherwise $v$ carries a signal and a noise edge at $w$, so is mixed, and
$\binom{m_{v,w}}{\mathbf s_{v,w}}\le m_{v,w}^{\mathbf s_{v,w}}\le D^{\mathbf s_{v,w}}$, as
$m_{v,w}\le|E_G|\le D$. Multiplying, and by Lemma~\ref{lem:complexity}(a),
\[
  \prod_{v,w}\binom{m_{v,w}}{\mathbf s_{v,w}}\ \le\
  D^{\sum_{v\in V_{G\mathbf s}\cap V_{G\mathbf n}}\sum_w\mathbf s_{v,w}}
  \ =\ D^{\sum_{v\in V_{G\mathbf s}\cap V_{G\mathbf n}}\deg_{G\mathbf s}(v)}
  \ \le\ D^{\,\comp+\lceil p\rceil\,\miCC} .
\]

\noindent\textbf{(b)} We follow the shadow argument of \cite[\S9.3]{carpentier2026low}. For $\mathbf M\in\mathcal M(\mathbf s,H)$ and
$\Pai\in\Ppins(\mathbf s,H,\mathbf M)$, call an edge of $G_{\mathbf s}$ \textbf{perfectly paired} when $\Pai$ pairs
it with an edge of $H_{\mathbf s^\dagger}$, and let $\underline\Pai\subseteq\Pai$ be the set of the other pairs,
all made of two edges of $G_{\mathbf s}$ (Lemma~\ref{lem:move}(a)). Let $\underline{PE}$ be the partition of the
perfectly paired edges according to the pinned node whose edge each of them is paired with. The part of a
pinned node $u$ has $\deg_H(u)$ elements (Lemma~\ref{lem:move}(a)), and $\deg_H(u)\ge1$. The \textbf{shadow} of $(\mathbf M,\Pai)$
is $(\underline\Pai,\underline{PE})$: it forgets $\mathbf M$ and which pinned node each part belongs to.

\noindent\uline{Step 1: the shadow determines $H$.} Let $H'$ be the graph $G_{\mathbf n}$ together with
one further individual per part of $\underline{PE}$, joined to each variable $w$ by as many edges as that part
has at $w$. Then $\mathbf M$, extended by sending the part of each pinned node $u$ to $u$, is an isomorphism
$\psi:H'\to H$, as in the proof of Lemma~\ref{lem:complexity}(c): $\Pai$ pairs the edges of the part of $u$ at
$w$ bijectively with those of $u$ at the variable merged with $w$, and every variable of $H'$, and every
variable of $H$ \eqref{eq:Mset}, is matched with a variable of the other graph, through
$\mathbf M|_{G_{\mathbf n}}$ if it carries a noise edge and through $\Pai$ otherwise. So $H$ is the member of
$\mathcal G_{\le D}$ isomorphic to $H'$.

\noindent\uline{Step 2: one shadow, few pairs $(\mathbf M,\Pai)$.} This is the analogue of \cite[Lemma~13]{carpentier2026low}.
Given the shadow, $(\mathbf M,\Pai)$ is determined by $\psi$ and, for each pinned node $u$ and each variable $w$
of $H'$, by a bijection from the edges of the part of $u$ at $w$ onto the $m^H_{u,\psi(w)}$ edges of $u$ at
$\psi(w)$: $\mathbf M$ is the restriction of $\psi$ to $V_{G\mathbf n}$ and to the variables of $H'$ (by Step~1,
no other variable is matched), and $\Pai$ is $\underline\Pai$ together with these bijections. As $\psi$ is one of
the $|\Aut(H)|$ isomorphisms $H'\to H$, the pairs $(\mathbf M,\Pai)$, $\mathbf M\in\mathcal M(\mathbf s,H)$,
$\Pai\in\Ppins(\mathbf s,H,\mathbf M)$, with shadow $(\underline\Pai,\underline{PE})$ number at most
\[
  |\Aut(H)|\,\prod_{u,w}m^H_{u,w}! .
\]

\noindent\uline{Step 3: the number of shadows.} Fix $\mathbf s$ and an integer $c\ge0$. We show that the shadows
of the pairs $(\mathbf M,\Pai)$, over all graphs $H$ with $\comp=c$, number at most $D^{6c}$.
\begin{itemize}[itemsep=1pt,topsep=2pt,leftmargin=*]
\item $\underline\Pai$ is a set of $j=(|E_G|-|E_H|)/2\le c/2$ pairs of edges of $G_{\mathbf s}$
(Lemma~\ref{lem:move}(a)), chosen among at most $|E_G|^2\le D^2$ pairs: at most
$\sum_{0\le j\le c/2}D^{2j}\le(1+\tfrac c2)D^{c}\le D^{2c}$ choices, as $1+\tfrac c2\le2^{c}\le D^{c}$.
\item $\underline{PE}$ is determined by $\underline\Pai$ and by the partition of $V_{G\mathbf s}$ into the sets
$V_C\cap V_{G\mathbf s}$, $C\in\CC_{\Gcup[\Pai]}$: the part of $u$ is the set of the edges at the individuals of
$V_C\cap V_{G\mathbf s}$, $C$ the component of $u$, that lie in no pair of $\underline\Pai$, while those of a
free component all lie in pairs of $\underline\Pai$ (Lemma~\ref{lem:move}(c)). By Lemma~\ref{lem:move}(b), this
partition has $|V_{H\mathbf s^\dagger}|+|\CC_{\Gcup[\Pai]}(V_{H\mathbf s^\dagger})|$ parts, so it is obtained
from the partition into singletons by
$|V_{G\mathbf s}|-|V_{H\mathbf s^\dagger}|-|\CC_{\Gcup[\Pai]}(V_{H\mathbf s^\dagger})|\le c$ merges of two
parts (Lemma~\ref{lem:complexity}(b)), hence by a sequence of exactly $c$ steps, each merging the parts that
hold the two individuals of an ordered pair of $V_{G\mathbf s}$ (a step changes nothing when they lie in one
part). As $|V_{G\mathbf s}|^2\le4D^2\le D^4$, there are at most $D^{4c}$ such sequences.
\end{itemize}
Multiplying, the shadows number at most $D^{2c}D^{4c}=D^{6c}$. With $c=\comp$ and Step~2,
\[
  \sum_{\mathbf M\in\mathcal M(\mathbf s,H)}\big|\Ppins(\mathbf s,H,\mathbf M)\big|\ \le\ D^{6\comp}\,|\Aut(H)|\,\prod_{u,w}m^H_{u,w}! ,
\]
which is \textbf{(b)}.

\noindent\textbf{(c)} Let $c,m\ge0$. For a fixed $\mathbf s$, every $H$ with $\comp=c$ is determined by a shadow of
Step~3 of \textup{(b)} (Step~1), so these $H$ number at most $D^{6c}$. Now let $(\mathbf s,H)$ have $\comp=c$
and $\miCC=m$, and fix $\mathbf M\in\mathcal M(\mathbf s,H)$ and $\Pai\in\Ppins(\mathbf s,H,\mathbf M)$. An
individual other than $v_1,v_2$ which is neither mixed nor purely signal carries noise edges only, so
$\mathbf M|_{G_{\mathbf n}}$ maps it to a noise individual of $H$ other than $v_1,v_2$, of the same degree, which is
at least $\lceil p\rceil$ \eqref{eq:blocksplit}. Hence an individual other than $v_1,v_2$ of degree at most $\lceil p\rceil-1$ is mixed or purely signal, and $v_1,v_2$ carry no signal edge unless they are mixed. So
$\mathbf s$ is determined by the three data below, counted with $|V_G|\le2D$.
\begin{enumerate}[label=\textup{(\arabic*)},itemsep=1pt,topsep=3pt,leftmargin=*]
\item The set of mixed individuals. It has at most $c$ elements (Lemma~\ref{lem:complexity}(a)), so there are
at most $\sum_{k\le c}\binom{2D}k\le(c+1)(2D)^{c}\le4^{c}D^{c}\le D^{3c}$ choices, as $c+1\le2^{c}$ and
$4\le D^2$.
\item The set of purely signal individuals of degree at least $\lceil p\rceil$. Such a $v$ has an edge of
$G_{\mathbf s}$ paired with an edge of $G_{\mathbf s}$, by Lemma~\ref{lem:move}(c) if $v$ lies in a free
component, and otherwise because the only edges of $H_{\mathbf s^\dagger}$ paired with edges of $v$ are the
$\deg_H(u)\le\lceil p\rceil-1<\deg_G(v)$ edges of the pinned node $u$ of its component. These edges being
distinct, the set has at most $|E_G|-|E_H|\le c$ elements (Lemma~\ref{lem:move}(a)), so there are at most
$D^{3c}$ choices, as in \textup{(1)}.
\item The values $\mathbf s_{v,w}$ at the at most $|E_G|\le D$ pairs $(v,w)$ with $v$ mixed and
$m_{v,w}\ge1$. They are non-negative integers of sum at most $t:=c+\lceil p\rceil m$
(Lemma~\ref{lem:complexity}(a)), so, by stars and bars, there are at most
$\binom{D+t}{t}=\prod_{i=1}^{t}\frac{D+i}{i}\le(D+1)^{t}\le D^{2t}=D^{\,2c+2\lceil p\rceil m}$ choices.
\end{enumerate}
Multiplying the four counts, the pairs $(\mathbf s,H)$ number at most
$D^{3c}D^{3c}D^{2c+2\lceil p\rceil m}D^{6c}=D^{\,14c+2\lceil p\rceil m}$.
\end{proof}

\subsection{The matrix of interest for balancing correlation and projection intensity on the blocks}\label{sec:matrix}

 Put
\begin{equation}\label{eq:eta}
\eta_G\ :=\ \hat s_G\,D^{50(|E_G|+|V_G|)} .
\end{equation}

\subsubsection{The entries, and the diagonal}\label{sec:entries}

For $G,H\in\mathcal G_{\le D}$ put
\begin{equation}\label{eq:Adef}
A_{H,G}\ :=\ \frac{\eta_H}{\eta_G}\,\big\langle\Psh_G,\ \bar Q_H\big\rangle
\ =\ \frac{\hat s_H}{\hat s_G}\,D^{-50\sdropb}\,\big\langle\Psh_G,\ \bar Q_H\big\rangle
\end{equation}
by \eqref{eq:eta}.
 For a splitting $\mathbf s$ of $G$
put $A^{\mathbf s}_{H,G}:=\frac{\eta_H}{\eta_G}\langle\Psh_{G,\mathbf s},\bar Q_H\rangle$, so that,
by $\Psh_G=\sum_{\mathbf s}\Psh_{G,\mathbf s}$,
$A_{H,G}=\sum_{\mathbf s}A^{\mathbf s}_{H,G}$.

As the edges at a variable $w$ have $(\deg_Gw-1)!!$ partitions into pairs,
$\#\Pcal(G)\le\prod_w(\deg_Gw-1)!!\le(|E_G|-1)!!\le D^{|E_G|}$,
so Proposition~\ref{prop:num} and \eqref{eq:eta} give
\begin{equation}\label{eq:rbound}
\big|\E\big[(x-\E x)\,\Psh_G\big]\big|\ \le\ \hat s_G\,D^{|E_G|}\ \le\ \eta_G .
\end{equation}
\begin{proposition}\label{prop:entry}
For $G,H\in\mathcal G_{\le D}$ and every splitting $\mathbf s$ of $G$, the entry $A^{\mathbf s}_{H,G}$ lies
between $\tfrac12$ and $2$ times
\begin{equation}\label{eq:entry}
\frac{D^{-50\sdropb}}{|\Aut(H)|\ \prod_{u,w}m^H_{u,w}!}
\Big(\prod_{v,w}\binom{m_{v,w}}{\mathbf s_{v,w}}\Big)
\!\!\!\!\sum_{\mathbf M\in\mathcal M(\mathbf s,H)}\ \sum_{\Pai\in\Ppins(\mathbf s,H,\mathbf M)}\!\!\!\!K^{\,|\CC_{\Gcup[\Pai]}(V_{H\mathbf s^\dagger})|-|V_{G\mathbf s}\cap V_{G\mathbf n}|} .
\end{equation}
\end{proposition}

\begin{proof}
Theorem~\ref{thm:cancel} gives the inner product of $\Psh_{G,\mathbf s}$ with one block of $H$,
and its factor $|\Pi_G(\mathbf M,\pi^0)|$ depends neither on $\mathbf M$ nor on the block. So the
inner product with the invariant block $\bar Q_H$ is an explicit multiple of it, and collecting the
factors of $A^{\mathbf s}_{H,G}$ leaves \eqref{eq:entry} times three ratios of counts whose product
lies between $\tfrac12$ and $2$.

If $\mathcal M(\mathbf s,H)=\varnothing$, both $A^{\mathbf s}_{H,G}$ (by \eqref{eq:cancelQ}) and
\eqref{eq:entry} vanish. Assume $\mathcal M(\mathbf s,H)\ne\varnothing$.

\noindent\uline{Step 1: the labelings following a matching.} Let $Q=Q_{H,\pi^0_1,\pi^0_2,\sigma}$ be a block
of $H$ and $\mathbf M\in\mathcal M(\mathbf s,H)$. A $\pi\in\Pi_G(\mathbf M,\pi^0)$ is prescribed on the matched
nodes, which are the individuals of $V_{G\mathbf n}$ and $|W_H|$ variables \eqref{eq:Mset}, and places the
others injectively in the remaining rows and columns \eqref{eq:PiM}, so
\begin{equation}\label{eq:entrylab}
  |\Pi_G(\mathbf M,\pi^0)|\ =\ \prod_{k=|V_{G\mathbf n}|}^{|V_G|-1}(n-k)\cdot
  \prod_{k=|W_H|}^{|W_G|-1}(d-k) ,
\end{equation}
the same number for every $\mathbf M$ and every block, between $(\tfrac78)^2$ and $1$ times
$n^{|V_G|-|V_{G\mathbf n}|}d^{|W_G|-|W_H|}$ by \eqref{eq:Picount}.

\noindent\uline{Step 2: the invariant block.} By \eqref{eq:cancelQ} and Step~1, for any $\mathbf M\in\mathcal M(\mathbf s,H)$,
\[
  \big\langle\Psh_{G,\mathbf s},Q\big\rangle\ =\ \tau_{\mathbf s,H}\,|\Pi_G(\mathbf M,\pi^0)|
  \!\!\sum_{\mathbf M'\in\mathcal M(\mathbf s,H)}\ \sum_{\Pai\in\Ppins(\mathbf s,H,\mathbf M')}\!\!
  K^{\,|\CC_{\Gcup[\Pai]}(V_{H\mathbf s^\dagger})|-|V_{G\mathbf s}|} ,
\]
which does not depend on the block $Q$. Hence, by \eqref{eq:Qbarnorm},
\[
  \big\langle\Psh_{G,\mathbf s},\bar Q_H\big\rangle\ =\
  \frac{\sum_{\pi^0_1,\pi^0_2,\sigma}\big\langle\Psh_{G,\mathbf s},Q_{H,\pi^0_1,\pi^0_2,\sigma}\big\rangle}
  {\big\|\sum_{\pi^0_1,\pi^0_2,\sigma}Q_{H,\pi^0_1,\pi^0_2,\sigma}\big\|}
  \ =\ \Big(\frac{\#\{(\pi^0_1,\pi^0_2,\sigma)\}}{|\Aut(H)|}\Big)^{1/2}
  \big\langle\Psh_{G,\mathbf s},Q\big\rangle .
\]
Finally, by \eqref{eq:Qbarnorm}, and by \eqref{eq:Pisize} and \eqref{eq:Picount},
\[
  \frac12\ \le\ \frac{\#\{(\pi^0_1,\pi^0_2,\sigma)\}}{n^{|V_{H\mathbf n^\dagger}|-2}d^{|W_H|}K^{|V_{H\mathbf s^\dagger}|}}\ \le\ 1
  \qquad\text{and}\qquad
  \Big(\frac78\Big)^2\ \le\ \frac{|\Pi_G|}{n^{|V_G|-2}d^{|W_G|}}\ \le\ 1 .
\]

\noindent\uline{Step 3: collecting the factors.} Multiply by
$\eta_H/\eta_G=\frac{\hat s_H}{\hat s_G}D^{-50\sdropb}$ \eqref{eq:Adef}, and insert
\eqref{eq:shat} for $\hat s_H$ and \eqref{eq:tau} for $\tau_{\mathbf s,H}$. We claim
\begin{equation}\label{eq:entryratio}
  A^{\mathbf s}_{H,G}\ =\
  \Big(\frac{\#\{(\pi^0_1,\pi^0_2,\sigma)\}}{n^{|V_{H\mathbf n^\dagger}|-2}d^{|W_H|}K^{|V_{H\mathbf s^\dagger}|}}\Big)^{\!1/2}
  \cdot\frac{|\Pi_G(\mathbf M,\pi^0)|}{n^{|V_G|-|V_{G\mathbf n}|}d^{|W_G|-|W_H|}}
  \cdot\frac{n^{|V_G|-2}d^{|W_G|}}{|\Pi_G|}\cdot\big(\text{the quantity \eqref{eq:entry}}\big) .
\end{equation}
Indeed, by \eqref{eq:noisematch}, and as the isomorphism $\mathbf M|_{G_{\mathbf n}}$ makes the factorials of
$\hat s_H$ at $V_{H\mathbf n^\dagger}\times W_H$ equal to the $\mathbf n_{v,w}!$, the factors $\hat s_G$, the powers
of $\Delta/\sqrt d$, $n$, $d$ and the noise factorials of $\hat s_H$ and $\tau_{\mathbf s,H}$ cancel, leaving
$|\Aut(H)|^{-1}\big(\prod_{u,w}m^H_{u,w}!\big)^{-1}$. By
$|V_G|=|V_{G\mathbf n}|+|V_{G\mathbf s}|-|V_{G\mathbf s}\cap V_{G\mathbf n}|$, the exponents of $K$ add up to
$|\CC_{\Gcup[\Pai]}(V_{H\mathbf s^\dagger})|-|V_{G\mathbf s}\cap V_{G\mathbf n}|$.

\noindent\uline{Step 4: the constants.} By Steps~1 and~2, the three ratios of \eqref{eq:entryratio} lie in
$[2^{-1/2},1]$, $[(\tfrac78)^2,1]$ and $[1,(\tfrac87)^2]$, so their product lies in
\[
  \big[2^{-1/2}\big(\tfrac78\big)^{2},\ \big(\tfrac87\big)^{2}\big]\ \subset\ \big[\tfrac12,2\big] .
  \qedhere
\]
\end{proof}

\begin{corollary}\label{cor:diag}
For every $G\in\mathcal G_{\le D}$, $A_{G,G}\ \ge\ \tfrac12$.
\end{corollary}

\begin{proof}
We bound $A_{G,G}$ below by the term of the block splitting $\mathbf s^\dagger(G)$, where the
automorphisms of $G$ alone give enough pairs $(\mathbf M,\Pai)$ to make \eqref{eq:entry} at least
$1$. By Proposition~\ref{prop:entry}, $A^{\mathbf s}_{G,G}\ge0$ for every splitting $\mathbf s$ of $G$, so
$A_{G,G}=\sum_{\mathbf s}A^{\mathbf s}_{G,G}\ge A^{\mathbf s^\dagger(G)}_{G,G}$, which is at least $\tfrac12$ times
\eqref{eq:entry} at $H=G$, $\mathbf s=\mathbf s^\dagger(G)$. There $\sdrop=0$, there is no mixed individual
\eqref{eq:blocksplit}, every binomial is $1$ and every summand is at least $1$. So it suffices to show
\[
  \sum_{\mathbf M\in\mathcal M(\mathbf s^\dagger(G),G)}\big|\Ppins(\mathbf s^\dagger(G),G,\mathbf M)\big|
  \ \ge\ |\Aut(G)|\prod_{u\in V_{G\mathbf s^\dagger},\,w\in W_G}m_{u,w}! .
\]
Each $\varphi\in\Aut(G)$ gives $\prod_{v\in V_{G\mathbf s^\dagger},\,w\in W_G}m_{v,w}!$ such pairs $(\mathbf M,\Pai)$: $\mathbf M$ is the restriction of
$\varphi$ to $V_{G\mathbf n^\dagger}\cup W_G$, which lies in $\mathcal M(\mathbf s^\dagger(G),G)$ as $\varphi$ preserves
$V_{G\mathbf n^\dagger}$ \eqref{eq:blocksplit}, and $\Pai$ pairs the edges of $G_{\mathbf s}$ joining $v$ to $w$
with those of $H_{\mathbf s^\dagger}$ joining $\varphi(v)$ to $\varphi(w)$, so that the components of
$\Gcup[\Pai]$ are the sets $\{v,\varphi(v)\}$ of one free and one pinned node. Distinct automorphisms give
disjoint sets of pairs, as $\varphi$ sends each free node to the pinned node of its component, which proves
the claim.
\end{proof}

\subsubsection{The column sum}\label{sec:colsum}

\begin{lemma}\label{lem:colsum}
For every $G\in\mathcal G_{\le D}$,
\[
  \sum_{(\mathbf s,H)\,\ne\,(\mathbf s^\dagger(G),G)}A^{\mathbf s}_{H,G}\ \le\ \frac15 ,
\]
the sum running over the pairs of a splitting $\mathbf s$ of $G$ and a graph $H\in\mathcal G_{\le D}$, other
than the diagonal one. In particular $\sum_{H\ne G}A_{H,G}\le\tfrac15$, every $A^{\mathbf s}_{H,G}$ being
non-negative (Proposition~\ref{prop:entry}).
\end{lemma}

\begin{proof}
We bound each entry by Proposition~\ref{prop:entry} and
Lemma~\ref{lem:counts}\textup{(a)},\textup{(b)}, and the number of pairs by
Lemma~\ref{lem:counts}\textup{(c)}, so that the pairs with $\comp=c$ and $\miCC=m$ contribute at
most $2D^{-29c-m}$ in total. As the diagonal pair is the only one with $\comp=0$, the sum over
$c\ge1$ and $m\ge0$ is at most $8D^{-29}$.

By Proposition~\ref{prop:entry}, only the pairs with $\Ppins(\mathbf s,H,\mathbf M)\ne\varnothing$ for some
$\mathbf M\in\mathcal M(\mathbf s,H)$ contribute. For such a pair, with $\comp$ and $\miCC$ as in
\eqref{eq:complexity}, the exponent of $K$ in \eqref{eq:entry} is at most $-\miCC$, and
Lemma~\ref{lem:complexity}(b) gives
\[
  \sdrop\ \ge\ \comp-2\miCC .
\]
So Proposition~\ref{prop:entry} and Lemma~\ref{lem:counts}(a),(b) give
\[
\begin{aligned}
  A^{\mathbf s}_{H,G}\ &\le\ 2\,D^{-50(\comp-2\miCC)}\,D^{\,\comp+\lceil p\rceil\miCC}\,D^{6\comp}\,K^{-\miCC}\\
  &=\ 2\,D^{-43\comp+(\lceil p\rceil+100)\miCC}\,K^{-\miCC} .
\end{aligned}
\]
For integers $c,m\ge0$, Lemma~\ref{lem:counts}(c) gives at most $D^{14c+2\lceil p\rceil m}$ pairs with
$\comp=c$ and $\miCC=m$. As $K\ge D^{800p\lceil p\rceil}$ \eqref{eq:dK} and
$800p\lceil p\rceil\ge3\lceil p\rceil+101$, these pairs contribute at most
\[
  2\,D^{-29c+(3\lceil p\rceil+100)m}\,K^{-m}\ \le\ 2\,D^{-29c-m} .
\]
The diagonal pair is the only one with $\comp=0$ (Lemma~\ref{lem:complexity}(c)), and $\miCC\ge0$
(Lemma~\ref{lem:complexity}(a)). Summing over $c\ge1$ and $m\ge0$, with $D\ge2$,
\[
  \sum_{(\mathbf s,H)\ne(\mathbf s^\dagger(G),G)}A^{\mathbf s}_{H,G}\ \le\
  2\sum_{c\ge1}D^{-29c}\sum_{m\ge0}D^{-m}\ \le\ 8D^{-29}\ \le\ \frac15 .
\]
\end{proof}

\subsection{Proof of Theorem~\ref{thm:LB}}\label{sec:solve}

\begin{proposition}\label{prop:reduce}
\[
  \Corr_{\le D}\ \le\ \frac{10}3\Big(\sum_{G\in\mathcal G_{\le D}}\eta_G^2\Big)^{1/2} .
\]
\end{proposition}

\begin{proof}
We expand $P$ along the centered polynomials $\Psh_G$, with coefficients $u_G$. The numerator
$\langle x-\E x,P\rangle$ is at most $\sum_G\eta_G|u_G|$, while
$\sum_H\eta_H|\langle P,\bar Q_H\rangle|$ is at least $\tfrac3{10}\sum_G\eta_G|u_G|$, as the matrix
$A$ has diagonal entries at least $\tfrac12$ and off-diagonal column sums at most $\tfrac15$.
Cauchy--Schwarz and Bessel's inequality for the orthonormal $\bar Q_H$ bound
$\sum_H\eta_H|\langle P,\bar Q_H\rangle|$ by $(\sum_H\eta_H^2)^{1/2}\|P\|$.

\noindent\uline{Step 1: the span.} By Proposition~\ref{prop:sym}, and as isomorphic graphs have the same $P_G$, the
supremum \eqref{eq:corr} may be taken
over $P\in\Span\big(1,(P_G)_{G\in\mathcal G_{\le D}}\big)$, $P\ne0$. This space equals
$\Span\big(1,(\Psh_G)_{G\in\mathcal G_{\le D}}\big)$: by \eqref{eq:Psh} and
Lemma~\ref{lem:deletable}(a), $\Psh_G-P_G$ is a linear combination of $\Psh_\varnothing=1$ and of
$\Psh_{G'}$ with $G'\in\mathcal G_{\le D}$, $|E_{G'}|+|V_{G'}|<|E_G|+|V_G|$, so by induction on
$|E_G|+|V_G|$ each $\Psh_G$ lies in the first space and each $P_G$ in the second. So let
\[
  P\ =\ u_0+\sum_{G\in\mathcal G_{\le D}}u_G\Psh_G\ \ne\ 0 ,
\]
with real numbers $u_0$, $u_G$. We show
$\langle x-\E x,P\rangle\le\tfrac{10}3\big(\sum_G\eta_G^2\big)^{1/2}\|P\|$.

\noindent\uline{Step 2: the numerator.} As $\E[x-\E x]=0$, \eqref{eq:rbound} gives
\begin{equation}\label{eq:numu}
  \big\langle x-\E x,\ P\big\rangle\ =\ \sum_G u_G\,\E\big[(x-\E x)\Psh_G\big]\ \le\
  \sum_G\eta_G|u_G| .
\end{equation}

\noindent\uline{Step 3: the coefficients along the invariant blocks.} Every block is non-constant, so has
mean $0$ \eqref{eq:hermprod}, and $\langle1,\bar Q_H\rangle=0$. So, by \eqref{eq:Adef},
for every $H\in\mathcal G_{\le D}$,
\[
  \eta_H\,\langle P,\bar Q_H\rangle\ =\ \sum_G\eta_Gu_G\,\frac{\eta_H}{\eta_G}\langle\Psh_G,\bar Q_H\rangle
  \ =\ \sum_GA_{H,G}\,\eta_Gu_G .
\]
As $A_{H,G}=\sum_{\mathbf s}A^{\mathbf s}_{H,G}\ge0$ for all $G,H$ (Proposition~\ref{prop:entry}),
\[
  \eta_H\,\big|\langle P,\bar Q_H\rangle\big|\ \ge\ A_{H,H}\,\eta_H|u_H|
  -\sum_{G\ne H}A_{H,G}\,\eta_G|u_G| .
\]

\noindent\uline{Step 4: the column sums.} Summing Step~3 over $H$ (all sums are finite, by
Lemma~\ref{lem:count}) and exchanging the order of summation in the double sum,
\[
  \sum_H\eta_H\big|\langle P,\bar Q_H\rangle\big|\ \ge\
  \sum_G\eta_G|u_G|\Big(A_{G,G}-\sum_{H\ne G}A_{H,G}\Big)\ \ge\
  \Big(\frac12-\frac15\Big)\sum_G\eta_G|u_G| ,
\]
by Corollary~\ref{cor:diag} and Lemma~\ref{lem:colsum}.

\noindent\uline{Step 5: conclusion.} The $\bar Q_H$ being orthonormal (Lemma~\ref{lem:ortho}), Cauchy--Schwarz
and Bessel's inequality give
\[
  \sum_H\eta_H\big|\langle P,\bar Q_H\rangle\big|\ \le\
  \Big(\sum_H\eta_H^2\Big)^{1/2}\Big(\sum_H\langle P,\bar Q_H\rangle^2\Big)^{1/2}\ \le\
  \Big(\sum_H\eta_H^2\Big)^{1/2}\|P\| .
\]
With Step~4, $\sum_G\eta_G|u_G|\le\tfrac{10}3\big(\sum_G\eta_G^2\big)^{1/2}\|P\|$, and with
\eqref{eq:numu} this is the claim.
\end{proof}

\begin{lemma}\label{lem:master}
For every $G\in\mathcal G_{\le D}$,
\[
\hat s_G\ \le\ \frac1K\,D^{-200|E_G|},\qquad\text{and therefore}\qquad
\eta_G\ \le\ \frac1K\,D^{-2|E_G|} .
\]
\end{lemma}

\begin{proof}
We bound $\hat s_G$ by $1/K$ times two factors, one for the signal individuals and one for the
noise individuals, and show that they are at most $D^{-200|E_{G\mathbf s^\dagger}|}$ and
$D^{-200|E_{G\mathbf n^\dagger}|}$. Signal individuals have degree at most $\lceil p\rceil-1$, and
Assumption~\ref{ass:A}\textup{(iii)} bounds the first factor. Noise individuals other than $v_1,v_2$
have degree at least $\lceil p\rceil$, and Assumption~\ref{ass:A}\textup{(i)} and \textup{(ii)} bound
the second.

All quantities below are those of the block splitting \eqref{eq:blocksplit} of $G$, as in \eqref{eq:shat}. By
\eqref{eq:evencons}, $|W|\le|E|/2=(|E_{G\mathbf s^\dagger}|+|E_{G\mathbf n^\dagger}|)/2$, and $d\ge1$. So,
dropping $|\Aut(G)|^{-1/2}\le1$ and $\big(\prod m_{v,w}!\big)^{-1/2}\le1$, \eqref{eq:shat} gives
\[
  \hat s_G\ \le\ \Big[d^{\frac{|E_{G\mathbf s^\dagger}|}4}K^{-\frac{|V_{G\mathbf s^\dagger}|}2}\Big]
  \cdot\Big[\Big(\frac{\Delta}{d^{1/4}}\Big)^{|E_{G\mathbf n^\dagger}|}
  n^{\frac{|V_{G\mathbf n^\dagger}|-2}2}K^{-(|V_{G\mathbf n^\dagger}|-2)}\Big]\cdot\frac1K .
\]
A signal individual has degree at most $\lceil p\rceil-1$ \eqref{eq:blocksplit}, so
$|E_{G\mathbf s^\dagger}|\le(\lceil p\rceil-1)|V_{G\mathbf s^\dagger}|$. With $K=d^{p/2}$ and
$d^{\,p-\lceil p\rceil+1}\ge D^{1600\lceil p\rceil}$
(Assumption~\ref{ass:A}\textup{(iii)}), the first bracket is
\[
  d^{\frac{|E_{G\mathbf s^\dagger}|-p|V_{G\mathbf s^\dagger}|}4}
  \ \le\ d^{-\frac{(p-\lceil p\rceil+1)|V_{G\mathbf s^\dagger}|}{4}}
  \ \le\ D^{-400\lceil p\rceil|V_{G\mathbf s^\dagger}|}
  \ \le\ D^{-200|E_{G\mathbf s^\dagger}|} .
\]
A noise individual other than $v_1,v_2$ has degree at least $\lceil p\rceil$, and all its edges are noise, so
$|E_{G\mathbf n^\dagger}|\ge\lceil p\rceil(|V_{G\mathbf n^\dagger}|-2)$. By Assumption~\ref{ass:A}\textup{(i)},
and by \textup{(ii)}, whose left-hand side is the square of the second base below, the second bracket is
\[
\begin{aligned}
  &\Big(\frac{\Delta}{d^{1/4}}\Big)^{|E_{G\mathbf n^\dagger}|-\lceil p\rceil(|V_{G\mathbf n^\dagger}|-2)}
  \bigg(\frac{\Delta^{\lceil p\rceil}\sqrt n}{d^{\lceil p\rceil/4}K}\bigg)^{|V_{G\mathbf n^\dagger}|-2}\\
  &\qquad\le\ D^{-200(|E_{G\mathbf n^\dagger}|-\lceil p\rceil(|V_{G\mathbf n^\dagger}|-2))}\,
  D^{-200\lceil p\rceil(|V_{G\mathbf n^\dagger}|-2)}\ =\ D^{-200|E_{G\mathbf n^\dagger}|} .
\end{aligned}
\]
Multiplying, $\hat s_G\le\tfrac1KD^{-200|E|}$. Finally $|V|\le2|E|$ (Lemma~\ref{lem:sizes}), so
$50(|E|+|V|)\le150|E|$, and \eqref{eq:eta} gives $\eta_G\le\tfrac1KD^{-50|E|}\le\tfrac1KD^{-2|E|}$.
\end{proof}

\begin{proof}[Proof of Theorem~\ref{thm:LB}]
By Lemma~\ref{lem:master}, $\eta_G^2\le K^{-2}D^{-4|E_G|}$. By Lemma~\ref{lem:count} the members of
$\mathcal G_{\le D}$ with $|E_G|=e$ number at most $(3D^2)^e$, and by \eqref{eq:evencons} only
even $e\ge2$ occur. As $D\ge2$ gives $3D^{-2}\le\tfrac34$,
\[
  \sum_{G\in\mathcal G_{\le D}}\eta_G^2\ \le\ \frac1{K^2}\sum_{\substack{e\ge2\\ e\ \text{even}}}(3D^2)^{e}D^{-4e}
  \ =\ \frac1{K^2}\sum_{\substack{e\ge2\\ e\ \text{even}}}\big(3D^{-2}\big)^{e}
  \ \le\ \frac1{K^2}\cdot\frac{(3/4)^2}{1-(3/4)^2}\ =\ \frac{9}{7K^2} .
\]
By Proposition~\ref{prop:reduce},
\[
  \Corr_{\le D}\ \le\ \frac{10}3\Big(\frac9{7K^2}\Big)^{1/2}\ =\ \frac{10}{\sqrt7\,K}\ \le\ \frac4K .
\]
\end{proof}

\section{Proof of the upper bound}\label{part:ub}
This section proves Theorem~\ref{thm:UB}, by applying the test of Lemma~\ref{lem:criterion} to the
statistic $P_{G}$ at a template $G$ (\S\ref{sec:template}).

\subsection{Overview of the proof}\label{sec:overview-ub}

The proof of Theorem~\ref{thm:UB} applies Chebyshev's inequality (Lemma~\ref{lem:criterion}) to the
invariant monomial $P_G$ of the template $G$ of Definition~\ref{def:template}.

\noindent\uline{Step 1: the mean.} $\E[P_G\mid x=0]=0$ and $\E[P_G\mid x=1]=\hat s_GK$ exactly
(Lemma~\ref{lem:gap}). Section~\ref{sec:moments} sketches this step.

\noindent\uline{Step 2: the second moment.} $\E[P_G^2\mid x]$ is a sum over the matchings
$\mathbf M$, which say which individual and variable nodes of two copies of $G$ are placed on the
same individuals and variables, and over the pairings $\Pai$ of the entries at each variable
\eqref{eq:partition}. Each term is at most $\hat s_G^2K^2K^{\mathcal E}$, up to a factor close to
$1$, with an explicit exponent $\mathcal E$ (Proposition~\ref{prop:exponent}). The single pair with
no matched node other than $v_1,v_2$ and no matched variable gives the squared mean. Section~\ref{sec:secondsketch} sketches this step.

\noindent\uline{Step 3: the exponent.} Theorem~\ref{thm:master} shows that every other pair has
$\mathcal E\le-\frac2{L_B}(|\mathbf M_V|-2+|\mathbf M_W|)$. The exponent is written exactly as a sum
of four groups \eqref{eq:Ekinds}, and the groups are bounded through the traces on one copy of the
components of the symmetric difference graph (Definition~\ref{def:trace}). A component missing the
matched individuals contributes a factor $K$, the free choice of its group, which is compensated by
the links leaving its trace. If the trace holds at most half of the body, the expansion of the body
gives at least $q$ such links, and $\ciso$ times the size of the trace
(Lemma~\ref{lem:small}), and $q$ links suffice because of the critical order condition
$d^{q/2}\ge K^{1+3201q/L}$ \eqref{eq:E}. If the trace holds more than half of
the body, the fasteners give at least $2L_\Phi$ such links (Lemmas~\ref{lem:fast}
and~\ref{lem:large}), and $2L_\Phi$ links suffice because $\Dcomp\ge1$ and $\Delta^2$ exceeds
$\Dcomp^2$ by the factor $K^{1/L_\Phi+4/L_B}$ at least, by \ref{H4}, where $\Dcomp$ is computed with
$q$ in place of $\lceil p\rceil$. Sections~\ref{sec:shapesketch} and~\ref{sec:overlapsketch} sketch this step.

\noindent\uline{Step 4: summation.} The pairs $(\mathbf M,\Pai)$ with $|\mathbf M_V|-2+|\mathbf M_W|$
equal to a given $\bar m$ number at most $(3|V|^2|W|^2)^{\bar m}$, so the variances are at most a
small multiple of the squared mean difference $\hat s_G^2K^2$ (Proposition~\ref{prop:summation}). Section~\ref{sec:secondsketch} sketches this step.

\subsection{The statistic used for the upper bound}\label{sec:ubsetting}

\subsubsection{Preliminaries}\label{sec:ubmodel}

The statistic is the polynomial $P_{G}$ of Definition~\ref{def:PG}, at a connected $G$ as in
\S\ref{sec:graphs} all of whose variables are links. As every multiplicity of such a $G$ is $1$
and $\He_{1}(z)=z$, every Hermite factor of \eqref{eq:PG} is $Y$ itself and every factorial
there is $1$, so that
\begin{equation}\label{eq:PGII}
  P_{G}\ =\ \hat s_{G}\,\frac{K^{\,|V|-1}}{|\Pi_{G}|}\Big(\frac{\sqrt d}{\Delta}\Big)^{2|W|}
   \sum_{\pi\in\Pi_{G}}\ \prod_{w=\{v,v'\}\in W}
      Y_{\pi_{1}(v),\pi_{2}(w)}\,Y_{\pi_{1}(v'),\pi_{2}(w)} .
\end{equation}

\paragraph{$G$ as an ordinary multigraph.} Such a $G$ --- the template of \S\ref{sec:template} is
one --- is the same thing as an ordinary multigraph on its individual nodes: a link
$w=\{v,v'\}$ is one edge of that multigraph, joining $v$ to $v'$, and the two elements of $E$
that $w$ carries are its two ends, one at $v$ and one at $v'$. The multigraph therefore has
$|W|$ edges where $G$ has $|E|=2|W|$. Section~\ref{part:ub} keeps the graph of \S\ref{sec:graphs} and its
words, writing \textbf{node} for an individual node since the variable nodes are the links.

\paragraph{The symmetric difference and the intersection.} Read that way, a pairing $\Pai$ of a
merged graph $\Gcup(G^{(1)},G^{(2)},\mathbf M)$ turns $\Gcup[\Pai]$ into two ordinary graphs on
the individual nodes of $\Gcup$: the \emph{symmetric difference graph}
$\Gomi(G^{(1)},G^{(2)},\mathbf M,\Pai)$, whose edges are its links, and the \emph{intersection
graph} $\Gcap(G^{(1)},G^{(2)},\mathbf M,\Pai)$, whose edges are its self-loops, a self-loop
giving a loop at the node it sits at. The edge set of $\Gcap$ is written $E_{\cap}$ and the set
of variables of $\Gcup$ carrying one of its edges $W_{\cap}$. The $|E_{\Gcup}|/2$ variables of
$\Gcup[\Pai]$ split between the two graphs, and the components of $\Gomi$ are those of
$\Gcup[\Pai]$ read on the individual nodes, so that the components of \S\ref{sec:graphs} are
read on $\Gomi$ below. As for $\Gcup$ in \S\ref{sec:graphs}, the
arguments of $\Gomi$ and $\Gcap$ are dropped.

\paragraph{The pairings of $\Gcup(G,G,\mathbf M)$.}\label{sec:linkpairings} Let $\mathbf M\in\mathcal M(G,G)$ and
$\Gcup:=\Gcup(G,G,\mathbf M)$, its pairings $\Pai$ ranging over $\PcalSL(\Gcup)$, not over $\Pcal(\Gcup)$. A
variable of $\Gcup$ carries one link of each copy of $G$ when it is matched, and one link when it is not, so it
has degree $4$, resp.\ $2$. An unmatched variable has a single pairing, joining its
link's two ends, and a matched one three. So only a matched variable can carry a self-loop. A self-loop at it
makes its two links share an individual node, which is then matched, and two self-loops make them share both.
So $W_{\cap}\subseteq\mathbf M_{W}$, and $\Gcup[\Pai]$, which has $|E_{\Gcup}|/2=2|W|$ variables, has
$|E_{\cap}|$ self-loops. A matched variable carries two of them exactly when it carries no edge of $\Gomi$, and
none exactly when it carries two, so $|E_{\cap}|-|W_{\cap}|$ matched variables carry two self-loops,
$2|W_{\cap}|-|E_{\cap}|$ carry one and $|\mathbf M_{W}|-|W_{\cap}|$ carry none.

\begin{lemma}\label{lem:mono}
Let $C\subseteq[n]$ contain the individuals $1$ and $2$, be partitioned into nonempty sets
$C_{1},\dots,C_{N}$, and let $K\ge2$. Let $\mathcal A$ be the event that $k^{*}$ is constant on
each $C_{a}$. Then, exactly,
\[
\begin{aligned}
  \Pr[\mathcal A\mid x=1]\ &=\ K^{\,N-|C|+\one\{1,2\in C_a\ \text{for some}\ a\}},\\
  \Pr[\mathcal A\mid x=0]\ &=\ \big(1-\one\{1,2\in C_a\ \text{for some}\ a\}\big)\,K^{\,N-|C|} .
\end{aligned}
\]
\end{lemma}

\begin{proof}
By Lemma~\ref{lem:colourcount}, $\Pr[\mathcal A]=K^{N-|C|}$. If $1,2\in C_a$ for some $a$, then
$\mathcal A\subseteq\{x=1\}$, so $\Pr[\mathcal A\mid x=1]=\Pr[\mathcal A]/\Pr[x=1]=K^{N-|C|+1}$ and
$\Pr[\mathcal A\mid x=0]=0$. Otherwise $\mathcal A\cap\{x=1\}$ is the event that $k^{*}$ is constant on each set
of the partition of $C$ obtained by merging the two sets holding $1$ and $2$, which has $N-1$ sets, so
$\Pr[\mathcal A\cap\{x=1\}]=K^{N-1-|C|}$ by Lemma~\ref{lem:colourcount}. Dividing by $\Pr[x=1]=K^{-1}$, and
subtracting from $\Pr[\mathcal A]$ and dividing by $\Pr[x=0]=1-K^{-1}>0$, gives
$\Pr[\mathcal A\mid x=1]=\Pr[\mathcal A\mid x=0]=K^{N-|C|}$.
\end{proof}

\begin{lemma}\label{lem:criterion}
Let $P$ be a statistic whose two conditional means differ. The test that declares $x=1$ when
$P$ is closer to $\E[P\mid x=1]$ than to $\E[P\mid x=0]$, and $x=0$ otherwise, declares the wrong
value of $x$, under either conditioning, with probability at most
\[
  \frac{4\big(\Var(P\mid x=1)+\Var(P\mid x=0)\big)}
       {\big(\E[P\mid x=1]-\E[P\mid x=0]\big)^{2}} .
\]
\end{lemma}

\begin{proof}
Let $i\in\{0,1\}$. Under $x=i$, the test declares $x=1-i$ only if
$|P-\E[P\mid x=1-i]|\le|P-\E[P\mid x=i]|$, and then, by the triangle inequality,
\[
  \big|\E[P\mid x=1]-\E[P\mid x=0]\big|\ \le\ 2\,\big|P-\E[P\mid x=i]\big| .
\]
By Chebyshev's inequality, this has conditional probability at most
\[
  \frac{4\Var(P\mid x=i)}{\big(\E[P\mid x=1]-\E[P\mid x=0]\big)^{2}} ,
\]
which is at most the bound of the statement.
\end{proof}

\subsubsection{The properties of the template used for the upper bound}\label{sec:template}

The graph used for the upper bound, the
template, is defined through the properties it has to satisfy. These properties are satisfied,
with positive probability, by a random expander defined in
Appendix~\ref{app:body}.
\begin{definition}\label{def:template}
Fix integers $L_{B},L_{\Phi}\ge1$ and let the \emph{body} $B$ be a set of $|B|=2L_{\Phi}L_{B}$ nodes carrying a
graph $G_{B}$ --- with parallel edges allowed but no self-loop --- in which every node meets
exactly $q$ edges and which satisfies, for some constant $\ciso>0$ and with
$\partial_{B}(S)$ denoting the number of edges of $G_{B}$ with exactly one endpoint in $S$,
\begin{equation}\label{eq:expand}
  \partial_{B}(S)\ \ge\ q\quad\big(\varnothing\ne S\subsetneq B\big),
  \qquad
  \partial_{B}(S)\ \ge\ \ciso\min\{|S|,|B|-|S|\}\quad\big(S\subseteq B\big) .
\end{equation}
Let the \emph{fasteners} $\Phi$ be $2L_{\Phi}$ nodes of $B$ at pairwise distance at least $4/\ciso$ in
$G_{B}$. The \emph{template} is the graph $G$ with individual nodes $V:=B\cup\{v_{1},v_{2}\}$,
the two extra ones being the distinguished nodes of \S\ref{sec:graphs}, with one link for each edge
of $G_{B}$ and one \emph{fastener link} from each of $L_{\Phi}$ of the fasteners to $v_{1}$ and from
each of the other $L_{\Phi}$ to $v_{2}$.
\end{definition}

In Appendix~\ref{app:body}, $G_{B}$ is a cycle through all the nodes of $B$ together with
$q-2$ random perfect matchings of $B$. Proposition~\ref{prop:body} shows that, with positive
probability, this gives a template satisfying Definition~\ref{def:template}, at every degree $q\ge3$, with
$\ciso=1/400$, once $L_{B}>q^{1601}$ and $|B|\ge(400q)^{1200q^{2}}$. Figure~\ref{fig:template} draws the shape of such a template at $q=3$, on far fewer nodes than Proposition~\ref{prop:body} asks.

Counting $|B|q/2$ body links and $2L_{\Phi}$ fastener links,
\[
  \deg_{G} v=q\ \ (v\in B\setminus\Phi),\qquad \deg_{G} v=q+1\ \ (v\in\Phi),\qquad
  \deg_{G} v_{1}=\deg_{G} v_{2}=L_{\Phi} ,
\]
\begin{equation}\label{eq:W-count}
  |W|=\frac{|B|q}{2}+2L_{\Phi}=L_{\Phi}(q L_{B}+2) ,
\end{equation}
and every variable of $G$ is a link, as $G_B$ has no self-loop. $G$ is also connected, as \eqref{eq:PGII} asks: the first
inequality of \eqref{eq:expand} gives $\partial_{B}(S)\geq>0$ for every nonempty
proper $S\subsetneq B$, so $G_{B}$ is connected, and $v_{1}$ and $v_{2}$ each meet $L_{\Phi}\ge1$
fastener links.

\begin{lemma}\label{lem:fast}
For every $S\subseteq B$ with $|S|<|B|/2$,
\[
  \partial_{B}(S)\ \ge\ 2\,|S\cap\Phi| .
\]
\end{lemma}

\begin{proof}
We reduce to a set $S$ inducing a connected subgraph of $G_{B}$. If $S$ holds one fastener, the
first inequality of \eqref{eq:expand} suffices. If it holds $k\ge2$ fasteners, pairwise at distance
at least $4/\ciso$, then $S$ has more than $2k/\ciso$ nodes, and the second inequality of
\eqref{eq:expand} concludes.

\noindent\uline{Step 1: reduction.} Both sides are additive over the node sets of the components of the subgraph of
$G_{B}$ induced on $S$, as no edge of $G_{B}$ joins two of them, and each of these sets has fewer than $|B|/2$
elements. So we may assume that this subgraph is connected. Put $k:=|S\cap\Phi|$. If $k=0$ there is nothing to
prove. If $k=1$, then $\varnothing\ne S\subsetneq B$, and the first inequality of \eqref{eq:expand} gives
$\partial_{B}(S)\geq\ge2=2k$.

\noindent\uline{Step 2: the case $k\ge2$.} The closed walk running along a spanning tree of the subgraph
induced on $S$, through each tree edge twice, has length $2(|S|-1)$ and visits every node of
$S$. Order the $k$ fasteners of $S$ cyclically by their first visit along it. The first visits
cut the closed walk into $k$ portions, each joining two consecutive fasteners of this cyclic
order. Two consecutive ones are distinct, as $k\ge2$, hence at distance at least $4/\ciso$ in
$G_{B}$ (Definition~\ref{def:template}), so each portion has length at least $4/\ciso$. Adding
over the $k$ portions,
\[
  2(|S|-1)\ \ge\ \frac{4k}{\ciso},\qquad\text{that is}\qquad |S|\ \ge\ \frac{2k}{\ciso}+1 .
\]
As $|S|<|B|/2$, $\min\{|S|,|B|-|S|\}=|S|$, and the second inequality of \eqref{eq:expand} gives
\[
  \partial_{B}(S)\ \ge\ \ciso|S|\ \ge\ \ciso\Big(\frac{2k}{\ciso}+1\Big)\ =\ 2k+\ciso\ >\ 2k .
\]
\end{proof}

\subsection{The standing hypotheses}\label{sec:exponents}

Fix an integer $q\ge3$ and write
\[
  n=K^{\nu}\ \text{points in }\mathbb R^{d},\qquad d=K^{2/q+\upp},
  \qquad \Delta^{2}=K^{\delta_{*}+\upd+\upp/2} ,
\]
where
\begin{equation}\label{eq:deltabar}
  \delta_{*}\ :=\ \frac{3-\nu}{q} .
\end{equation}
The combination that Assumption~\ref{ass:B}\textup{(ii)} bounds is
\begin{equation}\label{eq:invariant}
  \frac{\Delta^{2q}\,n}{K^{2}\,d^{\,q/2}}\ =\ K^{q\,\upd} .
\end{equation}

Besides $\delta_{*}$ and the parameters $q,\nu,\upd,\upp,L_{B},L_{\Phi}$, the argument
uses two combinations of exponents:
\begin{equation}\label{eq:defs}
  \zeta:=\frac{\nu-2}{q}-\upd+\frac{2}{L_{B}} ,\qquad
  \theta:=\delta_{*}+\upd+\frac{\upp}2-\frac3{L_{B}} .
\end{equation}
By \eqref{eq:deltabar}, $K^{\zeta-2/L_{B}}=\sqrt d/\Delta^{2}$ and $K^{\theta+3/L_{B}}=\Delta^{2}$,
while $K^{1/q+\upp/2}=\sqrt d$. Hence
\begin{equation}\label{eq:identities}
  \zeta+\theta\ =\ \frac1{q}+\frac{\upp}2-\frac1{L_{B}} .
\end{equation}
Until the proof of Theorem~\ref{thm:UB}, $L_{B}$ and $L_{\Phi}$ are subject only to the
hypotheses below.

\medskip\noindent\textbf{The hypotheses}, each statement naming the ones it uses.
\begin{enumerate}[label=\textup{(H\arabic*)},leftmargin=3em,itemsep=1pt]
\item\label{H1} $q\ge3$, $K\ge d$, $\nu>2$ and $n\le K^2d^{q/2}$, and $G$ is the template of \S\ref{sec:template}, with expansion $\ciso$.
\item\label{H2} $L_{B}\ge8/\ciso$.
\item\label{H3} $\upd<(\nu-2)/q$, equivalently $\delta_{*}+\upd<1/q$, equivalently
  $\Delta^{2}<\sqrt d$.
\item\label{H4} $\upd\ge4/L_{B}+1/L_{\Phi}$ and $\upp\ge2/L_{B}+16/(\ciso L_{B})$.
\item\label{H5} $d\ge K^{2/q}$ and $K^{1/q}\ge400\,|W|^{2}$.
\end{enumerate}
By \ref{H3} and \ref{H5}, $\Delta^{2}/d\le d^{-1/2}\le K^{-1/q}$, and $|V|\le|W|+1\le2|W|$ at a connected
graph all of whose variables are links, so \ref{H5} makes
\begin{equation}\label{eq:smallconseq}
  \frac{|V|^{2}}K,\qquad \frac{|W|^{2}}d,\qquad \frac{|W|\Delta^{2}}d
  \qquad\text{each at most}\qquad 4\,|W|^{2}K^{-1/q}\ \le\ \tfrac1{100} .
\end{equation}

Four consequences are used below. First, $\zeta>2/L_{B}\,>0$ by \ref{H3}. Second, as $K^2d^{q/2}=K^{3+q\upp/2}$, the bound $n\le K^2d^{q/2}$ of \ref{H1} reads $\nu\le3+q\upp/2$, that
is, $\delta_{*}+\upp/2\ge0$. With \ref{H4},
\begin{equation}\label{eq:H4bis}
  \theta\ \ge\ \frac{1}{L_{\Phi}}+\frac{1}{L_{B}} \ >\ 0 .
\end{equation}
Third, $q(\upd-2/L_{B})$ is at
least $2q/L_{B}\ge6/L_{B}$ by \ref{H4} and $q\ge3$, while
$\zeta<(\nu-2)/q\le\log d/(2\log K)\le1/2$ by \ref{H4} and \ref{H1}. So
\begin{equation}\label{eq:C}
  q\Big(\upd-\frac2{L_{B}}\Big)-\frac{2}{L_{B}}\ \ge\ \frac{4\zeta}{L_{B}} .
\end{equation}
Fourth, \ref{H4} gives $\upp/2-1/L_{B}\ge8/(\ciso L_{B})$, so that, by \eqref{eq:identities},
\begin{equation}\label{eq:E}
  q(\zeta+\theta)\ =\ 1+q\Big(\frac{\upp}2-\frac1{L_{B}}\Big)\ \ge\ 1+\frac{8}{\ciso L_{B}} .
\end{equation}

\subsection{The mean}\label{sec:mean}

\begin{lemma}\label{lem:gap}
Let $K\ge2$ and let $G$ be connected with every variable a link. Then $\E[P_{G}\mid x=0]=0$ and
$\E[P_{G}\mid x=1]=\hat s_{G}K$, both exactly.
\end{lemma}
\begin{proof}
Fix a labeling $(\pi_{1},\pi_{2})\in\Pi_{G}$. The variables $\pi_{2}(w)$, $w\in W$, are distinct, so the factors
at distinct links are independent given $k^{*}$ (\S\ref{sec:model}), and $\pi_{1}(v)\ne\pi_{1}(v')$ at every link
$w=\{v,v'\}$. Hence, by \eqref{eq:cov},
\begin{equation}\label{eq:gapone}
  \E\Big[\prod_{w=\{v,v'\}\in W}Y_{\pi_{1}(v),\pi_{2}(w)}Y_{\pi_{1}(v'),\pi_{2}(w)}\ \Big|\ k^{*}\Big]
  \ =\ \Big(\frac{\Delta^{2}}{d}\Big)^{|W|}\prod_{\{v,v'\}\in W}\one\{\pi_{1}(v)\sim\pi_{1}(v')\} ,
\end{equation}
and, $G$ being connected, the product of indicators is the indicator that $k^{*}$ is constant on
$\pi_{1}(V)\ni1,2$. Under $x=0$ it vanishes, so $\E[P_{G}\mid x=0]=0$ by \eqref{eq:PGII}. Under $x=1$ it has
probability $K^{-(|V|-2)}$ (Lemma~\ref{lem:mono} with $N=1$), so, by \eqref{eq:PGII},
\[
  \E[P_{G}\mid x=1]\ =\ \hat s_{G}\,\frac{K^{\,|V|-1}}{|\Pi_{G}|}\Big(\frac{\sqrt d}{\Delta}\Big)^{2|W|}
   \cdot|\Pi_{G}|\Big(\frac{\Delta^{2}}{d}\Big)^{|W|}K^{-(|V|-2)}\ =\ \hat s_{G}K .
\]
\end{proof}

\subsection{The second moment}\label{sec:second}

From now on $G$ is the template of \S\ref{sec:template}, so $|V|=|B|+2$ and $|W|$ is given by
\eqref{eq:W-count}. Squaring \eqref{eq:PGII} gives, with the square of its scalar in front,
\[
\begin{aligned}
  \E[P_{G}^{2}\mid x=1]
  =\bigg(\hat s_{G}\,\frac{K^{\,|V|-1}}{|\Pi_{G}|}\Big(\frac{\sqrt d}{\Delta}\Big)^{2|W|}\bigg)^{2}
  &\sum_{(\pi^{1}_{1},\pi^{1}_{2})}\ \sum_{(\pi^{2}_{1},\pi^{2}_{2})}\\[-2pt]
  &\times\E\Big[\prod_{r=1,2}\prod_{w=\{v,v'\}\in W}
    Y_{\pi^{r}_{1}(v),\pi^{r}_{2}(w)}Y_{\pi^{r}_{1}(v'),\pi^{r}_{2}(w)}\ \Big|\ x=1\Big] ,
\end{aligned}
\]
both outer sums running over the labelings of $G$. The two copies of $G$ they label are the
\emph{replicas}, $\pi^{r}$ labeling replica $r$. Which nodes of the two replicas land on one
individual, and which links on one variable, is a matching $\mathbf M$ of $G$ with $G$ in the
sense of \S\ref{sec:graphs}, the four maps being injective: by \eqref{eq:PiM} it is the unique
one with $\pi^{1}\in\Pi_{G}(\mathbf M,\pi^{2})$. It fixes $v_{1}$ and $v_{2}$, so it lies in
$\mathcal M(G,G)$ of \S\ref{sec:graphs}, and the pairs of labelings giving it are the $(\pi^{1},\pi^{2})\in\Pi_{G}\times\Pi_{G}$ with
$\pi^{1}\in\Pi_{G}(\mathbf M,\pi^{2})$.

Write $V^{(r)}(\mathbf M)\subseteq V$ for the nodes $\mathbf M$ matches in replica $r$,
and $\Gcup:=\Gcup(G,G,\mathbf M)$ for the merged graph of \S\ref{sec:linkpairings}
on the two replicas, with $2|V|-|\mathbf M_{V}|$ individual and $2|W|-|\mathbf M_{W}|$ variable
nodes. For a pair $(\pi^{1},\pi^{2})$ giving $\mathbf M$, call an individual \textbf{matched} when it lies in $\pi^{1}_{1}(V^{(1)}(\mathbf M))$ and
\textbf{unmatched} otherwise: there are $|\mathbf M_{V}|$ matched individuals, $1$ and $2$ among
them, so $|\mathbf M_{V}|-2$ counts the matched nodes of $B$ in each replica, and the unmatched
individuals of replica $r$ are the $\pi^{r}_{1}(v)$ with $v$ in the set
$B\setminus V^{(r)}(\mathbf M)$, of size $|B|-|\mathbf M_{V}|+2$.

\subsubsection{Decomposition as a sum over matchings and pairings $(\mathbf M,\Pai)$}\label{sec:wick}

Fix $(\pi^{1},\pi^{2})\in\Pi_{G}\times\Pi_{G}$ with $\pi^{1}\in\Pi_{G}(\mathbf M,\pi^{2})$, and send each node of
$\Gcup$ carrying a node $v$, resp.\ a variable $w$, of replica $r$ to $\pi^{r}_{1}(v)$, resp.\ $\pi^{r}_{2}(w)$. By
\eqref{eq:PiM} the coincidences of $\pi^{1}$ and $\pi^{2}$ are exactly the pairs of $\mathbf M$, so this map is well
defined and injective: the $2|V|-|\mathbf M_{V}|$ individual nodes of $\Gcup$ go to as many individuals, $1$ and $2$
among them, and its $2|W|-|\mathbf M_{W}|$ variable nodes to as many variables. Condition on $k^{*}$. The expectation factorizes over the variables of $\Gcup$, and at each of
them the entries are centered and jointly Gaussian, so by Isserlis's formula \cite[Theorem~1.28]{janson1997gaussian} the expectation of
the product of its entries is the sum, over the pairings of that variable, of the product of
$\E[Y\,Y]$ over its pairs. By
\eqref{eq:cov} that factor is $(\Delta^{2}/d)\one\{i\sim i'\}$ at a link of $\Gcup[\Pai]$ and
$1+\Delta^{2}/d$ at a self-loop, and $\Gcup[\Pai]$ has $2|W|$ variables of which
$|E_{\cap}|$ are self-loops
(\S\ref{sec:linkpairings}). Multiplying over the variables of $\Gcup$ and summing over the pairings,
and writing $\chi$ for the coloring of the individuals of $\Gcup$ by $k^{*}$,
\begin{equation}\label{eq:wick}
\begin{aligned}
 \E\Big[\prod_{r=1,2}\prod_{w=\{v,v'\}\in W}&Y_{\pi^{r}_{1}(v),\pi^{r}_{2}(w)}Y_{\pi^{r}_{1}(v'),\pi^{r}_{2}(w)}
   \ \Big|\ k^{*}\Big]\\
 &=\sum_{\Pai\in\PcalSL(\Gcup)}\Big(1+\frac{\Delta^{2}}{d}\Big)^{|E_{\cap}|}
   \Big(\frac{\Delta^{2}}{d}\Big)^{2|W|-|E_{\cap}|}f_{\Gomi}(\chi) ,
\end{aligned}
\end{equation}
where $f_{\Gomi}$ is the function $f$ of Proposition~\ref{prop:master} at
$\Gomi$, so that $f_{\Gomi}(\chi)$ is, by Corollary~\ref{cor:selfloop}, the indicator that $\chi$ is constant on each component of the
symmetric difference graph $\Gomi$ of \S\ref{sec:ubmodel}.
$\CC_{\Gomi}(\mathbf M)\subseteq\CC_{\Gomi}$
is the set of its \textbf{free} components, those holding no matched individual, the others being \textbf{pinned} (\S\ref{sec:graphs}, with $\bar V$ the set of matched individuals), and $1\sim_{\Gomi}2$
says
that the individuals $1$ and $2$ lie in one component --- a property of $(\mathbf M,\Pai)$, not a
random event.

The two exponents of \eqref{eq:wick} depend on
$(\mathbf M,\Pai)$ alone, and with the factor
$(\sqrt d/\Delta)^{4|W|}$ of the square of the scalar of \eqref{eq:PGII} they give:
\[
  \Big(\frac{\sqrt d}{\Delta}\Big)^{4|W|}\Big(1+\frac{\Delta^{2}}{d}\Big)^{|E_{\cap}|}
  \Big(\frac{\Delta^{2}}{d}\Big)^{2|W|-|E_{\cap}|}
  \ =\ \Big(1+\frac{d}{\Delta^{2}}\Big)^{|E_{\cap}|} .
\]
So put
\[
\begin{aligned}
  T(\mathbf M,\Pai)\ :=\ &\bigg(\hat s_{G}\,\frac{K^{\,|V|-1}}{|\Pi_{G}|}\bigg)^{2}
    \Big(1+\frac{d}{\Delta^{2}}\Big)^{|E_{\cap}|}\\
  &\times\sum_{\pi^{2}\in\Pi_{G}}\ \ \sum_{\pi^{1}\in\Pi_{G}(\mathbf M,\pi^{2})}
     \E\big[f_{\Gomi}(\chi)\mid x=1\big] .
\end{aligned}
\]
Taking in \eqref{eq:wick} the expectation given $x=1$ and summing over the pairs of labelings,
grouped by the matching they give,
\begin{equation}\label{eq:partition}
  \E[P_{G}^{2}\mid x=1]\ =\ \sum_{\mathbf M\in\mathcal M(G,G)}\ \ \sum_{\Pai\in\PcalSL(\Gcup)}T(\mathbf M,\Pai) ,
\end{equation}
an $\mathbf M\in\mathcal M(G,G)$ given by no pair of labelings contributing $0$. Everything above is pointwise in
$k^{*}$, so \eqref{eq:partition} holds with $x=0$ in place of $x=1$ throughout.

\subsubsection{Upper bound on $T(\mathbf M,\Pai)$}\label{sec:exppart}

\begin{proposition}\label{prop:exponent}
Assume \ref{H1}, \ref{H3} and \ref{H5}. For $\mathbf M\in\mathcal M(G,G)$ and $\Pai\in\PcalSL(\Gcup)$, put
\begin{equation}\label{eq:exponent}
\begin{aligned}
  \mathcal E\ :=\ &-(\nu-1)(|\mathbf M_{V}|-2)+|\CC_{\Gomi}|+\one\{1\sim_{\Gomi}2\}-2\\
   &+\frac{\nu-1}{q}\,|E_{\cap}|-\frac2{q}\,|\mathbf M_{W}|\\
   &+\Big(\frac{\upp}2-\upd\Big)|E_{\cap}|-\upp\,|\mathbf M_{W}| .
\end{aligned}
\end{equation}
Then
\[
  T(\mathbf M,\Pai)\ \le\ \hat s_{G}^{2}K^{2}\,K^{\mathcal E}\,\big(1+32|W|^{2}K^{-1/q}\big) ,
\]
and the contribution of $(\mathbf M,\Pai)$ to $\E[P_{G}^{2}\mid x=0]$ in \eqref{eq:partition} is $0$ if $1\sim_{\Gomi}2$ and at
most the same quantity otherwise.
\end{proposition}

\begin{proof}
Every term of the sum defining $T(\mathbf M,\Pai)$ is the same conditional expectation of
$f_{\Gomi}(\chi)$, which Lemma~\ref{lem:mono} gives exactly, so it remains to bound the number of
terms. The other factors are powers of $K,n,d,\Delta$ up to two correction factors with product at
most $1+32|W|^{2}K^{-1/q}$, and collecting the powers gives $K^{\mathcal E}$. Under $x=0$ the
conditional expectation is $0$ or the one under $x=1$.

\noindent\uline{Step 1: the conditional expectation of $f_{\Gomi}(\chi)$.} Fix $(\pi^{1},\pi^{2})\in\Pi_{G}\times\Pi_{G}$ with
$\pi^{1}\in\Pi_{G}(\mathbf M,\pi^{2})$. The map of \S\ref{sec:wick} sends the $2|V|-|\mathbf M_{V}|$ individual nodes of $\Gcup$
 onto a set of $2|V|-|\mathbf M_{V}|$ individuals containing $1$ and $2$.
The components of $\Gomi$ partition that set into $|\CC_{\Gomi}|$ nonempty sets, $1$ and $2$
lying in one of them exactly when $1\sim_{\Gomi}2$, and $f_{\Gomi}(\chi)$ is the indicator that
$k^{*}$ is constant on each of these sets (\S\ref{sec:wick}). Lemma~\ref{lem:mono}, with
$N=|\CC_{\Gomi}|$ and with $K\ge2$ by \ref{H5}, gives
\begin{equation}\label{eq:expprob}
\begin{aligned}
  \E\big[f_{\Gomi}(\chi)\mid x=1\big]
  \ &=\ K^{\,|\CC_{\Gomi}|+\one\{1\sim_{\Gomi}2\}-(2|V|-|\mathbf M_{V}|)} ,\\
  \E\big[f_{\Gomi}(\chi)\mid x=0\big]
  \ &=\ \big(1-\one\{1\sim_{\Gomi}2\}\big)\,\E\big[f_{\Gomi}(\chi)\mid x=1\big] ,
\end{aligned}
\end{equation}
and neither right-hand side depends on $(\pi^{1},\pi^{2})$.

\noindent\uline{Step 2: the number of sum elements.} By \eqref{eq:PiM}, $\pi^{1}\in\Pi_{G}(\mathbf M,\pi^{2})$ is a condition on
$(\pi^{1}_{1},\pi^{2}_{1})$ and one on $(\pi^{1}_{2},\pi^{2}_{2})$. A pair satisfying the first is determined by the
injection of \S\ref{sec:wick} of the individual nodes of $\Gcup$ into $[n]$, which sends $v_{1},v_{2}$ to $1,2$, and a pair
satisfying the second by an injection of the variable nodes of $\Gcup$ into $[d]$. Hence
\begin{equation}\label{eq:expcount}
  \big|\big\{(\pi^{1},\pi^{2})\in\Pi_{G}\times\Pi_{G}:\pi^{1}\in\Pi_{G}(\mathbf M,\pi^{2})\big\}\big|
  \ \le\ n^{\,2|V|-|\mathbf M_{V}|-2}\,d^{\,2|W|-|\mathbf M_{W}|} .
\end{equation}

\noindent\uline{Step 3: the two correction factors.} Write
$(1+d/\Delta^{2})^{|E_{\cap}|}=(d/\Delta^{2})^{|E_{\cap}|}(1+\Delta^{2}/d)^{|E_{\cap}|}$. As
$|E_{\cap}|\le2|W|$ (\S\ref{sec:linkpairings}), $\Delta^{2}/d\le K^{-1/q}$ (\S\ref{sec:exponents}) and
$2|W|\Delta^{2}/d\le1$ \eqref{eq:smallconseq}, and $e^{y}\le1+2y$ for $0\le y\le1$,
\begin{equation}\label{eq:exppower}
  \Big(1+\frac{\Delta^{2}}{d}\Big)^{|E_{\cap}|}\ \le\ e^{2|W|\Delta^{2}/d}
  \ \le\ 1+6|W|^{2}K^{-1/q} .
\end{equation}
By \eqref{eq:Pisize}, the Weierstrass product inequality (as in Lemma~\ref{lem:embcount}), \eqref{eq:smallconseq} and $n\ge K$,
$|\Pi_{G}|\ge n^{|V|-2}d^{|W|}\big(1-|V|^{2}/n-|W|^{2}/d\big)\ge n^{|V|-2}d^{|W|}\big(1-8|W|^{2}K^{-1/q}\big)$,
so, as $8|W|^{2}K^{-1/q}\le\tfrac16$ by \ref{H5} and $(1-y)^{-2}\le1+3y$ for $0\le y\le\tfrac16$,
\begin{equation}\label{eq:expscalar}
  |\Pi_{G}|^{-2}\ \le\ n^{-(2|V|-4)}d^{-2|W|}\big(1+24|W|^{2}K^{-1/q}\big) .
\end{equation}

\noindent\uline{Step 4: collecting.} By the definition of $T(\mathbf M,\Pai)$, by \eqref{eq:expprob}, whose
right-hand side does not depend on the labelings, and by \eqref{eq:expcount},
\[
\begin{aligned}
  T(\mathbf M,\Pai)\ \le\ &\hat s_{G}^{2}K^{2|V|-2}|\Pi_{G}|^{-2}
   \Big(1+\frac{d}{\Delta^{2}}\Big)^{|E_{\cap}|}\\
  &\times n^{\,2|V|-|\mathbf M_{V}|-2}\,d^{\,2|W|-|\mathbf M_{W}|}\,
   K^{\,|\CC_{\Gomi}|+\one\{1\sim_{\Gomi}2\}-(2|V|-|\mathbf M_{V}|)} .
\end{aligned}
\]
Inserting \eqref{eq:exppower} and \eqref{eq:expscalar} and collecting the powers of each of
$K,n,d,\Delta$,
\[
\begin{aligned}
  T(\mathbf M,\Pai)\ \le\ &\hat s_{G}^{2}K^{2}\,\big(1+6|W|^{2}K^{-1/q}\big)\big(1+24|W|^{2}K^{-1/q}\big)\\
  &\times K^{\,|\mathbf M_{V}|+|\CC_{\Gomi}|+\one\{1\sim_{\Gomi}2\}-4}\ n^{-(|\mathbf M_{V}|-2)}\
  d^{\,|E_{\cap}|-|\mathbf M_{W}|}\ \Delta^{-2|E_{\cap}|} .
\end{aligned}
\]
At $n=K^{\nu}$, $\Delta^{2}=K^{\delta_{*}+\upd+\upp/2}$ and $d=K^{2/q+\upp}$, the
second line is a power of $K$ with exponent
\[
\begin{aligned}
  &(|\mathbf M_{V}|-2)-2+|\CC_{\Gomi}|+\one\{1\sim_{\Gomi}2\}-\nu(|\mathbf M_{V}|-2)\\
  &\qquad+\Big(\frac2{q}+\upp\Big)\big(|E_{\cap}|-|\mathbf M_{W}|\big)
   -\Big(\delta_{*}+\upd+\frac{\upp}2\Big)|E_{\cap}|\\
  &\quad=\ -(\nu-1)(|\mathbf M_{V}|-2)+|\CC_{\Gomi}|+\one\{1\sim_{\Gomi}2\}-2\\
  &\qquad+\Big(\frac{\nu-1}{q}+\frac{\upp}2-\upd\Big)|E_{\cap}|
   -\Big(\frac2{q}+\upp\Big)|\mathbf M_{W}|
  \ =\ \mathcal E ,
\end{aligned}
\]
by $2/q-\delta_{*}=(\nu-1)/q$ \eqref{eq:deltabar} and by
\eqref{eq:exponent}. Finally
$(1+6|W|^{2}K^{-1/q})(1+24|W|^{2}K^{-1/q})\le1+32|W|^{2}K^{-1/q}$, as $144|W|^{2}K^{-1/q}\le2$ by \ref{H5}.

\noindent\uline{Step 5: under $x=0$.} By \eqref{eq:partition} with $x=0$ in place of $x=1$ (\S\ref{sec:wick}),
the contribution of $(\mathbf M,\Pai)$ to $\E[P_{G}^{2}\mid x=0]$ is the expression defining
$T(\mathbf M,\Pai)$ with $\E[f_{\Gomi}(\chi)\mid x=0]$ in place of
$\E[f_{\Gomi}(\chi)\mid x=1]$. By \eqref{eq:expprob} it is $0$ if $1\sim_{\Gomi}2$ and equals
$T(\mathbf M,\Pai)$ otherwise, so in both cases it is at most the bound on $T(\mathbf M,\Pai)$.
\end{proof}

\subsection{Bounding the exponent $\mathcal E$}\label{sec:master}

By \S\ref{sec:linkpairings}, $|E_{\cap}|-|W_{\cap}|$, $2|W_{\cap}|-|E_{\cap}|$ and
$|\mathbf M_{W}|-|W_{\cap}|$ matched variables carry two, one and no self-loop, so that
$|E_{\cap}|=2(|E_{\cap}|-|W_{\cap}|)+(2|W_{\cap}|-|E_{\cap}|)$ and $|\mathbf M_{W}|$ is the sum
of the three numbers. Written with them, and with
$|\CC_{\Gomi}|=|\CC_{\Gomi}(\mathbf M)|+|\CC_{\Gomi}\setminus\CC_{\Gomi}(\mathbf M)|$,
\eqref{eq:exponent} gives, by \eqref{eq:deltabar}, \eqref{eq:defs} and \eqref{eq:identities}, exactly
\begin{equation}\label{eq:Ekinds}
\begin{aligned}
  &\mathcal E+\frac{2}{L_{B}}\big(|\mathbf M_{V}|-2+|E_{\cap}|+|\mathbf M_{W}|-|W_{\cap}|\big)\\
  &\quad=\ \zeta\big(2(|E_{\cap}|-|W_{\cap}|)-q(|\mathbf M_{V}|-2)\big)
   &&\textup{(I)}\\
  &\qquad+|\CC_{\Gomi}(\mathbf M)|-\Big(\theta+\frac1{L_{B}}\Big)(2|W_{\cap}|-|E_{\cap}|)
   -2(\zeta+\theta)(|\mathbf M_{W}|-|W_{\cap}|)
   &&\textup{(II)}\\
  &\qquad-\big(|\mathbf M_{V}|-\one\{1\sim_{\Gomi}2\}-|\CC_{\Gomi}\setminus\CC_{\Gomi}(\mathbf M)|\big)
   &&\textup{(III)}\\
  &\qquad-\Big(q\Big(\upd-\frac2{L_{B}}\Big)-\frac2{L_{B}}\Big)(|\mathbf M_{V}|-2) .
   &&\textup{(IV)}
\end{aligned}
\end{equation}

Lemma~\ref{lem:parts} bounds the four groups, \textup{(III)} and \textup{(IV)} directly, and
\textup{(I)} and \textup{(II)} through the links and the traces of one replica (\S\ref{sec:linksets}).

\subsubsection{The sets of links, and the traces}\label{sec:linksets}

Fix a pair $(\pi^{1},\pi^{2})$ giving $\mathbf M$, which exists under \ref{H1} and \ref{H5}, as then $n\ge K\ge2|V|$ and
$d\ge K^{2/q}\ge2|W|$. The map of \S\ref{sec:wick} identifies the individuals it places with the individual nodes of
$\Gcup$, and nothing below depends on this choice.

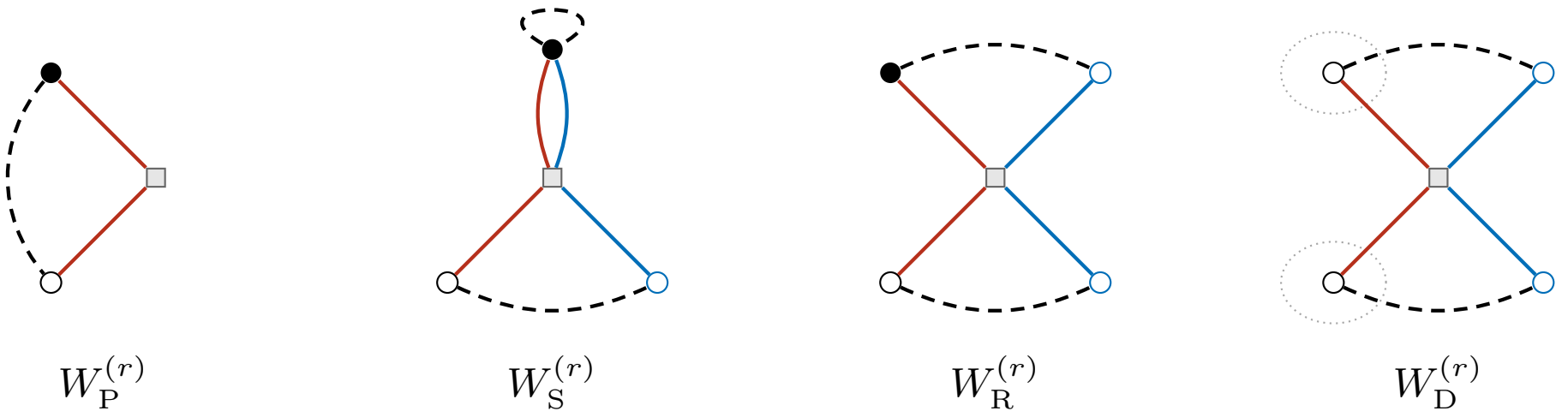
\begin{figure}[ht]
\centering
\begin{tikzpicture}[scale=0.9,
  mt/.style={circle,fill=black,inner sep=1.6pt},
  um/.style={circle,draw=black,fill=white,inner sep=1.6pt},
  va/.style={rectangle,draw=black!60,fill=black!10,inner sep=2pt},
  rr/.style={BrickRed,thick},
  ro/.style={NavyBlue,thick},
  pr/.style={black,thick,densely dashed}]
\begin{scope}
  \node[mt] (m) at (0,0.9) {};
  \node[um] (u) at (0,-0.9) {};
  \node[va] (w) at (0.9,0) {};
  \draw[rr] (m)--(w)--(u);
  \draw[pr] (m) to[bend right=40] (u);
  \node[font=\footnotesize,align=center] at (0.45,-1.75) {$\WP$};
\end{scope}
\begin{scope}[xshift=3.4cm]
  \node[mt] (m) at (0.9,1.1) {};
  \node[um] (u) at (0,-0.9) {};
  \node[va] (w) at (0.9,0) {};
  \node[um,draw=NavyBlue] (y) at (1.8,-0.9) {};
  \draw[rr] (m) to[bend right=20] (w); \draw[rr] (w)--(u);
  \draw[ro] (m) to[bend left=20] (w); \draw[ro] (w)--(y);
  \draw[pr] (m) to[out=150,in=30,looseness=12] (m);
  \draw[pr] (u) to[bend right=25] (y);
  \node[font=\footnotesize,align=center] at (0.9,-1.75) {$\WS$};
\end{scope}
\begin{scope}[xshift=7.2cm]
  \node[mt] (m) at (0,0.9) {};
  \node[um] (u) at (0,-0.9) {};
  \node[va] (w) at (0.9,0) {};
  \node[um,draw=NavyBlue] (x) at (1.8,0.9) {};
  \node[um,draw=NavyBlue] (y) at (1.8,-0.9) {};
  \draw[rr] (m)--(w)--(u); \draw[ro] (x)--(w)--(y);
  \draw[pr] (m) to[bend left=25] (x); \draw[pr] (u) to[bend right=25] (y);
  \node[font=\footnotesize] at (0.9,-1.75) {$\WR$};
\end{scope}
\begin{scope}[xshift=11cm]
  \draw[black!35,densely dotted] (0,0.9) ellipse (0.45 and 0.35);
  \draw[black!35,densely dotted] (0,-0.9) ellipse (0.45 and 0.35);
  \node[um] (u1) at (0,0.9) {};
  \node[um] (u2) at (0,-0.9) {};
  \node[va] (w) at (0.9,0) {};
  \node[um,draw=NavyBlue] (x) at (1.8,0.9) {};
  \node[um,draw=NavyBlue] (y) at (1.8,-0.9) {};
  \draw[rr] (u1)--(w)--(u2); \draw[ro] (x)--(w)--(y);
  \draw[pr] (u1) to[bend left=25] (x); \draw[pr] (u2) to[bend right=25] (y);
  \node[font=\footnotesize] at (0.9,-1.75) {$\WD$};
\end{scope}
\end{tikzpicture}
\caption{The link sets of replica $r$. Red: the link of replica $r$, through its variable
(box). Blue: the link of the other replica at the same variable, when the variable is matched.
Filled nodes are matched, open ones unmatched, and a node with a blue rim may be of
either kind. Dashed: the pairs of $\Pai$ at the variable, a
pair joining two individuals being an edge of $\Gomi$, one sitting at a single individual a
self-loop. In $\WP$, $\Pai$ pairs the two ends of the link together. In $\WS$, it
pairs them with the ends of the other link, the pair at the matched endpoint being a self-loop.
In $\WR$ it does so with no self-loop. In $\WD$ both endpoints are
unmatched and lie in different components of $\Gomi$ (dotted), and $\Pai$ then pairs the
two ends with those of the other link, with no self-loop
(Lemma~\ref{lem:links}\textup{(iii)}).}
\label{fig:linksets}
\end{figure}

Fix $r\in\{1,2\}$ and a link of $G$ with exactly one endpoint in $V^{(r)}(\mathbf M)$, the other then
lying in $B\setminus V^{(r)}(\mathbf M)$. The
link lies in $\WP$ when $\Pai$ pairs the two ends of its replica-$r$ copy with each other.
Otherwise its variable is matched, an unmatched one having the single pairing that
joins them (\S\ref{sec:linkpairings}), and $\Pai$ pairs each of these ends with an end of the other
replica's link at the same variable. The pair at the unmatched endpoint is not a
self-loop, the individual of that endpoint belonging to replica $r$ alone. The link lies in $\WS$
when the pair at its matched endpoint is a self-loop, and in $\WR$ when it is not. So $\WP$,
$\WS$, $\WR\subseteq W$ partition the links with exactly one endpoint in
$V^{(r)}(\mathbf M)$. Finally, $\WD\subseteq W$ is the
set of links whose two individuals in replica $r$ are both unmatched and lie in different
components of $\Gomi$. A link of each kind is drawn in Figure~\ref{fig:linksets}.

For $S\subseteq B$ let $\partial(S)$ be the number of links of $G$ with exactly one endpoint in
$S$. As $v_{1},v_{2}\notin S$, these are the fastener links at the fasteners in $S$
and the links of the edges of $G_{B}$ with exactly one endpoint in $S$, so that
\begin{equation}\label{eq:boundary}
  \partial(S)\ =\ \partial_{B}(S)+|S\cap\Phi| ,
\end{equation}
with $\partial_{B}(S)$ as in Definition~\ref{def:template}.

\begin{definition}\label{def:trace}
For $C\in\CC_{\Gomi}$ and a replica $r$ set
$\Vomi{r}(C):=\{v\in B\setminus V^{(r)}(\mathbf M):\pi^{r}_{1}(v)\in V_{C}\}$, the
\emph{trace} of $C$ on replica $r$. It is \emph{large} when it has more than $|B|/2$
elements.
\end{definition}

\subsubsection{Combinatorics for compensations between the different kinds of links and nodes}

\begin{lemma}\label{lem:links}
Fix a replica $r$. Distinct matched variables carry distinct links of replica $r$, and:
\begin{enumerate}[label=\textup{(\roman*)},leftmargin=2.4em,itemsep=2pt]
\item the link of replica $r$ at a variable carrying two self-loops has both endpoints in $V^{(r)}(\mathbf M)$,
  and
  \begin{equation}\label{eq:halfedges}
    2(|E_{\cap}|-|W_{\cap}|)+|\WP|+|\WS|+|\WR|\ \le\
    q(|\mathbf M_{V}|-2)+|V^{(r)}(\mathbf M)\cap\Phi|+2L_{\Phi} ,
  \end{equation}
\item a link of $\WS$ lies at a variable carrying exactly one self-loop, and
  $|\WS|\le2|W_{\cap}|-|E_{\cap}|$,
\item a link of $\WD$ or of $\WR$ lies at a matched variable carrying no self-loop, and
  $|\WD|+|\WR|\le|\mathbf M_{W}|-|W_{\cap}|$.
\end{enumerate}
\end{lemma}

\begin{proof}
Every bound of the lemma compares a set of links of replica $r$ with the number of matched
variables of one kind, carrying two, one or no self-loop, at which these links lie. Item
\textup{(i)} also counts the ends of links at $V^{(r)}(\mathbf M)$.

\noindent\textbf{The first claim.} A matched variable of $\Gcup$ is the variable obtained by identifying the
two variables of one pair of $\mathbf M_{W}$, one of replica $1$ and one of replica $2$, and every
variable of $G$ occurs in at most one pair of $\mathbf M_{W}$ (\S\ref{sec:graphs}). So a matched
variable carries exactly one link of replica $r$, and distinct matched variables carry distinct
ones. By \S\ref{sec:linkpairings}, the matched variables carrying two, one and no self-loop number
\begin{equation}\label{eq:IIb-kinds}
  |E_{\cap}|-|W_{\cap}|,\qquad 2|W_{\cap}|-|E_{\cap}|,\qquad |\mathbf M_{W}|-|W_{\cap}| .
\end{equation}
Hence the links of replica $r$ at the matched variables of each of these three kinds are exactly
as many as the corresponding number in \eqref{eq:IIb-kinds}, and a set of such links has at most
that many elements.

\noindent\textbf{(i)} At a variable carrying two self-loops the two links share both their
individuals (\S\ref{sec:linkpairings}), which are therefore matched. So the link of replica $r$ there
has both endpoints in $V^{(r)}(\mathbf M)$. The set $V^{(r)}(\mathbf M)$ consists of $v_{1}$ and
$v_{2}$, of degree $L_{\Phi}$, and of $|\mathbf M_{V}|-2$ nodes of $B$, of degree
$q+1$ in $\Phi$ and $q$ outside it (\S\ref{sec:template}), and a link
has one end at each of its two endpoints. Counting the ends of links at $V^{(r)}(\mathbf M)$ in
these two ways,
\[
\begin{aligned}
  q(|\mathbf M_{V}|-2)+|V^{(r)}(\mathbf M)\cap\Phi|+2L_{\Phi}
  &=\sum_{v\in V^{(r)}(\mathbf M)}\deg_{G}(v)\\
  &=2\,\big|\{w\in W:N(w)\subseteq V^{(r)}(\mathbf M)\}\big|\\
  &\qquad+\big|\{w\in W:|N(w)\cap V^{(r)}(\mathbf M)|=1\}\big| .
\end{aligned}
\]
The links with exactly one endpoint in $V^{(r)}(\mathbf M)$ are partitioned by $\WP$, $\WS$,
$\WR$ (\S\ref{sec:linksets}), and, by the first sentence of this item, the first claim and
\eqref{eq:IIb-kinds}, at least $|E_{\cap}|-|W_{\cap}|$ links have both endpoints in
$V^{(r)}(\mathbf M)$:
\[
\begin{aligned}
  \big|\{w\in W:|N(w)\cap V^{(r)}(\mathbf M)|=1\}\big|\ &=\ |\WP|+|\WS|+|\WR|,\\
  \big|\{w\in W:N(w)\subseteq V^{(r)}(\mathbf M)\}\big|\ &\ge\ |E_{\cap}|-|W_{\cap}| .
\end{aligned}
\]
Inserting these two relations into the previous display gives \eqref{eq:halfedges}.

\noindent\textbf{(ii)} The variable of a link of $\WS$ is matched, and $\Pai$ pairs each of the
two ends of the link with an end of the other replica's link (\S\ref{sec:linksets}). These are the
two pairs of $\Pai$ at that variable. The pair at the matched endpoint is a self-loop, by the
definition of $\WS$. The pair at the unmatched endpoint is not a self-loop (\S\ref{sec:linksets}). So the variable carries exactly one
self-loop, and the first claim and \eqref{eq:IIb-kinds} give
\[
  |\WS|\ \le\ 2|W_{\cap}|-|E_{\cap}| .
\]

\noindent\textbf{(iii)} Let $v,v'$ be the endpoints of a link of $\WD$. The pair formed by its
two ends in replica $r$ would be an edge of $\Gomi$ joining $\pi^{r}_{1}(v)$ to
$\pi^{r}_{1}(v')$, which lie in different components of $\Gomi$. So $\Pai$ does not contain that
pair. An unmatched variable has that pairing only (\S\ref{sec:linkpairings}). So the variable of the
link is matched, and $\Pai$ pairs each of the two ends with an end of the other replica's link.
Neither pair is a self-loop, as $\pi^{r}_{1}(v)$ and $\pi^{r}_{1}(v')$ are unmatched and so
belong to replica $r$ alone. A link of $\WR$ lies at a matched variable, and neither of the two
pairs there is a self-loop: the one at the matched endpoint by the definition of $\WR$, the one at
the unmatched endpoint as in item \textup{(ii)}. Finally $\WD\cap\WR=\varnothing$, a link of
$\WR$ having an endpoint in $V^{(r)}(\mathbf M)$ and a link of $\WD$ none. By the first claim and
\eqref{eq:IIb-kinds},
\[
  |\WD|+|\WR|\ =\ |\WD\cup\WR|\ \le\ |\mathbf M_{W}|-|W_{\cap}| .
\]
\end{proof}

\begin{lemma}\label{lem:traces}
Fix a replica $r$. The traces $\Vomi{r}(C)$, $C\in\CC_{\Gomi}$, are pairwise disjoint with union
$B\setminus V^{(r)}(\mathbf M)$, at most one of them is large, and:
\begin{enumerate}[label=\textup{(\roman*)},leftmargin=2.4em,itemsep=2pt]
\item summing over the components,
  \begin{equation}\label{eq:bookkeeping}
  \begin{gathered}
    \sum_{C\in\CC_{\Gomi}}\partial(\Vomi{r}(C))=|\WP|+|\WS|+|\WR|+2|\WD|,\\
    \sum_{C\in\CC_{\Gomi}}\big|\Vomi{r}(C)\cap\Phi\big|=2L_{\Phi}-|V^{(r)}(\mathbf M)\cap\Phi| ,
  \end{gathered}
  \end{equation}
\item $|\WP|\le\sum_{C\in\CC_{\Gomi}\setminus\CC_{\Gomi}(\mathbf M)}\partial(\Vomi{r}(C))$,
\item for every $C\in\CC_{\Gomi}(\mathbf M)$, both $\Vomi{1}(C)$ and $\Vomi{2}(C)$ are nonempty.
\end{enumerate}
\end{lemma}

\begin{proof}
Each unmatched individual of replica $r$ lies in exactly one component of $\Gomi$, which gives
the first claims and, counting each link at the traces holding its endpoints, item \textup{(i)}.
Items \textup{(ii)} and \textup{(iii)} follow a link from an endpoint in a trace $\Vomi{r}(C)$: the
pair of $\Pai$ at that endpoint is an edge of $\Gomi$, so its other individual lies in $C$ as well.

\noindent\textbf{The first claims.} For $v\in B\setminus V^{(r)}(\mathbf M)$ the individual
$\pi^{r}_{1}(v)$ is unmatched (\S\ref{sec:second}) and lies in exactly one component of $\Gomi$,
the components partitioning the individual nodes of $\Gcup$ (\S\ref{sec:graphs},
\S\ref{sec:ubmodel}). By Definition~\ref{def:trace}, $v$ therefore lies in exactly one trace:
\begin{equation}\label{eq:IIb-partition}
  \Vomi{r}(C)\cap\Vomi{r}(C')=\varnothing\quad(C\ne C'),\qquad
  \bigcup_{C\in\CC_{\Gomi}}\Vomi{r}(C)\ =\ B\setminus V^{(r)}(\mathbf M) .
\end{equation}
Two disjoint subsets of $B$ with more than $|B|/2$ elements each would have more than $|B|$
elements together. So at most one trace is large.

\noindent\textbf{(i)} By \eqref{eq:IIb-partition}, $\sum_{C\in\CC_{\Gomi}}\partial(\Vomi{r}(C))$ counts every
link of $G$ once for each trace holding exactly one of its endpoints. A link of $\WP\cup\WS\cup\WR$ has exactly
one endpoint in $B\setminus V^{(r)}(\mathbf M)$, so it is counted once. A link with both endpoints in
$V^{(r)}(\mathbf M)$ is not counted. Every other link has both endpoints in $B\setminus V^{(r)}(\mathbf M)$, as
$v_{1},v_{2}\in V^{(r)}(\mathbf M)$, and its two individuals in replica $r$ are unmatched (\S\ref{sec:second}). If
they lie in different components of $\Gomi$, the link lies in $\WD$ and its endpoints lie in two distinct traces,
so it is counted twice. Otherwise its two endpoints lie in one trace, and it is not counted. This gives the
first identity. As $\Phi\subseteq B$, \eqref{eq:IIb-partition} gives the second:
\[
  \sum_{C\in\CC_{\Gomi}}\big|\Vomi{r}(C)\cap\Phi\big|=\big|\Phi\setminus V^{(r)}(\mathbf M)\big|
  =2L_{\Phi}-|V^{(r)}(\mathbf M)\cap\Phi| .
\]

\noindent\textbf{(ii)} Let a link of $\WP$ have unmatched endpoint $v$ and matched endpoint $v'$, and let $C$ be the
component of $\Gomi$ holding $\pi^{r}_{1}(v)$, so that $v\in\Vomi{r}(C)$. The pair of $\Pai$ formed by the two ends
of the link in replica $r$ joins the distinct individuals $\pi^{r}_{1}(v)$ and $\pi^{r}_{1}(v')$, so it is an edge
of $\Gomi$, and $\pi^{r}_{1}(v')$ is matched. So $C$ is pinned. As $v'\in V^{(r)}(\mathbf M)$ lies outside
$\Vomi{r}(C)$, the link is counted by $\partial(\Vomi{r}(C))$. Distinct links being counted separately, this
gives \textup{(ii)}.

\noindent\textbf{(iii)} Let $C\in\CC_{\Gomi}(\mathbf M)$.

\noindent\uline{Step 1: a nonempty trace.} $C$ holds an individual, which is unmatched, hence of the form
$\pi^{1}_{1}(v)$ with $v\in B\setminus V^{(1)}(\mathbf M)$ or $\pi^{2}_{1}(v)$ with
$v\in B\setminus V^{(2)}(\mathbf M)$ (\S\ref{sec:second}). So $\Vomi{1}(C)$ or $\Vomi{2}(C)$ is
nonempty. As the replicas play symmetric roles, we may
assume $\Vomi{1}(C)\ne\varnothing$.

\noindent\uline{Step 2: a link leaving $\Vomi{1}(C)$.} The set $\Vomi{1}(C)$ is nonempty and contained in
$B\subsetneq V$, as $v_{1},v_{2}\notin B$. Hence, $G$ being connected (\S\ref{sec:template}), some link
$\{v,v'\}$ of $G$ has $v\in\Vomi{1}(C)$ and $v'\notin\Vomi{1}(C)$.

\noindent\uline{Step 3: its two ends in replica $1$ are not paired together.} Otherwise that pair would be an edge of
$\Gomi$ joining $\pi^{1}_{1}(v)\in V_{C}$ to $\pi^{1}_{1}(v')$, so $\pi^{1}_{1}(v')\in V_{C}$.
Then, $C$ being free, $\pi^{1}_{1}(v')$ would be unmatched, so $v'\in B\setminus V^{(1)}(\mathbf M)$,
and then $v'\in\Vomi{1}(C)$, a contradiction.

\noindent\uline{Step 4: a node of $\Vomi{2}(C)$.} So the variable of the link is matched, an unmatched one having only
the pairing just excluded (\S\ref{sec:linkpairings}), and $\Pai$ pairs the end at $\pi^{1}_{1}(v)$ with
an end of the replica-$2$ link at the same variable, at $\pi^{2}_{1}(v'')$ say. As
$\pi^{1}_{1}(v)$ is unmatched, $\pi^{2}_{1}(v'')\ne\pi^{1}_{1}(v)$, so this pair is not a
self-loop: it is an edge of $\Gomi$, and $\pi^{2}_{1}(v'')\in V_{C}$. $C$ being free,
$\pi^{2}_{1}(v'')$ is unmatched, so $v''\in B\setminus V^{(2)}(\mathbf M)$ and
$v''\in\Vomi{2}(C)$.
\end{proof}

\subsubsection{The bound on $\mathcal E$ as a trade-off between matched fasteners and frontier links}

\begin{lemma}\label{lem:parts}
Assume \ref{H1} and \ref{H3}--\ref{H5}, and fix a replica $r$. Every pair $(\mathbf M,\Pai)$ satisfies
\[
\begin{aligned}
  -\mathcal E-\frac{2}{L_{B}}\big(|\mathbf M_{V}|-2+|\mathbf M_{W}|\big)
  \ \ge\ &\frac{4\zeta}{L_{B}}(|\mathbf M_{V}|-2)-2\zeta\,|V^{(r)}(\mathbf M)\cap\Phi|\\
  &+\zeta\sum_{C\in\CC_{\Gomi}}\partial_{B}\big(\Vomi{r}(C)\big)
   +\sum_{C\in\CC_{\Gomi}(\mathbf M)}\Big(\theta\,\partial\big(\Vomi{r}(C)\big)-1\Big).
\end{aligned}
\]

\end{lemma}

\begin{proof}
We bound the four groups of \eqref{eq:Ekinds} one at a time, \textup{(I)} and \textup{(II)} by
Lemmas~\ref{lem:links} and~\ref{lem:traces}. Their bounds contain $2\zeta|\WD|$ and $-2\zeta|\WD|$,
which cancel in the sum.

\noindent\uline{Step 1: positivity of $\zeta,\theta$.} By \S\ref{sec:exponents}, $\zeta>2/L_{B}>0$ and $\theta\ge1/L_{\Phi}+1/L_{B}>0$.

\noindent\uline{Step 2: groups \textup{(III)} and \textup{(IV)}.} The components in
$\CC_{\Gomi}\setminus\CC_{\Gomi}(\mathbf M)$ are disjoint and each holds a matched individual.
There are $|\mathbf M_{V}|$ matched individuals (\S\ref{sec:second}), and if $1\sim_{\Gomi}2$, one
of these components holds the two matched individuals $1$ and $2$. So
$|\CC_{\Gomi}\setminus\CC_{\Gomi}(\mathbf M)|\le|\mathbf M_{V}|-\one\{1\sim_{\Gomi}2\}$. By
\eqref{eq:C}, the coefficient of $-(|\mathbf M_{V}|-2)$ in \textup{(IV)} is at least
$4\zeta/L_{B}$, and $|\mathbf M_{V}|-2\ge0$. Hence
\[
  \textup{(III)}\ \le\ 0,\qquad \textup{(IV)}\ \le\ -\frac{4\zeta}{L_{B}}(|\mathbf M_{V}|-2).
\]

\noindent\uline{Step 3: group \textup{(I)}.} By \eqref{eq:halfedges} and $\zeta>0$, then by \eqref{eq:bookkeeping},
then by \eqref{eq:boundary},
\[
\begin{aligned}
  \textup{(I)}\ &\le\ \zeta\big(|V^{(r)}(\mathbf M)\cap\Phi|+2L_{\Phi}-|\WP|-|\WS|-|\WR|\big)\\
  &=\ \zeta\Big(2|V^{(r)}(\mathbf M)\cap\Phi|+\sum_{C\in\CC_{\Gomi}}\big|\Vomi{r}(C)\cap\Phi\big|
   -\sum_{C\in\CC_{\Gomi}}\partial(\Vomi{r}(C))+2|\WD|\Big)\\
  &=\ 2\zeta\,|V^{(r)}(\mathbf M)\cap\Phi|
   -\zeta\sum_{C\in\CC_{\Gomi}}\partial_{B}(\Vomi{r}(C))+2\zeta\,|\WD| .
\end{aligned}
\]

\noindent\uline{Step 4: group \textup{(II)}.} By Lemma~\ref{lem:links}\textup{(ii)} and $\theta+1/L_{B}\ge\theta>0$,
and by Lemma~\ref{lem:links}\textup{(iii)} and $2(\zeta+\theta)\ge\theta>0$,
\[
\begin{aligned}
  -\Big(\theta+\frac1{L_{B}}\Big)(2|W_{\cap}|-|E_{\cap}|)\ &\le\ -\theta\,|\WS|,\\
  -2(\zeta+\theta)(|\mathbf M_{W}|-|W_{\cap}|)\ &\le\ -2(\zeta+\theta)|\WD|-\theta\,|\WR| .
\end{aligned}
\]
By \eqref{eq:bookkeeping}, then by Lemma~\ref{lem:traces}\textup{(ii)} and $\theta>0$,
\[
\begin{aligned}
  -\theta\big(|\WS|+|\WR|\big)\ &=\ -\theta\Big(\sum_{C\in\CC_{\Gomi}}\partial(\Vomi{r}(C))-|\WP|-2|\WD|\Big)\\
  &\le\ -\theta\sum_{C\in\CC_{\Gomi}(\mathbf M)}\partial(\Vomi{r}(C))+2\theta\,|\WD| .
\end{aligned}
\]
Inserting these bounds into \textup{(II)}, and writing $|\CC_{\Gomi}(\mathbf M)|$ as $\sum_{C\in\CC_{\Gomi}(\mathbf M)}1$,
\[
  \textup{(II)}\ \le\ -\sum_{C\in\CC_{\Gomi}(\mathbf M)}\big(\theta\,\partial(\Vomi{r}(C))-1\big)-2\zeta\,|\WD| .
\]

\noindent\uline{Step 5: the components.} Adding the bounds on \textup{(I)} and \textup{(II)}, the terms
$2\zeta|\WD|$ cancel:
\begin{equation}\label{eq:IIb-onetwo}
\begin{aligned}
  \textup{(I)}+\textup{(II)}\ \le\ &2\zeta\,|V^{(r)}(\mathbf M)\cap\Phi|
  -\zeta\sum_{C\in\CC_{\Gomi}}\partial_{B}(\Vomi{r}(C))\\
  &-\sum_{C\in\CC_{\Gomi}(\mathbf M)}\big(\theta\,\partial(\Vomi{r}(C))-1\big) .
\end{aligned}
\end{equation}

\noindent\uline{Step 6: conclusion.} By \eqref{eq:Ekinds}, adding \eqref{eq:IIb-onetwo} and the bounds on
\textup{(III)} and \textup{(IV)},
\[
\begin{aligned}
  &\mathcal E+\frac{2}{L_{B}}\big(|\mathbf M_{V}|-2+|E_{\cap}|+|\mathbf M_{W}|-|W_{\cap}|\big)
  \ \le\ 2\zeta\,|V^{(r)}(\mathbf M)\cap\Phi|-\frac{4\zeta}{L_{B}}(|\mathbf M_{V}|-2)\\
  &\qquad-\zeta\sum_{C\in\CC_{\Gomi}}\partial_{B}(\Vomi{r}(C))
   -\sum_{C\in\CC_{\Gomi}(\mathbf M)}\big(\theta\,\partial(\Vomi{r}(C))-1\big).
\end{aligned}
\]
As $|E_{\cap}|-|W_{\cap}|\ge0$, a number of variables (\S\ref{sec:linkpairings}), the left-hand side is
at least $\mathcal E+\frac{2}{L_{B}}\big(|\mathbf M_{V}|-2+|\mathbf M_{W}|\big)$. Multiplying by
$-1$ gives the lemma.
\end{proof}

\subsubsection{Bounds per trace}

In Lemmas~\ref{lem:large} and \ref{lem:small}, \ref{H1}--\ref{H5} are assumed, $r$ is a replica
and $C\in\CC_{\Gomi}$. By \S\ref{sec:exponents}, $\zeta>0$ and $\theta\ge1/L_{\Phi}+1/L_{B}>0$.

\begin{lemma}[a large trace]\label{lem:large}
Let $\Vomi{r}(C)$ be large. Then \textup{(i)}
$\zeta\,\partial_{B}(\Vomi{r}(C))\ge2\zeta|\Phi\setminus\Vomi{r}(C)|\ge2\zeta|V^{(r)}(\mathbf M)\cap\Phi|$,
and \textup{(ii)} $\theta\,\partial(\Vomi{r}(C))>1$.
\end{lemma}

\begin{proof}
Both parts apply Lemma~\ref{lem:fast} to the complement $B\setminus\Vomi{r}(C)$ of the large
trace, which has fewer than $|B|/2$ nodes and contains every fastener outside the trace.

\noindent\textbf{(i)} As $\Vomi{r}(C)\subseteq B\setminus V^{(r)}(\mathbf M)$ and
$\Phi\subseteq B$,
\[
  V^{(r)}(\mathbf M)\cap\Phi\ \subseteq\ \Phi\setminus\Vomi{r}(C) .
\]
As $\Vomi{r}(C)$ is large, $|B\setminus\Vomi{r}(C)|<|B|/2$, so Lemma~\ref{lem:fast} applies to
$B\setminus\Vomi{r}(C)$, and
\[
  \partial_{B}(\Vomi{r}(C))\ =\ \partial_{B}\big(B\setminus\Vomi{r}(C)\big)
  \ \ge\ 2\,\big|\Phi\setminus\Vomi{r}(C)\big|\ \ge\ 2\,|V^{(r)}(\mathbf M)\cap\Phi| :
\]
the equality because an edge of $G_{B}$ has exactly one endpoint in $\Vomi{r}(C)$ if and only if
it has exactly one in $B\setminus\Vomi{r}(C)$, the first inequality by Lemma~\ref{lem:fast} and
$\Phi\subseteq B$, the second by the inclusion above. Multiplying by $\zeta>0$ gives
\textup{(i)}.

\noindent\textbf{(ii)} By \eqref{eq:boundary} and the previous display,
\[
  \partial(\Vomi{r}(C))\ =\ \partial_{B}(\Vomi{r}(C))+|\Vomi{r}(C)\cap\Phi|
  \ \ge\ 2\,|\Phi\setminus\Vomi{r}(C)|+|\Vomi{r}(C)\cap\Phi|\ \ge\ |\Phi|\ =\ 2L_{\Phi} .
\]
Multiplying by $\theta>0$, and then by \eqref{eq:H4bis},
\[
  \theta\,\partial(\Vomi{r}(C))\ \ge\ 2L_{\Phi}\,\theta\ \ge\ 2L_{\Phi}\Big(\frac1{L_{\Phi}}+\frac1{L_{B}}\Big)
  \ =\ 2+\frac{2L_{\Phi}}{L_{B}}\ >\ 1 .
\]
\end{proof}

\begin{lemma}[a trace that is not large]\label{lem:small}
Let $\Vomi{r}(C)$ be not large. Then \textup{(i)}
$\zeta\,\partial_{B}(\Vomi{r}(C))\ge\frac{4\zeta}{L_{B}}|\Vomi{r}(C)|$, and \textup{(ii)} if
$C\in\CC_{\Gomi}(\mathbf M)$, also
$\zeta\,\partial_{B}(\Vomi{r}(C))+\theta\,\partial(\Vomi{r}(C))-1\ge\frac{4\zeta}{L_{B}}|\Vomi{r}(C)|$.
\end{lemma}

\begin{proof}
Part \textup{(i)} is the second inequality of \eqref{eq:expand} with \ref{H2}. In part
\textup{(ii)}, the first inequality of \eqref{eq:expand} and \eqref{eq:E} show that
$(\zeta+\theta)\,\partial_{B}(\Vomi{r}(C))-1$ is at least a fixed fraction of
$(\zeta+\theta)\,\partial_{B}(\Vomi{r}(C))$, which the second inequality of \eqref{eq:expand} bounds
below by a multiple of $|\Vomi{r}(C)|$.

\noindent\textbf{(i)} As $\Vomi{r}(C)$ is not large, $|\Vomi{r}(C)|\le|B|/2$, so
$\min\{|\Vomi{r}(C)|,|B|-|\Vomi{r}(C)|\}=|\Vomi{r}(C)|$. The second inequality of
\eqref{eq:expand}, then $\ciso L_{B}\ge8$ \textup{(}\ref{H2}\textup{)}, give
\begin{equation}\label{eq:IIb-small}
  \partial_{B}(\Vomi{r}(C))\ \ge\ \ciso\,|\Vomi{r}(C)|\ \ge\ \frac{8}{L_{B}}\,|\Vomi{r}(C)|
  \ \ge\ \frac{4}{L_{B}}\,|\Vomi{r}(C)| .
\end{equation}
Multiplying by $\zeta>0$ gives \textup{(i)}.

\noindent\textbf{(ii)} \uline{Step 1: a lower bound.} Let $C\in\CC_{\Gomi}(\mathbf M)$. By
Lemma~\ref{lem:traces}\textup{(iii)}, $\Vomi{r}(C)\ne\varnothing$, and
$|\Vomi{r}(C)|\le|B|/2<|B|$. So $\Vomi{r}(C)$ is a nonempty proper subset of $B$. The first
inequality of \eqref{eq:expand} and $\zeta+\theta>0$, then \eqref{eq:E}, give
\[
  (\zeta+\theta)\,\partial_{B}(\Vomi{r}(C))\ \ge\ q(\zeta+\theta)
  \ \ge\ 1+\frac{8}{\ciso L_{B}} .
\]

\noindent\uline{Step 2: conclusion.} By \eqref{eq:boundary},
$\partial(\Vomi{r}(C))\ge\partial_{B}(\Vomi{r}(C))$, so, as $\theta>0$,
\[
  \zeta\,\partial_{B}(\Vomi{r}(C))+\theta\,\partial(\Vomi{r}(C))-1
  \ \ge\ (\zeta+\theta)\,\partial_{B}(\Vomi{r}(C))-1 .
\]
Dividing the lower bound of Step~1 by $1+8/(\ciso L_{B})$ gives
$1\le(\zeta+\theta)\,\partial_{B}(\Vomi{r}(C))/\big(1+8/(\ciso L_{B})\big)$. Hence
\[
\begin{aligned}
  (\zeta+\theta)\,\partial_{B}(\Vomi{r}(C))-1\ &\ge\ \frac{8/(\ciso L_{B})}{1+8/(\ciso L_{B})}
   (\zeta+\theta)\,\partial_{B}(\Vomi{r}(C))\\
  &\ge\ \frac{4\zeta}{\ciso L_{B}}\,\partial_{B}(\Vomi{r}(C))\ \ge\ \frac{4\zeta}{L_{B}}\,|\Vomi{r}(C)| ,
\end{aligned}
\]
the second inequality because $8/(\ciso L_{B})\le1$ by \ref{H2} and $\theta>0$, the third by the
first inequality of \eqref{eq:IIb-small} and $\zeta>0$. The last two displays give \textup{(ii)}.
\end{proof}

\subsubsection{The final bound on $\mathcal E$}

\begin{theorem}\label{thm:master}
Assume \ref{H1}--\ref{H5}. Then \emph{every} pair $(\mathbf M,\Pai)$ satisfies
\[
  \mathcal E\ \le\ -\frac{2}{L_{B}}\,\big(|\mathbf M_{V}|-2+|\mathbf M_{W}|\big) .
\]
\end{theorem}

\begin{proof}
By Lemma~\ref{lem:parts}, it suffices to show that, in one replica, the terms of the traces and
of the matched nodes of $B$ add up to at least $2\zeta$ times the number of matched fasteners. A large
trace does this alone, by Lemma~\ref{lem:large}. When no trace is large, Lemma~\ref{lem:small}
bounds the term of each trace below by a multiple of its size, and the traces together with the
matched nodes of $B$ fill $B$, which has $L_{B}$ nodes per fastener.

\noindent\uline{Step 1: reduction.} Fix a replica $r$. By Lemma~\ref{lem:parts}, it suffices to show
\begin{equation}\label{eq:masterbracket}
  \frac{4\zeta}{L_{B}}(|\mathbf M_{V}|-2)+\zeta\sum_{C\in\CC_{\Gomi}}\partial_{B}(\Vomi{r}(C))
  +\sum_{C\in\CC_{\Gomi}(\mathbf M)}\big(\theta\,\partial(\Vomi{r}(C))-1\big)
  \ \ge\ 2\zeta\,|V^{(r)}(\mathbf M)\cap\Phi| .
\end{equation}

\noindent\uline{Step 2: putting together results on connected fastener sets and frontier links.}
The traces of replica $r$ are disjoint with union $B\setminus V^{(r)}(\mathbf M)$, at most one of
them is large (Lemma~\ref{lem:traces}), and, by Lemma~\ref{lem:small}, every $C$ whose trace is not large
satisfies
\begin{equation}\label{eq:masterb}
  \zeta\,\partial_{B}(\Vomi{r}(C))+\one\{C\in\CC_{\Gomi}(\mathbf M)\}\big(\theta\,\partial(\Vomi{r}(C))-1\big)\ \ge\ \frac{4\zeta}{L_{B}}\,|\Vomi{r}(C)|\ \ge\ 0 .
\end{equation}
The left-hand side of \eqref{eq:masterbracket} is $(4\zeta/L_{B})(|\mathbf M_{V}|-2)$ plus the sum over
$C\in\CC_{\Gomi}$ of the left-hand side of \eqref{eq:masterb}.

\noindent\uline{Step 3: final proof, distinguishing on whether or not there is a large trace.}\\
\textbf{Case 1: a trace of replica $r$, say $\Vomi{r}(C_{0})$, is large.} By Lemma~\ref{lem:large},
\[
  \zeta\,\partial_{B}(\Vomi{r}(C_{0}))\ \ge\ 2\zeta\,|V^{(r)}(\mathbf M)\cap\Phi| ,\qquad
  \one\{C_{0}\in\CC_{\Gomi}(\mathbf M)\}\big(\theta\,\partial(\Vomi{r}(C_{0}))-1\big)\ \ge\ 0 .
\]
The other traces are not large, so \eqref{eq:masterb} holds for every $C\ne C_{0}$, and
$(4\zeta/L_{B})(|\mathbf M_{V}|-2)\ge0$. Adding these inequalities gives \eqref{eq:masterbracket}.

\textbf{Case 2: no trace of replica $r$ is large.} Then \eqref{eq:masterb} holds for every
$C$, and summing it,
\[
\begin{aligned}
  &\frac{4\zeta}{L_{B}}(|\mathbf M_{V}|-2)+\zeta\sum_{C\in\CC_{\Gomi}}\partial_{B}(\Vomi{r}(C))
   +\sum_{C\in\CC_{\Gomi}(\mathbf M)}\big(\theta\,\partial(\Vomi{r}(C))-1\big)\\
  &\quad\ge\ \frac{4\zeta}{L_{B}}(|\mathbf M_{V}|-2)+\frac{4\zeta}{L_{B}}\sum_{C\in\CC_{\Gomi}}|\Vomi{r}(C)|\\
  &\quad=\ \frac{4\zeta}{L_{B}}\big(|\mathbf M_{V}|-2+|B|-|\mathbf M_{V}|+2\big)
  \ =\ \frac{4\zeta}{L_{B}}|B|\ =\ 4\zeta|\Phi| ,
\end{aligned}
\]
as $|B|=2L_{\Phi}L_{B}=L_{B}|\Phi|$. As $V^{(r)}(\mathbf M)\cap\Phi\subseteq\Phi$ and $\zeta>0$, this is at least
$2\zeta|V^{(r)}(\mathbf M)\cap\Phi|$, which is \eqref{eq:masterbracket}.
\end{proof}

\subsection{The summation, and the proof of Theorem~\ref{thm:UB}}\label{sec:sum}

\subsubsection{The summation}

\begin{proposition}\label{prop:summation}
Assume \ref{H1}, \ref{H3} and \ref{H5}, that $3|V|^{2}|W|^{2}K^{-2/L_{B}}\le\tfrac18$, and that every pair
$(\mathbf M,\Pai)$ satisfies $\mathcal E\ \le\ -\frac2{L_{B}}\big(|\mathbf M_{V}|-2+|\mathbf M_{W}|\big)$. Then
\[
  \Var(P_{G}\mid x=1)+\Var(P_{G}\mid x=0)\ \le\ \big(128\,|W|^{2}K^{-1/q}+24\,|V|^{2}|W|^{2}K^{-2/L_{B}}\big)\,
  \hat s_{G}^{2}K^{2} .
\]
\end{proposition}

\begin{proof}
We group the pairs $(\mathbf M,\Pai)$ by the number $\bar m$ of pairs of $\mathbf M$ other than
$(v_{1},v_{1})$ and $(v_{2},v_{2})$. There are at most $(3|V|^{2}|W|^{2})^{\bar m}$ of them, and, by
Proposition~\ref{prop:exponent} and the hypothesis on $\mathcal E$, each contributes at most
$K^{-2\bar m/L_{B}}\hat s_{G}^{2}K^{2}$, up to a factor close to $1$. The single pair with $\bar m=0$
gives the squared mean $\hat s_{G}^{2}K^{2}$ of Lemma~\ref{lem:gap} under $x=1$, up to that factor,
and $0$ under $x=0$.

\noindent\uline{Step 1: counting the number of pairs leading to the same exponent.} For a pair $(\mathbf M,\Pai)$
let $\bar m_V:=|\mathbf M_{V}|-2$, $\bar m_W:=|\mathbf M_{W}|$ and $\bar m:=\bar m_V+\bar m_W$. Then $\mathbf M_{V}$
is $\{(v_{1},v_{1}),(v_{2},v_{2})\}$ plus $\bar m_V$ pairs in $B\times B$, $\mathbf M_{W}$ consists of $\bar m_W$
pairs in $W\times W$, and $|\PcalSL(\Gcup)|=3^{\bar m_W}$ (\S\ref{sec:linkpairings}). So the pairs with a given
$\bar m$ number at most
\[
  \sum_{\bar m_V+\bar m_W=\bar m}\big(|V|^{2}\big)^{\bar m_V}\big(3|W|^{2}\big)^{\bar m_W}
  \ \le\ \big(|V|^{2}+3|W|^{2}\big)^{\bar m}\ \le\ \big(3|V|^{2}|W|^{2}\big)^{\bar m} ,
\]
the last inequality as $|V|\ge2$ and $|W|\ge1$. In particular at most one pair has $\bar m=0$. As no variable is then matched, its $\Pai$ joins the two ends of
every link (\S\ref{sec:linkpairings}). So $|E_{\cap}|=0$, and $\Gomi$ is the union of the two replicas, which are
connected and share $1$ and $2$. Hence $|\CC_{\Gomi}|=1$, $1\sim_{\Gomi}2$, and $\mathcal E=0$ by
\eqref{eq:exponent}.

\noindent\uline{Step 2: bounding the moments.}  The hypothesis reads
$\mathcal E\le-2\bar m/L_{B}$, so, by Step~1 and the geometric series,
\begin{equation}\label{eq:IIb-geom}
  \sum_{(\mathbf M,\Pai):\ \bar m\ge1}K^{\mathcal E}
  \ \le\ \sum_{\bar m\ge1}\big(3|V|^{2}|W|^{2}K^{-2/L_{B}}\big)^{\bar m}
  \ \le\ 6\,|V|^{2}|W|^{2}\,K^{-2/L_{B}} .
\end{equation}
By \eqref{eq:partition} and Proposition~\ref{prop:exponent}, every pair contributes at most
$\hat s_{G}^{2}K^{2}(1+32|W|^{2}K^{-1/q})K^{\mathcal E}$ to $\E[P_{G}^{2}\mid x=i]$, $i\in\{0,1\}$, and the pair with
$\bar m=0$ contributes $0$ given $x=0$, as $1\sim_{\Gomi}2$. So, by Step~1 and \eqref{eq:IIb-geom},
$\E[P_{G}^{2}\mid x=1]$ is at most $\hat s_{G}^{2}K^{2}(1+32|W|^{2}K^{-1/q})(1+6|V|^{2}|W|^{2}K^{-2/L_{B}})$, and
$\E[P_{G}^{2}\mid x=0]$ at most this minus $\hat s_{G}^{2}K^{2}$. As $\E[P_{G}\mid x=1]=\hat s_{G}K$ and
$\E[P_{G}\mid x=0]=0$ (Lemma~\ref{lem:gap}), for $i\in\{0,1\}$,
\[
\begin{aligned}
  \Var(P_{G}\mid x=i)\ &\le\ \hat s_{G}^{2}K^{2}\Big[\big(1+32|W|^{2}K^{-1/q}\big)\big(1+6|V|^{2}|W|^{2}K^{-2/L_{B}}\big)-1\Big]\\
  &\le\ 2\big(32|W|^{2}K^{-1/q}+6|V|^{2}|W|^{2}K^{-2/L_{B}}\big)\,\hat s_{G}^{2}K^{2} ,
\end{aligned}
\]
as $(1+a)(1+b)-1=a+b+ab\le2(a+b)$ when $0\le a\le1$ and $b\ge0$, and $32|W|^{2}K^{-1/q}\le1$ by \ref{H5}. Adding the two variances gives the claim.
\end{proof}

\subsubsection{Proof of Theorem~\ref{thm:UB}}

\begin{proof}[Proof of Theorem~\ref{thm:UB}]
We take $L_{\Phi}=L_{B}=L$ and the body of Proposition~\ref{prop:body}, check
\ref{H1}--\ref{H5} from Assumption~\ref{ass:B}, and combine Lemma~\ref{lem:gap},
Theorem~\ref{thm:master}, Proposition~\ref{prop:summation} and Lemma~\ref{lem:criterion}. The two
smallness conditions \eqref{eq:IIb-assume} on $|V|$ and $|W|$ may be assumed, as the error bound
exceeds $1$ otherwise.

\noindent\uline{Step 1: the template.} Put $L_{\Phi}:=L_{B}:=L$. As $L\ge(400q)^{600q^{2}}$,
$L_{B}>q^{1601}$ and $|B|=2L^{2}\ge(400q)^{1200q^{2}}$, so
Proposition~\ref{prop:body} applies, as $q\ge3$. Let $G$ be the template built on its body and
fasteners, for which $\ciso=1/400$. By \eqref{eq:PGII} and \eqref{eq:W-count}, $P_{G}$ is a polynomial of degree
$2|W|=2L(q L+2)\le2(q+2)L^{2}$ in the entries $Y_{ij}$, and, as $|V|=2L^{2}+2$ and
$q\ge3$, $|V|$ and $|W|$ are at most $3q L^{2}$. So
\begin{equation}\label{eq:IIb-target}
  512\,|W|^{2}K^{-1/q}+96\,|V|^{2}|W|^{2}K^{-2/L}\ \le\ 4608q^{2}L^{4}K^{-1/q}+7776q^{4}L^{8}K^{-2/L} ,
\end{equation}
and it suffices to show that the error probability is at most the left-hand side.

\noindent\uline{Step 2: two inequalities that may be assumed.} An error probability is at most $1$. If
$K^{1/q}<400|W|^{2}$, then $512|W|^{2}K^{-1/q}>1$. If $3|V|^{2}|W|^{2}K^{-2/L}>\tfrac18$, then
$96|V|^{2}|W|^{2}K^{-2/L}>1$. In either case there is nothing to prove. So assume
\begin{equation}\label{eq:IIb-assume}
  K^{1/q}\ \ge\ 400\,|W|^{2},\qquad 3|V|^{2}|W|^{2}K^{-2/L}\ \le\ \tfrac18 .
\end{equation}

\noindent\uline{Step 3: \ref{H1}--\ref{H5} hold.} By \eqref{eq:invariant}, Assumption~\ref{ass:B}\textup{(ii)} reads
$K^{q\upd}\ge K^{5q/L}$, and, as $d=K^{2/q+\upp}$,
Assumption~\ref{ass:B}\textup{(i)} reads $K^{\upp}\ge K^{6402/L}$. So, as $K>1$ by \eqref{eq:IIb-assume},
\[
  \upd\ \ge\ \frac5L\ =\ \frac4{L_{B}}+\frac1{L_{\Phi}} ,\qquad
  \upp\ \ge\ \frac{6402}{L}\ =\ \frac2{L_{B}}+\frac{16}{\ciso L_{B}} .
\]
\ref{H1}: $q\ge3$, $K\ge d$ and $K^2<n\le K^2d^{q/2}$ are in Assumption~\ref{ass:B}, and $G$ is the template.
\ref{H2}: $L_{\Phi}=L_{B}=L\ge3200=8/\ciso$. \ref{H3}: it is Assumption~\ref{ass:B}\textup{(iii)}. \ref{H4}: it is the
display. \ref{H5}: $d\ge K^{2/q}$ by Assumption~\ref{ass:B}\textup{(i)}, and $K^{1/q}\ge400|W|^{2}$ is \eqref{eq:IIb-assume}.

\noindent\uline{Step 4: the error.} By Theorem~\ref{thm:master}, every pair $(\mathbf M,\Pai)$ has
$\mathcal E\le-\frac2{L}(|\mathbf M_{V}|-2+|\mathbf M_{W}|)$. With \ref{H1}, \ref{H3}, \ref{H5} and
\eqref{eq:IIb-assume}, these are the hypotheses of Proposition~\ref{prop:summation}, which bounds the sum of the
two variances, while $(\E[P_{G}\mid x=1]-\E[P_{G}\mid x=0])^{2}=\hat s_{G}^{2}K^{2}>0$ by Lemma~\ref{lem:gap} and
\eqref{eq:shat}. Lemma~\ref{lem:criterion} then bounds the error probability, under either conditioning, by
\[
  4\big(128\,|W|^{2}K^{-1/q}+24\,|V|^{2}|W|^{2}K^{-2/L}\big)\ =\ 512\,|W|^{2}K^{-1/q}+96\,|V|^{2}|W|^{2}K^{-2/L} .
  \qedhere
\]
\end{proof}

\section{Proofs of the corollaries}\label{ss:tech1}

\subsection{Proof of Corollary~\ref{cor:LBpoly}}\label{sss:LBpoly}

\begin{proof}
For $K$ large in terms of $p$ and $C_{\mathrm{poly}}$, $D=\lfloor C_{\mathrm{poly}}\log K\rfloor$
satisfies $2\le D\le C_{\mathrm{poly}}\log K$ and $16D^{2}\le(\log K)^{3}\le n$, and $d=K^{2/p}$ is
$K=d^{p/2}$. The two displayed hypotheses give Assumption~\ref{ass:A}\textup{(i)} and \textup{(ii)}, as $D\le C_{\mathrm{poly}}\log K$,
and Assumption~\ref{ass:A}\textup{(iii)} holds because $d^{\,p-\lceil p\rceil+1}=K^{2(p-\lceil p\rceil+1)/p}$ is a fixed
positive power of $K$, while $D^{1600\lceil p\rceil}$ is at most a fixed power of
$C_{\mathrm{poly}}\log K$. So Assumption~\ref{ass:A} holds, and Theorem~\ref{thm:LB} applies.
\end{proof}

\subsection{Proof of Corollary~\ref{cor:UBpoly}}\label{sss:UBpoly}

\begin{proof}
We take $L$ of order $(\log K)^{1/2}$, the largest value for which the polynomial of
Theorem~\ref{thm:UB} has degree at most $D$, and check Assumption~\ref{ass:B}.
Put $L:=\lfloor(D/(2q+4))^{1/2}\rfloor$, so that the degree $2(q+2)L^{2}$ of
Theorem~\ref{thm:UB} is at most $D$, and $c:=40000q^{2}(C_{\mathrm{poly}}^{1/2}+C_{\mathrm{poly}}^{-1/2})$.
For $K$ large, $\frac12C_{\mathrm{poly}}\log K\le D\le C_{\mathrm{poly}}\log K$, so
$L\ge(400q)^{600q^{2}}$ and
$(C_{\mathrm{poly}}\log K)^{1/2}/(4q)\le L\le(C_{\mathrm{poly}}\log K)^{1/2}$, as $(4q+8)^{1/2}<4q$. The lower bound gives
$K^{6402/L}$ and $K^{5q/L}$ at most $e^{c(\log K)^{1/2}}$, so the hypotheses give
Assumption~\ref{ass:B}. The upper bound gives $K^{-2/L}\le e^{-2C_{\mathrm{poly}}^{-1/2}(\log K)^{1/2}}$, so the
error probability of Theorem~\ref{thm:UB} is at most
\[
  4608q^{2}\big(C_{\mathrm{poly}}\log K\big)^{2}K^{-1/q}
  +7776q^{4}\big(C_{\mathrm{poly}}\log K\big)^{4}e^{-2C_{\mathrm{poly}}^{-1/2}(\log K)^{1/2}}
  \ \le\ e^{-(\log K)^{1/2}/c}
\]
for $K$ large, as $c\ge C_{\mathrm{poly}}^{1/2}$.
\end{proof}

\section*{Acknowledgments}
The work of A.~Carpentier is partially supported by the Deutsche Forschungsgemeinschaft (DFG) -
Project-ID 318763901 - SFB1294 ``Data Assimilation'', Project A03, by the DFG on the
Forschungsgruppe FOR5381 ``Mathematical Statistics in the Information Age - Statistical Efficiency
and Computational Tractability'', Project TP 02 (Project-ID 460867398), by the Universit\'e
franco-allemande (UFA) through the coll\`ege doctoral franco-allemand CDFA-02-25 ``Statistisches
Lernen f\"ur komplexe stochastische Prozesse'', by the European Research Council (ERC) through the
ERC Consolidator Grant SOCE (Nr: 101229569). Funded by the European Union (ERC, project SOCE,
101229569). Views and opinions expressed are however those of the author(s) only and do not
necessarily reflect those of the European Union or the European Research Council. Neither the
European Union nor the granting authority can be held responsible for them. The work of
A.~Carpentier and N.~Verzelen is also supported by the Universit\'e franco-allemande (UFA) through
the coll\`ege doctoral franco-allemand CDFA-02-25 ``Statistisches Lernen f\"ur komplexe
stochastische Prozesse''.

\section*{AI disclosure}
The authors used the AI assistant Claude (Anthropic), through the Claude Code interface, to help
write and edit this manuscript, including the proofs, and to audit the proofs. All the text produced
with this assistance was reviewed by the authors, who take full responsibility for the content of
the manuscript, including all its mathematical statements and proofs.

\appendix
\section{The body of the template}\label{app:body}

Definition~\ref{def:template} builds the template on a \emph{body}: a set $B$ of $|B|=2L_{\Phi}L_{B}$
nodes carrying a graph $G_{B}$, together with a set $\Phi$ of $2L_{\Phi}$ of its nodes.
Throughout, $\partial_{B}(S)$ is the number of edges of $G_{B}$ with exactly one endpoint in
$S\subseteq B$, an edge occurring twice being counted twice. As an edge has
exactly one endpoint in $S$ if and only if it has exactly one in $B\setminus S$,
$\partial_{B}(S)=\partial_{B}(B\setminus S)$.

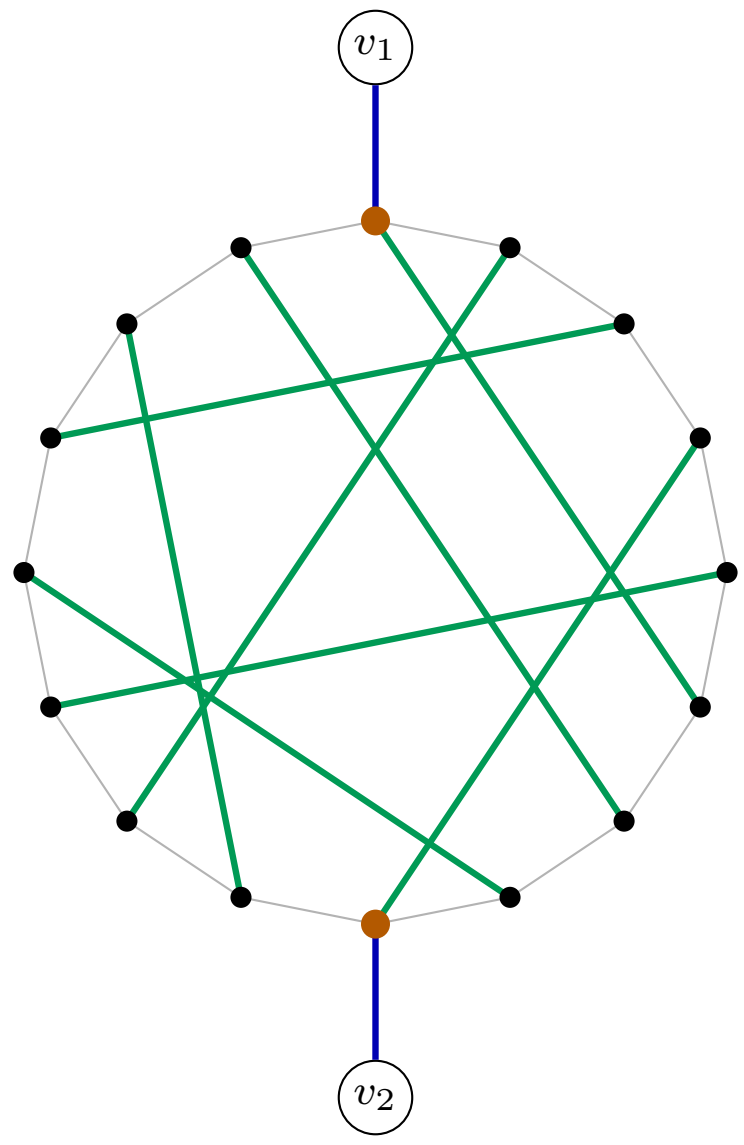
\begin{figure}[ht]
\centering
\begin{tikzpicture}[scale=0.62,
   fast/.style={fill=orange!70!black},
   hub/.style={draw,circle,fill=white,inner sep=1.6pt,font=\scriptsize}]
  \foreach \a in {1,...,16} {
    \pgfmathsetmacro{\ang}{90-(\a-1)*22.5}
    \coordinate (c\a) at (\ang:3.75);
  }
  \foreach \a in {1,...,16} {
    \pgfmathtruncatemacro{\b}{mod(\a,16)+1}
    \draw[gray!60] (c\a) -- (c\b);
  }
  \foreach \a/\b in {1/6,2/11,3/14,4/9,5/12,7/16,8/13,10/15} {
    \draw[very thick, ForestGreen] (c\a) -- (c\b);
  }
  \node[hub] (k1) at (90:5.6) {$v_1$};
  \node[hub] (k2) at (270:5.6) {$v_2$};
  \draw[very thick, blue!70!black] (k1) -- (c1);
  \draw[very thick, blue!70!black] (k2) -- (c9);
  \foreach \a in {1,9} {\fill[fast] (c\a) circle (4.4pt);}
  \foreach \a in {2,3,4,5,6,7,8,10,11,12,13,14,15,16} {\fill[black] (c\a) circle (3.2pt);}
\end{tikzpicture}
\caption{The template at $q=3$, with $|B|=2L_{\Phi}L_{B}=16$, $L_{B}=8$ and $L_{\Phi}=1$. In the body graph $G_{B}$ every node meets
$q$ edges, and $G_{B}$ satisfies \eqref{eq:expand}. The fasteners (orange) are two
nodes far apart in $G_{B}$, and their fastener links (blue) run to the distinguished nodes
$v_{1}$ and $v_{2}$, drawn white. Here
$|W|=L_{\Phi}(q L_{B}+2)=26$. 
(The proof of Theorem~\ref{thm:UB} takes $L_{B}$ much larger. This is drawn at $L_{B}=8$ for legibility.)}
\label{fig:template}
\end{figure}

\subsection{The requirements}\label{app:what}

\begin{enumerate}[label=\textup{(\alph*)},leftmargin=2.4em,itemsep=3pt]
\item \textbf{Every node meets exactly $q$ edges, none of them a self-loop.}
\item 
\textbf{The two inequalities of} \eqref{eq:expand}: for every $S\subseteq B$,
\[
  \text{(i)}\ \ \partial_{B}(S)\ \ge\ q\quad(S\ne\varnothing,\ S\ne B),
  \qquad
  \text{(ii)}\ \ \partial_{B}(S)\ \ge\ \ciso\min\{|S|,|B|-|S|\} .
\]
\item \textbf{The $2L_{\Phi}$ fasteners are at pairwise distance at least $4/\ciso$ in $G_{B}$.}
\end{enumerate}
Here $\ciso>0$ is a constant for which (ii) holds. Proposition~\ref{prop:body} gives $\ciso=1/400$.

Requirement (a) gives the degrees and the link count \eqref{eq:W-count} of
\S\ref{sec:template} and is used in Lemma~\ref{lem:links}, and (i) makes $G$ connected (\S\ref{sec:template}). Past
that, (i) and (ii) are used in Lemma~\ref{lem:fast} --- (i) at a connected subset of $B$ holding one
fastener, (ii) at one holding two or more --- and in Lemma~\ref{lem:small}. Finally, (c) is used only in
Lemma~\ref{lem:fast}, which Lemma~\ref{lem:large} uses.

Three demands are made of $\ciso$: \ref{H2} asks $\ciso L_{B}\ge8$, \ref{H4} asks
$\upp\ge2/L_{B}+16/(\ciso L_{B})$, and (c) asks the fasteners to be at pairwise distance at least
$4/\ciso$. At $\ciso=1/400$ the first two read $L_{B}\ge3200$ and $\upp\ge6402/L_{B}$.

For an even $q$,  the $(q/2)$-fold blow-up
of a cycle --- the graph of \cite{carpentier2025phaseII}, whose fastener density $a$ is our $1/L_{B}$ and whose
Proposition~3.1 plays the role of our Theorem~\ref{thm:master} --- satisfies the first inequality of
\eqref{eq:expand} but not the second: a set $S$ of consecutive layers has
$\partial_{B}(S)$ bounded independently of $|S|$.

\subsection{A body, for every degree}\label{app:grid}

\begin{proposition}\label{prop:body}
Let $q\ge3$, and let $L_{\Phi},L_{B}\ge1$ be integers with $L_{B}>q^{1601}$ and
$|B|=2L_{\Phi}L_{B}\ge(400q)^{1200q^{2}}$. Then there are a graph $G_{B}$ on
$B$ satisfying \textup{(a)} and \eqref{eq:expand} with $\ciso=1/400$, and a set $\Phi$ of $2L_{\Phi}$
nodes satisfying \textup{(c)}.
\end{proposition}

\begin{proof}
We draw $G_{B}$ at random, a cycle through $B$ together with $q-2$ random perfect
matchings, show that \eqref{eq:expand} holds with positive probability, and then choose the fasteners
pairwise far apart (Step~4). By \eqref{eq:appsplit}, a set at which \eqref{eq:expand} fails holds many
edges of the matchings compared with the number of runs of consecutive nodes of the cycle that it
occupies, and Steps~1 to~3 show, by a union bound, that such sets are unlikely.

\noindent\uline{The construction.} Identify $B$ with $\mathbb Z_{|B|}$ and let $G_{B}$ be the cycle joining
each $i$ to $i+1$ together with $q-2$ perfect matchings of $B$, at least one as $q\ge3$, drawn one
after another, each uniform among the perfect matchings of $B$ using no pair already joined.
Before each draw every node is joined to at most $q-1$ others, so the pairs not yet joined form a
graph on the even set $B$ of minimum degree at least $|B|-q>|B|/2$, which has a perfect matching
(Dirac). So each draw is possible, every node meets $q$ distinct edges, \textup{(a)} holds and
$G_{B}$ is simple. It remains to show that \eqref{eq:expand} holds with positive probability. For
$S\subseteq B$, $U(S)$ denotes the number of runs of consecutive positions of
$\mathbb Z_{|B|}$ that $S$ occupies.

\noindent\uline{The identity \eqref{eq:appsplit}.} As $G_{B}$ is simple and every node meets
$q$ edges,
\[
  \partial_{B}(S)\ =\ q|S|-2\,\#\{\text{edges inside }S\}\qquad(S\subseteq B) .
\]
As $\partial_{B}(S)=\partial_{B}(B\setminus S)$, both inequalities of \eqref{eq:expand} are
unchanged when $S$ is replaced by $B\setminus S$, so we may assume $|S|\le|B|/2$. Then
$S\ne B$, $S$ carries $|S|-U(S)$ edges of the cycle, and
\begin{equation}\label{eq:appsplit}
  \partial_{B}(S)\ =\ 2\,U(S)\ +\ \sum\Big(|S|-2\,\#\{\text{its edges inside }S\}\Big) ,
\end{equation}
the sum running over the $q-2$ matchings. Each of its terms is non-negative, a matching
having at most $|S|/2$ edges inside $S$.

\noindent\uline{Step 1: prescribed pairs are unlikely.} Consider the draw of one of the matchings, the pairs
already joined forming a graph of maximum degree at most $q-1$, and fix $j$ disjoint
pairs, none already joined, with $2j\le|B|/2$, and let $|B|\ge80q$. We show that the matching contains all
$j$ of them with probability at most
\begin{equation}\label{eq:appstep1}
  \Big(\frac{4.5}{|B|}\Big)^{j} .
\end{equation}
Let $\{a,b\}$ be the first of the $j$ pairs, and call a pair $\{a',b'\}$ of a matching admissible
when it meets no prescribed pair and neither of the at most $2(q-1)$ nodes already
joined to $a$ or to $b$. Replacing $\{a,b\}$ and an admissible $\{a',b'\}$ by $\{a,a'\}$ and
$\{b,b'\}$, $a'$ being the smaller of the two in $\{0,\dots,|B|-1\}$, sends a matching containing
the $j$ prescribed pairs to a matching that uses no already joined pair, contains the other $j-1$
prescribed pairs and does not contain $\{a,b\}$. The image determines $\{a',b'\}$, since in it $a$
is matched to $a'$ and $b$ to $b'$, hence determines the original matching. A matching containing
the $j$ prescribed pairs has $|B|/2$ pairs, of which $j$ are prescribed and at most
$2(q-1)$ meet a node already joined to $a$ or to $b$, so it has at least
\[
  \frac{|B|}2-j-2(q-1)\ \ge\ \frac{|B|}2-\frac{|B|}4-\frac{|B|}{40}\ =\ \frac{9|B|}{40}
\]
admissible pairs, as $j\le|B|/4$ and $2q\le|B|/40$. Hence the matchings using no already joined pair and containing the $j$
prescribed pairs are at most $40/(9|B|)\le4.5/|B|$ times as many as those containing the
other $j-1$. The matching being uniform among those using no already joined pair, induction on
$j$ gives \eqref{eq:appstep1}, the hypotheses holding for $j-1$. 

\noindent\uline{Step 2: \textup{(ii)} with $\ciso=1/400$.} Let $|S|\le|B|/2$, so that
$\min\{|S|,|B|-|S|\}=|S|$. For $S=\varnothing$ both sides of \textup{(ii)} vanish. For
$S\ne\varnothing$, $U(S)\ge1$, and \eqref{eq:appsplit} gives
\[
  \partial_{B}(S)\ \ge\ 2\ >\ \frac{|S|}{400}\qquad\text{when }|S|<800 .
\]
So let $800\le|S|\le|B|/2$ and suppose $\partial_{B}(S)<|S|/400$. Dropping
from \eqref{eq:appsplit} the terms of the other matchings, which are non-negative,
\[
  2\,\#\{\text{edges of the first matching inside }S\}\ \ge\ 2U(S)+|S|-\partial_{B}(S)
  \ >\ 2U(S)+\frac{399}{400}|S| ,
\]
so the first matching has at least $j:=\lceil\tfrac{399}{800}|S|\rceil+U(S)$ edges inside $S$. As
at most $|S|/2$ of its edges lie inside $S$, $2j\le|S|$, whence
\[
  2U(S)\ \le\ |S|-2\Big\lceil\frac{399}{800}|S|\Big\rceil\ \le\ \frac{|S|}{400},
  \qquad\text{that is}\qquad U(S)\ \le\ \frac{|S|}{800} .
\]
Now count, writing $|S|=k$ and $U(S)=U$. The sets $S$ with these two values number at most
$\binom{|B|}{U}\binom{k-1}{U-1}$, a set being determined by the first positions of its $U$ runs
and by their lengths, $U$ positive integers adding up to $k$. The families of $j$ disjoint pairs
inside such a set number at most $\binom{k}{2j}(2j-1)!!$, and each family is contained in the
first matching with probability at most $(4.5/|B|)^{j}$, by \eqref{eq:appstep1}, as $2j\le k\le|B|/2$, and with probability $0$ if it holds an edge of the cycle. Hence the expected number of sets of size $k$ with
$U$ runs at which \textup{(ii)} fails is at most
\[
  \binom{|B|}{U}\binom{k-1}{U-1}\binom{k}{2j}(2j-1)!!\,\Big(\frac{4.5}{|B|}\Big)^{j} .
\]
Bound the three factors by
\[
\begin{gathered}
  \binom{|B|}{U}\binom{k-1}{U-1}\ \le\ \exp\Big\{U\ln\frac{e^{2}|B|k}{U^{2}}\Big\},\\[3pt]
  \binom{k}{2j}=\binom{k}{k-2j}\ \le\ \Big(\frac{ek}{k-2j}\Big)^{k-2j}\ \le\ (400e)^{k/400},\\[3pt]
  (2j-1)!!\,\Big(\frac{4.5}{|B|}\Big)^{j}\ \le\ e\sqrt j\,\Big(\frac{9j}{e|B|}\Big)^{j}\\
  \le\ e\sqrt k\,\Big(\frac{4.5\,k}{e|B|}\Big)^{\frac{399}{800}k} ,
\end{gathered}
\]
the first by $\binom{k-1}{U-1}\le\binom{k}{U}$ and $\binom{a}{b}\le(ea/b)^{b}$, the second using
$k-2j\le k/400$ and $t\mapsto(ek/t)^{t}$ increasing for $t\le k$, and the third using
$(2j-1)!!\le2^{j}j!\le e\sqrt j\,(2j/e)^{j}$ together with $\tfrac{399}{800}k\le j\le k/2$ and
$4.5k/(e|B|)\le2.25/e<1$. The exponent $U\ln(e^{2}|B|k/U^{2})$ increases in $U$ up to
$U=\sqrt{|B|k}$, and $\sqrt{|B|k}\ge k\ge k/800\ge U$, so
\[
  U\ln\frac{e^{2}|B|k}{U^{2}}\ \le\ \frac{k}{800}\ln\frac{800^{2}e^{2}|B|}{k} .
\]
Multiplying the three bounds, and using $\ln(|B|/k)\ge\ln2$ as $k\le|B|/2$,
\[
\begin{aligned}
  &\binom{|B|}{U}\binom{k-1}{U-1}\binom{k}{2j}(2j-1)!!\,\Big(\frac{4.5}{|B|}\Big)^{j}\\
  &\quad\le\ e\sqrt k\,\exp\Big\{\frac{k}{800}\ln\frac{800^{2}e^{2}|B|}{k}
   +\frac{k}{400}\ln(400e)+\frac{399k}{800}\ln\frac{4.5k}{e|B|}\Big\}\\
  &\quad=\ e\sqrt k\,\exp\Big\{k\Big(\frac{\ln(800^{2}e^{2})}{800}+\frac{\ln(400e)}{400}
   +\frac{399}{800}\ln\frac{4.5}{e}\Big)-\frac{398k}{800}\ln\frac{|B|}{k}\Big\}\\
  &\quad\le\ e\sqrt k\,e^{-k/20} ,
\end{aligned}
\]
as
\[
  \frac{\ln\big(800^{2}e^{2}\big)}{800}+\frac{\ln(400e)}{400}
      +\frac{399}{800}\ln\frac{4.5}{e}-\frac{398}{800}\ln 2
  \ \le\ 0.0193+0.0175+0.2515-0.3448\ <\ -\frac1{20} .
\]
Summing over the at most $k$ values of $U$ and over $k\ge800$, the probability that
\textup{(ii)} fails for some $S$ is at most $\sum_{k\ge800}e\,k^{3/2}e^{-k/20}<\tfrac14$.

\noindent\uline{Step 3: \textup{(i)}.} Let $S\ne\varnothing$ and $|S|\le|B|/2$. If $|S|\ge400q$ and
\textup{(ii)} holds at $S$, then $\partial_{B}(S)\ge|S|/400\geq$. If $|S|\le2$, then, $G_{B}$ being
simple, $\partial_{B}(S)\ge|S|q-2(|S|-1)\geq$. Let now $3\le|S|<400q$ and suppose $\partial_{B}(S)<q$. As
$\partial_{B}(S)=q|S|-2\,\#\{\text{edges inside }S\}$, this means
\[
  \#\{\text{edges inside }S\}\ >\ \tfrac12q(|S|-1) .
\]
The cycle contributes $|S|-U(S)$ of these edges, so the matchings contribute at least
$j:=\big\lfloor\tfrac12q(|S|-1)\big\rfloor+1-|S|+U(S)$ edges inside $S$, and, as
$\lfloor\tfrac12q(|S|-1)\rfloor\ge\tfrac12q(|S|-1)-\tfrac12$,
$q\ge3$ and $|S|\ge3$,
\[
  j-U(S)\ \ge\ \tfrac12\big(q-2\big)\big(|S|-1\big)-\tfrac12\ >\ 0 ,
  \qquad\text{so}\qquad j-U(S)\ \ge\ 1 .
\]
Writing $|S|=k$ and $U(S)=U$, which determine $j$, the sets number at most $|B|^{U}k^{U}$, a run
being named by its first position and its length. The families of $j$ pairs inside such a set,
each assigned to one of the matchings, number at most $(q k^{2})^{j}$, and each is
present with probability at most $(4.5/|B|)^{j}$, by \eqref{eq:appstep1} applied to each matching in turn, conditionally on the earlier ones, as
$k<400q\le|B|/2$, and with probability $0$ if two of its pairs assigned to one matching meet, or
if one of its pairs is already joined when its matching is drawn. So, by $U-j\le-1$, the expected number of sets $S$
with $3\le|S|<400q$ at which \textup{(i)} fails is at most
\[
  \sum_{k=3}^{400q-1}\ \sum_{U=1}^{k}k^{U}\big(4.5q k^{2}\big)^{j}\,|B|^{U-j}
  \ \le\ \frac1{|B|}\sum_{k=3}^{400q-1}\ \sum_{U=1}^{k}k^{U}\big(4.5q k^{2}\big)^{j}\ .
\]
The double sum has at most $(400q)^{2}$ terms. As $k<400q$, $U\le k$,
$4.5q k^{2}\le(400q)^{3}$ and $j\leq k/2$, each term is at most
$(400q)^{400q+600q^{2}}$. So the double sum is at most
$(400q)^{1200q^{2}}/4\le|B|/4$, and the expected number is at most $\tfrac14$.

\noindent\uline{Conclusion of Steps 2 and 3.} With probability at least $1-\tfrac14-\tfrac14>0$,
\textup{(ii)} holds at every $S$ and \textup{(i)} at every $S$ with $3\le|S|<400q$ and
$|S|\le|B|/2$. By the other cases of Step~3 and the reduction to $|S|\le|B|/2$, \eqref{eq:expand} then holds
with $\ciso=1/400$. Fix such a realization.

\noindent\uline{Step 4: \textup{(c)}.} Here $4/\ciso=1600$. A ball of radius $1600$ in $G_{B}$ holds at most
$\sum_{i\le1600}q^{i}\leq^{1601}<L_{B}$ nodes. So the fasteners can be chosen one at a
time, each outside the balls of radius $1600$ around those already chosen, as these hold fewer than
$2L_{\Phi}L_{B}=|B|$ nodes. They are then at pairwise distance more than $4/\ciso$, which is \textup{(c)}.
\end{proof}

\subsection{Another body: a torus}\label{app:torus}

The body of Proposition~\ref{prop:body} is random. An explicit one, a torus, also works, under 
different requirements. We describe it, and say what it satisfies and what changes in Section~\ref{part:ub}, without
proofs.

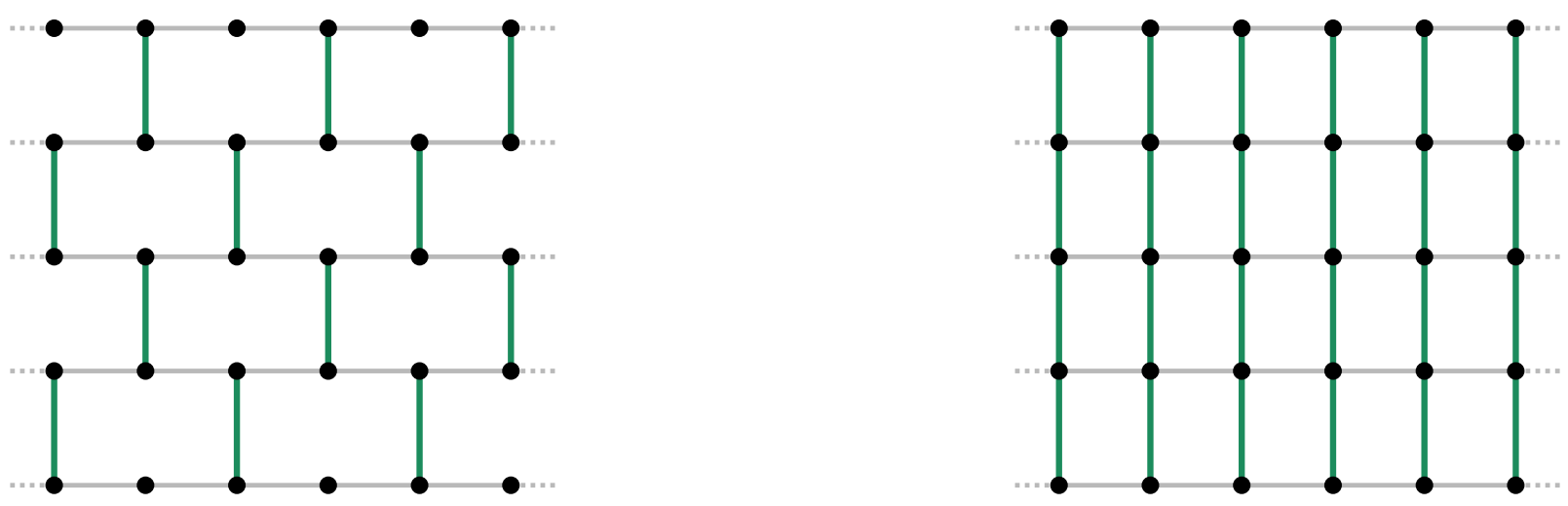
\begin{figure}[ht]
\centering
\begin{tikzpicture}[scale=0.62,
  nd/.style={circle,fill=black,inner sep=1.2pt},
  rw/.style={gray!55,line width=0.9pt},
  cl/.style={ForestGreen!80!black,line width=1.2pt},
  wr/.style={gray!55,line width=0.9pt,densely dotted}]
 \begin{scope}
  \foreach \i in {1,...,5}{
    \foreach \j in {1,...,5}{\draw[rw] (\j,1.25*\i) -- (\j+1,1.25*\i);}
    \draw[wr] (1,1.25*\i) -- ++(-0.5,0); \draw[wr] (6,1.25*\i) -- ++(0.5,0);}
  \foreach \i in {1,...,4}{\foreach \j in {1,...,6}{
    \pgfmathtruncatemacro{\t}{mod(\i+\j,2)}
    \ifnum\t=0 \draw[cl] (\j,1.25*\i) -- (\j,1.25*\i+1.25);\fi}}
  \foreach \i in {1,...,5}{\foreach \j in {1,...,6}{\node[nd] at (\j,1.25*\i){};}}
  \node[font=\footnotesize] at (3.5,-0.55) {$q$ odd: alternating rungs};
 \end{scope}
 \begin{scope}[xshift=11cm]
  \foreach \i in {1,...,5}{
    \foreach \j in {1,...,5}{\draw[rw] (\j,1.25*\i) -- (\j+1,1.25*\i);}
    \draw[wr] (1,1.25*\i) -- ++(-0.5,0); \draw[wr] (6,1.25*\i) -- ++(0.5,0);}
  \foreach \i in {1,...,4}{\foreach \j in {1,...,6}{
    \draw[cl] (\j,1.25*\i) -- (\j,1.25*\i+1.25);}}
  \foreach \i in {1,...,5}{\foreach \j in {1,...,6}{\node[nd] at (\j,1.25*\i){};}}
  \node[font=\footnotesize] at (3.5,-0.55) {$q$ even: every rung};
 \end{scope}
\end{tikzpicture}
\caption{A portion of the torus body at $r=1$, that is at $q=3$ (left) and $q=4$ (right). The rows wrap
round (dotted), and so do the columns. At a larger degree the edges along the rows are repeated $r$ times.}
\label{fig:torus}
\end{figure}

Let $a$ be an even integer with $a\ge2q$, let $B:=\mathbb Z_{a}\times\mathbb Z_{a}$, so that $|B|=a^{2}$,
and let $r:=\lfloor(q-1)/2\rfloor\ge1$. The graph $G_{B}$ has $r$ parallel edges joining $(i,j)$ to
$(i,j+1)$, for all $i,j$, and one edge, a \emph{rung}, joining $(i,j)$ to $(i+1,j)$, for every $j$ when
$q$ is even and for the $j$ with $i+j$ even when $q$ is odd. Every node then meets
$2r+2$, resp.\ $2r+1$, that is $q$, edges. This is the square lattice, resp.\ the honeycomb lattice,
on a torus (Figure~\ref{fig:torus}).

\noindent\uline{What it satisfies.} It satisfies \textup{(a)} and \textup{(i)}. It does not satisfy
\textup{(ii)} with $\ciso$ independent of $|B|$, as a band of $a/2$ consecutive rows has $|B|/2$ nodes and at most
$2a$ edges leaving it. It satisfies instead, for every $S\subseteq B$,
\begin{equation}\label{eq:torusiso}
  \partial_{B}(S)\ \ge\ \frac25\sqrt{\min\{|S|,|B|-|S|\}} .
\end{equation}

\noindent\uline{What changes in Section~\ref{part:ub}.} The inequality \textup{(ii)} is used only in Lemmas~\ref{lem:fast}
and~\ref{lem:small}. With \eqref{eq:torusiso} in its place, and $\ciso:=2/5$, they hold under two changes.
\begin{itemize}[leftmargin=1.6em,itemsep=2pt,topsep=3pt]
\item Lemma~\ref{lem:small} holds if $L_{\Phi}\le\ciso^{2}L_{B}/(16q^{2})$, and if the lower bound
$16/(\ciso L_{B})$ on $\upp$ in \ref{H4} is replaced by $16q/(\ciso^{2}L_{B})$.
\item Lemma~\ref{lem:fast} holds if the $2L_{\Phi}$ fasteners are at pairwise distance at least $100L_{\Phi}$ in
$G_{B}$. For instance, with $m$ even, $2L_{\Phi}=m^{2}$ and $a=mh$, the fasteners can be the nodes $(hk,hl)$,
$k,l\in\{0,\dots,m-1\}$, at pairwise distance at least $h$, and $h\ge50m^{2}$ amounts to
$L_{\Phi}\le\sqrt{L_{B}}/100$, as $L_{B}=h^{2}$.
\end{itemize}
The other results of Section~\ref{part:ub} up to Proposition~\ref{prop:summation} are unchanged. As \ref{H4} asks $\upd\ge4/L_{B}+1/L_{\Phi}$, and $L_{\Phi}$ is now much smaller than
$L_{B}$, this body needs at degree $D$ a value of
$\upd$ at least of order $D^{-1/3}$, larger than the order
$D^{-1/2}$ of Proposition~\ref{prop:body}.

\section{Elementary technical lemmas}\label{app:tech}

\begin{lemma}\label{lem:sizes}
Every graph of \S\ref{sec:graphs} has $|V|\le|E|+2$ and $|W|\le|E|$. In particular $|V|\le2|E|$
when $|E|\ge2$.
\end{lemma}

\begin{proof}
Every variable, and every individual other than $v_1,v_2$, carries an edge, and an edge sits at
one variable and one individual. So $|W|\le|E|$ and $|V|-2\le|E|$. Finally $|E|+2\le2|E|$ when
$|E|\ge2$.
\end{proof}

\begin{lemma}\label{lem:dK}
Under Assumption~\ref{ass:A},
\begin{equation}\label{eq:dK}
d\ \ge\ d^{\,p-\lceil p\rceil+1}\ \ge\ D^{1600\lceil p\rceil}\ \ge\ D^{3200},\qquad
K\ =\ d^{p/2}\ \ge\ D^{800p\lceil p\rceil}\ \ge\ D^{3200} .
\end{equation}
\end{lemma}

\begin{proof}
As $d\ge1$ and $p-\lceil p\rceil+1\le1$, Assumption~\ref{ass:A}\textup{(iii)} gives
$d\ge d^{\,p-\lceil p\rceil+1}\ge D^{1600\lceil p\rceil}$. So
$K=d^{p/2}\ge D^{800p\lceil p\rceil}$, and $\lceil p\rceil\ge p\ge2$ gives the bounds by
$D^{3200}$.
\end{proof}

\begin{lemma}\label{lem:embcount}
Under Assumption~\ref{ass:A}, for $N$ any of $n,d,K$ and every $0\le a\le a+t\le2D$,
\begin{equation}\label{eq:Picount}
  \tfrac78\,N^{t}\ \le\ (N-a)(N-a-1)\cdots(N-a-t+1)\ \le\ N^{t} .
\end{equation}
\end{lemma}

\begin{proof}
By $n\ge16D^{2}$ and \eqref{eq:dK}, $4D^{2}\le N/4$, so $a+t\le2D<N$, every factor is positive and at most $N$, and the upper bound follows. By the Weierstrass product inequality,
\[
  (N-a)\cdots(N-a-t+1)\ =\ N^{t}\prod_{k=0}^{t-1}\Big(1-\frac{a+k}{N}\Big)
  \ \ge\ N^{t}\Big(1-\frac{(a+t)^{2}}{2N}\Big)\ \ge\ N^{t}\Big(1-\frac{4D^{2}}{2N}\Big)
  \ \ge\ \frac78\,N^{t} .
\]
\end{proof}

\section{Technical lemmas taken from the literature}\label{app:lit}

\subsection{Proof of Lemma~\ref{lem:split}}\label{sss:split}

\begin{proof}
This is the addition theorem for Hermite polynomials \cite[(18.18.11)]{DLMF}, written there for
$H_m(z)=2^{m/2}\He_m(\sqrt2\,z)$ and holding for $\He_m$, with $(a_1,a_2)=(\Delta/\sqrt d,1)$ and
$(x_1,x_2)=(X_{k^*(i),j},\epsilon_{ij})$, so that $a_1x_1+a_2x_2=Y_{ij}$ \eqref{eq:model} and
$a_1^2+a_2^2=1+\Delta^2/d$. Multiply it by $\sqrt{m!}$ and use
$\sqrt{m!}/(\mathbf s!\,\mathbf n!)=\sqrt{\binom m{\mathbf s}}/\sqrt{\mathbf s!\,\mathbf n!}$. It is
also the comparison of the coefficients of $t^m$ in
$e^{tY_{ij}-t^2(1+\Delta^2/d)/2}=e^{t\frac{\Delta}{\sqrt d}X_{k^*(i),j}-t^2\frac{\Delta^2}{2d}}\,e^{t\epsilon_{ij}-t^2/2}$,
each factor being expanded by $\sum_{m\ge0}\frac{t^m}{m!}\He_m(z)=e^{tz-t^2/2}$.
\end{proof}

\subsection{Proof of Proposition~\ref{prop:sym}}\label{sss:sym}

\begin{proof}
We follow \cite[Lemmas~1 and~4]{carpentier2026low}, where the same reduction is made in a related model.
Permuting the rows $Y_{i,\cdot}$ other than $1,2$, permuting the columns $Y_{\cdot,j}$, and
flipping the signs of a set of columns are linear maps which, by \eqref{eq:model}, leave the law
of $(Y,x)$ unchanged. So averaging a polynomial $P$ of degree at most $D$ over them keeps
$\E[(x-\E x)P]$, keeps the degree at most $D$, and does not increase $\|P\|$. Expand $P$ in the basis
of Definition~\ref{def:hermite}. As $\He_m(-Z)=(-1)^m\He_m(Z)$, averaging over the $2^d$ sign
flips keeps exactly the products with $\sum_im_{i,j}$ even for every column $j$, that is, those
indexed by graphs with every $\deg_G(w)$ even, so with $|E_G|=0$ or $2\le|E_G|\le D$. The
permutations act transitively on $\Pi_G$, so averaging over them turns the product of $(G,\pi)$
into a positive multiple of $P_G$ \eqref{eq:PG}, and leaves the constant $1$ unchanged. Hence, if
$\E[(x-\E x)P]>0$, the average of $P$ is a non-zero element of $\Span(1,(P_G)_G)$ whose ratio in
\eqref{eq:corr} is at least that of $P$. If $\E[(x-\E x)P]\le0$, the ratio at $P$ is at most $0$,
its value at the constant $1$. Conversely, $\Span(1,(P_G)_G)$ consists of polynomials of degree
at most $D$.
\end{proof}

\subsection{Proof of Lemma~\ref{lem:wick}}\label{sss:wick}

\begin{proof}
As $\He_a(Z)$ is the Wick power $:\!Z^a\!:$, this is the formula for the expectation of a product
of Wick products as a sum over the Feynman diagrams with no edge inside a factor
\cite[Chapter~3]{janson1997gaussian}. Directly: by $\sum_{a\ge0}\frac{t^a}{a!}\He_a(Z)=e^{tZ-t^2/2}$ and
$\E e^{(t_1+\dots+t_r)Z}=e^{(t_1+\dots+t_r)^2/2}$,
\[
  \sum_{a_1,\dots,a_r\ge0}\frac{t_1^{a_1}\cdots t_r^{a_r}}{a_1!\cdots a_r!}\,
  \E\Big[\prod_{i=1}^r\He_{a_i}(Z)\Big]\ =\ \prod_{i<j}e^{t_it_j} ,
\]
and expanding the right-hand side gives
$\mu(a_1,\dots,a_r)=\sum\prod_ia_i!/\prod_{i<j}\varrho_{ij}!$ over the symmetric arrays of
non-negative integers with $\sum_{j\ne i}\varrho_{ij}=a_i$, the summand being the number of
Wick pairings with $\varrho_{ij}$ pairs joining the $i$th factor to the $j$th (the $a_i$ points of each factor $i$ are split among the other factors in $a_i!/\prod_{j\ne i}\varrho_{ij}!$ ways, and the $\varrho_{ij}$ points sent from $i$ to $j$ are matched with those sent from $j$ to $i$ in $\varrho_{ij}!$ ways). Items
\textup{(i)} and (ii) are read off the count.
\end{proof}

\section{Proof of Corollary~\ref{cor:partition}}\label{app:MOM}

\begin{lemma}[Median of means]\label{lem:MoM}
Assume that $n\ge \beta+2$ and $d\ge \beta$, and that, in the model~\eqref{eq:model} with parameters
$(K,n_0,d_0,\Delta_0)$ in place of $(K,n,d,\Delta)$, the test $\psi$ satisfies
$\Pr[\psi(Y)\ne x\mid x=\xi]\le1/4$ for $\xi\in\{0,1\}$. Then $\Pr[\hat x=x]\ge1-n^{-3}$.
\end{lemma}

\begin{proof}[Proof of Lemma~\ref{lem:MoM}]
\emph{Law of a batch.} As $\Delta/\sqrt d=\Delta_0/\sqrt{d_0}$, the entries of $Y^{(b)}$ are
$\frac{\Delta_0}{\sqrt{d_0}}X_{k^*(u),j}+\epsilon_{u,j}$, for $u\in\{1,2\}\cup I_b$ and
$j\in J_b$. The $n_0$ labels $k^*(u)$ are i.i.d.\ uniform on $[K]$, and the $X_{k,j}$ and
$\epsilon_{u,j}$, $k\in[K]$, $j\in J_b$, are i.i.d.\ $\mathcal N(0,1)$. So $Y^{(b)}$ follows the
model~\eqref{eq:model} with parameters $(K,n_0,d_0,\Delta_0)$, in which $x$ plays the same role as
in that model.

\emph{Conditioning on the labels of the pair.} In the model~\eqref{eq:model}, let $\sigma$ be a
permutation of $[K]$, and put $k'(u):=\sigma(k^*(u))$ and $X'_{\sigma(k),j}:=X_{k,j}$. Then
$(k',X',\epsilon)$ has the same law as $(k^*,X,\epsilon)$, and it gives the same data $Y$, as
$X'_{k'(u),j}=X_{k^*(u),j}$. Hence the conditional law of $Y$ given $(k^*(1),k^*(2))=(a,a')$ is
also its conditional law given $(k^*(1),k^*(2))=(\sigma(a),\sigma(a'))$. As the permutations of
$[K]$ map any pair $(a,a)$ to any other pair $(a'',a'')$, and any pair $(a,a')$ with $a\ne a'$ to
any other such pair, this conditional law only depends on $\one\{a=a'\}$, and it is therefore the
conditional law of $Y$ given $x=\one\{a=a'\}$. Applied to the batches, this gives, for all
$a,a'\in[K]$ and $b\in[\beta]$,
\[
  \Pr\big[\psi(Y^{(b)})\ne x\ \big|\ k^*(1)=a,\ k^*(2)=a'\big]\ \le\ \frac14 .
\]

\emph{Independence.} Given $(k^*(1),k^*(2))$, the matrix $Y^{(b)}$ is a function of the labels
$k^*(u)$, $u\in I_b$, of the $X_{k,j}$, $k\in[K]$, $j\in J_b$, and of the $\epsilon_{u,j}$,
$u\in\{1,2\}\cup I_b$, $j\in J_b$. As the sets $I_b$ are disjoint, as the sets $J_b$ are
disjoint, and as all these variables are independent and independent of $(k^*(1),k^*(2))$, the
matrices $Y^{(1)},\dots,Y^{(\beta)}$ are independent conditionally on $(k^*(1),k^*(2))$.

\emph{Majority vote.} Let $E_b:=\one\{\psi(Y^{(b)})\ne x\}$. If $\hat x\ne x$, then at least half of
the $\beta$ tests are wrong, that is $\sum_bE_b\ge \beta/2$. Conditionally on $(k^*(1),k^*(2))$, the $E_b$
are independent, with values in $\{0,1\}$ and means at most $1/4$, so Hoeffding's inequality gives
\[
  \Pr\Big[\sum_{b=1}^{\beta}E_b\ge\frac \beta2\ \Big|\ k^*(1),k^*(2)\Big]\ \le\
  \exp\Big(-\frac2\beta\Big(\frac \beta4\Big)^2\Big)\ =\ e^{-\beta/8}\ \le\ n^{-3} ,
\]
as $\beta\ge24\log n$. Integrating over $(k^*(1),k^*(2))$ gives $\Pr[\hat x\ne x]\le n^{-3}$.
\end{proof}

As $n\le K^2d^{q/2}\le K^{2+q/2}$, we have $\beta\le24(2+q/2)\log K+1$, so that, for $K$
large, $n\ge \beta+2$ and $2\beta\le\min\{50\log n,\ e^{c(\log K)^{1/2}}(2\beta)^{-q/2}\}$; in particular
$d\ge2\beta$. Then $n_0\le n$, and
\[
  n_0\ \ge\ \frac{n-2}\beta+1\ \ge\ \frac n\beta,\qquad d_0\ \ge\ \frac d\beta-1\ \ge\ \frac d{2\beta} .
\]
We check the conditions of Corollary~\ref{cor:UBpoly} for the parameters $(K,n_0,d_0,\Delta_0)$; the threshold on $K$ in
Corollary~\ref{cor:UBpoly} only depends on $q$ and $C_{\mathrm{poly}}$.
\begin{itemize}
\item $K\ge d\ge d_0$, and $d_0\ge d/(2\beta)\ge K^{2/q}e^{c(\log K)^{1/2}}$, as $2\beta\le e^{c(\log K)^{1/2}}$.
\item As $\Delta_0^2=\Delta^2d_0/d$,
\[
  \frac{\Delta_0^{2q}\,n_0}{K^{2}d_0^{\,q/2}}
  \ =\ \frac{\Delta^{2q}\,n}{K^{2}d^{\,q/2}}\cdot\frac{n_0}n\Big(\frac{d_0}d\Big)^{q/2}
  \ \ge\ e^{2c(\log K)^{1/2}}\,\beta^{-1}(2\beta)^{-q/2}\ \ge\ e^{c(\log K)^{1/2}} .
\]
\item $\Delta_0^2=\Delta^2d_0/d<d_0/\sqrt d\le\sqrt{d_0}$, as $d_0\le d$.
\item $n_0\le n\le K^2(d/(50\log n))^{q/2}\le K^2d_0^{q/2}$, as $2\beta\le50\log n$. By the two previous items,
$K^2d_0^{q/2}<\Delta_0^{2q}n_0<d_0^{q/2}n_0$, so that $n_0>K^2$.
\end{itemize}
So Corollary~\ref{cor:UBpoly} gives a test $\psi$ which, in the model with parameters
$(K,n_0,d_0,\Delta_0)$, decides $x$ with error probability at most $e^{-(\log K)^{1/2}/c}\le1/4$
under either conditioning, for $K$ large. By Lemma~\ref{lem:MoM}, $\Pr[\hat x=x]\ge1-n^{-3}$.

\section{Proof of Theorem~\ref{thm:LBsym}}\label{app:LBsym}
We follow the proof of Theorem~\ref{thm:LB} (Section~\ref{part:lb}) and only describe the changes.
Write $\Corr^{\mathrm{sym}}_{\le D}$ for the correlation of the statement, and
$\mathcal G^{\mathrm{even}}_{\le D}$ for the set of the graphs $G\in\mathcal G_{\le D}$ all of whose
individual nodes, $v_1$ and $v_2$ included, have even degree.

\noindent\uline{Step 1: reduction to the model \eqref{eq:model} and to the even graphs.} Put
$\epsilon'_{i,j}:=z_i\epsilon_{i,j}$ and
$Y'_{i,j}:=z_iY_{i,j}=\frac{\Delta}{\sqrt d}X_{k^*(i),j}+\epsilon'_{i,j}$. Given $z$, the
$\epsilon'_{i,j}$ are i.i.d.\ $\mathcal N(0,1)$ and independent of $(k^*,X)$. So $(k^*,X,\epsilon')$
has the law of Section~\ref{sec:model} and is independent of $z$: the pair $(x,Y')$ follows the
model \eqref{eq:model}, and $z$ is independent of $(x,Y')$.

Let $P$ be a polynomial of degree at most $D$, and let
$\bar P(y):=2^{-n}\sum_{s\in\{-1,1\}^n}P\big((s_iy_{i,j})_{i,j}\big)$ be its average over the sign
changes of the rows. As $Y=(z_iY'_{i,j})_{i,j}$, with $z$ uniform on $\{-1,1\}^n$ and independent of
$(x,Y')$, we have $\E[P(Y)\mid x,Y']=\bar P(Y')$, so that
\[
  \E\big[(x-\E x)P(Y)\big]\ =\ \E\big[(x-\E x)\bar P(Y')\big]
  \qquad\text{and}\qquad
  \E\big[P(Y)^2\big]\ \ge\ \E\big[\bar P(Y')^2\big] .
\]
Expand $P$ in the basis of Definition~\ref{def:hermite}. As $\He_m(-y)=(-1)^m\He_m(y)$, the sign
change of the row $i$ multiplies the product of multiplicities $(m_{i,j})$ by
$(-1)^{\sum_jm_{i,j}}$. So $\bar P$ keeps exactly the products with $\sum_jm_{i,j}$ even for every
$i\in[n]$, that is, those indexed by graphs all of whose individual nodes have even degree. The
averaging of the proof of Proposition~\ref{prop:sym}, applied to $\bar P$ in the model
\eqref{eq:model}, sends each of these products to $0$ or to a positive multiple of $P_G$, $G$ being
its graph. Hence, as in that proof,
\[
  \Corr^{\mathrm{sym}}_{\le D}\ \le\ \sup\bigg\{\frac{\E\big[(x-\E x)P(Y')\big]}{\sqrt{\E\big[P(Y')^2\big]}}\ :\
  P\in\Span\big(1,(P_G)_{G\in\mathcal G^{\mathrm{even}}_{\le D}}\big),\ P\ne0\bigg\} .
\]
From now on, we work in the model \eqref{eq:model}, with observations $Y'$.

\noindent\uline{Step 2: the results of Section~\ref{part:lb}.} The signal and noise individuals are
still defined by \eqref{eq:blocksplit}, with the threshold $p$. The results of
Sections~\ref{sec:hermite} to~\ref{sec:matrix} are about the model \eqref{eq:model}, and use
Assumption~\ref{ass:A} only through \eqref{eq:dK}: in Lemma~\ref{lem:embcount}, which needs
$d,K\ge16D^2$, and in Lemma~\ref{lem:colsum}, which needs $K\ge D^{3\lceil p\rceil+101}$. As $d\ge1$
and $p-\lceil p\rceil_{\mathrm{even}}+2\le2$, condition \textup{(iii)} gives
$d^2\ge D^{1600\lceil p\rceil_{\mathrm{even}}}$, so that
$K=d^{p/2}\ge d\ge D^{800\lceil p\rceil_{\mathrm{even}}}\ge D^{800\lceil p\rceil}$. This gives both
requirements, as $D\ge2$ and $\lceil p\rceil\ge2$. So these results hold under the conditions of
Theorem~\ref{thm:LBsym}.

\noindent\uline{Step 3: Proposition~\ref{prop:reduce} on the even graphs.} We claim that
\[
  \Corr^{\mathrm{sym}}_{\le D}\ \le\ \frac{10}3\Big(\sum_{G\in\mathcal G^{\mathrm{even}}_{\le D}}\eta_G^2\Big)^{1/2} .
\]
The proof of Proposition~\ref{prop:reduce} applies with $\mathcal G_{\le D}$ replaced by
$\mathcal G^{\mathrm{even}}_{\le D}$, both for $G$ and for $H$, and with Step~1 in place of
Proposition~\ref{prop:sym}. Three points need checking. In Step~1 of that proof,
$\Span\big(1,(P_G)_{G\in\mathcal G^{\mathrm{even}}_{\le D}}\big)$ is also
$\Span\big(1,(\Psh_G)_{G\in\mathcal G^{\mathrm{even}}_{\le D}}\big)$: removing a deletable set $V'$
does not change the degrees of the individuals outside $V'$ (Step~3 of the proof of
Lemma~\ref{lem:deletable}\textup{(a)}), so that $\Gs\in\mathcal G^{\mathrm{even}}_{\le
D}\cup\{\varnothing\}$ for every $G\in\mathcal G^{\mathrm{even}}_{\le D}$ and $V'\in\mathcal D(G)$.
In Step~4, Corollary~\ref{cor:diag} still gives $A_{G,G}\ge\frac12$ and, the entries $A_{H,G}$ being
non-negative, the off-diagonal column sums restricted to $H\in\mathcal G^{\mathrm{even}}_{\le D}$
are at most those of Lemma~\ref{lem:colsum}, hence at most $\frac15$. In Step~5, the $\bar Q_H$,
$H\in\mathcal G^{\mathrm{even}}_{\le D}$, are orthonormal by Lemma~\ref{lem:ortho}, and Bessel's
inequality applies to them.

\noindent\uline{Step 4: Lemma~\ref{lem:master} on the even graphs.} We claim that
$\eta_G\le K^{-1}D^{-2|E_G|}$ for every $G\in\mathcal G^{\mathrm{even}}_{\le D}$. This is the only
step where $\lceil p\rceil_{\mathrm{even}}$ enters. In the proof of Lemma~\ref{lem:master}, a signal
individual of $G$ has an even degree smaller than $p$, hence at most
$\lceil p\rceil_{\mathrm{even}}-2$, and a noise individual other than $v_1,v_2$ has an even degree
at least $p$, hence at least $\lceil p\rceil_{\mathrm{even}}$. So
$|E_{G\mathbf s^\dagger}|\le(\lceil p\rceil_{\mathrm{even}}-2)|V_{G\mathbf s^\dagger}|$ and, by
\textup{(iii)}, the first bracket of that proof is at most
\[
  d^{-\frac{(p-\lceil p\rceil_{\mathrm{even}}+2)|V_{G\mathbf s^\dagger}|}{4}}
  \ \le\ D^{-400\lceil p\rceil_{\mathrm{even}}|V_{G\mathbf s^\dagger}|}
  \ \le\ D^{-200|E_{G\mathbf s^\dagger}|} .
\]
Likewise, $|E_{G\mathbf n^\dagger}|\ge\lceil p\rceil_{\mathrm{even}}(|V_{G\mathbf n^\dagger}|-2)$,
and the bound on the second bracket holds with $\lceil p\rceil_{\mathrm{even}}$ in place of
$\lceil p\rceil$, by \textup{(i)} and \textup{(ii)}. The rest of the proof is unchanged.

\noindent\uline{Step 5: conclusion.} As $\mathcal G^{\mathrm{even}}_{\le D}\subseteq\mathcal
G_{\le D}$, Step~4 and the summation in the proof of Theorem~\ref{thm:LB} give
$\sum_{G\in\mathcal G^{\mathrm{even}}_{\le D}}\eta_G^2\le\frac9{7K^2}$, and Step~3 gives
$\Corr^{\mathrm{sym}}_{\le D}\le\frac{10}{\sqrt7\,K}\le\frac4K$.

\bibliographystyle{alpha}
\bibliography{biblio}

\end{document}